\documentclass[11pt, a4paper]{article}

\title{Cameron--Martin theorem for SLE welding measures}
\date{\today}
\usepackage[T1]{fontenc}

\usepackage[utf8]{inputenc}

\usepackage{mathtools}
\usepackage{lmodern}
\usepackage{amsthm,amsmath,amsfonts,amssymb}
\usepackage{graphicx}
\usepackage{enumerate,enumitem}
\usepackage{authblk}
\usepackage{cite,url}
\usepackage[normalem]{ulem}
\usepackage{mathrsfs}
\usepackage{bbm}

\usepackage{multirow}
\usepackage{array}
\newcolumntype{P}[1]{>{\centering\arraybackslash}p{#1}}
\usepackage{longtable,booktabs}

\usepackage{hyperref}
\hypersetup{
    colorlinks=true, 
    linkcolor=blue, 
    urlcolor=black, 
    citecolor=black,
    linktoc=all 
}

\usepackage[font=normal,labelfont=bf,labelsep=quad]{caption}

\usepackage[activate={true,nocompatibility},final,tracking=false,kerning=false,spacing=false,expansion=false,factor=1100,stretch=20,shrink=20]{microtype}
\microtypecontext{spacing=nonfrench}

\allowdisplaybreaks

\setlist[enumerate]{topsep = 1ex, leftmargin=1cm, itemsep= -2pt}
\setlist[itemize]{topsep = 1ex, leftmargin=1cm, itemsep= -2pt}

\let\OLDthebibliography\thebibliography
\renewcommand\thebibliography[1]{
  \OLDthebibliography{#1}
  \setlength{\parskip}{1pt}
  \setlength{\itemsep}{2pt}
}

\newtheorem{thm}{Theorem}[section]
\newtheorem{cor}[thm]{Corollary}
\newtheorem{conj}[thm]{Conjecture}
\newtheorem{lem}[thm]{Lemma}
\newtheorem{prop}[thm]{Proposition}

\theoremstyle{definition} 
\newtheorem{df}[thm]{Definition}
\newtheorem*{definition*}{Definition}

\newtheorem{remark}[thm]{Remark}

\numberwithin{equation}{section}

\usepackage{floatrow}

\usepackage{mathtools}

\usepackage[dvipsnames]{xcolor}
\usepackage{tikz}
\usetikzlibrary{arrows.meta,fit,positioning}
\usepackage{placeins}

\global\long\def\ii{\mathfrak{i}}

\renewcommand{\liminf}{\varliminf}
\renewcommand{\limsup}{\varlimsup}

\newcommand{\abs}[1]{\left\lvert #1 \right \rvert}

\newcommand{\brac}[1]{\left \langle #1 \right \rangle}
\newcommand{\norm}[1]{\lVert #1 \rVert}

\newcommand{\mc}[1]{\mathcal{#1}}
\newcommand{\m}[1]{\mathbb{#1}}

\makeatletter
\newcommand{\Rm}[1]{\expandafter\@slowromancap\romannumeral #1@}
\makeatother

\renewcommand\Re{\operatorname{Re}}

\def\W{\m R/ 2\pi \m Z}
\def\KK{\mc K}

\def\QIL{\mathrm{QIL}}
\def\QIR{\mathrm{QIR}}

\def\QI{\mathrm{QI}}
\def\QIC{\mathrm{QIC}}
\def\Rot{\mathrm{Rot}}
\def\1{\mathbf{1}}

\usepackage{comment}

\def\U{\operatorname{U}}

\def\SLE{\operatorname{SLE}}

\def \Homeo {\operatorname{Homeo}}
\def\mob{\mathrm{M\ddot ob}}
\def\Diff{\operatorname{Diff}}

\def\QS{\operatorname{QS}}
\def\Sy{\operatorname{S}}

\def\tr{\operatorname{Tr}}

\def\WP{\operatorname{WP}}

\def \WP {\operatorname{WP}}
\def \HS {\operatorname{HS}}
\def \QS {\operatorname{QS}}

\def\a{\alpha}

\def\g{\gamma}
\def\G{\Gamma}
\def\d{\delta}
\def\D{\Delta}
\def\e{\epsilon}

\def\t{\theta}

\def\l{\lambda}
\def\k{\kappa}

\def\L{\Lambda}
\def\s{\sigma}

\def\o{\omega}
\def\O{\Omega}
\def\vare{\varepsilon}

\def\dist{{\rm dist}}

\def\dd{\mathrm{d}}

\def\var{\mathbf{Var}}
\def\cov{\mathbf{Cov}}

\def\1{\mathbf{1}}

\author{Shuo Fan \thanks{Tsinghua University,shuofan.math@gmail.com}}

\begin{document}

\maketitle
\begin{abstract}
For $\kappa\in(0,4]$, we prove a Cameron--Martin theorem for the
welding law induced by SLE$_\kappa$ loop measure. Under simultaneous
left and right composition $h\mapsto\varphi\circ h\circ\psi^{-1}$,
the transformed law is mutually absolutely continuous with respect
to the original law exactly when both maps belong to the
Weil--Petersson group, and is mutually singular otherwise. Thus the
Weil--Petersson group is the quasi-invariance group of the welding law
under composition on either side. We also prove that it is its own
normalizer in the group of orientation-preserving circle homeomorphisms.
For Weil--Petersson pairs, we derive an explicit Radon--Nikodym
derivative in terms of a relative action well defined for SLE
weldings and prove that it belongs to $L^{1+\varepsilon}$ for some
$\varepsilon>0$ depending on the pair.
\end{abstract}
\tableofcontents
\section{Introduction}

We study the law of SLE welding homeomorphisms under deterministic left
and right composition. Our main result is a Cameron--Martin theorem:
simultaneous composition preserves the measure class precisely when
both maps belong to the Weil--Petersson group, and every other pair
gives a mutually singular law. We also determine the Radon--Nikodym
derivative in the Weil--Petersson directions.

\subsection{Background}

The classical Cameron--Martin theorem \cite{CM44} identifies the
translations preserving the measure class of a Gaussian measure and
gives their Radon--Nikodym derivatives. Translations outside the
Cameron--Martin space produce mutually singular measures. For a random
circle welding, left and right composition replace translation, and
the Weil--Petersson group takes the place of the Cameron--Martin space.

The Schramm--Loewner evolution (SLE$_\kappa$) loop measure is a
one-parameter family of $\sigma$-finite measures on loops on the Riemann
sphere. For $\kappa\in(0,4]$, it is supported on Jordan curves. SLE
first arose as a description of interfaces in the scaling limits of
critical lattice models
\cite{Sch00,SW01,LSW04LERWUST,SS06-6,SS09-4,CDCH}; its loop measures
were constructed and studied in
\cite{KS07,werner_measure,MR3,Benoist_loop,zhan2020sleloop}.
Conformal welding gives a second route to these measures.
Sheffield's quantum zipper \cite{Quantum_zipper} produces SLE
interfaces by gluing Liouville quantum gravity surfaces, and
\cite{AHS23} obtains SLE loop measure by welding quantum disks.
The construction uses Gaussian fields through a nonlinear welding map.
There is no a priori Gaussian structure on the loop measure or on the
induced law on circle homeomorphisms, so the classical Cameron--Martin
theorem does not apply directly to these laws.

For a Jordan loop separating $0$ from $\infty$, let $f:\m D\to\Omega$
and $g:\m D^*\to\Omega^*$ be conformal maps onto its complementary domains.
Their boundary values define the welding homeomorphism $h=g^{-1}\circ f$.
The group $\WP(\m S^1)$ consists of the absolutely continuous
orientation-preserving circle homeomorphisms $\varphi$ with
$\varphi'>0$ almost everywhere and $\log\varphi'\in H^{1/2}(\m S^1)$.
Equivalently, they admit quasiconformal extensions with hyperbolically
square-integrable Beltrami coefficient \cite{TT06,Shen18}.
Their welded curves are exactly the finite-Loewner-energy curves.
The identity between Loewner energy and the universal Liouville action
\cite{TT06,W2}, together with its role as the small-$\kappa$ rate
function for SLE \cite{W1,RW,PW24_LDP_multichordal}, motivates the
Cameron--Martin analogy.

Fan and Sung \cite{F1} prove WP quasi-invariance for GMC welding laws,
including SLE welding laws for $\kappa<4$. They use the action on the
Cameron--Martin space of a log-correlated field and a Hilbert--Schmidt
criterion for mutual absolute continuity of Gaussian laws.
Baverez and Jego \cite[Corollary~1.3]{BJ25} obtain a Radon--Nikodym
formula for analytic deformations. Section~\ref{sec:rn-wp-sle}
extends Fan and Sung's quasi-invariance result to $\kappa=4$ for SLE
welding and extends the analytic density formula to all fixed WP
deformations for $0<\kappa\le4$. The infinitesimal identities needed for our density
formula come from the conformal field theory approach to SLE loop
measure: Baverez and Jego \cite{BJ24} establish the Kontsevich--Suhov
conjecture for $\kappa\in(0,4]$, while Gordina, Qian and Wang
\cite{GQW24} construct a Virasoro representation from infinitesimal
conformal restriction. The finite-mode integration-by-parts formula,
its horizontal $H^{3/2}$ covariance and densities for finite products
of horizontal Fourier flows also appear in
\cite[Theorems~1.2--1.3, Propositions~4.4 and~6.2]{Cai26}.

Our principal contribution is the classification of singular
directions, including simultaneous left and right composition. We give
a proof using logarithmic adjacent-arc observables for $\kappa<4$ and
a local shape argument at $\kappa=4$; the latter also extends to each
fixed subcritical parameter. In the WP directions, we extend the
analytic density to general WP maps and prove moment bounds.

\subsection{The welding law and the Cameron--Martin theorem}\label{introsec:main-results}

We first specify the normalization of the welding law. Let $m$ be Haar
probability measure on $\m S^1$, and write $\Rot_z(w)=zw$ for rotation by
$z\in\m S^1$. Fan and Sung \cite[Definitions~4.4 and~4.7]{F1} work on
the quotient by left rotations,
\[
 \m S^1\backslash\Homeo_+(\m S^1),
\]
identified with the homeomorphisms fixing $1$. The representative of the
class of $h$ is
\[
 N(h)=\Rot_{h(1)^{-1}}\circ h.
\]
We retain the rotation coordinate, so that both composition actions can
be considered on the full space $\Homeo_+(\m S^1)$.

\begin{df}[SLE welding law]\label{introdef:sle-haar-lift}
Fix $\k\in(0,4]$ and sample a loop $\eta$ from the SLE$_\k$ shape
probability measure of \cite[Definition~4.7]{F1}. Thus $\eta$ separates
$0$ from $\infty$ and has conformal radius $1$ viewed from $0$. Let
$f:\m D\to\Omega$ and $g_0:\m D^*\to\Omega^*$ be the Riemann maps
onto the complementary domains, normalized by
\[
 f(0)=0,\qquad f'(0)=1,\qquad g_0(\infty)=\infty,
 \qquad f(1)=g_0(1).
\]
Independently of $\eta$, sample $Z$ with law $m$, and set
\begin{equation}\label{introeq:sle-haar-lift}
 g_Z=g_0\circ\Rot_{Z^{-1}},\qquad
 \widehat H=g_Z^{-1}\circ f=\Rot_Z\circ(g_0^{-1}\circ f).
\end{equation}
We denote the law of $\widehat H$ by
$\SLE_\k^{\mathrm{weld}}=\m P^{\mathrm{weld}}$.
\end{df}

Equivalently, after sampling $Z$, we normalize the exterior map by
$f(1)=g_Z(Z)$. The difference from the fixed-$1$ convention is precisely
this independent uniform rotation. Our law is the unique lift of the
Fan--Sung quotient law invariant under left rotations. The same
randomization is used in \cite[Corollary~4.12]{F1}; it also agrees with
\eqref{eq:law}, since averaging over rotations removes the initial choice
of phase of the exterior map.

For a circle LGF $X$ and $\gamma\in(0,2)$, let
$\widehat{\mc M}_X^\gamma$ be its normalized GMC measure and put
\[
 F_X(z)=\exp\bigl(2\pi\ii\widehat{\mc M}_X^\gamma([1,z])\bigr),
 \qquad z\in\m S^1,
\]
where $[1,z]$ is the counterclockwise arc. We also use $F_X$ for its
degree-one lift fixing $0$ when working in real coordinates.

\begin{lem}[Comparison with GMC welding]\label{introlem:sle-gmc-lift}
Let $\k\in(0,4)$ and $\gamma=\sqrt\k$. For independent circle LGFs
$X_1,X_2$ and an independent uniform rotation $U$, set
\[
 H=F_{X_2}^{-1}\circ U\circ F_{X_1},\qquad
 \m P=\operatorname{Law}(H).
\]
Then $\m P\sim\m P^{\mathrm{weld}}$ as probability measures on
$\Homeo_+(\m S^1)$.
\end{lem}
The Haar-coordinate argument proving this lifting statement is given in
Section~\ref{subsec:subcritical-model}.

\begin{remark}[Why we keep the rotation coordinate]\label{introrem:fixed-point-normalizer}
The passage to simultaneous left and right composition uses
$\mathcal N(\WP(\m S^1))=\WP(\m S^1)$ in $\Homeo_+(\m S^1)$.
For the subgroup $\mathcal K_1=\{\varphi\in\WP(\m S^1):\varphi(1)=1\}$,
the normalizer is strictly larger: it even contains bi-Lipschitz maps
outside WP fixing $1$. Theorem~\ref{thm:fixed-point-normalizer} gives
an explicit example and proof, also for any finite set of fixed points.
Keeping the rotation coordinate allows us to use the full WP group in
the two-sided argument.
\end{remark}

Definition~\ref{introdef:sle-haar-lift} applies throughout
$\kappa\in(0,4]$. The global GMC comparison in
Lemma~\ref{introlem:sle-gmc-lift} is used only for $\kappa<4$.
At $\kappa=4$, we instead compare the laws on proper arcs by
Theorem~\ref{thm:crit-local-bridge}.

\begin{thm}[Cameron--Martin theorem]\label{introthm:main}
Let $\kappa\in(0,4]$, let $\varphi,\psi\in\Homeo_+(\m S^1)$, and set
$T_{\psi,\varphi}(h)=\psi^{-1}\circ h\circ\varphi$.
\begin{enumerate}[label=(\roman*)]
\item If $\varphi,\psi\in\WP(\m S^1)$, then
$(T_{\psi,\varphi})_*\m P^{\mathrm{weld}}$ and
$\m P^{\mathrm{weld}}$ are mutually absolutely continuous.
Their Radon--Nikodym derivative is given $\m P^{\mathrm{weld}}$-almost surely by
\[
 \rho_{\psi,\varphi}(h)
 :=\frac{\dd(T_{\psi,\varphi})_*\m P^{\mathrm{weld}}}
        {\dd\m P^{\mathrm{weld}}}(h)
 =\exp\{-\mathscr J_{\psi,\varphi}(h)\}.
\]
The relative action $\mathscr J_{\psi,\varphi}$ is defined in
\eqref{introeq:two-sided-action} of Section~\ref{introsec:rn}. For some
$\varepsilon>0$ depending on the pair,
$\rho_{\psi,\varphi}\in L^{1+\varepsilon}(\m P^{\mathrm{weld}})$
and its negative moment of order $\varepsilon$ is finite.
\item If at least one of $\varphi,\psi$ lies outside $\WP(\m S^1)$,
then
\[
 (T_{\psi,\varphi})_*\m P^{\mathrm{weld}}
 \perp\m P^{\mathrm{weld}}.
\]
\end{enumerate}
Thus the pairs giving a nonzero common part are exactly
$\WP(\m S^1)\times\WP(\m S^1)$.
\end{thm}

The singularity assertion is the converse to WP quasi-invariance, for
arbitrary simultaneous changes on the two sides. Its proof first
classifies the one-sided non-singular directions, then uses the
normalizer theorem to pass to pairs. The explicit Radon--Nikodym
derivative in the WP case, its geometric interpretation and its moment
bounds are stated in Section~\ref{introsec:rn} and proved in
Section~\ref{sec:rn-wp-sle}.

\begin{figure}[!htbp]
\centering
\begin{tikzpicture}[
  font=\footnotesize,
  part/.style={draw=black!50,rounded corners=2pt,align=center,
                text width=.257\linewidth,inner sep=5pt},
  result/.style={draw=black!55,rounded corners=2pt,align=center,
                  text width=.88\linewidth,inner sep=6pt},
  arr/.style={-{Stealth[length=2mm]},semithick,draw=black!65},
  note/.style={align=center,font=\footnotesize,inner sep=2pt}]
\node[part,fill=black!3] (recovered) at (4.58,0) {
  \textbf{Earlier results recovered}\\
  \cite{F1}: subcritical SLE WP quasi-invariance\\
  \cite{BJ25}: analytic density};
\node[part,fill=NavyBlue!5,anchor=south] (sub) at (-4.58,0 |- recovered.south) {
  \textbf{Section~\ref{sec:subcritical-estimates}}\\
  One-sided quasi-invariance\\
  implies WP for $0<\kappa<4$};
\coordinate[below=22pt of recovered.south] (methods);
\node[part,fill=NavyBlue!5,anchor=north,minimum height=2.0cm]
  (log) at (-4.58,0 |- methods) {
  \textbf{Section~\ref{sec:subcritical-singularity}}\\
  Logarithmic observables\\
  $\kappa<4$};
\node[part,fill=NavyBlue!5,anchor=north,minimum height=2.0cm]
  (shape) at (0,0 |- methods) {
  \textbf{Section~\ref{sec:critical-local}}\\
  Local comparison and shape\\
  Each fixed $0<\kappa\le4$};
\node[part,fill=BrickRed!6,anchor=north,minimum height=2.0cm]
  (rn) at (4.58,0 |- methods) {
  \textbf{Section~\ref{sec:rn-wp-sle}}\\
  WP quasi-invariance\\
  WP densities and moments\\
  $0<\kappa\le4$};
\node[fit=(log)(shape)(rn),inner sep=0pt] (methodrow) {};
\coordinate[below=40pt of methodrow.south] (conclusions);
\node[result,text width=.53\linewidth,anchor=north]
  (both) at (-2.29,0 |- conclusions) {
  \textbf{Section~\ref{sec:normalizer}}\\
  WP normalizer and conjugation give the two-sided dichotomy};
\node[part,fill=BrickRed!6,anchor=north]
  (density) at (4.58,0 |- conclusions) {
  \textbf{Two-sided density}\\
  $\rho_{\psi,\varphi}=e^{-\mathscr J_{\psi,\varphi}}$};
\node[fit=(both)(density),inner sep=0pt] (resultrow) {};
\coordinate[below=18pt of resultrow.south] (summary);
\node[result,anchor=north] (main) at (0,0 |- summary) {
  \textbf{Theorem~\ref{introthm:main}}\\
  Two-sided dichotomy and the Radon--Nikodym derivative for WP pairs};
\draw[arr] (sub.south) -- (log.north);
\coordinate[below=8pt of sub.south] (reuse);
\draw[arr,dashed] (reuse) -| (shape.north);
\draw[arr] (rn.north) -- (recovered.south);
\coordinate[below=19pt of methodrow.south] (wp);
\node[note] (singularity) at (both.north |- wp) {One-sided singularity};
\draw[arr] (log.south) |- (singularity.west);
\draw[arr] (shape.south) |- (singularity.east);
\draw[arr] (singularity.south) -- (both.north);
\draw[arr] (rn.south) -- (density.north);
\draw[arr] (rn.south |- wp) -- node[note,above] {Quasi-invariance}
  (1.37,0 |- wp) -- (1.37,0 |- both.north);
\draw[arr] (both.south) -- (both.south |- main.north);
\draw[arr] (density.south) -- (density.south |- main.north);
\end{tikzpicture}
\caption{Organization of the proof. Section~\ref{sec:subcritical-estimates}
shows that one-sided quasi-invariance implies WP regularity when $0<\kappa<4$.
Sections~\ref{sec:subcritical-singularity} and~\ref{sec:critical-local}
give two routes to one-sided singularity. The dashed arrow indicates
the estimates reused below criticality (Remark~\ref{rem:subcritical-shape-proof}).
Section~\ref{sec:rn-wp-sle} proves WP quasi-invariance and the density
formula for $0<\kappa\le4$. These results recover subcritical SLE welding
quasi-invariance from \cite{F1} and the analytic density from \cite{BJ25}.
Technical estimates are collected in Section~\ref{sec:appendix}.}
\label{fig:intro-paper-structure}
\end{figure}
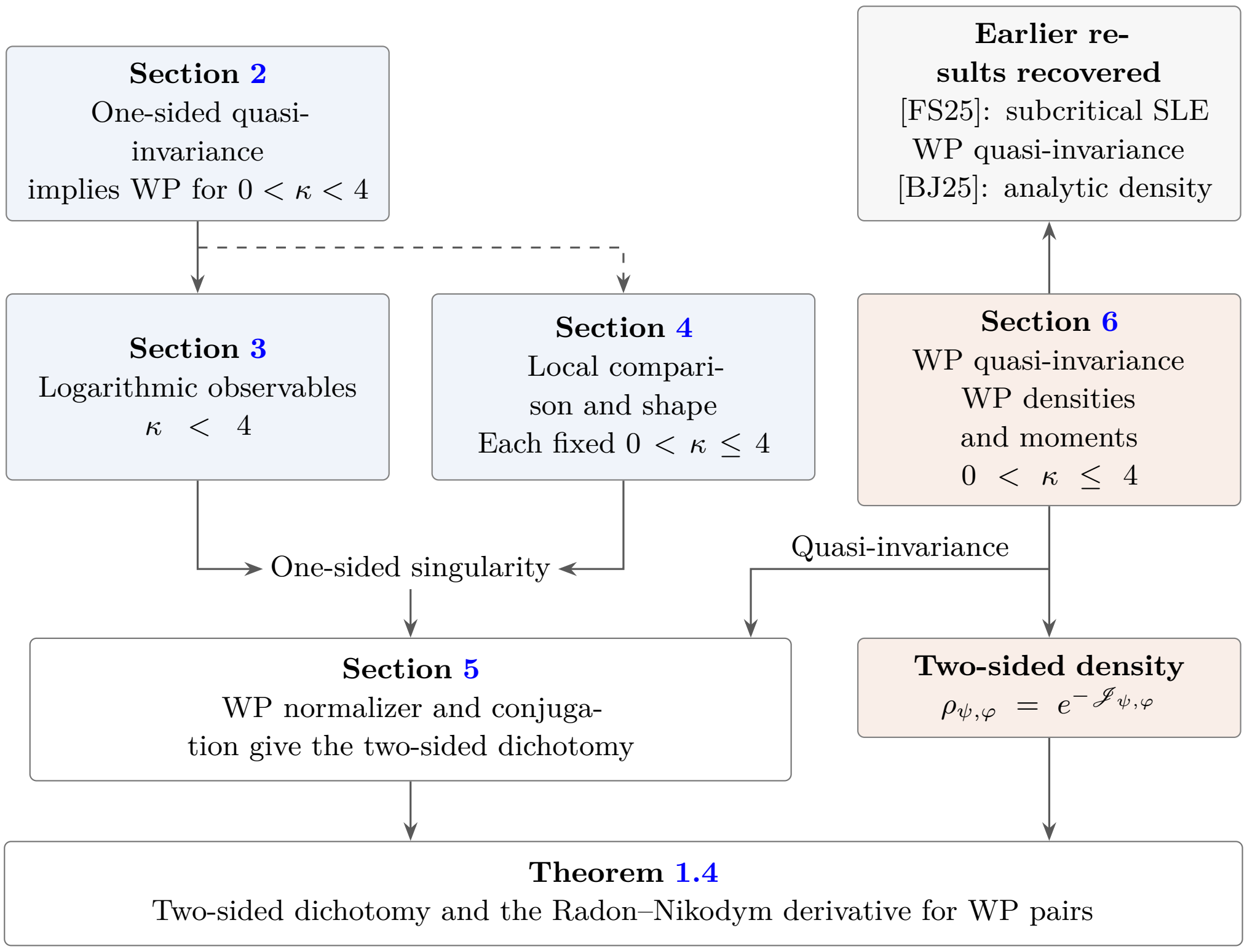

Sections~\ref{sec:subcritical-estimates}--\ref{sec:subcritical-singularity}
give the logarithmic-observable proof for $\kappa<4$.
Section~\ref{sec:critical-local} adapts the strategy to $\kappa=4$
using local shape statistics. This second proof also applies at each
fixed $\kappa<4$ (Remark~\ref{rem:subcritical-shape-proof}).
Section~\ref{sec:normalizer} gives the two-sided
reduction. Figure~\ref{fig:intro-paper-structure} shows how the
proofs fit together.

\subsection{Singularity in the subcritical regime}\label{introsec:subcritical}

Fix $\kappa\in(0,4)$ and $\gamma=\sqrt\kappa$. We use the auxiliary
GMC welding law $\m P$ of Lemma~\ref{introlem:sle-gmc-lift} and write
$\m E$ for its expectation. Identify $\m S^1$ with $\m R/\m Z$ via
$x\mapsto e^{2\pi\ii x}$. For $h\in\Homeo_+(\m S^1)$, choose a
degree-one lift $H$ and let $\nu_h([a,b])=H(b)-H(a)$. Define
\[
 Z_{x,r}(h)=\log\frac{\nu_h([x,x+r])}{\nu_h([x-r,x])},
 \qquad x\in\m R,\quad 0<r\le1/2.
\]
For a deterministic map $\varphi$ with lift $\Phi$, write
$\delta_\varphi(x,r)=Z_{x,r}(\varphi)$. The size of this one
observable distinguishes the three regularity classes
\[
 \WP(\m S^1)\subsetneqq\Sy(\m S^1)
 \subsetneqq\QS(\m S^1)\subsetneqq\Homeo_+(\m S^1),
\]
where $\Sy(\m S^1)$ denotes the symmetric homeomorphisms and
$\QS(\m S^1)$ the quasisymmetric homeomorphisms.

\begin{lem}\label{introlem:discrete WP}
For $\varphi\in\Homeo_+(\m S^1)$:
\begin{enumerate}[label=(\roman*)]
\item $\varphi\in\QS(\m S^1)$ if and only if
$\sup_{x,\,0<r\le1/2}|\delta_\varphi(x,r)|<\infty$.
\item $\varphi\in\Sy(\m S^1)$ if and only if
$\sup_x|\delta_\varphi(x,r)|\to0$ as $r\downarrow0$.
\item $\varphi\in\WP(\m S^1)$ if and only if
$\varphi\in\Sy(\m S^1)$ and
\[
 E(\varphi)=\int_{\m S^1\times(0,1/2]}
 \delta_\varphi(x,r)^2\,\dd\eta(x,r)<\infty,
 \qquad \dd\eta=\frac{\dd x\,\dd r}{r^2}.
\]
\end{enumerate}
\end{lem}
The first two characterizations and the sufficiency in (iii) are proved
in Lemma~\ref{lem:discrete WP}; the necessity in (iii) is
Theorem~\ref{thm:wp-finite-increment-energy}. We equip the circle
homeomorphism group with its uniform group metric
\[
 d_H(h_1,h_2)=\sup_{z\in\m S^1}d(h_1(z),h_2(z))
 +\sup_{z\in\m S^1}d(h_1^{-1}(z),h_2^{-1}(z)),
\]
where $d$ is circle distance.

The three ways to fail these conditions lead to three choices of
statistic. Outside $\QS(\m S^1)$, choose shrinking adjacent arcs with
$d_n:=|\delta_\varphi(x_n,r_n)|\to\infty$ and divide $Z_{x_n,r_n}$
by $d_n$. Lemma~\ref{lem:moment of Z} shows that the normalized
statistic tends to zero under the original law and stays bounded in
$L^2$ under the transformed law. In
$\QS(\m S^1)\setminus\Sy(\m S^1)$, the distortion stays bounded
away from zero along a shrinking sequence. Here
Lemma~\ref{lem:subcritical-weak-cancellation} gives weak convergence
to zero without normalization. For
$\varphi\in\Sy(\m S^1)\setminus\WP(\m S^1)$, every individual
distortion tends to zero, so we accumulate observations over positions
and scales.

\begin{figure}[!htbp]
\centering
\begin{tikzpicture}[
  font=\footnotesize,
  csplit/.style={draw=NavyBlue},
  msplit/.style={draw=BrickRed},
  dot/.style={circle,inner sep=1.55pt,outer sep=0pt},
  arr/.style={-{Stealth[length=2mm]},semithick}]
\node[font=\footnotesize\bfseries] at (3.55,3.0) {(a) Source and image arcs};
\begin{scope}[shift={(1.48,1.22)}]
  \draw[black!45] (0,0) circle (1.02);
  \draw[semithick,-{Stealth[length=1.35mm]}]
    (24:1.02) arc[start angle=24,end angle=156,radius=1.02];
  \foreach \angle in {24,90,156}
    \fill (\angle:1.02) circle (1.45pt);
  \node[right] at (24:1.05) {$x-r$};
  \node[above] at (90:1.04) {$x$};
  \node[left] at (156:1.05) {$x+r$};
  \node at (0,0) {$\m S^1$};
  \node[below] at (0,-1.10) {Equal halves of $I$};
\end{scope}
\draw[arr] (2.85,1.22) -- node[above] {$\varphi$} (4.31,1.22);
\begin{scope}[shift={(5.64,1.22)}]
  \draw[black!45] (0,0) circle (1.02);
  \draw[semithick,-{Stealth[length=1.35mm]}]
    (24:1.02) arc[start angle=24,end angle=156,radius=1.02];
  \foreach \angle in {24,156}
    \fill (\angle:1.02) circle (1.45pt);
  \node[dot,fill=NavyBlue] at (63:1.02) {};
  \node[dot,msplit,fill=white] at (90:1.02) {};
  \node[right] at (24:1.06) {$a$};
  \node[above right,text=NavyBlue] at (63:1.04) {$c$};
  \node[above,text=BrickRed] at (90:1.04) {$M$};
  \node[left] at (156:1.06) {$b$};
  \node at (0,0) {$\m S^1$};
  \node[below] at (0,-1.10) {Same endpoints for both splits};
\end{scope}
\node[font=\footnotesize\bfseries] at (3.55,-.52) {(b) Two mass ratios on the same image arc};
\draw[black!30,densely dotted] (.14,-1.02) -- (.14,-2.29)
                             (4.45,-1.02) -- (4.45,-2.29);
\draw[csplit,line width=1pt,opacity=.47] (.14,-1.23) -- (1.414,-1.23);
\draw[csplit,line width=1.7pt] (1.414,-1.23) -- (4.45,-1.23);
\draw[msplit,line width=1pt,opacity=.47] (.14,-2.14) -- (2.295,-2.14);
\draw[msplit,line width=1.7pt] (2.295,-2.14) -- (4.45,-2.14);
\foreach \yy in {-1.23,-2.14} {
  \draw (.14,\yy-.065) -- (.14,\yy+.065);
  \draw (4.45,\yy-.065) -- (4.45,\yy+.065);
  \node[above=2pt] at (.14,\yy) {$a$};
  \node[above=2pt] at (4.45,\yy) {$b$};
}
\node[dot,fill=NavyBlue] at (1.414,-1.23) {};
\node[dot,msplit,fill=white] at (2.295,-2.14) {};
\node[above=2pt,text=NavyBlue] at (1.414,-1.23) {$c$};
\node[above=2pt,text=BrickRed] at (2.295,-2.14) {$M$};
\node[right,text=NavyBlue] at (4.68,-1.23)
  {$\displaystyle\log\frac{\nu_h([c,b])}{\nu_h([a,c])}$};
\node[right,text=BrickRed] at (4.68,-2.14)
  {$\displaystyle\log\frac{\nu_h([M,b])}{\nu_h([a,M])}$};
\node[font=\footnotesize\bfseries] at (10.72,3.0) {(c) Normalized shape at $\kappa=4$};
\begin{scope}[shift={(9.12,-1.5)},x=3.13cm,y=3.67cm]
  \draw[black!23] (0,0) rectangle (1,1);
  \draw[arr] (0,0) -- (1.10,0) node[right] {$u$};
  \draw[arr] (0,0) -- (0,1.09);
  \draw[black!50,dashed] (0,0) -- (1,1);
  \node[above left,black!65] at (.8,.8) {$\operatorname{id}$};
  \draw[msplit,densely dotted] (0,.5) -- (.5,.5);
  \draw[csplit,densely dotted] (0,.29545) -- (.5,.29545) -- (.5,0);
  \draw[csplit,line width=.85pt]
    (0,0) .. controls (.15,.012) and (.32,.13) .. (.5,.29545)
          .. controls (.68,.46) and (.88,.74) .. (1,1);
  \node[right,text=NavyBlue] at (.80,.62) {$k_I$};
  \node[dot,msplit,fill=white] at (.5,.5) {};
  \node[dot,fill=NavyBlue] at (.5,.29545) {};
  \node[below left] at (0,0) {$0$};
  \node[below] at (.5,0) {$1/2$};
  \node[below] at (1,0) {$1$};
  \node[left] at (0,1) {$1$};
  \node[left,text=BrickRed] at (0,.5) {$1/2$};
  \node[left,text=NavyBlue] at (0,.29545) {$k_I(1/2)$};
\end{scope}
\node[align=center] at (10.72,-2.16)
  {$\displaystyle k_I(1/2)=\frac{c-a}{b-a}$};
\node[draw=black!45,rounded corners=2pt,align=center,inner sep=6pt]
  at (6.30,-3.09) {
   Multiply by $\delta_\varphi(x,r)$ and integrate over $\Omega_N$\\[3pt]
   $T_N(h\circ\varphi)=B_N(h)+A_N(h),\qquad A_N(h)\ge0$};
\end{tikzpicture}
\caption{Two splits of the same image arc. In increasing lifts,
$(a,c,b)=(\Phi(x-r),\Phi(x),\Phi(x+r))$ and $M=(a+b)/2$.
The blue ratio evaluates the split at $c$ and the red ratio evaluates
the split at $M$. Their difference contributes $A_N$ after multiplication
by $\delta_\varphi$ and integration. The critical observable uses
normalized restrictions. Rescaling the two arcs to $[0,1]$ sends
$x$ to $1/2$, $c$ to $k_I(1/2)$, and $M$ to $1/2$ in the image.}
\label{fig:intro-splittings}
\end{figure}
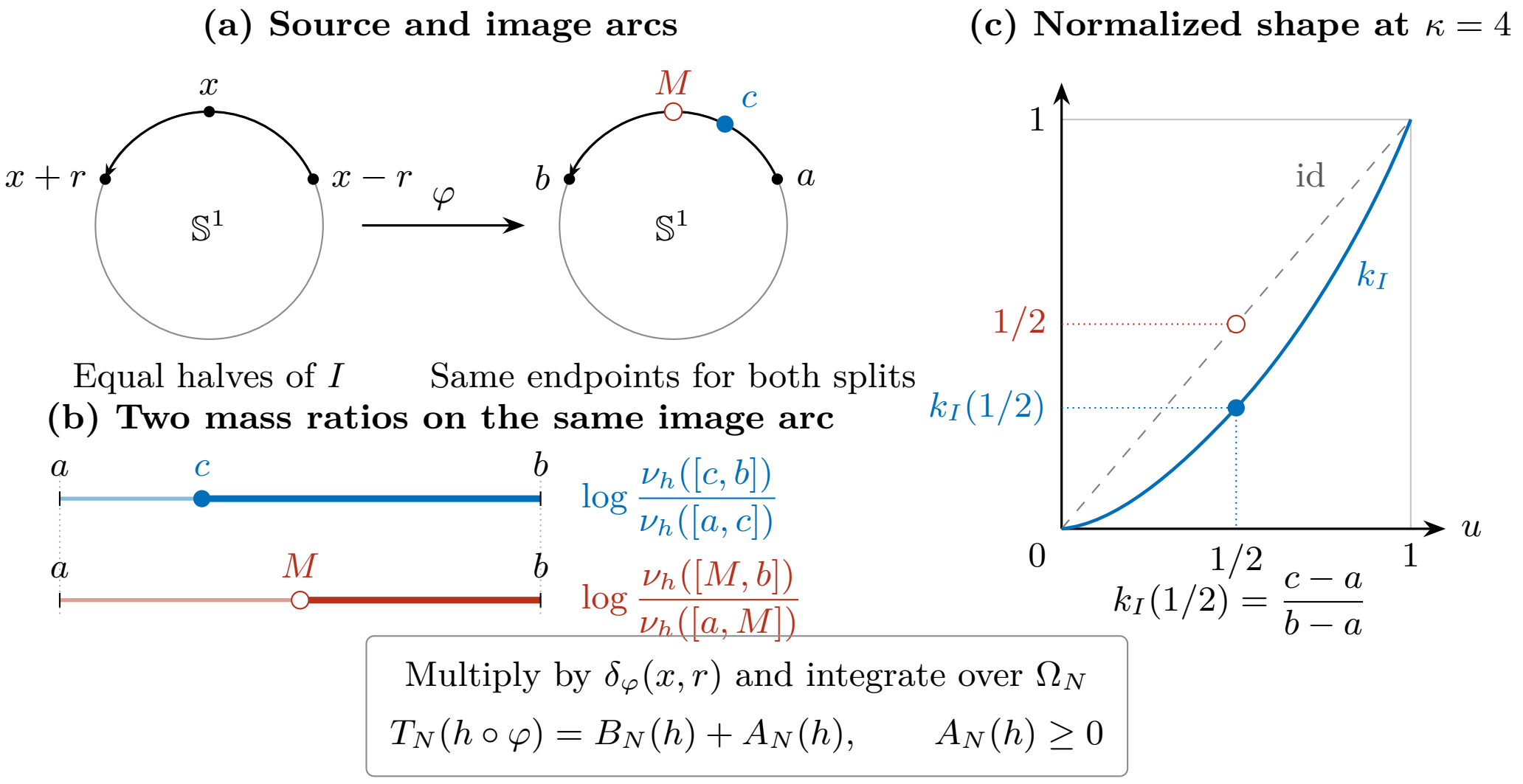

All three cases use the same decomposition. Put
\[
 M_r(x)=\frac{\Phi(x+r)+\Phi(x-r)}2,\qquad
 L_r(x)=\frac{\Phi(x+r)-\Phi(x-r)}2,
\]
and let
$R_\varphi(x,r;h)=|Z_{x,r}(h\circ\varphi)-Z_{M_r(x),L_r(x)}(h)|$.
Monotonicity of the two image masses as the splitting point moves gives
\begin{equation}\label{introeq:sing-decomposition}
 Z_{x,r}(h\circ\varphi)
 =Z_{M_r(x),L_r(x)}(h)
 +\operatorname{sgn}(\delta_\varphi(x,r))R_\varphi(x,r;h).
\end{equation}
The first term uses the geometric midpoint of the image interval.
The second records the change from that midpoint to $\Phi(x)$.
Figure~\ref{fig:intro-splittings} shows the two splits on the same
image arc.

For the multiscale case, set
\[
 \Omega_N=\m S^1\times[2^{-N-1},1/2),\qquad
 E_N=1+\int_{\Omega_N}\delta_\varphi^2\,\dd\eta,
\]
and define
\[
 \begin{aligned}
 T_N(h)&=\int_{\Omega_N}\delta_\varphi Z_{x,r}(h)\,\dd\eta,\\
 B_N(h)&=\int_{\Omega_N}\delta_\varphi
                       Z_{M_r(x),L_r(x)}(h)\,\dd\eta,\\
 A_N(h)&=\int_{\Omega_N}|\delta_\varphi|R_\varphi(x,r;h)\,\dd\eta.
 \end{aligned}
\]
Then $E_N\to\infty$ and $T_N(h\circ\varphi)=B_N(h)+A_N(h)$,
with $A_N\ge0$. Lemma~\ref{lem:Bessel} gives
\[
 \left\|\frac{T_N}{E_N}\right\|_2+
 \left\|\frac{B_N}{E_N}\right\|_2
 \le C_{\gamma,\varphi}E_N^{-1/2}\longrightarrow0.
\]
Lemma~\ref{lem:checkpoint} supplies a deterministic subsequence on
which $A_N/E_N$ is bounded in $L^2$; the decomposition then gives the
required bound under the transformed law.

Testing against a possible common part requires a lower bound weighted
by its density. Lemma~\ref{lem:weighted-lower-bound} averages over the
independent rotation with the phase weight still present. This bounds
the welding displacement below by a source displacement minus an error
whose contribution vanishes after normalization
(Lemma~\ref{lem:nonwp-sym-integrable}). Source positivity and deletion
of finitely many Fourier modes give the weighted estimate in
Proposition~\ref{prop:nonwp-weighted}. In the cumulative case,
\[
 \liminf_k\m E\!\left[W\frac{A_{N_k}}{E_{N_k}}\right]
 \ge c\m E[W]\qquad\text{for every bounded measurable }W\ge0.
\]
The source $L^2$ bound on the same subsequence follows by taking
$\Psi\equiv1$ in the rotation lower bound and applying conditional
Jensen's inequality, as in \eqref{eq:nonwp-cumulative-input}.

The measure-theoretic step is the following criterion, proved in
Lemma~\ref{lem:sing-criterion}.

\begin{lem}\label{introlem:sing-criterion}
Let $P$ be a probability measure on a standard Borel space $\mathcal X$,
let $T$ be a Borel bijection with Borel inverse, and set $Q=T_*P$.
Suppose that measurable $D_n$ satisfy
\[
 D_n\rightharpoonup0\quad\text{in }L^2(P),\qquad
 \sup_n\|D_n\|_{L^2(Q)}<\infty,
\]
and that $D_n\circ T=N_n+Y_n$, where
$N_n\rightharpoonup0$ in $L^2(P)$, $Y_n\ge0$ and, for some $c>0$,
\[
 \liminf_n\int wY_n\,\dd P\ge c\int w\,\dd P
 \qquad\text{for every bounded measurable }w\ge0.
\]
Then $Q\perp P$.
\end{lem}
Indeed, a nonzero common part would pull back under $T$ to a bounded
nonnegative density against which the positive term cannot vanish.
The uniform $L^2$ bound under $Q$ makes the original statistic vanish
against that same common part.

Write $\mu\not\perp\nu$ when two measures have a nonzero common part,
and define
\[
 \tilde\QIR=\{\varphi\in\Homeo_+(\m S^1):
 (h\mapsto h\circ\varphi^{-1})_*\m P^{\mathrm{weld}}
 \not\perp\m P^{\mathrm{weld}}\}.
\]
The regularity classes are closed under inversion, so the argument for
$h\circ\varphi$ also applies to this convention. The three steps are
assembled in Theorem~\ref{thm:sing-conditional} and
Corollary~\ref{cor:one-sided-dichotomy}.
\begin{thm}\label{introthm:one-sided singularity}
For $\kappa\in(0,4)$:
\begin{enumerate}[label=(\Alph*)]
\item If $\varphi\in\tilde\QIR$, then $\varphi\in\QS(\m S^1)$.
\item If $\varphi\in\tilde\QIR\cap\QS(\m S^1)$, then
$\varphi\in\Sy(\m S^1)$.
\item If $\varphi\in\tilde\QIR\cap\Sy(\m S^1)$, then
$\varphi\in\WP(\m S^1)$.
\end{enumerate}
Consequently, $\tilde\QIR=\WP(\m S^1)$.
\end{thm}
Remark~\ref{rem:subcritical-shape-proof} gives a second proof using
local shape statistics.

\FloatBarrier
\subsection{Adapting the singularity argument to the critical parameter}\label{introsec:critical}

The proof at $\kappa=4$ adapts the strategy of
Section~\ref{introsec:subcritical}: choose local observations, control
their centered fluctuations, and obtain a positive weighted mean after
composition to exclude a common part. The changes concern the
comparison of welding laws and the choice of observations. This local
construction works for every fixed $0<\kappa\le4$; its critical form
requires the additional estimates of Section~\ref{sec:critical-local}.

Let $F_1,F_2$ be the cumulative functions of two independent normalized
critical GMC measures, and let $U$ be independent and uniform in the
mass coordinate. Put
\[
 H_4=F_2^{-1}\circ R_U\circ F_1,\qquad
 \m P_4=\operatorname{Law}(H_4).
\]
The comparison with the SLE welding law is now on proper arcs:
\begin{equation}\label{introeq:critical-local-comparison}
 (\operatorname{Res}_I)_*\m P_4^{\mathrm{weld}}
 \sim(\operatorname{Res}_I)_*\m P_4
 \qquad\text{for every closed proper arc }I,
\end{equation}
where $\operatorname{Res}_I(h)=h|_I$ and $\sim$ means mutual absolute
continuity. This is Theorem~\ref{thm:crit-local-bridge}, proved by
critical wedge welding in
Proposition~\ref{prop:crit-geometric-branch-proved}.
Proposition~\ref{prop:crit-local-sing} then transfers a separating
event on one arc from the model to the SLE law, since
$(h\circ\varphi)|_I$ depends only on $h|_{\varphi(I)}$.

For $I=[a,a+r]$, the new observations use the normalized restrictions
\[
 h_I(u)=\frac{h(a+ru)-h(a)}{h(a+r)-h(a)},\qquad
 k_I=\varphi_I,\qquad 0\le u\le1,
\]
with differences computed in increasing lifts. Panel (c) of
Fig.~\ref{fig:intro-splittings} shows how $k_I$ records the shape of
$\varphi$ after its source and image arcs have been rescaled to
$[0,1]$. The bounded statistic
$F_I(h)$ and the distortion coefficient $w_I$ of
Section~\ref{subsec:crit-shape} replace $Z_{x,r}(h)$ and
$\delta_\varphi(x,r)$. Critical mass estimates and smooth source
deformations supply the single-scale steps
(Theorem~\ref{thm:crit-nonqs-total} and
Lemma~\ref{lem:crit-common-sym}). At accumulated scales,
Proposition~\ref{prop:crit-shape-layers} relates the energy of $w_I$
to $\|k_I-\operatorname{id}\|_\infty^2$, so it detects the same
failure of WP regularity.

For $\varphi\in\Sy(\m S^1)\setminus\WP(\m S^1)$, the corresponding
statistics are $\mathcal T_N^\beta,\mathcal B_N^\beta,\mathcal A_N^\beta$
of Lemma~\ref{lem:crit-record}. They already include the energy
normalization, and the deterministic cutoff $\beta$ localizes them to
a proper arc. They satisfy
\[
 \mathcal T_N^\beta(h\circ\varphi)
 =\mathcal B_N^\beta(h)+\mathcal A_N^\beta(h).
\]
The covariance estimate in Proposition~\ref{prop:crit-two-frames}
and Lemma~\ref{lem:crit-record} make $\mathcal T_N^\beta$ and
$\mathcal B_N^\beta$ tend to zero in $L^2(\m P_4)$ along suitable scales.
The signed remainder $\mathcal A_N^\beta$ has the positive weighted
lower bound of Proposition~\ref{prop:crit-weighted-cutoff}; its proof
uses the cancellation of the extra length term in
Lemma~\ref{lem:crit-square-cancellation}. The same test against a
possible common part then gives singularity.

\begin{thm}[Critical singularity]\label{introthm:critical-singularity}
For $\varphi\in\Homeo_+(\m S^1)\setminus\WP(\m S^1)$, both
\[
 (h\mapsto h\circ\varphi)_*\m P_4^{\mathrm{weld}}
 \perp\m P_4^{\mathrm{weld}},\qquad
 (h\mapsto\varphi\circ h)_*\m P_4^{\mathrm{weld}}
 \perp\m P_4^{\mathrm{weld}}.
\]
\end{thm}
This is Theorem~\ref{thm:crit-all-nonwp-total}.
Remark~\ref{rem:subcritical-shape-proof} explains the extension to each
fixed $\kappa<4$: the mass estimates of
Section~\ref{subsec:subcritical-gmc} and the subcritical GMC coordinate
rule replace their critical counterparts. We retain the independent
logarithmic proof because it follows the adjacent-arc observable
directly. No uniformity as $\kappa\uparrow4$ is asserted.

\FloatBarrier
\subsection{From one-sided to two-sided composition}\label{introsec:two-sided}

The full one-sided dichotomy is the input to the two-sided argument:
WP maps preserve the measure class, and every map outside
$\WP(\m S^1)$ gives a singular law. In addition to $\tilde\QIR$ above,
define
\[
 \begin{aligned}
 \tilde\QIL&=\{\varphi:(h\mapsto\varphi\circ h)_*
                  \m P^{\mathrm{weld}}\not\perp\m P^{\mathrm{weld}}\},\\
 \tilde\QI&=\{(\varphi,\psi):(h\mapsto\varphi\circ h\circ\psi^{-1})_*
                  \m P^{\mathrm{weld}}\not\perp\m P^{\mathrm{weld}}\}.
 \end{aligned}
\]
All maps here lie in $\Homeo_+(\m S^1)$. Inversion invariance gives
the left-sided statement from the right-sided one
(Lemma~\ref{lem:inversion}). Writing $\QIR$ for the right
quasi-invariance group, we obtain
\[
 \tilde\QIL=\tilde\QIR=\QIR=\WP(\m S^1).
\]
The deterministic ingredient is the normalizer theorem,
Theorem~\ref{thm:main}:
\[
 G:=\mathcal N_{\Homeo_+(\m S^1)}(\WP(\m S^1))=\WP(\m S^1).
\]
Its proof first establishes quasisymmetry of a normalizing map
(Theorem~\ref{thm:qs regularity}), then WP regularity through the
composition operator (Theorem~\ref{thm:from qs to wp}).
Theorem~\ref{thm:two-sided-singularity} uses the one-sided dichotomy to
obtain
\[
 \tilde\QI\subset\mc N(\QIR)\times\mc N(\QIR)
 =\WP(\m S^1)\times\WP(\m S^1).
\]
The first inclusion uses $\tilde\QIR=\QIR$ and connectedness in the
WP topology. The equality follows from the normalizer theorem
since $\QIR=\WP(\m S^1)$.

\subsection{Radon--Nikodym derivatives for WP deformations}\label{introsec:rn}

Fix $\kappa\in(0,4]$ and let
$c_{\mathrm L}=25+6(2/\sqrt\kappa-\sqrt\kappa/2)^2$.
Expectations $\m E^{\mathrm{weld}}$ and spaces $L^p$ of functions
of $h$ refer to $\m P^{\mathrm{weld}}$.
For a normalized welding pair $h=g_h^{-1}f_h$, put
\begin{equation}\label{introeq:conformal-radius}
 K(h)=\log|g_h'(\infty)|-\log|f_h'(0)|.
\end{equation}
For $h,\varphi\in\WP(\m S^1)$, our relative action satisfies
\begin{equation}\label{introeq:one-energy}
 \begin{aligned}
 \mathscr I(\varphi,h)
 &=\frac{c_{\mathrm L}}{24}
   [I^{\mathrm L}(h\varphi^{-1})-I^{\mathrm L}(h)]\\
 &\quad+2[K(h\varphi^{-1})-K(h)],
 \end{aligned}
\end{equation}
where $I^{\mathrm L}$ is Loewner energy
(Proposition~\ref{prop:rn-finite-energy}). For SLE weldings, a
limit construction must control both the relative action and its
exponential.

Write $G_1$ for the inner Grunsky operator of \eqref{eq:TT}.
For analytic $\varphi$, the action is
\begin{equation}\label{introeq:explicit-action}
 \begin{aligned}
 \mathscr I(\varphi,h)
 ={}&-\frac{c_{\mathrm L}}2
       \log\det(I-G_1(\varphi)^*G_1(\varphi))\\
 &+c_{\mathrm L}\operatorname{Re}\log\det
       (I-G_1(\varphi)\overline{G_1(h)})\\
 &+2[K(h\varphi^{-1})-K(h)].
 \end{aligned}
\end{equation}
The analytic formula follows from Lemma~\ref{lem:analytic-cocycle}.
For $\varphi\in\WP(\m S^1)$, $G_1(\varphi)$ is Hilbert--Schmidt,
whereas the Grunsky inequality only gives a bound in operator norm
for the SLE welding. The mixed product can therefore be outside
trace class. We use $\det_2$ for this product and combine the
linear part of the removed trace with the first variation of $K$.
For a Beltrami representative $\mu$ as in \eqref{eq:mu}, this gives
\[
 \begin{gathered}
 \mathscr I_\mu=C_\mu+\mathfrak d_{v_\mu}+E_\mu,\\
 \|E_\mu\|_\infty\le C_{\kappa,r_0}\frac{N(\mu)^2}{(1-k)^3},
 \qquad k=\|\mu\|_\infty<1.
 \end{gathered}
\]
The constant $C_\mu$ is deterministic. The combined linear term
extends continuously from $H^{3/2}$ to an $L^2$ random variable,
and the resolvent estimates bound the quadratic remainder uniformly
in the welding. See \eqref{body:eq:action} and
Lemmas~\ref{lem:divergence-estimate} and~\ref{lem:bounds}.

Radial cutoffs $\mu_j$ preserve $\|\mu_j\|_\infty\le k<1$ and
give $N(\mu_j-\mu)\to0$ in the hyperbolic $L^2$ norm of \eqref{eq:mu}.
Their analytic boundary maps $\varphi_j$ tend uniformly to $\varphi$.
Continuity of the regularized terms gives
\begin{equation}\label{introeq:action-limit}
 \mathscr I(\varphi_j,\cdot)\longrightarrow\mathscr I(\varphi,\cdot)
 \quad\text{in }L^2(\m P^{\mathrm{weld}}).
\end{equation}
To pass from actions to densities we need uniform integrability.
For real $r$ with $|r|k<1$, $r\mu_j$ defines another analytic
deformation and $v_{r\mu_j}=rv_{\mu_j}$. Thus the linear terms cancel
in $\mathscr I_{r\mu_j}-r\mathscr I_{\mu_j}$ and the remaining
difference is bounded uniformly in $j$. The normalization of the
analytic density gives
\[
 \begin{aligned}
 \m E^{\mathrm{weld}}\!\left[e^{-\mathscr I_{r\mu_j}}\right]&=1,\\
 \sup_j\m E^{\mathrm{weld}}\!\left[e^{-r\mathscr I_{\mu_j}}\right]
 &<\infty \qquad (|r|k<1).
 \end{aligned}
\]
With $r>1$, this bound upgrades density convergence to $L^1$.
Together with uniform convergence of $\varphi_j$, it identifies the
limit as the density of the WP pushforward (Lemma~\ref{lem:approx}
and Theorem~\ref{body:thm:one}).

\begin{thm}\label{thm:one}
For fixed $\varphi\in\WP(\m S^1)$, the action
\eqref{introeq:action-limit} is finite and independent of the
representative, and
\begin{equation}\label{eq:one}
 \rho_\varphi(h):=
 \frac{\dd(R_\varphi)_*\m P^{\mathrm{weld}}}{\dd\m P^{\mathrm{weld}}}(h)
 =e^{-\mathscr I(\varphi,h)}>0
 \quad\m P^{\mathrm{weld}}\text{-a.s.},
\end{equation}
where $R_\varphi(h)=h\circ\varphi$. For $k>0$ as above,
\begin{equation}\label{eq:one-moments}
 \m E^{\mathrm{weld}}[\rho_\varphi^r]<\infty\qquad(|r|k<1).
\end{equation}
Moreover, $\rho_\varphi\in L^{2+\varepsilon}$ for some
$\varepsilon>0$, and $\mathscr I(\varphi,\cdot)\in L^p$ for every
finite $p$. For a rotation, $\rho_\varphi=1$.
\end{thm}
The inverse-density identity gives the $L^{2+\varepsilon}$ bound
in Theorem~\ref{body:thm:one}.

For the simultaneous action, set
\begin{equation}\label{eq:two-action}
 T_{\psi,\varphi}(h)=\psi^{-1}\circ h\circ\varphi,
\end{equation}
\begin{equation}\label{introeq:two-sided-action}
 \mathscr J_{\psi,\varphi}(h)
 =\mathscr I(\varphi,h)+\mathscr I(\psi,\varphi h^{-1}).
\end{equation}
\begin{thm}\label{thm:two}
For fixed $\psi,\varphi\in\WP(\m S^1)$,
\begin{align}
 \rho_{\psi,\varphi}(h)
 &:=\frac{\dd(T_{\psi,\varphi})_*\m P^{\mathrm{weld}}}
              {\dd\m P^{\mathrm{weld}}}(h),
              \label{introeq:main-density}\\
 &=\rho_\varphi(h)\rho_\psi(\varphi h^{-1})
  =e^{-\mathscr J_{\psi,\varphi}(h)}>0
  \quad\m P^{\mathrm{weld}}\text{-a.s.}
  \label{eq:two}
\end{align}
The action is well defined modulo null sets, and, for some
$\varepsilon>0$ depending on the pair,
\begin{equation}\label{introeq:main-moments}
 \m E^{\mathrm{weld}}[\rho_{\psi,\varphi}^{\,1+\varepsilon}
                     +\rho_{\psi,\varphi}^{-\varepsilon}]<\infty.
\end{equation}
\end{thm}
This follows by inversion and the chain rule
(Theorem~\ref{body:thm:two}). For deterministic
$h,\varphi,\psi\in\WP(\m S^1)$, \eqref{introeq:one-energy} gives
\begin{equation}\label{introeq:main-energy}
 \begin{aligned}
 \mathscr J_{\psi,\varphi}(h)
 ={}&\frac{c_{\mathrm L}}{24}
 [I^{\mathrm L}(\psi h\varphi^{-1})-I^{\mathrm L}(h)]\\
 &+2[K(\psi h\varphi^{-1})-K(h)].
 \end{aligned}
\end{equation}
The comparison point is $T_{\psi,\varphi}^{-1}(h)=\psi h\varphi^{-1}$.

\FloatBarrier
\textbf{Acknowledgements.} 

The author thanks Hong-Bin Chen for discussion on Gaussian measures, Jinwoo Sung for helping him understand GMC measures in our previous work \cite{F1}, Dapeng Zhan for discussion of SLE loop measure in Venice, Fredrik Viklund for bringing up the $\log\det_2$ normalization in \cite{F3} that is helpful for this work, Hao Wu for discussion and guidance at Tsinghua University and Mingchang Liu for a conversation about an early version. The author thanks Yilin Wang for her helpful discussions and generous guidance at IHES, especially the key observation of the relation between Weil--Petersson class and SLE welding measures.

\textbf{Research statement}

After \cite{F1}, the author came to view the Weil--Petersson normalizer problem as a deterministic counterpart to quasi-invariance. The author discussed this connection widely and had only a few exchanges about a second observation, that the normalizer could reduce the two-sided problem to the one-sided case. Another approach finally led to exploring the Grunsky operator of a welding, which was one motivation for \cite{F3}. The author completed a short proof of the normalizer theorem in spring 2026 and kept it in his notes while pursuing the full Cameron--Martin theorem. By early September, the author had a preliminary blow-up argument for the welding law.

Before using AI on this project, the author completed the proofs in Sections~\ref{sec:subcritical-estimates} and~\ref{sec:normalizer} and drafted the local approach in Section~\ref{sec:critical-local}. Drawing on \cite{F1,GQW24,BJ24,BJ25,Cai26,F3}, the author also conjectured the density formula in Section~\ref{sec:rn-wp-sle} and outlined a route to prove it for general Weil--Petersson maps.

\textbf{AI statement}

The author began using ChatGPT Plus on 23 October 2026 through an offer that made the first month free. Before then the author had used ChatGPT in brief conversations. The author initially gave Sol individual sections to revise separately. During the development of the proofs, the author used Astra under his direction to explore extensions and work through technical arguments.

For the subcritical singularity theorem (Theorem~\ref{thm:sing-conditional}), the author developed the midpoint decomposition, the multiscale estimates, and the weighted argument for ruling out a common part of the two laws, including the rotation lower bound. Astra provided the decorrelation estimate obtained by removing finitely many Fourier modes (Lemma~\ref{lem:nonwp-mode-deletion}).

After completing all the proofs, the author gave ChatGPT the complete \LaTeX{} source to revise the exposition and organization of the paper as a whole. Sol drafted alternative passages, especially for the introduction, and the author selected and revised them using the exposition written by his coauthors in \cite{F3} as a model. When GPT-6.2 Astra became available, the author decided to finish the revision promptly. The author reviewed the manuscript and takes responsibility for its contents.

\section{Subcritical welding estimates and one-sided regularity}\label{sec:subcritical-estimates}

We prove that one-sided quasi-invariance forces WP regularity for
$\k\in(0,4)$. The argument proceeds in three steps: quasisymmetry,
symmetry, and finally finite multiscale distortion energy. We first
introduce a common adjacent-arc observable and the auxiliary GMC welding
law. Local moment estimates, weak cancellation, and cumulative bounds
then supply the respective steps. Section~\ref{sec:subcritical-singularity}
uses these estimates in a weighted argument to exclude a nonzero common
part of the original and transformed laws.

\subsection{The one-sided problem and the auxiliary law}\label{subsec:subcritical-model}

We express the three regularity conditions through the same logarithmic
ratio of adjacent image arcs. This identifies the deterministic
quantities that the probabilistic argument must control. We then
introduce the GMC representation in which the required estimates will
be proved.

Fix $\k\in(0,4)$ throughout this section. We use the welding law of
Definition~\ref{introdef:sle-haar-lift} and write $P\sim Q$ for mutual
absolute continuity. The right quasi-invariance class is
\[
 \QIR=\{\varphi\in\Homeo_+(\m S^1):
 \operatorname{Law}(\widehat H\circ\varphi^{-1})
 \sim\m P^{\mathrm{weld}}\}.
\]
We first establish the analytic implications below. The two-sided
consequences are proved in Section~\ref{subsec:two-sided-singularity},
after the one-sided singularity argument and the normalizer theorem.
\begin{thm}
    For $\varphi \in \Homeo_+(\m S^1)$, the following implications hold:
    \begin{itemize}
        \item[(A)] If $\varphi \in \QIR$, then $\varphi \in \QS(\m S^1)$;
        \item[(B)] If $\varphi \in \QIR \cap \QS(\m S^1)$, then $\varphi \in \Sy(\m S^1)$;
        \item[(C)] If $\varphi \in \QIR \cap \Sy(\m S^1)$, then $\varphi \in \WP(\m S^1)$.
    \end{itemize}
    Together with WP quasi-invariance, these implications give $\QIR = \WP(\m S^1)$.
\end{thm}
Let $\pi:\m R \rightarrow \m S^1=\m R/ \m Z$ be the quotient map
$\pi(x) = e^{2\pi \ii x}$. Every $h \in \Homeo_+(\m S^1)$ admits a lift
$H:\m R \rightarrow \m R$ satisfying $H(x+1) = H(x) +1$ and
$\pi \circ H = h \circ \pi$.
For $h \in \Homeo_+(\m S^1)$, $x \in \m R$ and $r \in (0,1/2]$, define
\[
    Z_{x,r}(h) := \log\frac{\nu_h([x,x+r])}{\nu_h([x-r,x])},
\]
where $\nu_h([a,b])=H(b)-H(a)$ is the normalized length of the image arc. Thus $\nu_h(\m S^1)=1$.
We use this observable for both random weldings and deterministic
homeomorphisms. For a deterministic map $\varphi$, write
$\delta_\varphi(x,r)= Z_{x,r}(\varphi)$ and denote its lift by $\varPhi$.
\begin{lem}\label{lem:discrete WP}
    For $\varphi \in  \Homeo_+(\m S^1)$,
    \begin{itemize}
        \item $\varphi \in \QS(\m S^1)$ if and only if $\sup_{x \in \m R, r\in (0,1/2]} |\d_\varphi(x,r)| < \infty$.
        \item $\varphi \in  \Sy(\m S^1)$ if and only if $\lim_{r \to 0} \sup_{x \in \m R} |\d_\varphi(x,r)|=0$.
        \item If $\varphi \in \Sy(\m S^1)$ and
            \[
                E(\varphi) = \int_{\m S^1 \times (0,1/2]} \d_{\varphi}(x,r)^2  \dd \eta(x,r) < \infty,
            \]
        where $\dd \eta(x,r) = \dd x\dd r /r^2$,
        then $\varphi \in \WP(\m S^1)$.
    \end{itemize}
\end{lem}
\begin{remark}
    For a smooth lift, the local expansion is
    \[
        \d_\varphi(x,r)
        =\log\frac{r\varPhi'(x)+r^2\varPhi''(x)/2+o(r^2)}
        {r\varPhi'(x)-r^2\varPhi''(x)/2+o(r^2)}
        =r\frac{\varPhi''(x)}{\varPhi'(x)}+o(r).
    \]
    The integral defining $E(\varphi)$ involves all scales. More precisely,
    for a smooth periodic $v$ and $\varPhi_\vare(x)=x+\vare v(x)$,
    \[
        \lim_{\vare\to0}\frac{E(\varphi_\vare)}{\vare^2}
        =\int_{\m S^1}\int_0^{1/2}
        \frac{|v(x+r)-2v(x)+v(x-r)|^2}{r^4}\,\dd r\,\dd x
        \asymp \|v'\|_{\dot H^{1/2}(\m S^1)}^2,
    \]
    where the last comparison follows by Fourier expansion.
\end{remark}
\begin{proof}
    The first two assertions restate the definitions. For the third, we
    use the WP extension criterion. Since $\varphi$ is symmetric, it is
    quasisymmetric, and the Beurling--Ahlfors extension~\cite{BA56}
    \[
        U(x,y)=\frac{1}{2y}\int_{x-y}^{x+y}\varPhi(t)\,\dd t,
        \qquad
        V(x,y)=\frac{1}{y}\int_0^y
        \bigl(\varPhi(x+t)-\varPhi(x-t)\bigr)\,\dd t
    \]
    defines a quasiconformal map $F=U+\ii V$ of $\m H$.

    For $I=[a,b]$ of length $L$, let $\ell_I$ be the affine function
    agreeing with $\varPhi$ at its endpoints, and set
    \[
        \mc D_n I=\{[a+kL2^{-n},a+(k+1)L2^{-n}]:0\le k<2^n\},
        \qquad n\ge0.
    \]
    In particular, $\mc D_0I=\{I\}$. For
    $K=[m_K-r_K,m_K+r_K]$, put
    \[
        \psi_K(x)=\1_{x\in K}(1-|x-m_K|/r_K),\qquad
        c_K=\varPhi(m_K)-\frac{\varPhi(m_K-r_K)+\varPhi(m_K+r_K)}2.
    \]
    Successive dyadic linear interpolants give the uniformly convergent expansion
    \[
        \varPhi-\ell_I=\sum_{n\ge0}\sum_{K\in\mc D_nI}c_K\psi_K,
        \qquad
        \frac{c_K}{|\varphi(K)|}
        =-\frac12\tanh\frac{\d_\varphi(m_K,r_K)}2.
    \]
    Thus the average error
    \[
        T(I):=\frac1{|I|}\int_I(\varPhi(t)-\ell_I(t))\,\dd t
    \]
    satisfies
    \[
        T(I)=-\sum_{n\ge0}2^{-n-2}
        \sum_{K\in\mc D_nI}|\varphi(K)|
        \tanh\frac{\d_\varphi(m_K,r_K)}2.
    \]
    Since
    $\sum_{n\ge0}2^{-n-2}\sum_{K\in\mc D_nI}|\varphi(K)|/|\varphi(I)|=1/2$,
    weighted Cauchy--Schwarz gives
    \[
        \frac{T(I)^2}{|\varphi(I)|^2}
        \le C\sum_{n\ge0}2^{-n}
        \sum_{K\in\mc D_nI}\d_\varphi(m_K,r_K)^2.
    \]
    For $I_+(x,y)=[x,x+y]$, translation invariance gives
    \[
        \int_{\m S^1}\sum_{K\in\mc D_nI_+(x,y)}
        \d_\varphi(m_K,r_K)^2\,\dd x
        =2^n\int_{\m S^1}\d_\varphi(x,2^{-n-1}y)^2\,\dd x.
    \]
    Integrating over $0<y\le1/2$ and setting $r=2^{-n-1}y$, we obtain
    \[
        \int_{\m S^1}\int_0^{1/2}
        \frac{T(I_+)^2}{|\varphi(I_+)|^2}\frac{\dd y\,\dd x}{y^2}
        \le C\sum_{n\ge0}2^{-n-1}E(\varphi)
        \le CE(\varphi).
    \]
    The same estimate holds for $I_-(x,y)=[x-y,x]$.

    To estimate the Beltrami coefficient, write
    \[
        S=\varPhi(x+y)-\varPhi(x-y),\qquad
        D=\varPhi(x+y)+\varPhi(x-y)-2\varPhi(x),
    \]
    and set
    \[
        e_+=T(I_+)=\frac1y\int_x^{x+y}\varPhi(t)\,\dd t
        -\frac{\varPhi(x)+\varPhi(x+y)}2,
    \]
    \[
        e_-=T(I_-)=\frac1y\int_{x-y}^x\varPhi(t)\,\dd t
        -\frac{\varPhi(x-y)+\varPhi(x)}2.
    \]
    Direct differentiation yields
    \begin{align*}
        2y\bar\partial F&=e_+-e_-+
        \ii\left(\frac{5D}{4}-\frac{e_++e_-}{2}\right),\\
        2y\Re\partial F&=S-(e_+-e_-).
    \end{align*}
    Monotonicity gives $|e_\pm|\le|\varphi(I_\pm)|/2$, so
    $2y\Re\partial F\ge S/2$. Since
    $D/S=\tanh(\d_\varphi(x,y)/2)$, it follows that
    \[
        |\mu_F|^2\le C\left(
        \d_\varphi(x,y)^2+
        \frac{T(I_+)^2}{|\varphi(I_+)|^2}+
        \frac{T(I_-)^2}{|\varphi(I_-)|^2}\right).
    \]
    Consequently,
    \[
        \int_{\m S^1}\int_0^{1/2}|\mu_F|^2
        \frac{\dd y\,\dd x}{y^2}\le CE(\varphi)<\infty.
    \]
    The identity $F(z+1)=F(z)+1$ allows $F$ to descend to a
    quasiconformal extension on the disk. Near its boundary the
    hyperbolic area element is comparable to $\dd x\,\dd y/y^2$;
    the remaining compact part has finite hyperbolic area.
    The WP extension criterion therefore gives $\varphi\in\WP(\m S^1)$.
\end{proof}
For an integrable function $h$ on $\m S^1$, let $\pi_0 h$ denote its
mean-zero part. Throughout this section, $\dd x$ denotes Haar probability
measure on $\m R/\m Z$, identified with $\m S^1$ by $x\mapsto e^{2\pi\ii x}$.
Let $X$ be the mean-zero circle log-correlated Gaussian field (LGF) with
covariance kernel
\begin{equation}\label{Int eq:LGF}
 K_c(x,y)=-2\log|e^{2\pi\ii x}-e^{2\pi\ii y}|
          =-2\log|2\sin(\pi(x-y))|.
\end{equation}
Set $\g =\sqrt \k$, $\l = \g/2$ and $\tau = \l^2 \in (0,1)$.
For $\gamma\in(0,2)$, the associated GMC measure is formally written as
\begin{equation}\label{eq gmc formal}
    \mc M_X^\gamma(\dd x) = e^{\frac{\gamma}{2}X(e^{2\pi\ii x}) - \frac{\gamma^2}{8}\m E[X(e^{2\pi\ii x})^2]}\,\dd x.
\end{equation}
When $\gamma$ is fixed, we write $\mc M_X = \mc M_X^\gamma$. The corresponding normalized GMC measure is 
\begin{equation}\label{eq normalized gmc}
    \widehat{\mc M}_X^\g := \frac{1}{\mc M_X(\m S^1)} \mc M_X^\g.
\end{equation}
We also write $\mu_X=\widehat{\mc M}_X^\g$. Let $F_X$ be the degree-one lift of its cumulative function, so that $F_X(x+1) = F_X(x)+1$, and set $Q_X = F_X^{-1}$.
We use the auxiliary welding
\[
 H=F_{X_2}^{-1}\circ U\circ F_{X_1},
\]
where $X_1,X_2$ and the uniform rotation $U$ are independent.
By Lemma~\ref{introlem:sle-gmc-lift}, the law of $H$ is mutually absolutely
continuous with respect to the SLE welding law on the full homeomorphism space.
Write $\m P$ for the law of $H$ and $\m E$ for expectation under this law. For $\varphi \in \Homeo_+(\m S^1)$, define $R_\varphi :\Homeo_+(\m S^1) \rightarrow \Homeo_+(\m S^1) $ by $R_\varphi (h)  = h\circ \varphi$. Let $\m P_\varphi = (R_\varphi)_* \m P$ and $\m E_\varphi$ denote its expectation.

\begin{proof}[Proof of Lemma~\ref{introlem:sle-gmc-lift}]
Write $F_i=F_{X_i}$. By \cite[Lemma~4.8]{F1}, the law of $N(H)$ is
mutually absolutely continuous with respect to the law of the fixed-$1$
representative of the SLE welding. To lift this mutual absolute continuity, we check that $\m P$ is invariant under left
rotations. For fixed $z\in\m S^1$, set
\[
 c_z=F_2(z^{-1}),\qquad
 F_2^{(z)}=\Rot_{c_z^{-1}}\circ F_2\circ\Rot_{z^{-1}}.
\]
Rotation invariance of the LGF and covariance of its GMC imply that
$F_2^{(z)}$ has the same law as $F_2$. Moreover,
\[
 \Rot_z\circ H
 =(F_2^{(z)})^{-1}\circ(\Rot_{c_z^{-1}}\circ U)\circ F_1.
\]
Conditional on the fields, the middle factor is still a uniform rotation,
and its conditional law does not depend on the fields. Thus
$\Rot_z\circ H$ has law $\m P$.

Under the measurable identification
\[
 h\longleftrightarrow(N(h),h(1)),\qquad
 (h_0,z)\longmapsto\Rot_z\circ h_0,
\]
left rotation invariance gives
\[
 \begin{aligned}
 (N,\operatorname{ev}_1)_*\m P
 &=(N_*\m P)\otimes m,\\
 (N,\operatorname{ev}_1)_*\m P^{\mathrm{weld}}
 &=(N_*\m P^{\mathrm{weld}})\otimes m.
 \end{aligned}
\]
where $\operatorname{ev}_1(h)=h(1)$. Mutual absolute continuity of the first marginals therefore implies
mutual absolute continuity of the laws on the full space.
\end{proof}

\subsection{Local fields and logarithmic GMC moments}\label{subsec:subcritical-gmc}
We establish the local GMC estimates used in the welding argument.
A change from circle coordinates to a real chart gives a radial
representation in terms of Brownian motion and a stationary field.
We use this representation to control logarithmic masses and increments
of the inverse cumulative function. The resulting bounds,
Lemmas~\ref{lem:inverse rotation} and~\ref{lem:log GMC moment}, are combined
in the next subsection to estimate adjacent-arc ratios for the welding.
Let $K_c$ be the circle kernel in \eqref{Int eq:LGF}, and set
$K_l(t,u)=-2\log|t-u|$. Under the coordinate change
$x(t)=\arctan(t)/\pi$ for $t\in[-1,1]$, a direct calculation
\begin{align*}
    |e^{2\pi\ii x(t)} - e^{2\pi \ii x(u)}| = &|\frac{(1+it)^2}{1+t^2}-\frac{(1+iu)^2}{1+u^2}|\\
    =&|\frac{2(t-u)\left( -(u+t)+\ii(1-tu)\right)}{(u+t)^2+(1-tu)^2}|\\
    =& \frac{2|t-u|}{\sqrt{(1+t^2)(1+u^2)}}
\end{align*}
shows that
\[
    K_c(x(t),x(u)) = K_l(t,u)+b(t)+b(u),
\]
where $b(t)=\log(1+t^2)-\log2$. Let $X=X_c$ and $X_l$ denote the fields with covariances $K_c$ and $K_l$, respectively. Let $\mc M_c$ and $\mc M_l$ denote the corresponding GMC measures.
The covariance difference vanishes on pairs of test functions of integral
zero. Thus the fields have the same law modulo constants; coupling this
common Gaussian part and their Gaussian averages gives
$X_c\circ x=X_l+G_0$, where $G_0$ is a Gaussian random constant. In general
$G_0$ is correlated with $X_l$, and the covariance identity reads
\[
    \var G_0+2\cov(X_l(t),G_0)=2b(t).
\]
Consequently,
\[
    \mc M_c (\dd x(t)) = \frac{2^\tau}{\pi} e^{\l G_0}(1+t^2)^{-1-\tau} \mc M_l(\dd t).
\]
For a jointly Gaussian pair consisting of a field $Y$ and a random variable $C$,
\[
    \mc M_{Y+C} (\dd t) = e^{\l C-\tau \var(C)/2-\tau \cov(Y(t),C)} \mc M_Y(\dd t).
\]
And $\mc M_c(\dd x(t)) = x'(t) \mc M_{X_c \circ x}(\dd t)$.
Moreover,
\[
    \mc M^\bullet (\dd x):= |e^{2\pi\ii x}-1|^{-2\tau} \mc M_c(\dd x)
\]
satisfying
\[
    \mc M^\bullet (\dd x(t))= R_0 \frac{|t|^{-2\tau}}{1+t^2} \mc M_l(\dd t), \qquad R_0 = \frac{2^{-\tau}}{\pi} e^{\l G_0}.
\]
We next couple the fields on the two sides of the origin.
For $s,v>0$,
\[
    -2\log|\pm e^{-s}-e^{-v}|= 2\min(s,v) + 2\sum_{n \geq 1} \frac{(\pm 1)^n  e^{-n|s-v|}}{n}.
\]
Choose the radial coupling
\[
    X_l(\pm e^{-s})=-\sqrt2 B_s+N_\pm(s),\qquad
    N_\pm(s)=\sum_{n\ge1}\sqrt{2/n}(\pm1)^nY_n(s),
\]
where $B$ is standard Brownian motion, independent of the stationary
Gaussian processes $Y_n$, and
\[
    \cov(Y_n(s),Y_m(v))=\1_{n=m}e^{-n|s-v|}.
\]
Let $\mc N_\pm$ be the GMC measures with parameter $\l$ associated with
$N_\pm$. They are jointly stationary, independent of $B$, and have mean
measure $\dd s$.

Define
\[
    J_\pm (s) = \int_s^\infty e^{-A_u} \mc N_\pm (\dd u), \qquad  A_s = (1-\tau ) s + \sqrt{2\tau} B_s, \qquad d :=1-\tau >0.
\]
For $\xi\ge0$, define the hitting time
\[
    \t(\xi)=\inf\{s\ge0:A_s=\xi\}.
\]
For every $\xi\in\m R$, define
\[
    S_\pm(\xi)=\inf\{s\ge0:J_\pm(s)\le e^{-\xi}\}.
\]
Note that
\[
    P( \mc N_\pm([n,n+1])> n^3) \leq E[\mc N_\pm([n,n+1])]/n^3 \leq 1/n^3,
\]
The Borel--Cantelli lemma therefore gives, almost surely, $\mc N_\pm([n,n+1])\leq  n^3$ in the tail ($n \geq n_0(\o)$).
There exists $u_0(\o)$ such that $A_u \geq du/2$ for $u \geq u_0(\o)$,
\[
    J_\pm(s) \leq \sum_{n=\lfloor s\rfloor}^\infty n^3 e^{-d n/2} < \infty,
\]
for $s \geq n_0(\o)\vee u_0(\o)$,
and therefore
\[
    J_\pm(s) \to 0,~\text{as}~s \to \infty,
\]
combined with $P(\mc N(I)< \infty, \forall I=[p,q]\in \m Q^2) = 1$ and the lower boundedness of $A$, we obtain
\[
    J_\pm(s)  < \infty
\]
for all $s\geq0$. Since subcritical GMC is non-atomic and has full support,
$J_\pm$ is continuous and strictly decreasing.

Note that if $R \geq 0$ is a $B$-measurable random variable, not necessarily a stopping time, we have
\[
    \mathrm{Law} (\mc N_\pm (R + \cdot )) = \mathrm{Law}(\mc N_\pm(\cdot)) 
\]
and the shifted measure remains independent of $B$.
\begin{lem}\label{lem:radial-positive-moments}
    For $t\ge0$ and $0<q<\min(1,d/\tau)$,
    \[
        E[J_\pm(t)^q]
        =e^{-q(d-\tau q)t}E[J_\pm(0)^q]<\infty.
    \]
\end{lem}
\begin{proof}
    By independent Brownian increments and stationarity of $\mc N_\pm$,
    \[
        J_\pm(t)=e^{-A_t}
        \int_0^\infty e^{-(A_{t+v}-A_t)}\mc N_\pm(t+\dd v)
        =:e^{-A_t}\widetilde J_\pm(t),
    \]
    where $\widetilde J_\pm(t)$ is independent of $A_t$ and has the law
    of $J_\pm(0)$. Since $E[e^{-qA_t}]=e^{-q(d-\tau q)t}$, this gives
    the asserted identity. To prove finiteness, set
    \[
        K_{n,\pm}=\int_n^{n+1}e^{-(A_u-A_n)}\mc N_\pm(\dd u).
    \]
    Each $K_{n,\pm}$ is independent of $A_n$ and has the law of $K_{0,\pm}$.
    Moreover, $E[K_{0,\pm}]=\int_0^1e^{-(d-\tau)u}\,\dd u<\infty$.
    Jensen's inequality and subadditivity of the $q$th power yield
    \[
        E[J_\pm(0)^q]
        \le\sum_{n\ge0}E[e^{-qA_n}]E[K_{n,\pm}^q]
        \le C_q\sum_{n\ge0}e^{-q(d-\tau q)n}<\infty.
    \]
\end{proof}
\begin{lem}\label{lem:stationary-negative-moments}
    For every $q>0$, there is $v_q<\infty$ such that
    \[
        \sup_{|I|=1}E[\mc N_\pm(I)^{-q}]\le v_q.
    \]
\end{lem}
\begin{proof}
    On a unit interval the covariance is
    \[
        2\sum_{n\ge1}\frac{e^{-n|s-v|}}n
        =-2\log|s-v|+g(s,v),\qquad
        g(s,v)=-2\log\frac{1-e^{-|s-v|}}{|s-v|},
    \]
    where $g$ extends continuously across the diagonal and is bounded.
    The negative-moment theorem for subcritical GMC applies, with
    parameter $\sqrt2\l<\sqrt2$ in the covariance-$-\log$ convention;
    see \cite[Theorems~3.19 and~3.42]{GFFbook}. Stationarity makes the
    bound uniform over unit intervals.
\end{proof}
The same arguments apply with any positive drift $a$ in place of $d$.
For later use, write
\[
    A_s^{(a)}=as+\sqrt{2\tau}B_s,\qquad
    J_\pm^{(a)}(s)=\int_s^\infty e^{-A_u^{(a)}}\mc N_\pm(\dd u).
\]
The preceding positive-moment estimate holds for $0<q<\min(1,a/\tau)$.
Also,
\[
    J_\pm^{(a)}(0)\ge
    e^{-\sup_{0\le u\le1}A_u^{(a)}}\mc N_\pm([0,1]).
\]
The Brownian supremum has all exponential moments and is independent
of $\mc N_\pm$. Thus $J_\pm^{(a)}(0)$ has all negative moments, and
\begin{equation}\label{eq:radial-log-moments}
    E[|\log J_\pm^{(a)}(0)|^p]<\infty,
    \qquad a>0,\quad 1\le p<\infty.
\end{equation}

\begin{lem}\label{lem:two-sided estimate}
    There exist $c$, $C>0$ such that for all $\xi \geq 0$ and $t \geq 2$:
    \[
        P(|S_\pm(\xi) -\t(\xi)| > t) \leq Ce^{-ct}.
    \]
\end{lem}
\begin{proof}
    Fix $0<q<\min(1,d/(2\tau))$. For the right tail, using the monotonicity of $J_\pm$,
    \[
        \{S_\pm(\xi) > \t(\xi) + t\} =\{J_\pm(\t(\xi)+t)> e^{-\xi}\} = \{ e^{\xi}J_\pm(\t(\xi)+t)>1\},
    \]
    By the strong Markov property, the law of
    \[
        e^{\xi}J_\pm(\t(\xi)+t) = \int_t^\infty e^{-(A_{\t(\xi)+v}-\xi)} \mc N_\pm(\t(\xi)+\dd v)
    \]
    is the law of $J_\pm(t)$. Hence
    \[
        P(S_\pm(\xi) > \t(\xi) + t) \leq P(J_\pm(t) >1) \leq  E [J_\pm(t)^q] = e^{-q(d-\tau q)t} E[J_\pm(0)^q].
    \]
    For the left tail, it suffices to work on $\{\t(\xi) \geq t\}$. Set $I_t = [\t(\xi) -t,\t(\xi) -t+1] \subset [0,\t(\xi) ]$ and write $Y_u : = \xi -A_u$ for $ u \leq \t(\xi) $. Then
    \[
        e^\xi J_\pm(\t(\xi) -t) \geq \int_{I_t} e^{Y_u} \mc N_\pm(\dd u) =:X_t > 0,
    \]
    \[
        \{S_\pm(\xi) < \t(\xi) - t\}  = \{ e^{\xi}J_\pm(\t(\xi)-t)<1\} \subset \{X_t < 1\},
    \]
    Therefore
     \[
        P(S_\pm(\xi)<\t(\xi)-t)
        \le E[X_t^{-q};\t(\xi)\ge t].
    \]
    Let $\mu_t =\mc  N_\pm|_{I_t}$, $M_t = \mc N_\pm(I_t)$ and $\nu_t = \mu_t /M_t$. 
    \[
        X_t^{-q}=(M_t \int e^Y \dd \nu_t)^{-q} \leq M_t^{-q-1}\int_{I_t}e^{-qY_u }\mc N_\pm(\dd u)\leq M_t^{-q} \sup_{u \in I_t}  e^{-qY_u }.
    \]
    Hence
    \[
        E[X_t^{-q};\t(\xi)\ge t]
        \le v_q E_B[\sup_{u\in I_t}e^{-qY_u};\t(\xi)\ge t]
        \le Ce^{-q(d-\tau q)(t-1)}.
    \]
    To prove the last estimate, use the transition and hitting-time
    densities of the drift-$-d$ Brownian motion $Y$, killed at $0$.
    The left-tail event is empty when $\xi=0$, so assume $\xi>0$.
    Define the free transition density
    \[
        p^\text{free} (y_0,y_1;s):= (4\pi\tau s)^{-1/2} \exp(-\frac{(y_1-y_0+ds)^2}{4\tau s})
    \]
    and the killed transition density
    \[
        p^K (y_0,y_1;s)=p^\text{free} (y_0,y_1;s)-e^{dy_0/\tau} p^\text{free} (-y_0,y_1;s)=p^\text{free} (y_0,y_1;s)(1-e^{-y_0 y_1/\tau s}).
    \]
    In particular,
    \[
        0 \leq p^K (y_0,y_1;s) \leq p^\text{free} (y_0,y_1;s), \qquad p^K (y_0,0;s) =0.
    \]
    For $A_u = a +du+\sqrt{2\tau} B_u$, the Green kernel is
    \[
        G_\text{free} (a,x):=\int_0^\infty (4\pi\tau u)^{-1/2} \exp\left(-\frac{(x-a-du)^2}{4\tau u}\right)\dd u = \exp(d\min(x-a,0)/\tau)/d.
    \]
    For $A$ killed on hitting $\xi$, the Green kernel, for $x<\xi$, is
    \[
        G_\xi(0,x) := G_\text{free} (0,x)-G_\text{free} (\xi,x) = \frac{\exp(d\min(x,0)/\tau)-\exp(d(x-\xi)/\tau)}{d}.
    \]
    Define
    \[
        G^\text{def}(\xi,y)=G_\xi(0,\xi-y)\in[0,1/d],\qquad y>0.
    \]
    Let $\hat Y$ be a free Brownian motion of variance rate $2\tau$ and
    drift $-d$, started from $y_0>0$, and define
    \[
        K(y_0,\dd y_1):=E_{y_0}^{\mathrm{free}}
        [e^{-q\min_{[0,1]}\hat Y}
        \1_{\min_{[0,1]}\hat Y>0}\1_{\hat Y_1\in\dd y_1}].
    \]
    The density of the hitting time of $0$, starting from $y>0$, is
    \[
        f_y(s) = \frac{y}{\sqrt{4\pi \tau s^3}} \exp(-\frac{(y-ds)^2}{4\tau s}).
    \]
    Disintegrating at the hitting time and using the Markov property,
    for $s>t$ we have
    \[
        E_B[\sup_{u\in I_t}e^{-qY_u};\t(\xi)\in\dd s]
        =\left(\int_0^\infty p^K(\xi,y_0;s-t)\,\dd y_0
        \int_0^\infty K(y_0,\dd y_1)f_{y_1}(t-1)\right)\dd s.
    \]
    Conditioned on $\hat Y_1 = y_1$,
    \[
        \hat Y_u = (1-u)y_0  + u y_1 + \sqrt{2\tau} \beta_u,
    \]
    where $\beta_u$ is a Brownian bridge and for $\mc B = \sup_{[0,1]}|\beta|$, we have
    \[
        P(\mc B > x) \leq 2e^{-2x^2},
    \]
    therefore
    \[
        c_q:=E[e^{q\sqrt{2\tau}\mc B}]<\infty.
    \]
    Under this free bridge law, let $m=\min_{[0,1]}\hat Y$. Then
    $m\ge\min(y_0,y_1)-\sqrt{2\tau}\mc B$, so
    \[
        E[e^{-qm}\mid\hat Y_1=y_1]
        \le c_q(e^{-qy_0}+e^{-qy_1}).
    \]
    Conditioning further on $m>0$ decreases the expectation of
    $e^{-qm}$, since this is a decreasing function of $m$. Therefore
    \[
        K(y_0,\dd y_1)
        \le c_q(e^{-qy_0}+e^{-qy_1})p^K(y_0,y_1;1)\,\dd y_1.
    \]
    Note the identity
    \[
        f_{y_0}(t) = \int_0^\infty p^K(y_0,y_1;1)f_{y_1}(t-1) \dd y_1,
    \]
    we have
    \[
        E_B[ \sup_{u \in I_t}  e^{-qY_u } \1_{\t(\xi) \geq t}] \leq c_q(A+B),
    \]
    where
    \[
        A := \int_0^\infty e^{-qy_0}f_{y_0}(t)G^\text{def}(\xi,y_0) \dd y_0,
    \]
    and
    \[
        B:=\int_0^\infty \int_0^\infty e^{-qy_1} p^K(y_0,y_1;1) f_{y_1}(t-1)G^\text{def}(\xi,y_0) \dd y_1 \dd y_0.
    \]
    \[
        \begin{aligned}
        A &\leq \frac1d\int_0^\infty e^{-qy_0} f_{y_0}(t)\,\dd y_0\\
          &=\frac{e^{-q(d-\tau q)t}}{d}
            \int_0^\infty \frac{y_0}{t}(4\pi\tau t)^{-1/2}
                 e^{-(y_0-at)^2/(4\tau t)}\,\dd y_0\\
          &\leq \frac{e^{-q(d-\tau q)t}}{d}
                  \bigl(a+\sqrt{2\tau}E[|Z|]\bigr),
        \end{aligned}
    \]
    where
    \[
        a = d-2\tau q > 0,
    \]
    and $Z$ is a standard Gaussian random variable.
    \begin{align}
        B \leq & \int_0^\infty \int_0^\infty e^{-qy_1} p^\text{free}(y_0,y_1;1) f_{y_1}(t-1) \dd y_1 \dd y_0 /d \\
        \leq &  \int_0^\infty  e^{-qy_1} f_{y_1}(t-1) \dd y_1  /d\\
        \leq & \frac{e^{-q(d-\tau q)(t-1)}}{d} (a+\sqrt{2\tau}E[|Z|]).
    \end{align}

\end{proof}
\begin{cor}\label{cor:inverse-level-moments}
    For $p\in[1,\infty)$ and $\xi_1,\xi_2\in\m R$,
    \[
        E[|S_+(\xi_1)-S_-(\xi_2)|^p]
        \le C_p(1+|\xi_1-\xi_2|)^p.
    \]
\end{cor}
\begin{proof}
    First let $\xi_1,\xi_2\ge0$. Write
    \begin{align*}
        S_+(\xi_1)-S_-(\xi_2)
        ={}&(S_+(\xi_1)-\t(\xi_1))
        +(\t(\xi_1)-\t(\xi_2))\\
        &-(S_-(\xi_2)-\t(\xi_2)).
    \end{align*}
    Lemma~\ref{lem:two-sided estimate} bounds the first and last terms
    uniformly in $L^p$. By the strong Markov property, the absolute value
    of the middle term has the law of $\t(\D)$, where
    $\D=|\xi_1-\xi_2|$. The hitting-time density displayed above shows
    that $\t(1)$ has moments of every order. For $\D>0$,
    $\t(\D)\le\t(\lceil\D\rceil)$, and the latter is a sum of
    $\lceil\D\rceil$ independent increments with the law of $\t(1)$.
    Minkowski's inequality gives
    \[
        \|\t(\D)\|_{L^p}\le C_p(1+\D).
    \]
    This proves the assertion for nonnegative levels.
    If both levels are negative, use $0\le S_\pm(\xi)\le S_\pm(0)$.
    If one is negative and the other nonnegative, compare the negative
    level to $0$; the nonnegative level is at most $|\xi_1-\xi_2|$.
    These observations prove the same bound for all real levels.
\end{proof}
Recall that
\[
    \mc M^\bullet (\dd x):= |e^{2\pi\ii x}-1|^{-2\tau} \mc M_c(\dd x)
\]
satisfying
\[
    \mc M^\bullet (\dd x(t))= R_0 \frac{|t|^{-2\tau}}{1+t^2} \mc M_l(\dd t), \qquad R_0 = \frac{2^{-\tau}}{\pi} e^{\l G_0}.
\]
Let $T^\bullet = \mc M^\bullet (\m S^1)$ and $D_\bullet = 1/ T^\bullet$.
\begin{lem}\label{lem:relation of bullet and J}
    For $s\ge0$,
    \[
        \mc M^\bullet([0,x(e^{-s})]) = R_0 \int_s ^\infty \frac{1}{1+e^{-2u}}e^{-A_u} \mc N_+(\dd u),
    \]
    \[
        \mc M^\bullet([-x(e^{-s}),0]) = R_0 \int_s ^\infty \frac{1}{1+e^{-2u}}e^{-A_u} \mc N_-(\dd u).
    \]
\end{lem}
\begin{proof}
    Under $t=\pm e^{-s}$, the chosen radial coupling gives
    \[
        \mc M_l^\pm(\dd t)
        =e^{-\l\sqrt2 B_s-\tau s}e^{-s}\mc N_\pm(\dd s).
    \]
    Multiplying by $R_0|t|^{-2\tau}/(1+t^2)$ gives
    \[
        \mc M^\bullet_\pm(\dd x(t))
        =R_0\frac{e^{-A_s}}{1+e^{-2s}}\mc N_\pm(\dd s),
    \]
    proving both identities on the same probability space for $s\ge0$.
\end{proof}
Write $T=\mc M_c(\m S^1)$. The usual GMC moment bounds give
$E[T^q]<\infty$ for $0<q<1/\tau$, and $E[T^{-q}]<\infty$ for every
$q>0$; see \cite[Theorems~3.30 and~3.42]{GFFbook}.
Since $|e^{2\pi\ii x}-1|\le2$,
$T^\bullet\ge2^{-2\tau}T$. On the two local arcs,
Lemma~\ref{lem:relation of bullet and J} bounds the rooted mass by
$R_0(J_+(0)+J_-(0))$; on their complement the insertion weight is bounded.
Consequently,
\[
    T^\bullet\le C T+R_0(J_+(0)+J_-(0)).
\]
As $R_0$ has moments of every real order, H\"older's inequality and
Lemma~\ref{lem:radial-positive-moments} give a positive moment of
$T^\bullet$. With $C^\bullet=\log(T^\bullet/R_0)$, another application
of H\"older's inequality gives
\begin{equation}\label{eq:rooted-normalization-moments}
    E[D_\bullet^q]<\infty\quad(q>0),\qquad
    E[e^{a|C^\bullet|}]<\infty\quad\text{for some }a>0.
\end{equation}
No independence between $R_0$ and the radial fields is used here.

For a measure $\mu$, write $-z\mu$ for its image under rotation by $-z$.
\begin{lem}\label{lem:rooted-palm}
    Let $\Phi$ be a nonnegative measurable functional on finite nonzero
    positive measures satisfying $\Phi(c\mu)=\Phi(\mu)$ for $c>0$. Then
    \[
        E\left[\int_{\m S^1}\Phi(-z\mu_X)\,\mu_X(\dd z)\right]
        =E\left[D_\bullet\Phi\left(\frac{\mc M^\bullet}{T^\bullet}\right)\right].
    \]
\end{lem}
\begin{proof}
    Use a rotation-invariant Poisson regularization $X_\vare$ with covariance
    \[
        K_\vare(z,x)=-2\log|1-e^{-\vare}e^{2\pi\ii(x-z)}|,
    \]
    and let
    \[
        \mc M^\vare(\dd x)=
        e^{\l X_\vare(x)-\tau\var(X_\vare(x))/2}\,\dd x,
        \qquad m_\vare=\mc M^\vare(\m S^1),\qquad
        \mu_\vare=\mc M^\vare/m_\vare.
    \]
    For each $z$, set
    \[
        \mc M^{\bullet,z,\vare}(\dd x)
        =e^{\tau K_\vare(z,x)}\mc M^\vare(\dd x),
        \qquad T^{\bullet,z,\vare}=\mc M^{\bullet,z,\vare}(\m S^1).
    \]
    Gaussian Girsanov, followed by rotation invariance, gives
    \begin{align*}
        E\int_{\m S^1}\Phi(-z\mu_\vare)\,\mu_\vare(\dd z)
        &=\int_{\m S^1}E\left[
        \frac1{T^{\bullet,z,\vare}}
        \Phi\left(\frac{-z\mc M^{\bullet,z,\vare}}{T^{\bullet,z,\vare}}\right)
        \right]\dd z\\
        &=E\left[D^{\bullet,\vare}
        \Phi\left(\frac{\mc M^{\bullet,\vare}}{T^{\bullet,\vare}}\right)\right],
    \end{align*}
    where the quantities on the last line have insertion at $0$ and
    $D^{\bullet,\vare}=1/T^{\bullet,\vare}$.

    We justify the limit first for bounded continuous $\Phi$ on probability
    measures. Away from $0$, the insertion weights converge uniformly.
    Near $0$, the usual GMC scaling estimate and Kahane comparison
    \cite[Theorem~3.19 and Section~3.9]{GFFbook} give, for
    $0<q<\min(1,d/\tau)$,
    \[
        \sup_{0<\vare\le r}
        E[\mc M^{\bullet,\vare}([-r,r])^q]
        \le C_q r^{q(d-\tau q)},\qquad 0<r\le1/4.
    \]
    Indeed, on an annulus of radius $\rho\ge\vare$ the insertion is
    bounded by $C\rho^{-2\tau}$ and the ordinary mass has $q$th moment
    bounded by $C\rho^{(1+\tau)q-\tau q^2}$.
    The core of radius $\vare$ has the same bound at scale $\vare$.
    Summing over dyadic annuli, using subadditivity, proves the display.
    Thus $\mc M^{\bullet,\vare}\to\mc M^\bullet$ weakly in probability.
    Moreover, $K_\vare\ge-2\log2$, so
    $D^{\bullet,\vare}\le2^{2\tau}/m_\vare$.
    The uniform negative-moment bounds for these regularizations,
    obtained from Kahane comparison and
    \cite[Theorem~3.42]{GFFbook}, give uniform integrability.
    We may therefore pass to the limit. Equality of the resulting finite
    measures on the space of probability measures extends the identity
    to nonnegative measurable $\Phi$ by a monotone-class argument.
\end{proof}
\begin{lem}\label{lem:inverse rotation}
For any $1 \leq p < \infty$, $d_\pm>0$ with $d_++d_- \leq 1$, there exists $C_{\g,p}$
    \[
        E\left[\left|\log \frac{Q(V+d_+)-Q(V)}{Q(V)-Q(V-d_-)}\right|^p\right] \leq C_{\g,p} (1+ |\log\frac{d_+}{d_-}|)^p,
    \]
    where $U$ and $V$ are independent uniform rotations, independent of $X$, and $Q = F_\mu^{-1}$ with $\mu = \mu_X$ or $\mu = \mu^\bullet =U\mc M^\bullet/T^\bullet$ from $X = X_c$.
\end{lem}
\begin{proof}
    Put $Q=F_\mu^{-1}$ and
    \[
        G(\mu,v)=\log\frac{Q(v+d_+)-Q(v)}{Q(v)-Q(v-d_-)},
        \qquad \Phi_p(\mu)=|G(\mu,0)|^p,
    \]
    extending $\Phi_p$ to finite nonzero measures by normalization.
    The change of variables $v=F_X(z)$ and
    Lemma~\ref{lem:rooted-palm} give
    \begin{align*}
        E[|G(\mu_X,V)|^p]
        &=E\int_{\m S^1}\Phi_p(-z\mu_X)\,\mu_X(\dd z)\\
        &=E\left[D_\bullet\Phi_p\left(\frac{\mc M^\bullet}{T^\bullet}\right)\right].
    \end{align*}
    Define, for $s\ge0$ and $h\in\m R$,
    \[
        J_\pm^\o(s)=\int_s^\infty
        \frac{e^{-A_u}}{1+e^{-2u}}\mc N_\pm(\dd u),
        \qquad S_\pm^\o(h)=\inf\{s\ge0:J_\pm^\o(s)\le e^{-h}\}.
    \]
    Since $J_\pm/2\le J_\pm^\o\le J_\pm$,
    \[
        S_\pm(h-\log2)\le S_\pm^\o(h)\le S_\pm(h).
    \]
    Corollary~\ref{cor:inverse-level-moments} and these bounds imply
    \begin{equation}\label{eq:weighted-inverse-levels}
        \|S_+^\o(u)-S_-^\o(v)\|_{L^p}
        \le C_{\g,p}(1+|u-v|),\qquad u,v\in\m R.
    \end{equation}

    Let $\ell_\pm$ be the lengths of the forward and backward arcs from
    $0$ having respective masses $d_\pm$ under $\mc M^\bullet/T^\bullet$.
    Clip these lengths to the coordinate chart by setting
    $\bar\ell_\pm=\min(\ell_\pm,1/4)$. With
    \[
        h_\pm=-\log d_\pm-C^\bullet,
    \]
    Lemma~\ref{lem:relation of bullet and J} gives
    \[
        \bar\ell_\pm=x(e^{-S_\pm^\o(h_\pm)}).
    \]
    If the prescribed mass exceeds the mass of the local arc,
    $S_\pm^\o(h_\pm)=0$, so this formula still holds.
    Since $e^{-s}/4\le x(e^{-s})\le e^{-s}/\pi$ and
    $1\le\ell_\pm/\bar\ell_\pm\le4$, we obtain
    \begin{equation}\label{eq:inverse-arcs-clipped}
        \log\frac{\ell_+}{\ell_-}
        =-\bigl(S_+^\o(h_+)-S_-^\o(h_-)\bigr)+O(1),
    \end{equation}
    with a deterministic uniform bound on the error.

    The common shift $C^\bullet$ is correlated with $S_\pm^\o$,
    so we use a partition into deterministic levels.
    Put $a_\pm=-\log d_\pm$ and $E_j=\{j\le C^\bullet<j+1\}$,
    $j\in\m Z$. By monotonicity, on $E_j$,
    \[
        |S_+^\o(a_+-C^\bullet)-S_-^\o(a_--C^\bullet)|\le M_j,
    \]
    where
    \begin{align*}
        M_j={}&|S_+^\o(a_+-j)-S_-^\o(a_--j-1)|\\
        &+|S_-^\o(a_--j)-S_+^\o(a_+-j-1)|.
    \end{align*}
    By \eqref{eq:weighted-inverse-levels},
    $\|M_j\|_{L^{4p}}\le C_{\g,p}(1+|a_+-a_-|)$, uniformly in $j$.
    H\"older's inequality therefore gives
    \begin{align*}
        &E[D_\bullet|S_+^\o(h_+)-S_-^\o(h_-)|^p]\\
        &\quad\le\sum_{j\in\m Z}
        \|D_\bullet\|_{L^4}\|M_j\|_{L^{4p}}^pP(E_j)^{1/2}\\
        &\quad\le C_{\g,p}(1+|a_+-a_-|)^p.
    \end{align*}
    Here \eqref{eq:rooted-normalization-moments} ensures both
    $\|D_\bullet\|_{L^4}<\infty$ and
    $\sum_jP(E_j)^{1/2}<\infty$. Since
    $|a_+-a_-|=|\log(d_+/d_-)|$, the ordinary case follows from
    \eqref{eq:inverse-arcs-clipped}.

    For $\mu^\bullet=U\mc M^\bullet/T^\bullet$, the same Girsanov
    calculation as in Lemma~\ref{lem:rooted-palm}, without the factor
    $m_\vare^{-1}$, gives the different identity
    \[
        E[\Psi(\mu^\bullet)]=E[T\Psi(\mu_X)]
    \]
    for nonnegative measurable $\Psi$. Indeed, before taking the limit,
    $E[m_\vare\Psi(\mu_\vare)]$ equals
    $\int_{\m S^1}E[\Psi(\mc M^{\bullet,z,\vare}/T^{\bullet,z,\vare})]\,\dd z$;
    the uniform positive moments of $m_\vare$ justify the limit for bounded
    continuous $\Psi$, and a monotone-class argument gives the stated identity.
    Apply it to $\Psi(\mu)=\int_0^1|G(\mu,v)|^p\,\dd v$.
    Choose $1<r<1/\tau$ and let $r'=r/(r-1)$. The ordinary case already proved
    and H\"older's inequality yield
    \begin{align*}
        E[|G(\mu^\bullet,V)|^p]
        &=E[T|G(\mu_X,V)|^p]\\
        &\le E[T^r]^{1/r}
        E[|G(\mu_X,V)|^{pr'}]^{1/r'}\\
        &\le C_{\g,p}(1+|\log(d_+/d_-)|)^p.
    \end{align*}
\end{proof}
\begin{lem}\label{lem:log GMC moment}
    For $1\leq p <\infty$ and adjacent arcs $A$ and $B$ sharing one endpoint, we have
    \[
        E\left[\left|\log\frac{\mu(A)}{\mu(B)}\right|^p\right] \leq C_{\g,p} \left(1+ |\log\frac{|A|}{|B|}|\right)^p,
    \]
    where $\mu = \mu_X$ is the standard GMC measure.
\end{lem}
\begin{proof}
    Since normalization cancels in a mass ratio, we may work with
    $\mc M_c$. By rotation invariance, take the common endpoint to be $0$,
    with $A=[0,a]$ and $B=[-b,0]$.
    First assume $a,b\le1/4$, and set
    \[
        s_A=-\log\tan(\pi a),\qquad
        s_B=-\log\tan(\pi b),\qquad a_0=1+\tau.
    \]
    The ordinary GMC has drift $a_0$, rather than the rooted drift $1-\tau$:
    \[
        \mc M_l^\pm(\dd t)
        =e^{-(1+\tau)s-\sqrt{2\tau}B_s}\mc N_\pm(\dd s),
        \qquad t=\pm e^{-s}.
    \]
    The chart weight $(1+t^2)^{-1-\tau}$ is bounded above and below for
    $|t|\le1$, and the common factor $e^{\l G_0}$ cancels. Thus
    \[
        \log\frac{\mu(A)}{\mu(B)}
        =\log J_+^{(a_0)}(s_A)-\log J_-^{(a_0)}(s_B)+O(1),
    \]
    with a deterministic uniform error. For each deterministic $s$,
    \[
        J_\pm^{(a_0)}(s)=e^{-A_s^{(a_0)}}\widetilde J_\pm^{(a_0)}(s),
        \qquad
        \widetilde J_\pm^{(a_0)}(s)\overset{\mathrm{law}}=
        J_\pm^{(a_0)}(0).
    \]
    By \eqref{eq:radial-log-moments}, the logarithms of the last factors
    have uniformly bounded $p$th moments. Brownian increment moments give
    \begin{align*}
        E\left[\left|\log\frac{\mu(A)}{\mu(B)}\right|^p\right]
        &\le C_{\g,p}
        \bigl(1+|s_A-s_B|^p+|s_A-s_B|^{p/2}\bigr)\\
        &\le C_{\g,p}(1+|\log(a/b)|)^p,
    \end{align*}
    since $\tan(\pi a)/a$ and $\tan(\pi b)/b$ lie between $\pi$ and $4$.

    For larger arcs, shorten each at its far endpoint to length at most
    $1/4$, keeping the common endpoint. If $A'$ is the shortened part of
    an arc $A$ of length greater than $1/4$, then
    \[
        0\le\log\frac{\mc M_c(A)}{\mc M_c(A')}
        \le |\log T|+|\log\mc M_c(A')|.
    \]
    The right side has all finite moments by the positive and negative
    GMC moment bounds, uniformly over arcs $A'$ of length $1/4$.
    The same applies to $B$, and shortening changes the logarithmic length
    ratio by a bounded amount. This proves the general case.
\end{proof}
\subsection{Adjacent-arc moments and quasisymmetry}\label{subsec:subcritical-qs}

The preceding source and inverse-GMC estimates give moment bounds for
the welding observable. We pair these bounds with a deterministic
lower bound obtained by averaging logarithmic masses over rotations.
If a deformation is not quasisymmetric, its distortion diverges along
a sequence of adjacent arcs. After normalization, the observable tends
to zero under the original law but retains a positive mean under the
transformed law. Absolute continuity and uniform integrability give
the contradiction in Theorem~\ref{thm:A}.

\begin{lem}\label{lem:moment of Z}
    For $1\leq p <\infty$ and adjacent arcs $A$ and $B$ sharing one endpoint, we have
    \[
        \m E\left[\left|\log\frac{\nu_H(A)}{\nu_H(B)}\right|^p\right] \leq C_{\g,p} \left(1+ |\log\frac{|A|}{|B|}|\right)^p,
    \]
    in particular,
    \[
        \sup_{x \in \m S^1,r \in (0,1/2]} \m E|Z_{x,r}|^p \leq C_{\g,p}.
    \]
    For arcs $I\subset J$ with $0<|I|<|J|$ and $J\backslash I$ connected,
    \[
        \m E\left[\left|\log\frac{\nu_H(I)}{\nu_H(J)}\right|^p\right] \leq C_{\g,p} \left(1+ |\log\frac{|I|}{|J|-|I|}|\right)^p.
    \]
\end{lem}
\begin{proof}
    Let $\mu_i=\mu_{X_i}$. If $x$ is the common endpoint, label the arcs
    so that $A$ follows $x$ and $B$ precedes it. Conditional on $\mu_1$,
    the masses $d_+=\mu_1(A)$ and $d_-=\mu_1(B)$ are fixed, while
    $V=U+F_{X_1}(x)\pmod1$ is uniform and independent of $\mu_2$.
    Lemma~\ref{lem:inverse rotation} therefore gives
    \[
        \m E\left[\left|\log\frac{\nu_H(A)}{\nu_H(B)}\right|^p {\vert}  \mu_1 \right] \leq C_{\g,p} \left(1+ |\log\frac{\mu_1(A)}{\mu_1(B)}|\right)^p,
    \]
    then from Lemma~\ref{lem:log GMC moment},
    \[
        \m E\left[\left|\log\frac{\nu_H(A)}{\nu_H(B)}\right|^p\right] \leq C_{\g,p} \left(1+ |\log\frac{|A|}{|B|}|\right)^p.
    \]
    The case of $Z_{x,r}$ is a special case when $|A|=|B|$.
    For the last assertion, use $\log(1+e^u)\le\log2+|u|$,
    Minkowski's inequality, and the identity
    \[
        \log\frac{\nu_H(I)}{\nu_H(J)} =-\log(1+ \frac{\nu_H(J\backslash I)}{\nu_H(I)}),
    \]
    \[
        \m E\left[\left|\log\frac{\nu_H(I)}{\nu_H(J)}\right|^p\right]\leq \left(\log2+C_{\g,p} \left(1+ |\log\frac{|I|}{|J|-|I|}|\right)\right)^p.
    \]
\end{proof}
For any non-atomic and fully-supported probability measure $\nu$ on the unit circle, define for any $0 < s \leq 1$ and $0 < u < 1$
\[
    F(s) = F^\nu(s)= \int_{\m S^1} \log \nu([x,x+s]) \dd x,
\]
\[
    K(u) = K^\nu(u)= \int_{\m S^1} \frac{\nu(\dd y)}{\nu([y-u,y])}.
\]
\begin{lem}\label{lem:c_0}
    We have the following:
    \begin{itemize}
        \item For $ 0 < a \leq b \leq 1$,
        \[
            F(b)-F(a) = \int_a^b K(u)\dd u.
        \]
        \item For $0 < u <1/2$,
        \[
            K(u) \geq \frac{\log2}{2u}.
        \]
        \item For $ 0 < a \leq b \leq 1$
        \[
            F(b)-F(a) \geq c_0 \log\frac{b}{a}, \qquad c_0 = \frac{\log2}{2}.
        \]
    \end{itemize}
\end{lem}
\begin{proof}
    By Fubini's theorem,
    \[
         F(b)-F(a) = \int_{\m S^1} \dd x \int_a^b\frac{\nu(x+\dd s)}{\nu([x,x+s])} =\int_a^b \dd s \int_{\m S^1} \frac{\nu(\dd y)}{\nu([y-s,y])} =\int_a^b K(s) \dd s.
    \]
    Let $n = \lfloor1/u\rfloor$ and divide the circle into $n$ equal arcs $(A_j)_{j=1}^n$ with $(p_j = \nu(A_j)>0)_{j=1}^n$,
    then
    \[
         \int_{A_j} \frac{\nu(\dd y)}{\nu([y-u,y])} \geq \int_0^{p_j} \frac{\dd t}{p_{j-1}+t} = \log\left(1+\frac{p_j}{p_{j-1}}\right),
    \]
    therefore by Jensen's inequality with the function $v \mapsto \log(1+e^v)$
    \[
        K(u) \geq \sum_{j=1}^n  \log\left(1+\frac{p_j}{p_{j-1}}\right) \geq n\log2 \geq \frac{\log 2}{2u},
    \]
    where we use
    \[
        \sum_{j=1}^n  \log\left(\frac{p_j}{p_{j-1}}\right) = 0.
    \]
    When $b \leq 1/2$, we obtain the third statement by integration from the first two statements. When $a \geq 1/2$, let $u_b(x) = \nu([x+b,x+1])$,
    \begin{align*}
        F(b)-F(a) =& \int_{\m S^1}\log\frac{1-u_b(x)}{1-u_a(x)} \dd x = \int_{\m S^1} \dd x \int_{u_b}^{u_a} \frac{\dd v}{1-v} \\
        \geq&  \int_{\m S^1} \dd x (u_a(x)-u_b(x)) = b-a\geq c_0 \log\frac{b}{a}.
    \end{align*}
        
    For $a \leq 1/2 \leq b$, split the interval at $1/2$ and add the two estimates.
\end{proof}
\begin{lem}\label{lem:lowerbound}
    For $\varphi \in \Homeo_+(\m S^1)$, $0<r\le 1/2$,
    $I_+=\varphi([x,x+r])$ and $I_-=\varphi([x-r,x])$, we have
    \[
        \d_\varphi(x,r) s_\varphi(x,r) \geq c_0^2  \d_\varphi(x,r)^2,
    \]
    where $c_0 = \log 2/2$ and
    \[
        s_\varphi(x,r) = \m E_\varphi [Z_{x,r}].
    \]
\end{lem}
\begin{proof}
    It suffices to assume $|I_+| \geq |I_-|$.
    Recall that
    \[
        H=Q \circ V  \circ F,
    \]
    let $\mu_1$, $\mu_2$ denote the measure of $F$, $Q$, respectively.
    Recall the definition
    \[
        F(s) = F^\nu(s)= \int_{\m S^1} \log \nu([x,x+s]) \dd x,
    \]
    now we define a random distribution function
    \[
        \hat F_2(w) : = \m E_{\mu_2} [F^{\mu_2}(w)],\qquad 0< w \leq 1,
    \]
    and a deterministic distribution function
    \[
        \tilde F_2(l):= \m E_{\mu_1} [\hat F_2(\mu_1(J_l))],
    \]
    where $J_l$ is any fixed arc with length $l$. This definition is independent of the position of the arc, since the expectation is rotation invariant. Thus
    \[
        \m E[\log \nu_H(J)] =   \m E_{\mu_1} [\hat F_2(\mu_1(J))] = \tilde F_2(|J|),
    \]
    Applying the preceding lemma twice gives
    \begin{align*}
        s_\varphi(x,r) =& \tilde F_2(|I_+|)-\tilde F_2(|I_-|)\\ 
        \geq& c_0 \m E_{\mu_1} [\log\frac{\mu_1(I_+)}{\mu_1(I_-)}]= c_0 \m E_{\mu_1}[F^{\mu_1}(|I_+|)-F^{\mu_1}(|I_-|)] \\
        \geq & c_0^2 \log \frac{|I_+|}{|I_-|}.
    \end{align*}

\end{proof}
\begin{thm}\label{thm:A}
If $\varphi \in \QIR$, then $\varphi \in \QS(\m S^1)$.
\end{thm}
\begin{proof}
    Suppose $\varphi \notin \QS(\m S^1)$. Then there exist
    \[
        r_n \downarrow 0, x_n\rightarrow x \in \m S^1, d_n:= |\d_\varphi(x_n,r_n) |\geq n \rightarrow \infty, l_n:= \nu_\varphi([x_n-r_n,x_n+r_n]) \rightarrow 0.
    \]
    Let $I_+=\varphi([x_n,x_n+r_n])$, $I_-=\varphi([x_n-r_n,x_n])$,
    $a_n=|I_+|$, $b_n=|I_-|$, $\d_n=\d_\varphi(x_n,r_n)$, and $Z_n=Z_{x_n,r_n}(\cdot)$.
    Lemmas~\ref{lem:moment of Z} and~\ref{lem:lowerbound} give:
    \begin{itemize}
        \item $\norm{Z_n}_{L^2(\m P)}\leq C_{\g,2}$, so
        $\norm{Z_n/d_n}_{L^2(\m P)}\to0$.
        \item 
            \[
                \norm{Z_n}_{L^2(\m P_\varphi)}
                =\norm{\log \frac{\nu_\cdot(I_+)}{\nu_\cdot(I_-)}}_{L^2(\m P)}
                \leq C_{\g,2}(1+|\log(a_n/b_n)|),
            \]
            with $|\log(a_n/b_n)|=d_n$.
        \item $\mathrm{sign}(\d_n) \m E_{\varphi} [Z_n] \geq c_0^2 d_n$.
    \end{itemize}
    Define the measurable function $X_n = \mathrm{sign}(\d_n)  Z_n /d_n$. Let $D = \dd \m P_\varphi /\dd \m P \in L^1(\m P)$.
    From the first claim, $X_n \rightarrow 0$ in $\m P$-probability, hence also in $\m P_\varphi$-probability by the following estimate, letting $n \to \infty$ and $M \to \infty$
    \[
        \m P_\varphi(|X_n|>\e) \leq M\m P(|X_n|>\e)+ \int_{D>M} D \dd \m P,
    \]
    for any $\e>0$. The second claim bounds $X_n$ in $L^2(\m P_\varphi)$.
    Uniform integrability gives $\m E_\varphi X_n\to0$, contradicting
    $\m E_\varphi[X_n]\geq c_0^2$.
\end{proof}
\subsection{Small-scale cancellation and symmetry}\label{subsec:subcritical-sym}

Once quasisymmetry is known, the distortion is bounded, so the previous
normalization no longer detects a failure of symmetry. We instead show
that the adjacent-arc observable converges weakly to zero along shrinking
arcs. Localized smooth flows give the required cancellation while fixing
finitely many observations away from the limiting point. Combining this
weak convergence with the positive lower bound from the preceding
subsection proves Theorem~\ref{thm:B}.

Fix a nonnegative $\zeta\in C_c^\infty((-1,1))$ equal to $1$ on
$[-1/2,1/2]$. For a center $x_0$ and $0<l<1/2$, let $\k_t$ be the
flow of the periodized vector field $w_l(x)=l\zeta((x-x_0)/l)$.
It fixes the complement of $[x_0-l,x_0+l]$ and translates at speed $l$
while the trajectory stays in $[x_0-l/2,x_0+l/2]$.

\begin{lem}\label{lem:flow estimate}
For every $1\le p<\infty$ there are $t_p,C_p>0$, uniform in the center
and scale, such that for $|t|\le t_p$ the density
$D_t=\dd\m P_{\k_t}/\dd\m P$ exists and
$\|D_t-1\|_{L^p(\m P)}\le C_p|t|$.
If $\lambda_s$ is another such flow at the same scale with disjoint
support and $D_{t,s}=\dd\m P_{\k_t\circ\lambda_s}/\dd\m P$, then
\begin{equation}\label{eq:subcritical-mixed-density}
 \|D_{t,s}-D_{t,0}-D_{0,s}+1\|_{L^p(\m P)}\le C_p|ts|
 \qquad (|t|,|s|\le t_p).
\end{equation}
\end{lem}
\begin{proof}
Let $\Phi=\k_t\circ\lambda_s$, with $\lambda_s=\operatorname{id}$ in the
one-flow case. The GMC coordinate rule \cite[Lemma~3.8]{F1} gives the
normalized source field and its cumulative function as
\[
 \begin{gathered}
 X_1^\Phi=\pi_0(X_1\circ\Phi+Q_\gamma\log\Phi'),\qquad
 Q_\gamma=\frac2\gamma+\frac\gamma2,\\
 F_{X_1^\Phi}(x)=F_{X_1}(\Phi(x))-F_{X_1}(\Phi(0)).
 \end{gathered}
\]
Here $F_{X_1^\Phi}$ is formed using the same GMC functional as $F_{X_1}$.
At regularization scale $\epsilon$, the subtracted variance is that of
the regularized reference field $X_{1,\epsilon}$.
Replacing $U$ by $U+F_{X_1}(\Phi(0))$ modulo $1$ preserves the
conditional Haar law given both fields. Hence, if $L_{t,s}$ is the
Gaussian density of $X_1^\Phi$
relative to $X_1$, then $D_{t,s}(H)=E[L_{t,s}(X_1)\mid H]$.

Use the Cameron--Martin coordinates of
Appendix~\ref{subsec:crit-gaussian-proof}. The covariance and mean are
$S_{t,s}=C_{t,s}C_{t,s}^*$ and $b_{t,s}=Q_\gamma\pi_0\log\Phi'$,
where $C_{t,s}f=\pi_0(f\circ\Phi)$. Because the flows commute,
$C_{t,s}=C_{t,0}C_{0,s}$, and disjointness gives
$b_{t,s}=b_{t,0}+b_{0,s}$. Let $L_i$ be the two composition generators
and $A_i=L_i+L_i^*$. As in the appendix,
\[
 \partial_s S_{t,s}=C_{t,s}A_2C_{t,s}^*,\qquad
 \partial_t\partial_s S_{t,s}
 =C_{t,s}\bigl(L_1A_2+(L_1A_2)^*\bigr)C_{t,s}^*.
\]
The estimate \eqref{eq:crit-flow-LA} applies to $L_1A_2$: its proof
uses only the Fourier decay of the two vector fields at their common
scale, and translation leaves the absolute Fourier coefficients
unchanged. Thus the first and mixed covariance derivatives are
uniformly Hilbert--Schmidt on bounded time intervals. The corresponding
mean derivatives are uniformly bounded by \eqref{eq:crit-flow-bounds}
and the mixed derivative vanishes. In particular,
$\|S_{t,s}-I\|_{\mathrm{HS}}\le C(|t|+|s|)$.

For the Fourier compressions of these Gaussian laws, write $L_{t,s,N}$
for the finite-dimensional density. Choose $q>p$ and then $t_p$ small
enough that $qI-(q-1)S_{t,s}$ is uniformly bounded below by a positive
constant on $[-t_p,t_p]^2$. The determinant moment calculation in
Appendix~\ref{subsec:crit-gaussian-proof} now gives
$\sup_{t,s,N}\|L_{t,s,N}\|_{L^q}<\infty$. Differentiating
\eqref{eq:crit-log-gaussian}, the first and mixed derivatives of
$\log L_{t,s,N}$ are Gaussian polynomials of degree at most two with
uniform bounds in every finite $L^a$, by the covariance and mean estimates.
The product rule and H\"older's inequality therefore give uniform
$L^p$ bounds on $\partial_tL_{t,s,N}$ and
$\partial_t\partial_sL_{t,s,N}$. Integration on the time interval
and on the rectangle yields the bounds $C_p|t|$ and $C_p|ts|$ for the
first and mixed density increments. For each $t,s$, the compressions
form a martingale under the reference law. The uniform $L^q$ bound
gives an $L^p$ limit, which is $L_{t,s}$ by testing against bounded
cylinder functions. Conditional expectation given $H$ contracts
$L^p$ and gives the stated estimates.
\end{proof}
\begin{lem}\label{lem:subcritical-weak-cancellation}
    If $r_n \downarrow 0$, $x_n\rightarrow x \in \m S^1$, then for any $B \in L^2(\m P)$, let $Z_n =  Z_{x_n,r_n}$
    \[
        \m E[BZ_n  ] \to 0.
    \]
\end{lem}
\begin{proof}
    A circle homeomorphism is determined by its values on a countable dense
set such as $D = \m Q/\m Z \backslash \{x\}$. Bounded functions of
finitely many evaluations on $D$ are dense in $L^2(\m P)$. Given $B$ and
$\vare$, choose such a function $A$, depending on the values at
$(y_j)_{j=1}^m$, with $\norm{A-B}_{L^2(\m P)} \leq  \vare$:
    \[
        |\m E[BZ_n]| \leq |\m E[AZ_n]| + \norm{A-B}_{L^2(\m P)}\sup_n \norm{Z_n}_{L^2(\m P)} \leq  |\m E[AZ_n]| + \vare C_{\g,2}.
    \]
    It therefore suffices to prove $\m E[AZ_n] \to 0$ for each $A$.
    Use the flow of Lemma~\ref{lem:flow estimate} centered at $x$, with
    $l$ so small that $J=[x-l,x+l]$ avoids $\{y_1,\ldots,y_m\}$.
    For large $n$, the arc $[x_n-r_n,x_n+r_n]$ lies in
    $[x-l/2,x+l/2]$. Set $\k_n=\k_{t_n}$ for $t_n=r_n/l$; then
    \[
        \k_n([x_n-r_n,x_n]) = [x_n, x_n + r_n], \qquad \k_n|_{\m S^1 \backslash J} = \text{id}.
    \]
    Let
    \[
        G_n(H) = \log \nu_H([x_n-r_n,x_n]) \in L^2(\m P),
    \]
    as we have
    \[
        \norm{G_n}_{L^2(\m P)}  \leq C_{\g,2}(1+ \log(1/r_n)).
    \]
    using Lemma~\ref{lem:moment of Z} with $J = \m S^1$.
    then
    \[
        Z_n(H) = G_n(H \circ \k_n) - G_n(H),
    \]
    and
    \[
        A (H \circ \k_n) = A(H).
    \]
    Therefore, from Lemma~\ref{lem:flow estimate} with $D_n = \dd \m P_{\k_n}/\dd \m P$,
    \[
        \m E[AG_n(\cdot \circ \k_n)] = \m E[D_n AG_n].
    \]
    Therefore, from Lemma~\ref{lem:flow estimate}
    \begin{align*}
        |\m E[AZ_n] | =&|\m E[AG_n(D_n-1)]| \\
        \leq& \norm{A}_{L^\infty(\m P)} \norm{G_n}_{L^2(\m P)} \norm{D_n-1}_{L^2(\m P)} \\
        \leq & C_{J,\g} \norm{A}_{L^\infty(\m P)} \frac{r_n}{l}(1+ \log(1/r_n)) \to 0.
    \end{align*}
    
\end{proof}
\begin{thm}\label{thm:B}
    If $\varphi \in \QIR \cap \QS(\m S^1)$, then $\varphi \in \Sy(\m S^1)$.
\end{thm}
\begin{proof}
    Suppose $\varphi \notin \Sy(\m S^1)$. Then there exist $\vare > 0 $ and
    \[
        r_n \downarrow 0, x_n\rightarrow x \in \m S^1, d_n:= |\d_\varphi(x_n,r_n) |\geq \vare, l_n:= \nu_\varphi([x_n-r_n,x_n+r_n]) \rightarrow 0.
    \]
    For any $x \in \m S^1$ and $r \in (0,1/2]$, let $I_+ =\varphi([x,x+r])$ and $I_- =\varphi([x-r,x])$.
    \[
        Z_{x,r}(H \circ \varphi)  = \log \frac{\nu_H(I_+)}{\nu_H(I_-)},
    \]
    With $K_\varphi=\exp(\sup_{x,r}|\d_\varphi(x,r)|)<\infty$,
    \[
        \frac{|I_+|}{|I_-|} = \exp(\d_\varphi(x,r)) \in [K_\varphi^{-1},K_\varphi].
    \]
    Therefore,
    \[
        \sup_{x \in \m S^1,r \in (0,1/2]}\norm{Z_{x,r}(H \circ \varphi )}_{L^2(\m P)}
        \leq C_{\g,2}(1+\log K_\varphi) =: C_{\g,K_\varphi}.
    \]
    In particular, for any $x \in \m S^1$ and $r \in (0,1/2]$
    \[
        \m E_\varphi[Z_{x,r}^2] \leq C_{\g,K_\varphi}^2.
    \]
    Let $Z_n=Z_{x_n,r_n}(\cdot)$.
    With $D = \dd \m P_\varphi / \dd \m P$,
    \[
        \m E_\varphi [Z_n] = \m E [D Z_n] = \m E[D_M Z_n] + \m E[(D-D_M)Z_n],
    \]
    where $M > 0$ and
    \[
        D_M = D \wedge M \in L^\infty(\m P) \subset L^2(\m P).
    \]
    \[
        |\m E[(D-D_M)Z_n] | \leq \m E[D\1_{\{D > M\}}|Z_n|] = \m E_\varphi[\1_{\{D > M\}}|Z_n|] \leq   C_{\g,K_\varphi} \m P_\varphi(D > M)^{1/2}.
    \]
    Therefore,
    \[
        \limsup_{n \to \infty} |\m E_\varphi[Z_n]| \leq 0+C_{\g,K_\varphi} \m P_\varphi(D > M)^{1/2} \xrightarrow{M \to \infty } 0
    \]
    which contradicts
    \[
        |\m E_\varphi[Z_n]| \geq c_0^2 \vare
    \]
    for each $n$ by Lemma~\ref{lem:lowerbound}.
    
\end{proof}
\subsection{Multiscale estimates and WP regularity}\label{subsec:subcritical-wp}

Symmetry makes the distortion small at each individual scale, but WP
regularity requires its square to be integrable over positions and
scales. We therefore integrate the observable against the deterministic
distortion. The midpoint decomposition separates the resulting statistic
into a centered term and a nonnegative displacement. Lemma~\ref{lem:Bessel}
controls the centered terms, and Lemma~\ref{lem:checkpoint} supplies a
subsequence on which the displacement is uniformly integrable after
normalization. These estimates yield the contradiction proving
Theorem~\ref{thm:C}. We introduce the statistics first and then establish
the geometric and covariance estimates needed for these two lemmas.

For $j \in \m N_{+}$, let
\[
    t_j =2^{-j},\quad  C_j = \m S^1 \times [t_{j+1}, t_j), \quad e_j = \int_{C_j} \d_\varphi^2 \dd \eta.
\]
For $N \in  \m N_{+}$, let
\[
    \O_N =\cup_{j=1}^N C_j =\m S^1 \times [t_{N+1},1/2), \quad E_N = 1+ \sum_{j=1}^N e_j,
\]
and
\[
    T_N(H) = \int_{\O_N} \d_\varphi Z_{x,r}(H ) \dd \eta.
\]
Given $x \in \m S^1$ and $r \in (0,1/2]$, let $\varPhi$ denote the lift of $\varphi$
\[
    M = \frac{\varPhi(x-r)+\varPhi(x+r)}{2}, \quad L= \frac{\varPhi(x+r)-\varPhi(x-r)}{2},\quad C= \varPhi(x),
\]
recall by definition
\[
    Z_{M,L}(H) =   \log\frac{\nu_H([M,M+L])}{\nu_H[M-L,M]},
\]
let
\[
    B_N(H) = \int_{\O_N} \d_\varphi Z_{M,L}(H) \dd \eta,
\]
and
\[
    A_N(H)=T_N(H\circ \varphi)-B_N(H).
\]

\begin{lem}\label{lem:decomposition estimate}
    With $R_\varphi(H)=R_{\varphi,x,r}(H) \geq 0$ we have
    \[
        Z_{x,r}(H \circ \varphi) = Z_{M,L}(H) + \text{sign}(\d_\varphi) R_\varphi(H),
    \]
    \[
        A_N = \int_{\O_N} |\d_\varphi| R_\varphi \dd \eta \geq 0.
    \]
    Moreover, we have the estimate
    \[
        \m E[A_N] \geq c_0^2(E_N-1),
    \]
    and
    \[
        \m E[B_N^2] \leq C_{\g,\varphi} E_N.
    \]
\end{lem}
 \begin{proof}
    Recall the definition of $R_\varphi$ with $\nu= \nu_H$
    \[
        R_\varphi(x,r) = \int_{\min(C,M)}^{\max(C,M)} \left(\frac{1}{\nu([A,z])}+ \frac{1}{\nu([z,B])}\right) \nu(\dd z),
    \]
    where
    \[
        A = \varPhi(x-r),\quad B=\varPhi(x+r),
    \]
    \[
        M = \frac{\varPhi(x-r)+\varPhi(x+r)}{2}, \quad L= \frac{\varPhi(x+r)-\varPhi(x-r)}{2},\quad C= \varPhi(x).
    \]
    Since $\m E[Z_{M,L}]=0$ by rotation invariance, integrating
    Lemma~\ref{lem:lowerbound} gives the lower bound for $\m E[A_N]$.
    The last inequality follows from Lemma~\ref{lem:Bessel}.
 \end{proof}
For a proper arc $J$, let $\ell_J$ be the affine interpolant of $\Phi$
on $J$, let $L_J=|\varphi(J)|$, and set
\[
 \varepsilon_J=\frac{\|\Phi-\ell_J\|_{L^\infty(J)}}{L_J}.
\]
Dyadic bisection and monotonicity show that symmetry implies
$\varepsilon_J\to0$ uniformly as $|J|\to0$. Choose $j_0\ge4$ so that
$\varepsilon_J\le1/80$ for $|J|\le4t_{j_0}$ and
$9t_{j_0}\le r_0$ in Lemma~\ref{lem:sym ab} with
$\alpha=7/8$ and $\beta=9/8$.

\begin{lem}
For $\varphi\in\Sy(\m S^1)$ there are deterministic $S_N$ with
$\|A_N\|_{L^2(\m P)}\le C_{\gamma,\varphi}S_N$.
\end{lem}
\begin{proof}
For $j\ge j_0$, set $t=t_j$ and $h=t/100$. For $0\le k<h^{-1}$, put
\begin{equation}\label{eq:Jjk}
 \begin{aligned}
 x_k&=kh,\qquad B_k=[x_k-h/2,x_k+h/2]\times[t/2,t],\\
 J&=J_{j,k}=[x_k-2t,x_k+2t].
 \end{aligned}
\end{equation}
In the circle coordinate the blocks $B_k$ cover $C_j$ with disjoint
interiors. Set $s_j=\sum_k\varepsilon_{J_{j,k}}^2$ and
$S_N=1+\sum_{j=j_0}^Ns_j$.

Fix a block and write $\nu=\nu_H$. For $(x,r)\in B_k$, put
$A=\Phi(x-r)$, $B=\Phi(x+r)$ and use $M,C$ as above. Both $M$ and $C$
are within $\varepsilon_J L_J$ of $\ell_J(x)$. Since
$\ell_J'=L_J/(4t)$, their intervening segment is contained in
$[z_-,z_+]$, where
\[
 q=\ell_J(x_k),\qquad
 z_\pm=q\pm(\varepsilon_J+1/800)L_J.
\]
Let $A_0=\Phi(x_k+h/2-t/2)$ and $B_0=\Phi(x_k-h/2+t/2)$,
so $A\le A_0$ and $B_0\le B$. The distances $z_--A_0$ and $B_0-z_+$ are
both at least $(1/8-1/400-2\varepsilon_J)L_J\ge L_J/12$. The length
$s=z_+-z_-$ is less than $L_J/12$.

On $[z_-,z_+]$ define
$G_J(z)=\nu([A_0,z])^{-1}+\nu([z,B_0])^{-1}$. Its integral is
\[
 W_J=\int_{z_-}^{z_+}G_J(z)\,\nu(\dd z)
 =\log\frac{\nu([A_0,z_+])}{\nu([A_0,z_-])}
  +\log\frac{\nu([z_-,B_0])}{\nu([z_+,B_0])}.
\]
For the first logarithm, choose the terminal subarc $I'$ of
$[A_0,z_-]$ of length $s$. The logarithm is bounded above by
$\log 2+|\log(\nu([z_-,z_+])/\nu(I'))|$. The adjacent equal-length arcs
give a uniform $L^2$ bound by Lemma~\ref{lem:moment of Z}.
Using the initial subarc of $[z_+,B_0]$ gives the same bound for the
second logarithm. Hence $\|W_J\|_{L^2(\m P)}\le C_\gamma$.

The lengths $h_\pm=|\Phi(x\pm r)-\Phi(x)|$ differ from
$rL_J/(4t)$ by at most $2\varepsilon_JL_J$. Thus
$\min(h_+,h_-)\ge L_J/10$ and $|\delta_\varphi|\le40\varepsilon_J$.
For fixed $r,z$, the condition that $z$ lies between $C$ and $M$
implies $|z-\ell_J(x)|\le\varepsilon_J L_J$ and hence restricts $x$
to a set of length at most $8\varepsilon_Jt$.
The Stieltjes formula for $R_\varphi$ and the bounds $A\le A_0$,
$B_0\le B$ allow us to dominate the kernel by $G_J$ before changing
the order of integration. Tonelli's theorem now gives
\[
 \begin{aligned}
 \int_{B_k}|\delta_\varphi|R_\varphi\,\dd\eta
 &\le (40\varepsilon_J)(8\varepsilon_Jt)
       W_J\int_{t/2}^{t}\frac{\dd r}{r^2}\\
 &=320\varepsilon_J^2 W_J.
 \end{aligned}
\]
Sum over $k,j$ and use the triangle inequality in $L^2$. On the
finitely many coarse layers, the two image arms have lengths
uniformly bounded away from zero. Lemma~\ref{lem:moment of Z} and
Minkowski's inequality bound their total contribution in $L^2$.
This proves the assertion.
\end{proof}

\begin{lem}\label{lem:checkpoint}
If $\varphi\in\Sy(\m S^1)$, there is a deterministic subsequence
$(N_j)_{j\ge1}$ on which $S_N\le C_\varphi E_N$ and
$\|A_N\|_{L^2(\m P)}\le C_{\gamma,\varphi}E_N$.
\end{lem}
\begin{proof}
For $I$ of length $L\le9t_{j_0}$, let $\mc D_n I$ consist of its $2^n$
dyadic subintervals of length $2^{-n}L$, starting with $\mc D_0 I=\{I\}$.
For $K\in\mc D_n I$ write $m_K$ for its midpoint, $r_K$ for its
half-length, and $\psi_K$ for the hat function on $K$ with value one
at $m_K$. Successive dyadic linear interpolants give
\[
 \Phi-\ell_I=\sum_{n\ge0}\sum_{K\in\mc D_n I}c_K\psi_K,
 \quad
 c_K=\Phi(m_K)-\tfrac12\bigl(\Phi(m_K-r_K)+\Phi(m_K+r_K)\bigr).
\]
The series converges uniformly by uniform continuity of $\Phi$.
The adjacent length ratio gives
\begin{equation}\label{eq:estimate of cK}
 \frac{|c_K|}{|\varphi(K)|}
 =\frac12\left|\tanh\!\left(\frac{\delta_\varphi(m_K,r_K)}2\right)\right|
 \le\frac{|\delta_\varphi(m_K,r_K)|}4 .
\end{equation}
At every point at most one hat per generation is nonzero.
Cauchy--Schwarz with weights $2^{-n/4}$ and
Lemma~\ref{lem:sym ab} with $\alpha=7/8$ give
\[
 \varepsilon_I^2
 \le C_\varphi\sum_{n\ge0}2^{-3n/2}
       \sum_{K\in\mc D_n I}\delta_\varphi(m_K,r_K)^2.
\]

Fix $j\ge j_0$ and $t=t_j$. For $J=J_{j,k}$ and
$J'=[a,a+L]\supset J$ with $L\in[8t,9t]$, the affine difference
$\ell_{J'}-\ell_J$ is bounded on $J$ by its endpoint values.
Consequently $\varepsilon_J\le2\varepsilon_{J'}L_{J'}/L_J$,
where $L_{J'}/L_J\le C_\varphi$ by Lemma~\ref{lem:sym ab}.
For each $J$ the admissible parameters $(a,L)$ have measure
$\int_{8t}^{9t}(L-4t)\,\dd L=9t^2/2$, and a given $J'$ contains at
most $501$ of the arcs $J_{j,k}$. Averaging first and then summing over
$k$ shows that
\[
 s_j\le C_\varphi t^{-2}
        \int_{8t}^{9t}\int_{\m S^1}\varepsilon_{[a,a+L]}^2
        \,\dd a\,\dd L.
\]
Inserting the dyadic bound, the integral over $a$ of the sum at level
$n$ equals $2^n\int_{\m S^1}\delta_\varphi(x,r)^2\,\dd x$ with
$r=L/2^{n+1}$. The change of variables $\dd L=2^{n+1}\dd r$ gives
a factor $2^{2n+1}t^{-2}\asymp r^{-2}$. The resulting $r$ interval
$[4t2^{-n},(9/2)t2^{-n}]$ belongs to the scale band $C_{j+n-3}$.
Thus
\[
 s_j\le C_\varphi\sum_{n\ge0}2^{-3n/2}e_{j+n-3},
 \qquad
 S_N\le C_\varphi\sum_{n\ge0}2^{-3n/2}E_{N+n}.
\]
Since $\delta_\varphi$ is bounded, $E_N=O(2^N)$. The positive
sequence $E_N2^{-5N/4}$ tends to zero, so it attains its maximum on
every tail. Choose successively increasing indices $N_j$ at such
maxima. Then $E_{N_j+n}\le2^{5n/4}E_{N_j}$ for $n\ge0$. Substitution
in the last bound gives $S_{N_j}\le C_\varphi E_{N_j}$ and the
preceding lemma gives the desired $L^2$ bound.
\end{proof}

\begin{lem}\label{lem:Bessel}
    For $f \in L^2(\m S^1 \times (0,1/2], \dd \eta)$, there exists $C_\g$,
    \[
        \norm{\int_{\m S^1 \times (0,1/2]} f(x,r)Z_{x,r}(\cdot) \dd \eta(x,r)}_{L^2(\m P)}^2 \leq C_\g \int_{\m S^1 \times (0,1/2]} |f(x,r)| ^2\dd \eta(x,r).
    \]
    And for $\varphi \in \Sy(\m S^1)$, there exists $C_{\g,\varphi}$,
    \[
        \norm{\int_{\m S^1 \times (0,1/2]} f(x,r)Z_{M(x,r),L(x,r)}(\cdot) \dd \eta(x,r)}_{L^2(\m P)}^2 \leq C_{\g,\varphi} \int_{\m S^1 \times (0,1/2]} |f(x,r)| ^2\dd \eta(x,r).
    \]
\end{lem}
\begin{proof}
    Define
    \[
        \G((x,r),(y,R)) = \m E[Z_{x,r} Z_{y,R}],
    \]
    it suffices to consider bounded and compactly supported $f$,
    \[
        \iint_{(\m S^1 \times (0,1/2])^2} f(x,r)  \G((x,r),(y,R)) f(y,R) \dd \eta(x,r) \dd \eta(y,R)= <f,\mc C f>_{L^2(\eta)},
    \]
    where
    \[
        \mc C f (x,r)=    \int_{\m S^1 \times (0,1/2]}   \G((x,r),(y,R)) f(y,R) \dd \eta(y,R).
    \]
    Decompose $\mc C$ on $C_j \times C_k$,
    \[
        \mc C_{jk} f (x,r)=    \1_{(x,r)\in C_j}\int_{C_k}   \G((x,r),(y,R)) f(y,R) \dd \eta(y,R) ,
    \]
    and similarly decompose $f$
    \[
         f_j =f\1_{\cdot \in C_j},
    \]
    then we have
    \begin{align*}
        <f,\mc C f>_{L^2(\eta)} =& \sum_{jk} \iint_{C_j \times C_k} f(x,r)  \G((x,r),(y,R)) f(y,R) \dd \eta(x,r) \dd \eta(y,R)\\
        = & \sum_{jk}  <f_j,\mc C_{jk} f_k>_{L^2(\eta)} 
    \end{align*}
    Define
    \[
        A_{jk}=A(\mc C_{jk}) = \sup_{(x,r) \in C_j} \int_{C_k} | \G((x,r),(y,R))| \dd \eta(y,R),
    \]
    \[
         B_{jk}=B(\mc C_{jk}) =  \sup_{(y,R) \in C_k} \int_{C_j}  | \G((x,r),(y,R))| \dd \eta(x,r) ,
    \]
    note that by Cauchy--Schwarz inequality,
    \[
        <f_j,\mc C_{jk} f_k>_{L^2(\eta)} \leq \sqrt{A_{jk}B_{jk}}\norm{f_j}_{L^2(\eta)}\norm{f_k}_{L^2(\eta)},
    \]
    we claim
    \[
        Q_{jk}=\sqrt{A_{jk}B_{jk}} \leq C_\g  Q(|j-k|),
    \]
    with
    \[
        \sum_{l}|Q(l)| < \infty,
    \]
    then
    \begin{align*}
        <f,\mc C f>_{L^2(\eta)} \leq & \sum_{jk} Q_{jk}\norm{f_j}_{L^2(\eta)}\norm{f_k}_{L^2(\eta)}\\
        \leq & 2\sum_{l} Q(l) \sum_{j-k=l} \norm{f_j}_{L^2(\eta)}\norm{f_k}_{L^2(\eta)}\\
        \leq &  2\sum_{l} Q(l)\sum_j \norm{f_j}_{L^2(\eta)}^2 \\
        \leq & C_{\g}  \norm{f}_{L^2(\eta)}^2.
    \end{align*}
    By symmetry, take $j\ge k$. For $j>k$ we have $r\le R$. On the
    diagonal the radii differ by at most a factor of two, so the
    following integral bounds hold after exchanging the radii if needed.
    Set
    \[
        L = 1+ |\log R/r|.
    \]
    To prove the claim, apply Lemma~\ref{lem:Gamma estimate}
    \[
        \sup_{x \in \m S^1}\int_{\m S^1}  | \G((x,r),(y,R))| \dd y \leq C_\g rL^2,
    \]
    \[
        \sup_{y \in \m S^1}\int_{\m S^1}  | \G((x,r),(y,R))| \dd x \leq C_\g rL^2.
    \]
    Then we have
    \[
        A_{jk} =   \sup_{(x,r) \in C_j} \int_{C_k} | \G((x,r),(y,R))| \dd \eta(y,R) \leq C_\g \int_{t_{k+1}}^{t_k} \frac{\dd R}{R^2} rL^2 \leq C_\g 2^{-(j-k)} (1+j-k)^2,
    \]
    and
    \[
        B_{jk} =   \sup_{(y,R) \in C_k} \int_{C_j} | \G((x,r),(y,R))| \dd \eta(x,r) \leq C_\g \int_{t_{j+1}}^{t_j} \frac{\dd r}{r^2} rL^2 \leq C_\g  (1+j-k)^2,
    \]
    therefore,
    \[
        Q_{jk}=\sqrt{A_{jk}B_{jk}} \leq C_\g  Q(|j-k|)=C_\g 2^{-|j-k|/2}(1+|j-k|^2).
    \]
    If we replace $\G$ by
    \[
        \tilde \G((x,r),(y,R)) = \m E[Z_{M(x,r),L(x,r)} Z_{M(y,R),L(y,R)}],
    \]
    by the same argument, it suffices to estimate
    \[
        \tilde Q_{jk} =  \sqrt{\tilde A_{jk}\tilde B_{jk}},
    \]
     \[
        \tilde A_{jk}=A(\tilde {\mc C}_{jk}) = \sup_{(x,r) \in C_j} \int_{C_k} | \tilde \G((x,r),(y,R))| \dd \eta(y,R),
    \]
    \[
         \tilde B_{jk}=B(\tilde {\mc C}_{jk}) =  \sup_{(y,R) \in C_k} \int_{C_j}  | \tilde \G((x,r),(y,R))| \dd \eta(x,r) .
    \]
    We claim
    \[
        \tilde Q_{jk}=\sqrt{\tilde A_{jk}\tilde B_{jk}}
        \leq C_{\g,\varphi} 2^{-|j-k|/6},
    \]
    then we have the desired result.
    To prove the claim, apply Lemma~\ref{lem:Tilde Gamma estimate} with $\s=2/3$,
    \[
        \sup_{x \in \m S^1}\int_{\m S^1}  | \tilde \G((x,r),(y,R))| \dd y \leq C_{\g,\varphi} R(r/R)^\s,
    \]
    \[
        \sup_{y \in \m S^1}\int_{\m S^1}  |\tilde  \G((x,r),(y,R))| \dd x \leq C_{\g,\varphi} R(r/R)^\s.
    \]
    Then we have
    \[
        \tilde A_{jk} =   \sup_{(x,r) \in C_j} \int_{C_k} |\tilde  \G((x,r),(y,R))| \dd \eta(y,R) \leq C_{\g,\varphi} \int_{t_{k+1}}^{t_k} \frac{\dd R}{R^2}  R(r/R)^\s\leq C_{\g,\varphi} 2^{(k-j)\s} ,
    \]
    and
    \begin{align*}
        \tilde B_{jk} =  & \sup_{(y,R) \in C_k} \int_{C_j} | \tilde \G((x,r),(y,R))| \dd \eta(x,r) \\
        \leq &  C_{\g,\varphi} \int_{t_{j+1}}^{t_j} \frac{\dd r}{r^2} R(r/R)^\s \\
        \leq &  C_{\g,\varphi} 2^{(k-j)(\s-1)} ,
    \end{align*}
    therefore,
    \[
        \tilde Q_{jk}=\sqrt{\tilde A_{jk}\tilde B_{jk}}
        \leq C_{\g,\varphi} 2^{-|j-k|(2\s-1)/2}.
    \]
\end{proof}
We next estimate the covariance kernel in Lemma~\ref{lem:Bessel}.
Write $\dist$ for circular distance and set
\[
 D=\dist(x,y),\qquad
 d=\min_{z\in\{y-R,y,y+R\}}\dist(x,z).
\]

\begin{lem}\label{lem:Gamma estimate}
For $0<r,R\le1/2$ the kernel satisfies $|\Gamma((x,r),(y,R))|\le C_\gamma$.
There is a constant $C_*$ such that
\begin{align}
 |\Gamma((x,r),(y,R))|
 &\le C_\gamma\frac r d\bigl(1+\log(d/r)\bigr)
 &&\text{if }d\ge C_*r,\label{eq:subcritical-cov-local}\\
 |\Gamma((x,r),(y,R))|
 &\le C_\gamma\frac{rR}{D^2}\bigl(1+\log(D/r)\bigr)
                       \bigl(1+\log(D/R)\bigr)
 &&\text{if }D\ge C_*\max(r,R).\label{eq:subcritical-cov-far}
\end{align}
If $r\le R$, both integrals of $|\Gamma((x,r),(y,R))|$, over $x$ with
$y$ fixed and over $y$ with $x$ fixed, are bounded by
$C_\gamma r(1+\log(R/r))^2$.
\end{lem}
\begin{proof}
Cauchy--Schwarz and Lemma~\ref{lem:moment of Z} give the uniform bound.
For a flow $\kappa_t$ centered at $x$ with scale $\ell$ and support
in $J_x=[x-\ell,x+\ell]$, set
\[
 G_x(h)=\log\frac{\nu_h([x-r,x])}{\nu_h(J_x)}.
\]
When $r\le\ell/2$, Lemma~\ref{lem:moment of Z} gives
$\|G_x\|_{L^4(\m P)}\le C_\gamma(1+\log(\ell/r))$.
Indeed, $J_x\setminus[x-r,x]$ consists of two arcs adjacent to
$[x-r,x]$, each of length comparable to $\ell$. If $u,v$ are their
masses divided by $\nu_h([x-r,x])$, use
$\log(1+u+v)\le\log 3+|\log u|+|\log v|$.
For the same $r\le\ell/2$, the flow at time $t=r/\ell$ carries
$[x-r,x]$ to $[x,x+r]$ within the translation region and preserves $J_x$.
Hence $Z_{x,r}(h)=G_x(h\circ\kappa_t)-G_x(h)$.

For \eqref{eq:subcritical-cov-local}, take $\ell=d/2$. The flow fixes
$y-R,y,y+R$ and therefore preserves $Z_{y,R}$. Changing variables gives
\[
 \Gamma((x,r),(y,R))=\m E[(D_t-1)G_xZ_{y,R}],
\]
where $D_t$ is the density in Lemma~\ref{lem:flow estimate}.
Choose $C_*$ so that $t$ lies in the time interval of that lemma for
$p=4$. H\"older's inequality with exponents $4,4,2$ proves
\eqref{eq:subcritical-cov-local}.

For \eqref{eq:subcritical-cov-far}, use two flows $\kappa_t,\lambda_s$
of scale $\ell=D/4$ and centers $x,y$, with $t=r/\ell$ and $s=R/\ell$.
Their supports are disjoint. Define $G_y$ by replacing $(x,r)$
in $G_x$ with $(y,R)$. Each flow preserves the other logarithmic ratio.
Expanding the product of the two differences and changing variables gives
\[
 \Gamma((x,r),(y,R))
 =\m E[(D_{t,s}-D_{t,0}-D_{0,s}+1)G_xG_y].
\]
Increasing $C_*$ if needed, Lemma~\ref{lem:flow estimate} bounds the mixed
density difference in $L^2(\m P)$ by $C_\gamma |ts|$. Together with
the two $L^4$ bounds this gives \eqref{eq:subcritical-cov-far}.

To integrate, suppose $r\le R$. With either $x$ or $y$ fixed, the
distribution of $d$ under integration in the other variable has
Lebesgue density at most $6$. This follows by partitioning according
to the nearest of the three centers and bounding each contribution
by the density of distance to one point. Thus $d<C_*r$ contributes
$O(r)$. On $d\ge C_*r$ and $D<C_*R$, use $d\le D$ and bound the
contribution by
\[
 C_\gamma\int_{C_*r}^{C_*R}\frac{r}{u}
                      (1+\log(u/r))\,\dd u
 \le C_\gamma r(1+\log(R/r))^2.
\]
On $D\ge C_*R$, use \eqref{eq:subcritical-cov-far} and substitute
$v=D/R$. The bound is at most
\[
 C_\gamma r\int_{C_*}^\infty
       \frac{(1+\log(R/r)+\log v)(1+\log v)}{v^2}\,\dd v
 \le C_\gamma r(1+\log(R/r)).
\]
\end{proof}
We now prove the corresponding estimates after composition by a symmetric map. We begin with deterministic bounds on its distortion of relative arc lengths.
\begin{lem}\label{lem:sym ab}
Let $\varphi\in\Sy(\m S^1)$. For every $0<\alpha<1<\beta$, there
are $c,C,r_0>0$, depending on $\varphi,\alpha,\beta$, such that
\[
 c\left(\frac{|I|}{|J|}\right)^\beta
 \le\frac{|\varphi(I)|}{|\varphi(J)|}
 \le C\left(\frac{|I|}{|J|}\right)^\alpha
 \qquad(0<|I|,\ I\subset J,\ |J|\le r_0).
\]
After adjusting $c,C$, the bounds hold for all $|J|\le1$.
\end{lem}
\begin{proof}
Symmetry allows us to choose $r_0$ so that every half $J'$ of an arc
$J$ with $|J|\le r_0$ satisfies
\[
 2^{-\beta}\le\frac{|\varphi(J')|}{|\varphi(J)|}\le2^{-\alpha}.
\]
Indeed, both half-image ratios tend uniformly to $1/2$ as $|J|\to0$.
Iteration gives the corresponding power bounds on each dyadic
descendant of $J$. Any subarc $I\subset J$ contains such a descendant
of length at least $|I|/4$ and is covered by at most two of length at
most $2|I|$. Monotonicity therefore gives the claimed inequalities,
with $c=4^{-\beta}$ and $C=2^{1+\alpha}$.

If $|J|>r_0$ and $|I|\le r_0/2$, choose $I\subset J'\subset J$ with
$|J'|=r_0/2$ and apply the estimate on $J'$. Compactness bounds
$|\varphi(J')|/|\varphi(J)|$ below uniformly, and $|J'|$ and $|J|$
are comparable. If $|I|>r_0/2$, compactness bounds $|\varphi(I)|$
below and $|I|/|J|$ is bounded away from zero.
\end{proof}
Fix $\alpha=7/8$, $\beta=9/8$ and put $\upsilon=\alpha/\beta=7/9$.
For $0<R\le 1/2$, define
\[
 M_R(y)=\frac{\Phi(y-R)+\Phi(y+R)}2,
 \qquad q_R=\Phi^{-1}\circ M_R.
\]
The increasing maps $M_R,q_R$ commute with translation by one. For
$I=[a,b]$, write $|M_R(I)|=M_R(b)-M_R(a)$. Since
\[
 |M_R(I)|=\tfrac12\bigl(|\varphi(I-R)|+|\varphi(I+R)|\bigr),
\]
the two power bounds of Lemma~\ref{lem:sym ab} hold for the relative
interval lengths under $M_R$, uniformly in $R$. Applying the lower bound
to $M_R$, the upper bound to $\varphi$, and using $M_R\circ q_R^{-1}=\Phi$
gives
\[
 \frac{|q_R^{-1}(I)|}{|q_R^{-1}(J)|}
 \le C_\varphi \left(\frac{|I|}{|J|}\right)^\upsilon
 \qquad(I\subset J,\ |J|\le1).
\]
Also $|q_R(y)-y|\le R$ by monotonicity of $\Phi$, and therefore
$|q_R^{-1}(y)-y|\le R$.

\begin{lem}\label{lem:Tilde Gamma estimate}
Let $\varphi\in\Sy(\m S^1)$ and $0<r\le R\le1/2$. For $s=2/3$, the
kernel
\[
 \widetilde\Gamma((x,r),(y,R))
 :=\Gamma\bigl((M(x,r),L(x,r)),(M(y,R),L(y,R))\bigr)
\]
satisfies
\begin{align*}
 \sup_x\int_{\m S^1}|\widetilde\Gamma((x,r),(y,R))|\,\dd y
 &\le C_{\gamma,\varphi}r^sR^{1-s},\\
 \sup_y\int_{\m S^1}|\widetilde\Gamma((x,r),(y,R))|\,\dd x
 &\le C_{\gamma,\varphi}r^sR^{1-s}.
\end{align*}
\end{lem}
\begin{proof}
Write $(m,l)=(M(x,r),L(x,r))$ and $(M,L)=(M(y,R),L(y,R))$.
All distances below are circular. Let
\[
 D=\dist(x,y),\qquad d=\min_{z\in\{y-R,q_R(y),y+R\}}\dist(x,z).
\]
The images under $\varphi$ of the three points in this minimum are
$M-L,M,M+L$ modulo one. We first show that
\begin{equation}\label{eq:tilde-cov-point}
 |\widetilde\Gamma((x,r),(y,R))|
 \le C_{\gamma,\varphi}\min\{1,(r/d)^s\},
\end{equation}
where the minimum is taken to be $1$ at $d=0$.

If $t=\dist(x,z)>2r$, pull back a shortest arc from $m$ to
$\varphi(z)$ and extend it at its initial end to include $[x-r,x+r]$.
The resulting arc $J$ contains $x,z$, so $|J|\ge t$. The shortest
image arc has length at least $l$ since it must leave
$\varphi([x-r,x+r])$. The extension adds $l$ to its image length, so
$|\varphi(J)|\le2\dist(m,\varphi(z))$. The upper power bound gives
\[
 \frac{l}{\dist(m,\varphi(z))}\le C_\varphi(r/t)^\alpha.
\]
Consequently, for $d$ larger than a sufficiently large multiple of
$r$, \eqref{eq:subcritical-cov-local} applies at $(m,l)$ for either
ordering of $l,L$. Absorbing its logarithm with
exponent $\tau=s/\alpha<1$ proves \eqref{eq:tilde-cov-point} there.
For the other values of $d$, the uniform covariance bound gives
\eqref{eq:tilde-cov-point} after increasing its constant.

For sufficiently large $A$, depending on $\gamma,\varphi$, we also have
\begin{equation}\label{eq:tilde-cov-far}
 D>AR\quad\Longrightarrow\quad
 |\widetilde\Gamma((x,r),(y,R))|
 \le C_{\gamma,\varphi}(rR/D^2)^s .
\end{equation}
For $A>2$ the source intervals $[x-r,x+r]$ and $[y-R,y+R]$ are
disjoint. Pull back a shortest arc from $m$ to $M$ and extend it to
include both intervals, forming an arc $Q$ of length at least $D$.
The shortest image arc contains a half of each image interval, so
its length is at least $l+L$. The extensions add $l+L$ to that
length. By the upper power bound,
\[
 \frac{l}{\dist(m,M)}\le C_\varphi(r/D)^\alpha,\qquad
 \frac{L}{\dist(m,M)}\le C_\varphi(R/D)^\alpha .
\]
Increasing $A$ ensures $\dist(m,M)\ge C_*\max(l,L)$.
Then \eqref{eq:subcritical-cov-far}, with both logarithms absorbed
using $\tau$, proves \eqref{eq:tilde-cov-far}.

To integrate \eqref{eq:tilde-cov-point} on $D\le AR$, let $w$ denote either
integration variable with the other fixed. Write $F(u)$ for the length
of the set of $w$ where $D\le AR$ and $d\le u$. For fixed $y$, the three
points defining $d$
are fixed, so $F(u)\le6u$. For fixed $x$, the points $y\pm R$ contribute
at most $4u$. For $0<u\le R$, the remaining contribution is bounded by
\[
 |q_R^{-1}([x-u,x+u])|
 \le C_\varphi |q_R^{-1}([x-R,x+R])|(u/R)^\upsilon
 \le C_\varphi R(u/R)^\upsilon .
\]
The last inequality uses the displacement bound for $q_R^{-1}$.
In either variable, $F(u)\le C_\varphi R(u/R)^\upsilon$ for $0<u\le R$,
and $F(u)\le 2AR$ for $u>R$. Therefore
\[
 \int_{\{D\le AR\}}\min\{1,(r/d)^s\}\,\dd w
 = s r^s\int_r^\infty F(u)u^{-s-1}\,\dd u
 \le C_{\gamma,\varphi}r^sR^{1-s},
\]
since $s<\upsilon$. In the region $D>AR$, \eqref{eq:tilde-cov-far} integrates
in either variable to at most
$2C_{\gamma,\varphi}(rR)^s\int_{AR}^\infty D^{-2s}\,\dd D$,
which is bounded by $C_{\gamma,\varphi}r^sR^{1-s}$ since $2s>1$.
\end{proof}
\begin{thm}\label{thm:C}
    If $\varphi \in \QIR \cap \Sy(\m S^1)$, then $\varphi \in \WP(\m S^1)$.
\end{thm}
\begin{proof}
    Suppose $\varphi \in \Sy(\m S^1) \backslash \WP(\m S^1)$. By Lemma~\ref{lem:discrete WP},
    \[
        E_N \to \infty, (N \to \infty).
    \]
    By Lemma~\ref{lem:Bessel}, we have
    \[
        \m E[T_N^2] \leq C_{\g,\varphi} \int_{\O_N} \d_\varphi^2 \dd \eta = C_{\g,\varphi}(E_N-1),
    \]
    therefore,
    \[
        \frac{T_N}{E_N} \to 0 ,~\text{in}~\m P\text{-probability}.
    \]
    Using the same argument as in the proof of Theorem~\ref{thm:A}, we have
    \[
        \frac{T_N}{E_N} \to 0 ,~\text{in}~\m P_\varphi\text{-probability},
    \]
    that is
    \[
        \frac{T_N(\cdot \circ \varphi)}{E_N} \to 0 ,~\text{in}~\m P\text{-probability}.
    \]
    Using Lemma~\ref{lem:decomposition estimate}, we have
    \[
        \m E[B_N^2] \leq C_{\g,\varphi}^2 E_N,
    \]
    \[
         \frac{B_N}{E_N} \to 0 ,~\text{in}~L^2(\m P)~\text{and}~\m P\text{-probability},
    \]
    therefore,
    \[
         \frac{A_N}{E_N}=\frac{T_N(\cdot \circ \varphi)}{E_N}- \frac{B_N}{E_N}\to 0 ,~\text{in}~\m P\text{-probability}.
    \]
    By Lemma~\ref{lem:checkpoint}, along a subsequence,
    \[
        \m E[A_N^2] \leq C_{\g,\varphi}^2 E_N^2,
    \]
    Uniform integrability then gives
    \[
        \m E\left[\frac{A_N}{E_N}\right]\longrightarrow0,
    \]
    which contradicts
    \[
        \m E\left[\frac{A_N}{E_N}\right] \geq c_0^2\left(1-\frac{1}{E_N}\right) \to c_0^2
    \]
    by Lemma~\ref{lem:decomposition estimate}. 
\end{proof}
\section{Singularity in the subcritical regime}\label{sec:subcritical-singularity}

Fix $\kappa\in(0,4)$, let $\gamma=\sqrt\kappa$ and
$\lambda=\gamma/2$. We use $\m S^1=\m R/\m Z$, with Haar probability
measure $m$, and the subcritical GMC model of
Section~\ref{subsec:subcritical-model}:
\[
 \begin{gathered}
 H=F_2^{-1}\circ R_U\circ F_1,\qquad \m P=\operatorname{Law}(H),\\
 R_\varphi(h)=h\circ\varphi,\qquad \m P_\varphi=(R_\varphi)_*\m P.
 \end{gathered}
\]
Here $F_i$ is the cumulative function of the normalized measure
$\mu_i$ of $X_i$, and $X_1,X_2,U$ are independent, with $U$ Haar in
the mass coordinate. We write $P,E$ for probability and expectation on
this product space. Unspecified $L^p$ norms refer to this probability;
for functions of $H$ they agree with the norms under $\m P$.
Fix $\varphi\in\Homeo_+(\m S^1)$ and a degree-one lift $\Phi$.
The observables $Z_{x,r}$, $\delta_\varphi$ and the measure
$\dd\eta=\dd x\,\dd r/r^2$ are those of
Section~\ref{subsec:subcritical-model}.

This section proves one-sided singularity using logarithmic ratios
on adjacent arcs. The rotation lower bound transfers the cumulative
$L^2$ estimate from the welding to the source measure.

\begin{thm}[One-sided singularity]\label{thm:sing-conditional}
Let $\varphi\in\Homeo_+(\m S^1)\setminus\WP(\m S^1)$. Then
\[
 (R_\varphi)_*\m P\perp\m P,\qquad
 (R_\varphi)_*\m P^{\mathrm{weld}}\perp\m P^{\mathrm{weld}}.
\]
\end{thm}
The proof first separates the effect of composition into a midpoint
term and a nonnegative displacement. Averaging over $U$ gives a lower
bound for the displacement. Source-measure estimates make this bound
valid against arbitrary nonnegative weights, so that the singularity
criterion excludes a common part of the two laws. The last step uses
Lemma~\ref{introlem:sle-gmc-lift} to pass from the GMC law to SLE.

\subsection{Midpoint decomposition and accumulated statistics}\label{subsec:sing-decomposition}

Composition changes both the image interval and its splitting point.
The midpoint term below records the first change; the nonnegative
displacement records the second. We define the single-scale quantities
before accumulating them over positions and scales.

For $h\in\Homeo_+(\m S^1)$, let $\nu_h$ be the probability measure
defined by $\nu_h([u,v])=h(v)-h(u)$, using any degree-one lift of $h$.
For a non-atomic full-support finite positive measure $\nu$
and $a<c<b$ with $b-a<1$, set
\[
 R_\nu(a,c,M,b)=\left|
 \log\frac{\nu([c,b])}{\nu([a,c])}-
 \log\frac{\nu([M,b])}{\nu([a,M])}\right|,
 \qquad M=\frac{a+b}{2}.
\]
For $x\in\m S^1$ and $0<r<1/2$, put
\[
 a=\Phi(x-r),\quad c=\Phi(x),\quad b=\Phi(x+r),\qquad
 M_r(x)=\frac{a+b}{2},\quad L_r(x)=\frac{b-a}{2}.
\]
Write $R_\varphi(x,r;h)=R_{\nu_h}(a,c,M_r(x),b)$.
The function $z\mapsto\log\{\nu([z,b])/\nu([a,z])\}$ is decreasing.
It follows that
\begin{equation}\label{eq:sing-decomposition}
 Z_{x,r}(h\circ \varphi)=Z_{M_r(x),L_r(x)}(h)
   +\operatorname{sgn}(\delta_\varphi(x,r))R_\varphi(x,r;h).
\end{equation}

For the multiscale argument, write
\[
 \Omega_N=\m S^1\times[2^{-N-1},1/2),\qquad
 E_N=1+\int_{\Omega_N}\delta_\varphi(x,r)^2\,\dd\eta(x,r),
\]
and set
\[
 T_N(h)=\int_{\Omega_N}\delta_\varphi Z_{x,r}(h)\,\dd\eta,\qquad
 B_N(h)=\int_{\Omega_N}\delta_\varphi Z_{M_r(x),L_r(x)}(h)\,\dd\eta.
\]
For a non-atomic full-support probability measure $\nu$, also set
\[
 A_N^\nu=\int_{\Omega_N}|\delta_\varphi(x,r)|
 R_\nu\bigl(\Phi(x-r),\Phi(x),M_r(x),\Phi(x+r)\bigr)\,\dd\eta(x,r).
\]
The superscript in $A_N^\nu$ specifies the measure used to form the
image-length ratios. In particular, $A_N^{\mu_1}$ is a source-field
statistic, while $A_N(h):=A_N^{\nu_h}$ is a welding statistic.
Integrating \eqref{eq:sing-decomposition} gives
\begin{equation}\label{eq:sing-integrated}
 T_N(h\circ\varphi)=B_N(h)+A_N(h),\qquad A_N(h)\ge0.
\end{equation}

Let $\rho_2=(F_2)_*m$. Then $\rho_2([u,v])=F_2^{-1}(v)-F_2^{-1}(u)$,
and for every arc $I$,
\[
 {\nu_{H}(I)}=\rho_2(R_UF_1(I)),
 \qquad |F_1(I)|=\mu_1(I).
\]

\subsection{A measure-theoretic singularity criterion}\label{subsec:sing-criterion}

We will apply the following criterion to the decomposition just defined.
Its weighted lower bound tests the positive term against the density of
a possible common part. Uniform control under the transformed law is
needed to pair the original statistic with that common part.

\begin{lem}[A weighted singularity criterion]\label{lem:sing-criterion}
Let $P$ be a probability measure on a standard Borel space $\mathcal X$,
let $T$ be a Borel bijection with Borel inverse, and set $Q=T_*P$.
Suppose that measurable functions $D_n$ satisfy
\[
 D_n \rightharpoonup 0~\text{in }~L^2(P),\qquad
 \sup_n\|D_n\|_{L^2(Q)}<\infty.
\]
Suppose also that $D_n\circ T=N_n+Y_n$ $P$-a.s., where
$N_n\rightharpoonup0$ in $L^2(P)$, $Y_n\ge0$, and, for some $c>0$,
\begin{equation}\label{eq:sing-weighted}
 \liminf_{n\to\infty}\int wY_n\,\dd P\ge c\int w\,\dd P
 \quad\text{for every bounded measurable }w\ge0.
\end{equation}
Then $Q\perp P$.
\end{lem}

\begin{proof}
Let $Q=fP+Q^s$ be the Lebesgue decomposition and put $\mu=fP$.
Then $\mu\le Q$. Let $C=\sup_n\|D_n\|_{L^2(Q)}$. For $M>0$,
weak convergence gives
\[
 \limsup_n\left|\int D_n\,\dd\mu\right|
 \le C\,\mu\{f>M\}^{1/2}.
\]
Indeed, $f\mathbf1_{\{f\le M\}}\in L^2(P)$, and the remaining term
is bounded by Cauchy--Schwarz with respect to $\mu$. Letting $M\to\infty$
shows that $\int D_n\,\dd\mu\to0$.

Now $(T^{-1})_*\mu\le P$. Let $w$ be its density with respect to $P$.
Then $0\le w\le1$ and
\[
 0=\lim_n\int D_n\,\dd\mu
 =\lim_n\int w(N_n+Y_n)\,\dd P.
\]
Since $w\in L^2(P)$, the term with $N_n$ tends to zero.
By \eqref{eq:sing-weighted}, $0\ge c\int w\,\dd P=c\mu(\mathcal X)$.
Thus $\mu=0$, as required.
\end{proof}

\subsection{Averaging over the rotation}\label{subsec:sing-rotation}

For fixed fields, the independent phase $U$ moves an interval around the
target mass coordinate. We first prove a deterministic estimate for this
average. Applying it to the two arms of the splitting-point displacement
gives a source term minus an explicit error. The source estimates in the
next subsection will control both quantities.

Let $\nu$ be a non-atomic probability measure of full support on $\m S^1$.
For $0<s<1$ and a Lipschitz function $\Psi$, define
\[
 F^\nu(s)=\int_{\m S^1}\log\nu([x,x+s])\,\dd x,\qquad
 K^\nu_\Psi(s)=\int_{\m S^1}
 \frac{\Psi(y)\nu(\dd y)}{\nu([y-s,y])}.
\]
We also write $L_\nu(s)=-F^\nu(s)$ and $c_0=(\log2)/2$.
\begin{lem}[Averaging the reciprocal mass]\label{lem:weighted-kernel}
Let $\nu$ be a non-atomic probability measure of full support on $\m S^1$ and let $\Psi\ge0$ be Lipschitz. For $0<s<1$,
\begin{equation}\label{eq:weighted-kernel}
 K^\nu_\Psi(s)\ge\frac{c_0}{s}\int_{\m S^1}\Psi(x)\,\dd x
 -\operatorname{Lip}(\Psi)\{1-F^\nu(s)\}.
\end{equation}
The same bound holds with $\nu([y-s,y])$ replaced by $\nu([y,y+s])$.
\end{lem}
\begin{proof}
For fixed $s$, the denominator is continuous and has a positive minimum. Write $L=\operatorname{Lip}(\Psi)$. First assume $s\le1/2$. Set $n=\lfloor1/s\rfloor$ and $\ell=1/n$. For $\theta\in[0,\ell)$, let $I_j=[\theta+j\ell,\theta+(j+1)\ell]$, with indices modulo $n$, and put
\[
 p_j=\nu(I_j),\qquad a_j=\inf_{I_j}\Psi.
\]
For $y\in I_j$ let $u(y)$ be the $\nu$-mass of the part of $I_j$ to the left of $y$. Since $s\le\ell$,
\[
 \nu([y-s,y])\le p_{j-1}+u(y).
\]
The Stieltjes change of variables then gives
\[
 \int\frac{\Psi(y)\,\nu(\dd y)}{\nu([y-s,y])}
 \ge\sum_j a_j\log\left(1+\frac{p_j}{p_{j-1}}\right).
\]
Using $\log(1+e^v)\ge\log2+v/2$ and summation by parts on the cycle, the right side is at least
\[
 (\log2)\sum_j a_j+\frac12\sum_j(a_j-a_{j+1})\log p_j.
\]
Since $\Psi\le a_j+L\ell$ on $I_j$,
\[
 \sum_j a_j\ge n\int_{\m S^1}\Psi(x)\,\dd x-L.
\]
The intervals $I_j,I_{j+1}$ have a common endpoint, so $|a_j-a_{j+1}|\le L\ell$. Hence
\[
 \frac12\sum_j(a_j-a_{j+1})\log p_j
 \ge-\frac{L\ell}{2}\sum_j|\log p_j|.
\]
Averaging over $\theta\in[0,\ell)$ gives
\[
 \frac1\ell\int_0^\ell\ell\sum_j|\log p_j|\,\dd\theta
 =L_\nu(\ell)\le L_\nu(s).
\]
As $n\ge1/(2s)$, this proves \eqref{eq:weighted-kernel} for $s\le1/2$.

For $s>1/2$, the denominator is at most one. Thus
\[
 \int\frac{\Psi(y)\,\nu(\dd y)}{\nu([y-s,y])}
 \ge\int\Psi\,\dd\nu\ge\inf\Psi\ge\int_{\m S^1}\Psi(x)\,\dd x-L.
\]
Since $(\log2)/(2s)\le\log2<1$ and $L_\nu(s)\ge0$, the required bound follows. Reflection of the circle proves the other case.
\end{proof}

\begin{lem}[The rotation lower bound]\label{lem:weighted-lower-bound}
Let $F_1,F_2$ be lifts of orientation-preserving circle homeomorphisms, let $\mu_1$ be the measure with $\mu_1([a,b])=F_1(b)-F_1(a)$, and put $\rho_2=(F_2)_*m$. Let $H_U=F_2^{-1}\circ R_U\circ F_1$, where $R_U$ is rotation by $U$. For $a<c<b$, $b-a<1$, set $M=(a+b)/2$,
\[
 q=\mu_1([a,b]),\qquad
 s_*=\min\{\mu_1([a,\min(c,M)]),\mu_1([\max(c,M),b])\}.
\]
For every nonnegative Lipschitz $\Psi$,
\begin{equation}\label{eq:weighted-lower-bound}
 \begin{aligned}
 \int_{\m S^1}\Psi(U)R_{\nu_{H_U}}(a,c,M,b)\,\dd U
 &\ge\frac{\log2}{2}\left(\int_{\m S^1}\Psi(x)\,\dd x\right)R_{\mu_1}(a,c,M,b)\\
 &\quad-2\operatorname{Lip}(\Psi)q\{1+L_{\rho_2}(s_*)\}.
 \end{aligned}
\end{equation}
\end{lem}
\begin{proof}
Write $a_1=F_1(a)$, $c_1=F_1(c)$, $m_1=F_1(M)$ and $b_1=F_1(b)$,
and put $u_1=\min(c_1,m_1)$, $v_1=\max(c_1,m_1)$. Since
$\nu_{H_U}([z,w])=\rho_2([U+F_1(z),U+F_1(w)])$,
\[
 \begin{aligned}
 R_{\nu_{H_U}}(a,c,M,b)
 &=\log\frac{\rho_2([U+u_1,U+b_1])}{\rho_2([U+v_1,U+b_1])}\\
 &\quad+\log\frac{\rho_2([U+a_1,U+v_1])}{\rho_2([U+a_1,U+u_1])}.
 \end{aligned}
\]
Both terms are nonnegative. The Stieltjes formula for the logarithm
and Tonelli's theorem express their weighted sum as
\begin{align*}
 &\int_{b_1-v_1}^{b_1-u_1}\int_{\m S^1}
 \frac{\Psi(y-b_1+s)\,\rho_2(\dd y)}{\rho_2([y,y+s])}\,\dd s\\
 &\qquad+\int_{u_1-a_1}^{v_1-a_1}\int_{\m S^1}
 \frac{\Psi(y-a_1-s)\,\rho_2(\dd y)}{\rho_2([y-s,y])}\,\dd s.
\end{align*}
The translated weights have the same integral and Lipschitz constant as
$\Psi$. Lemma~\ref{lem:weighted-kernel} therefore gives the main term
\[
 \begin{aligned}
 &c_0\left(\int\Psi\,\dd m\right)
 \left(\log\frac{b_1-u_1}{b_1-v_1}
       +\log\frac{v_1-a_1}{u_1-a_1}\right)\\
 &\qquad=c_0\left(\int\Psi\,\dd m\right)R_{\mu_1}(a,c,M,b).
 \end{aligned}
\]
The two integration intervals have total length $2(v_1-u_1)\le2q$.
Throughout them,
\[
 s_* =\min\{u_1-a_1,b_1-v_1\}\le s\le q<1.
\]
Since $L_{\rho_2}$ is nonincreasing, their total error is at most
$2\operatorname{Lip}(\Psi)q\{1+L_{\rho_2}(s_*)\}$.
This proves both orders of $c$ and $M$ at once; when $c=M$, the two
integration intervals are empty and the desired bound is immediate.
\end{proof}

\subsection{Source-measure estimates and weighted decorrelation}\label{subsec:sing-source}

The main term in Lemma~\ref{lem:weighted-lower-bound} is the source
displacement $R_{\mu_1}$; its error is controlled by a small source mass
and logarithms of the two arm masses. The estimates below give a positive
mean for the first term and show that the second vanishes on shrinking
arcs. We then delete finitely many Fourier modes of $X_1$ to obtain the
weighted version of the mean lower bound.

For an arc $I$, write $|I|=m(I)$. The following consequence of
Lemma~\ref{lem:log GMC moment} will be used throughout this subsection.
\begin{lem}[Nested source-mass ratios]\label{lem:sing-source-log-bounds}
For every $1\le p<\infty$, there is $C_p<\infty$ such that, for
nondegenerate arcs $I\subset J$ with a common endpoint,
\begin{equation}\label{eq:nonwp-source-moments}
 \left\|\log\frac{\mu_1(J)}{\mu_1(I)}\right\|_p
 \le C_p\left(1+\log\frac{|J|}{|I|}\right),\qquad
 \|\log\mu_1(I)\|_p\le C_p\left(1+\log\frac1{|I|}\right).
\end{equation}
\end{lem}
\begin{proof}
If $I'$ is adjacent to $I$ and $|I'|=|I|$, then
\[
 0\le\log\frac{\mu_1(I\cup I')}{\mu_1(I)}
 \le\log2+\left|\log\frac{\mu_1(I')}{\mu_1(I)}\right|.
\]
Lemma~\ref{lem:log GMC moment} bounds the right side in $L^p$ by a
constant independent of $I$. Starting from $I$, double its length while
keeping the common endpoint fixed, stopping when the resulting arc covers
$J$. There are at most $1+\lceil\log_2(|J|/|I|)\rceil$ steps.
For the final step at lengths greater than $1/2$, use a half-circle
contained in the current arc and its complementary half-circle. The same
bound applies. Summing these bounds by Minkowski's inequality proves the
first estimate. Taking $J=\m S^1$ and using $\mu_1(\m S^1)=1$ proves
the second.
\end{proof}

The positive moments, negative moments and scaling relation for
unnormalized GMC used next are given in
\cite[Theorems~3.27, 3.30 and~3.42]{GFFbook}. To match that convention,
apply these results to $X_i/\sqrt2$, with dimension $d=1$ and GMC
parameter $\gamma/\sqrt2$. On a fixed circle chart the remaining
covariance term is bounded, so Kahane's comparison inequality gives
uniform constants; see \cite[Theorem~3.19]{GFFbook}.

\begin{lem}[Small-arc masses and inverse coordinates]\label{lem:nonwp-mass}
There are $p\in(1,2)$, $b>1$ and $C<\infty$ such that
\begin{equation}\label{eq:nonwp-mass}
 E[\mu_i(I)^p]\le C|I|^b
\end{equation}
for every arc $I$. There are $a_0>0$ and $K_2\ge1$, depending only on $X_2$, such that $E[\log K_2]<\infty$ and
\begin{equation}\label{eq:nonwp-quantile}
 \mu_2(I)\le K_2|I|^{a_0},\qquad
 L_{\rho_2}(s)\le a_0^{-1}\bigl(\log K_2+\log(1/s)\bigr).
\end{equation}
\end{lem}
\begin{proof}
Let $M_i$ be the unnormalized GMC and put $\tau=\gamma^2/4$.
The scaling and negative-moment estimates cited above give
\[
 E[M_i(I)^q]\le C_q|I|^{\zeta(q)},\qquad
 \zeta(q)=(1+\tau)q-\tau q^2,\quad1<q<1/\tau,
\]
and, by H\"older, for $a>1$ and $pa<1/\tau$,
\[
 E[\mu_i(I)^p]\le C_{p,a}|I|^{\zeta(pa)/a}.
\]
Since $\zeta(1)=1$ and $\zeta'(1)=1-\tau>0$, choose $p\in(1,2)$
and then $a>1$ sufficiently close to $1$ so that $b=\zeta(pa)/a>1$.
This proves \eqref{eq:nonwp-mass}.

Choose $0<a_0<(b-1)/p$ and let $I_{j,k}$ be the dyadic arcs of length
$2^{-j}$, with indices $0\le k<2^j$. Set
$Z=\sup_{j,k}2^{ja_0}\mu_2(I_{j,k})$. Then $Z\ge1$ and
\[
 E[Z^p]\le\sum_{j,k}2^{jpa_0}E[\mu_2(I_{j,k})^p]
 \le C\sum_{j\ge0}2^{-j(b-1-pa_0)}<\infty.
\]
In particular,
\begin{equation}\label{eq:dyadic-mass-tail}
 P(Z>t)\le Ct^{-p}\quad(t\ge1),\qquad E[\log Z]<\infty.
\end{equation}
Every arc $I$ with $|I|\le2^{-j}<2|I|$ is covered by at most two
level-$j$ dyadic arcs, including across the circle seam. Non-atomicity
therefore gives
\[
 \mu_2(I)\le2Z2^{-ja_0}\le K_2|I|^{a_0},\qquad
 K_2=2^{1+a_0}Z.
\]
Thus $E[\log K_2]<\infty$.

Finally, for $Q_2=F_2^{-1}$ and
$\ell=Q_2(u+s)-Q_2(u)=\rho_2([u,u+s])$, the preceding estimate gives
\[
 s=\mu_2([Q_2(u),Q_2(u+s)])\le K_2\ell^{a_0},\qquad
 -\log\ell\le a_0^{-1}(\log K_2+\log(1/s)).
\]
Integrating over $u$ proves the bound for $L_{\rho_2}$.
\end{proof}

\begin{lem}[Source displacement]\label{lem:nonwp-source-log-moment}
Let $a<c<b$, $b-a<1$, $M=(a+b)/2$ and $d=|\log((b-c)/(c-a))|$. For each $p\ge1$,
\begin{equation}\label{eq:nonwp-source-log-moment}
 E[R_{\mu_1}(a,c,M,b)]\ge c_0d,\qquad
 \|R_{\mu_1}(a,c,M,b)\|_p\le C_p(1+d).
\end{equation}
\end{lem}
\begin{proof}
For any non-atomic full-support probability measure $\nu$, the
Stieltjes formula, Tonelli's theorem and
Lemma~\ref{lem:weighted-kernel} with $\Psi=1$ give
\begin{equation}\label{eq:nonwp-F-increments}
 \begin{aligned}
 F^\nu(s_+)-F^\nu(s_-)
 &=\int_{s_-}^{s_+}\int_{\m S^1}
   \frac{\nu(\dd y)}{\nu([y-s,y])}\,\dd s\\
 &\ge c_0\log(s_+/s_-)\qquad(0<s_-<s_+<1).
 \end{aligned}
\end{equation}
Rotation invariance and \eqref{eq:nonwp-source-moments} imply that
$f(s)=E[\log\mu_1([x,x+s])]=E[F^{\mu_1}(s)]$ is finite and independent
of $x$. Thus
\begin{equation}\label{eq:nonwp-f-increments}
 f(s_+)-f(s_-)\ge c_0\log(s_+/s_-).
\end{equation}
Write $\ell_-=c-a$ and $\ell_+=b-c$. Monotonicity of
$z\mapsto\log\{\mu_1([z,b])/\mu_1([a,z])\}$ fixes the sign of the
difference defining $R_{\mu_1}$. The two halves cut at $M$ have equal
expected logarithmic mass, so
\[
 E[R_{\mu_1}(a,c,M,b)]
 =f(\max(\ell_-,\ell_+))-f(\min(\ell_-,\ell_+))\ge c_0d.
\]
For the upper bound, apply Lemma~\ref{lem:log GMC moment} directly
to the adjacent arcs $[a,c],[c,b]$ and to the two equal halves cut at
$M$. The triangle inequality gives $\|R_{\mu_1}\|_p\le C_p(1+d)$.
The same formulas include $c=M$, when $R_{\mu_1}=d=0$.
\end{proof}

\begin{lem}[The small-arc error]\label{lem:nonwp-error}
For $\ell=b-a\le1/2$, define $q,s_*$ as in Lemma~\ref{lem:weighted-lower-bound} and $d=|\log((b-c)/(c-a))|$. There exist $b_1>1$ and $C<\infty$ such that
\begin{equation}\label{eq:nonwp-error}
 E\left[q^2\{1+\log(1/s_*)\}^2\right]
 \le C\ell^{b_1}\{1+\log(1/\ell)+d\}^2.
\end{equation}
In particular, if $\ell_n\to0$ and $d_n$ is bounded, then
\[
 E[q_n\{1+\log(1/s_{*,n})\}]\longrightarrow0.
\]
If $\ell_n\to0$ and $d_n\ge1$, the same expectation divided by $d_n$
tends to zero.
\end{lem}
\begin{proof}
Take $p,b$ from Lemma~\ref{lem:nonwp-mass} and choose $1<h<b$. Since $q\le1$ and $2h>p$,
\[
 (E[q^{2h}])^{1/h}\le C\ell^{b/h}.
\]
The shorter of the two physical arms has length $\ell/(1+e^d)$. Both arcs in the definition of $s_*$ have at least this length. By \eqref{eq:nonwp-source-moments}, for every finite $k$,
\[
 \|1+\log(1/s_*)\|_k
 \le C_k\{1+\log(1/\ell)+d\}.
\]
Let $h'=h/(h-1)$. H\"older's inequality gives
\[
 \begin{aligned}
 E[q^2\{1+\log(1/s_*)\}^2]
 &\le E[q^{2h}]^{1/h}
 E[\{1+\log(1/s_*)\}^{2h'}]^{1/h'}\\
 &\le C\ell^{b/h}\{1+\log(1/\ell)+d\}^2.
 \end{aligned}
\]
Taking $b_1=b/h>1$ proves \eqref{eq:nonwp-error}. By Cauchy--Schwarz,
\[
 E[q\{1+\log(1/s_*)\}]
 \le C\ell^{b_1/2}\{1+\log(1/\ell)+d\}.
\]
If $d_n$ is bounded, the right side tends to zero. If $d_n\ge1$, divide by
$d_n$; the resulting bound is at most
$C\ell_n^{b_1/2}\{2+\log(1/\ell_n)\}$, which also tends to zero.
\end{proof}

\begin{lem}[Deletion of finitely many Fourier modes]\label{lem:nonwp-mode-deletion}
Let $a_n<c_n<b_n<a_n+1$, $M_n=(a_n+b_n)/2$, and
\[
 \ell_n=b_n-a_n,\qquad
 d_n=\left|\log\frac{b_n-c_n}{c_n-a_n}\right|.
\]
Assume $\ell_n\to0$. Let $V_n$ be either
\[
 \frac{R_{\mu_1}(a_n,c_n,M_n,b_n)}{d_n},\quad d_n\ge1,
 \qquad\text{or}\qquad
 R_{\mu_1}(a_n,c_n,M_n,b_n),\quad \sup_n d_n<\infty.
\]
Then, for every bounded $X_1$-measurable random variable $G_1$,
\begin{equation}\label{eq:nonwp-point-covariance}
 \operatorname{Cov}(G_1,V_n)\longrightarrow0.
\end{equation}
Let $N_k\to\infty$ be a deterministic sequence such that
$E_{N_k}\to\infty$ and $\sup_k\|A_{N_k}^{\mu_1}/E_{N_k}\|_2<\infty$.
Assume also that, for every $r_0>0$, the unnormalized integral defining
$A_{N_k}^{\mu_1}$ over $r\ge r_0$ has expectation bounded uniformly in $k$.
Then
\begin{equation}\label{eq:nonwp-integrated-covariance}
 \operatorname{Cov}\!\left(G_1,\frac{A_{N_k}^{\mu_1}}{E_{N_k}}\right)
 \longrightarrow0.
\end{equation}
\end{lem}
\begin{proof}
We use the moment bounds in \eqref{eq:nonwp-source-log-moment}. They imply
$\sup_n\|V_n\|_2<\infty$ in both point cases. The corresponding bound for
$A_{N_k}^{\mu_1}/E_{N_k}$ is an assumption of the lemma.

Fix $m\ge1$ and let $\mathcal F_m$ be generated by the first $m$
Fourier pairs of $X_1$. Decompose $X_1=X_{\le m}+X_{>m}$ into these
modes and the independent remaining field. Define
\[
 s_m(x)=\lambda X_{\le m}(x)
 -\frac{\lambda^2}{2}E[X_{\le m}(x)^2],
 \qquad L_m=\sup_x|s_m'(x)|.
\]
Let $M_{>m}$ be the subcritical GMC of $X_{>m}$, and put
$\widehat\mu_{>m}=M_{>m}/M_{>m}(\m S^1)$. The random measure
$\widehat\mu_{>m}$ is independent of $\mathcal F_m$. Wick multiplication
for the independent fields gives
\begin{equation}\label{eq:nonwp-finite-mode-wick}
 M_{X_1}(\dd x)=e^{s_m(x)}M_{>m}(\dd x).
\end{equation}
Since $s_m$ is a finite Gaussian Fourier sum, $E[e^{tL_m}]<\infty$ for
all $t>0$, and $p_m:=P(L_m\le1)>0$.

For a continuous function $s$ and an arc $I$, define
$\operatorname{osc}_I s=\sup_{x,y\in I}|s(x)-s(y)|$. Let $I=[a,b]$ and
$\nu^s(\dd x)=e^{s(x)}\nu(\dd x)$. For $a<z<b$, set
\[
 \ell_\nu(z)=\log\frac{\nu([z,b])}{\nu([a,z])}.
\]
The Stieltjes chain rule gives
\begin{equation}\label{eq:nonwp-stieltjes}
 -\dd\ell_\nu(z)=
 \left(\frac1{\nu([a,z])}+\frac1{\nu([z,b])}\right)\nu(\dd z).
\end{equation}
For $J=[a,z]$ or $J=[z,b]$, let
$\overline{e^s}_J=\nu(J)^{-1}\int_J e^{s(u)}\nu(\dd u)$. Then
\[
 e^{-\operatorname{osc}_I s}
 \le\frac{e^{s(z)}}{\overline{e^s}_J}
 \le e^{\operatorname{osc}_I s}.
\]
Both terms in \eqref{eq:nonwp-stieltjes} are nonnegative. Integrating
between $c$ and $M$, in either order, gives
\begin{equation}\label{eq:nonwp-delete}
 e^{-\operatorname{osc}_I s}R_\nu(a,c,M,b)
 \le R_{\nu^s}(a,c,M,b)
 \le e^{\operatorname{osc}_I s}R_\nu(a,c,M,b).
\end{equation}
The quantity $R_\nu$ does not change when $\nu$ is multiplied by a positive
constant. Thus \eqref{eq:nonwp-finite-mode-wick} and
$\operatorname{osc}_I s_m\le |I|L_m$ give the same comparison between
$\mu_1$ and $\widehat\mu_{>m}$.

In either point case let $\widehat V_n$ be $V_n$ with $\mu_1$ replaced by
$\widehat\mu_{>m}$. In the integral case set
$V_k=A_{N_k}^{\mu_1}/E_{N_k}$ and
$\widehat V_k=A_{N_k}^{\widehat\mu_{>m}}/E_{N_k}$. Applying
\eqref{eq:nonwp-delete} pointwise and using nonnegative weights, on
$\{L_m\le1\}$ we have $\widehat V\le eV$. By independence,
\begin{equation}\label{eq:nonwp-fine-moment}
 p_m E[\widehat V^2]
 =E[\mathbf1_{\{L_m\le1\}}\widehat V^2]
 \le e^2 E[V^2].
\end{equation}
Consequently, for fixed $m$, the fine-field variables have uniformly bounded
second moments. In particular their first moments are bounded by a constant
$C_m$ independent of the sequence index.

In the point cases, independence and \eqref{eq:nonwp-delete} give
\[
 E[|V_n-\widehat V_n|]
 \le E[e^{\ell_nL_m}-1]E[\widehat V_n].
\]
For all large $n$, $\ell_n\le1$. Since $E[e^{L_m}]<\infty$, dominated
convergence gives
\begin{equation}\label{eq:nonwp-point-deletion}
 E[|V_n-\widehat V_n|]\longrightarrow0.
\end{equation}

For the integral case, let
$\omega_\varphi(t)=\sup_{|x-y|\le t}|\Phi(x)-\Phi(y)|$.
Then $\omega_\varphi(t)\to0$ as $t\downarrow0$. Fix $r_0>0$ and split
both integrals at $r=r_0$. Let $C_k$ and $\widehat C_k$ be the respective
unnormalized contributions on $r\ge r_0$. By assumption,
$\sup_k E[C_k]<\infty$. Since every image parent arc has length less than
$1$, \eqref{eq:nonwp-delete} gives $\widehat C_k\le eC_k$ on
$\{L_m\le1\}$. Independence gives
\[
 \sup_k E[\widehat C_k]\le e p_m^{-1}\sup_k E[C_k]<\infty.
\]
For $r<r_0$, the image parent arc has length at most
$\omega_\varphi(2r_0)$. Using nonnegativity and independence, we obtain
\[
 \begin{aligned}
 E[|V_k-\widehat V_k|]
 &\le\frac{E[C_k]+E[\widehat C_k]}{E_{N_k}}\\
 &\quad+E[e^{\omega_\varphi(2r_0)L_m}-1]E[\widehat V_k]\\
 &\le\frac{C_m(r_0)}{E_{N_k}}
 +C_mE[e^{\omega_\varphi(2r_0)L_m}-1].
 \end{aligned}
\]
First let $k\to\infty$ with $r_0$ fixed. Then let $r_0\downarrow0$.
Dominated convergence gives
\begin{equation}\label{eq:nonwp-integrated-deletion}
 E[|V_k-\widehat V_k|]\longrightarrow0.
\end{equation}

It remains to pass from finite-mode weights to arbitrary source weights.
If $G_1$ is bounded and $\mathcal F_m$-measurable, then
$\widehat V$ is independent of $\mathcal F_m$, so
\[
 \operatorname{Cov}(G_1,V)
 =E[(G_1-E[G_1])(V-\widehat V)].
\]
Equations \eqref{eq:nonwp-point-deletion} and
\eqref{eq:nonwp-integrated-deletion} show that this covariance tends to zero.
For general bounded $G_1$, put $G_1^{(m)}=E[G_1\mid\mathcal F_m]$.
Since the Fourier coefficients generate the field sigma-field,
$G_1^{(m)}\to G_1$ in $L^2$. The uniform bound $\sup\|V\|_2<\infty$ gives
\[
 \limsup\left|\operatorname{Cov}(G_1,V)\right|
 \le C\|G_1-G_1^{(m)}\|_2.
\]
Letting $m\to\infty$ proves both assertions.
\end{proof}
\subsection{Weighted estimates in the three regularity regimes}\label{subsec:sing-estimates}

We combine the rotation lower bound with the source estimates. The
single-scale moment and cancellation estimates come directly from
Section~\ref{sec:subcritical-estimates}. At accumulated scales, we
use the rotation lower bound with a constant weight to obtain the
source bound needed by the mode-deletion argument.

For adjacent arcs $I_-,I_+$ of lengths less than $1/2$,
Lemma~\ref{lem:moment of Z}, with $p=2$, gives
\begin{equation}\label{eq:nonwp-weld-moment}
 \left\|\log\frac{\nu_{H}(I_+)}{\nu_{H}(I_-)}\right\|_2
 \le C_\gamma\left(1+\left|\log\frac{|I_+|}{|I_-|}\right|\right).
\end{equation}
Lemma~\ref{lem:subcritical-weak-cancellation} gives, whenever
$x_n\to x$ and $r_n\downarrow0$,
\begin{equation}\label{eq:nonwp-weak}
 Z_{x_n,r_n}\rightharpoonup0\quad\text{in }L^2(\m P).
\end{equation}
For $\varphi\in\Sy(\m S^1)$, apply Lemma~\ref{lem:Bessel} with
$f(x,r)=\delta_\varphi(x,r)\mathbf1_{\Omega_N}(x,r)$. It gives
\begin{equation}\label{eq:sing-bessel-bounds}
 \|T_N\|_2+\|B_N\|_2\le C_{\gamma,\varphi}\sqrt{E_N}.
\end{equation}
Lemma~\ref{lem:checkpoint} supplies a deterministic subsequence along
which $\|A_{N_k}(H)\|_2\le C_{\gamma,\varphi}E_{N_k}$.
For fixed $X_1,X_2$, take $\Psi\equiv1$ in
\eqref{eq:weighted-lower-bound}. The error term vanishes. Integrating
the resulting inequality against $|\delta_\varphi|\,\dd\eta$ on
$\Omega_N$ and using Tonelli's theorem gives
\[
 c_0 A_N^{\mu_1}\le E[A_N(H)\mid X_1,X_2],
 \qquad c_0=(\log 2)/2.
\]
Both sides are nonnegative. Conditional Jensen's inequality gives
$\|A_N^{\mu_1}\|_2\le c_0^{-1}\|A_N(H)\|_2$, so the same subsequence
gives the two bounds used in the multiscale argument:
\begin{equation}\label{eq:nonwp-cumulative-input}
 \begin{aligned}
 \|T_{N_k}(H)\|_2+\|B_{N_k}(H)\|_2
 &\le C_{\gamma,\varphi}\sqrt{E_{N_k}},\\
 \|A_{N_k}(H)\|_2+\|A_{N_k}^{\mu_1}\|_2
 &\le C_{\gamma,\varphi}E_{N_k}.
 \end{aligned}
\end{equation}
The three regimes use these estimates as follows. Outside $\QS(\m S^1)$, the
distortion diverges and the moment bound \eqref{eq:nonwp-weld-moment}
controls a single normalized observation. In $\QS(\m S^1)\setminus\Sy(\m S^1)$,
\eqref{eq:nonwp-weak} supplies cancellation along shrinking arcs.
In $\Sy(\m S^1)\setminus\WP(\m S^1)$, the divergence $E_{N_k}\to\infty$ and
\eqref{eq:nonwp-cumulative-input} allow the distortion to be accumulated
over scales. In all three cases, Lemma~\ref{lem:weighted-lower-bound}
and the source estimates give the positive term against a nonnegative
weight.

For symmetric maps, the rotation error for nonconstant weights must
also be summable over scales.

\begin{lem}[Integrability of the accumulated error]\label{lem:nonwp-sym-integrable}
Let $\varphi\in\Sy(\m S^1)$. For $x\in\m S^1$ and $0<r<1/2$, put
\[
 a=\Phi(x-r),\quad c=\Phi(x),\quad b=\Phi(x+r),\quad
 M=\frac{a+b}{2},\quad q(x,r)=\mu_1([a,b]),
\]
and
\[
 s_*(x,r)=\min\{\mu_1([a,\min(c,M)]),
 \mu_1([\max(c,M),b])\}.
\]
Then
\begin{equation}\label{eq:nonwp-error-integrable}
 \int_{\m S^1\times(0,1/2)}
 E\!\left[q(x,r)^2\{1+\log(1/s_*(x,r))\}^2\right]\,\dd\eta<\infty.
\end{equation}
Consequently, if $E_N\to\infty$ along a sequence, then
\begin{equation}\label{eq:nonwp-integrated-error}
 \frac1{E_N}\int_{\Omega_N}|\delta_\varphi|
 E\!\left[q\{1+\log(1/s_*)\}\right]\,\dd\eta\longrightarrow0.
\end{equation}
\end{lem}
\begin{proof}
Fix $0<\alpha<1<\beta$. Apply Lemma~\ref{lem:sym ab} inside a
finite collection of fixed parent arcs whose interiors cover the circle.
Absorbing their image lengths into the constants gives, uniformly in
$x$ and for sufficiently small $r$,
\begin{equation}\label{eq:nonwp-sym-powers}
 cr^\beta\le\ell(x,r):=\Phi(x+r)-\Phi(x-r)\le Cr^\alpha.
\end{equation}
Since $\varphi\in\QS(\m S^1)$, there is $D<\infty$ such that
$|\delta_\varphi(x,r)|\le D$. Let $b_1>1$ be the exponent in
Lemma~\ref{lem:nonwp-error}. Choose $\alpha$ close enough to $1$ that
$\alpha b_1>1$. Equations \eqref{eq:nonwp-error} and
\eqref{eq:nonwp-sym-powers} yield
\[
 E\!\left[q^2\{1+\log(1/s_*)\}^2\right]
 \le C r^{\alpha b_1}(1+\log(1/r))^2.
\]
The right side is integrable near $r=0$ against $\dd r/r^2$.
On a fixed band $r\ge r_0>0$, the two arcs defining $s_*$ have lengths
bounded below by a positive constant. The logarithmic moment bound
\eqref{eq:nonwp-source-moments} gives a uniform bound there. This proves
\eqref{eq:nonwp-error-integrable}.

Set $a(x,r)=E[q\{1+\log(1/s_*)\}]$. Jensen's inequality and
\eqref{eq:nonwp-error-integrable} give $\int a^2\,\dd\eta<\infty$.
Cauchy--Schwarz then gives
\[
 \frac1{E_N}\int_{\Omega_N}|\delta_\varphi|a\,\dd\eta
 \le\frac{\sqrt{E_N-1}}{E_N}
 \left(\int a^2\,\dd\eta\right)^{1/2}.
\]
The right side tends to zero when $E_N\to\infty$.
\end{proof}

\begin{prop}[A weighted positive lower bound]\label{prop:nonwp-weighted}
For every bounded measurable $W=W(X_1,X_2,U)\ge0$, the following hold.
\begin{enumerate}[label=(\roman*)]
\item If $\ell_n=\Phi(x_n+r_n)-\Phi(x_n-r_n)\to0$ and $d_n=|\delta_\varphi(x_n,r_n)|\to\infty$, then
\[\liminf_nE\!\left[W\frac{R_\varphi(x_n,r_n;H)}{d_n}\right]\ge c_0^2E[W].\]
\item If $\ell_n\to0$ and $0<\varepsilon\le d_n\le D<\infty$, then
\[\liminf_nE[W R_\varphi(x_n,r_n;H)]\ge c_0^2\varepsilon E[W].\]
\item If $\varphi\in\Sy(\m S^1)\setminus\WP(\m S^1)$ and $N_k$ is as in \eqref{eq:nonwp-cumulative-input}, then
\[\liminf_kE\!\left[W\frac{A_{N_k}(H)}{E_{N_k}}\right]\ge c_0^2E[W].\]
\end{enumerate}
\end{prop}
\begin{proof}
Let $W=G_1(X_1)G_2(X_2)\Psi(U)$, where $G_1,G_2$ are bounded and
nonnegative, and $\Psi$ is nonnegative and Lipschitz. For fixed $X_1,X_2$,
Lemma~\ref{lem:weighted-lower-bound} and the bound on $L_{\rho_2}$ in
\eqref{eq:nonwp-quantile} give
\[
 \begin{aligned}
 &\int\Psi(U)R_{\nu_{H_U}}(a,c,M,b)\,\dd U\\
 &\quad\ge c_0\left(\int\Psi(u)\,\dd u\right)R_{\mu_1}(a,c,M,b)\\
 &\qquad-C\operatorname{Lip}(\Psi)q
 \{1+\log(1/s_*)+\log K_2\}.
 \end{aligned}
\]
Here $q,s_*$ depend on $X_1$ and $K_2$ depends on $X_2$. Since the two
fields are independent and $E[\log K_2]<\infty$, the expected error after
multiplication by $G_1G_2$ is at most
\[
 C_{G_1,G_2,\Psi}E[q\{1+\log(1/s_*)\}].
\]
No independence between $q$ and $s_*$ is used.

In case (i), divide by $d_n$. In case (ii), use the inequality without
normalization. Lemma~\ref{lem:nonwp-error} gives an error tending to zero
in both cases. Define
\[
 V_n=\frac{R_{\mu_1}(a_n,c_n,M_n,b_n)}{d_n}
 \quad\text{in (i)},\qquad
 V_n=R_{\mu_1}(a_n,c_n,M_n,b_n)\quad\text{in (ii)}.
\]
Here $a_n=\Phi(x_n-r_n)$, $c_n=\Phi(x_n)$,
$b_n=\Phi(x_n+r_n)$ and $M_n=(a_n+b_n)/2$.
For either case, the main term is
\[
 c_0\left(\int\Psi(u)\,\dd u\right)E[G_2]E[G_1V_n].
\]
By Lemma~\ref{lem:nonwp-mode-deletion},
\[
 E[G_1V_n]=E[G_1]E[V_n]+o(1).
\]
Lemma~\ref{lem:nonwp-source-log-moment} gives
\[
 \liminf_n E[V_n]\ge c_0\quad\text{in (i)},\qquad
 \liminf_n E[V_n]\ge c_0\varepsilon\quad\text{in (ii)}.
\]
Since $E[W]=E[G_1]E[G_2]\int\Psi(u)\,\dd u$, the required bounds follow
for these product weights.

In case (iii), multiply the first displayed inequality by
$|\delta_\varphi(x,r)|/E_N$ and integrate over $\Omega_N$. By
Lemma~\ref{lem:nonwp-sym-integrable}, the resulting error tends to zero.
The main term is
\[
 c_0\left(\int\Psi(u)\,\dd u\right)E[G_2]
 E\!\left[G_1\frac{A_N^{\mu_1}}{E_N}\right].
\]
Tonelli's theorem and Lemma~\ref{lem:nonwp-source-log-moment} give
\[
 E[A_N^{\mu_1}]
 \ge c_0\int_{\Omega_N}\delta_\varphi^2\,\dd\eta
 =c_0(E_N-1).
\]
The second-moment bound in \eqref{eq:nonwp-cumulative-input} applies to
$A_N^{\mu_1}/E_N$ along $N=N_k$. On every band $r\ge r_0>0$, boundedness
of $\delta_\varphi$ and \eqref{eq:nonwp-source-log-moment} give a uniform
bound on the expectation of the unnormalized integral. Thus
Lemma~\ref{lem:nonwp-mode-deletion} gives
\[
 E\!\left[G_1\frac{A_N^{\mu_1}}{E_N}\right]
 =E[G_1]E\!\left[\frac{A_N^{\mu_1}}{E_N}\right]+o(1).
\]
As $E_N\to\infty$, the lower limit of the left side is at least
$c_0E[G_1]$. This proves (iii) for product weights.

For the extension to general $W$, let $Y_n$ denote the nonnegative variable
on the left of the relevant assertion, without $W$. In (iii), the index $n$
stands for $N_k$. By \eqref{eq:nonwp-weld-moment} and
\eqref{eq:sing-decomposition}, $\sup_n\|Y_n\|_2<\infty$ in (i) and (ii).
The bound in \eqref{eq:nonwp-cumulative-input} gives the same in (iii).
Denote this bound by $C$ and the corresponding lower-bound constant by
$c>0$.

Conditional expectations on increasing finite product partitions of
the field Fourier coefficients and the phase approximate $W$ in $L^2$.
They are finite nonnegative sums of products of field indicators and
phase indicators. Approximating the latter by nonnegative Lipschitz
functions gives product-weight sums $W_j$ with
$\|W-W_j\|_2\to0$. The uniform $L^2$ bound on $Y_n$ gives
\[
 \begin{aligned}
 \liminf_nE[WY_n]
 &\ge\liminf_nE[W_jY_n]-C\|W-W_j\|_2\\
 &\ge cE[W]-(C+c)\|W-W_j\|_2.
 \end{aligned}
\]
Letting $j\to\infty$ proves all three statements.
\end{proof}
\subsection{Proof of one-sided singularity}

We apply Lemma~\ref{lem:sing-criterion} with $T=R_\varphi$.
The first two cases use a shrinking sequence of arcs. The symmetric
non-WP case uses the normalized cumulative statistic.

\begin{proof}[Proof of Theorem~\ref{thm:sing-conditional}]
Apply Lemma~\ref{lem:sing-criterion} with
$\mathcal X=\Homeo_+(\m S^1)$, $P=\m P$ and $T=R_\varphi$.
For a bounded measurable $w\ge0$ on $\mathcal X$, the estimates in
Proposition~\ref{prop:nonwp-weighted} apply to $W=w(H)$.

Suppose first that $\varphi\notin\QS(\m S^1)$. There are $x_n$ and
$r_n\to0$ such that $d_n=|\delta_\varphi(x_n,r_n)|\to\infty$. If the
radii stayed bounded away from zero, continuity and strict monotonicity of
$\Phi$ would give a uniform bound on $d_n$. Set
\[
 D_n(h)=\frac{\operatorname{sgn}(\delta_\varphi(x_n,r_n))}{d_n}
 Z_{x_n,r_n}(h).
\]
The adjacent-arc estimate \eqref{eq:nonwp-weld-moment} gives
\[
 \|D_n\|_{L^2(P)}\le\frac{C}{d_n}\longrightarrow0,\qquad
 \sup_n\|D_n\|_{L^2(T_*P)}<\infty.
\]
By \eqref{eq:sing-decomposition}, $D_n\circ T=N_n+Y_n$, where
\[
 \begin{aligned}
 N_n(h)&=\frac{\operatorname{sgn}(\delta_\varphi(x_n,r_n))}{d_n}
 Z_{M_{r_n}(x_n),L_{r_n}(x_n)}(h),\\
 Y_n(h)&=\frac{R_\varphi(x_n,r_n;h)}{d_n}.
 \end{aligned}
\]
Again \eqref{eq:nonwp-weld-moment} gives $\|N_n\|_{L^2(P)}\le C/d_n$.
By uniform continuity, the lengths of the image parent arcs tend to zero.
Proposition~\ref{prop:nonwp-weighted}(i) therefore gives
\[
 \liminf_n E[w(H)Y_n(H)]\ge c_0^2E[w(H)].
\]
The conditions of Lemma~\ref{lem:sing-criterion} hold, so $T_*P\perp P$.

Next suppose $\varphi\in\QS(\m S^1)\setminus\Sy(\m S^1)$. There are
$x_n$, $r_n\to0$, and $0<\varepsilon\le D<\infty$ such that
$\varepsilon\le d_n:=|\delta_\varphi(x_n,r_n)|\le D$.
After passing to a subsequence, $x_n\to x$ and
\[
 M_{r_n}(x_n)\to\varphi(x),\qquad L_{r_n}(x_n)\to0.
\]
Pass to a further subsequence on which both $r_n$ and $L_{r_n}(x_n)$
decrease to zero, so that \eqref{eq:nonwp-weak} applies to both sequences.
Set $D_n(h)=\operatorname{sgn}(\delta_\varphi(x_n,r_n))Z_{x_n,r_n}(h)$.
By \eqref{eq:nonwp-weld-moment} and \eqref{eq:nonwp-weak},
\[
 \sup_n\|D_n\|_{L^2(T_*P)}<\infty,\qquad
 D_n\rightharpoonup0\quad\text{in }L^2(P).
\]
Using \eqref{eq:sing-decomposition}, put
\[
 N_n(h)=\operatorname{sgn}(\delta_\varphi(x_n,r_n))
 Z_{M_{r_n}(x_n),L_{r_n}(x_n)}(h),\qquad
 Y_n(h)=R_\varphi(x_n,r_n;h).
\]
Equation \eqref{eq:nonwp-weak} also gives $N_n\rightharpoonup0$ in $L^2(P)$.
Proposition~\ref{prop:nonwp-weighted}(ii) gives
\[
 \liminf_n E[w(H)Y_n(H)]\ge c_0^2\varepsilon E[w(H)].
\]
Lemma~\ref{lem:sing-criterion} again yields $T_*P\perp P$.

Finally suppose $\varphi\in\Sy(\m S^1)\setminus\WP(\m S^1)$. Use the
sequence $N_k$ in \eqref{eq:nonwp-cumulative-input}, and put
$D_k=T_{N_k}/E_{N_k}$. The first bound in that equation gives
\[
 \|D_k\|_{L^2(P)}+
 \left\|\frac{B_{N_k}(H)}{E_{N_k}}\right\|_2
 \le \frac{C}{\sqrt{E_{N_k}}}\longrightarrow0.
\]
By \eqref{eq:sing-integrated},
\[
 D_k\circ T=\frac{B_{N_k}}{E_{N_k}}+\frac{A_{N_k}}{E_{N_k}}.
\]
The two bounds in \eqref{eq:nonwp-cumulative-input} therefore give
\[
 \sup_k\|D_k\|_{L^2(T_*P)}<\infty.
\]
Proposition~\ref{prop:nonwp-weighted}(iii) gives, for bounded $w\ge0$,
\[
 \liminf_k E\!\left[w(H)\frac{A_{N_k}(H)}{E_{N_k}}\right]
 \ge c_0^2E[w(H)].
\]
Lemma~\ref{lem:sing-criterion} applies with the weakly convergent term
$B_{N_k}/E_{N_k}$ and the nonnegative term $A_{N_k}/E_{N_k}$.
Thus $T_*P\perp P$.

Lemma~\ref{introlem:sle-gmc-lift} gives
$\m P^{\mathrm{weld}}\sim\m P$. Mutual absolute continuity is preserved by measurable
pushforward, so $(R_\varphi)_*\m P^{\mathrm{weld}}\sim(R_\varphi)_*\m P$.
The singularity therefore transfers to the SLE welding law.
\end{proof}

\begin{cor}\label{cor:one-sided-dichotomy}
For $0<\kappa<4$, the one-sided non-singularity class of
Section~\ref{subsec:two-sided-singularity} satisfies
\[
 \widetilde{\QIR}=\WP(\m S^1).
\]
\end{cor}
\begin{proof}
Since $\WP(\m S^1)$ is a group, $\varphi^{-1}\notin\WP(\m S^1)$ whenever
$\varphi\notin\WP(\m S^1)$. The inclusion $\widetilde{\QIR}\subset\WP(\m S^1)$ follows from
Theorem~\ref{thm:sing-conditional}.

The reverse inclusion follows from WP quasi-invariance
\cite[Corollary~4.12]{F1}, with the Haar lift of
Definition~\ref{introdef:sle-haar-lift}.
\end{proof}

\section{A local approach to the SLE$_4$ welding law}\label{sec:critical-local}

At $\kappa=4$, we prove singularity by comparing the SLE welding law
with a critical GMC model on closed proper arcs. The two restriction
laws are mutually absolutely continuous, so a separating event for the model also separates
the SLE laws. The proof has three cases: maps outside $\QS(\m S^1)$, maps in
$\QS(\m S^1)\setminus\Sy(\m S^1)$, and maps in $\Sy(\m S^1)\setminus\WP(\m S^1)$. Local distortion
detects the first two cases. In the last case, normalized restrictions
approach the identity, and we instead use the distortion accumulated
over scales. The same local statistics also give a proof for each fixed
$0<\kappa<4$, using the mass estimates from Section~\ref{subsec:subcritical-gmc};
see Remark~\ref{rem:subcritical-shape-proof}.

We identify $\m S^1$ with $\m R/\m Z$ and use normalized Lebesgue
length and degree-one lifts. All arcs are oriented; a closed proper arc
has length in $(0,1)$. Write $R_t(x)=x+t$, $R_\varphi(h)=h\circ\varphi$,
and $\operatorname{Res}_I(h)=h|_I$. Let $\mathcal X_I$ be the space of
orientation-preserving continuous embeddings $I\to\m S^1$, with the
uniform topology; restriction laws are Borel measures on this space.
Distances on the circle are measured by normalized arc length.
In conformal coordinates, $R_t$
also denotes $z\mapsto e^{2\pi\mathrm{i}t}z$. The fields and estimates
below are defined at criticality, independently of the subcritical
constructions. We use $e_j^{(4)},s_j^{(4)},\Omega_N^{(4)}$ for the critical
layer energies and truncated parameter domain, and
$\mathcal E_N,\mathcal S_N,\mathcal T_N^\beta,\mathcal B_N^\beta$
for the corresponding cumulative quantities and statistics. These differ
from $e_j,s_j,\Omega_N,E_N,S_N,T_N,B_N$ in
Sections~\ref{subsec:subcritical-wp} and~\ref{sec:subcritical-singularity}.

\subsection{The laws and the local comparison theorem}\label{subsec:crit-laws}

We define the two laws and prove their local comparison in
Section~\ref{subsec:crit-comparison}. The estimates in
Sections~\ref{subsec:crit-shape}--\ref{subsec:crit-mass} yield local
regularity in Section~\ref{subsec:crit-local-regularity}. The accumulated
distortion estimates of Section~\ref{subsec:crit-energy} then complete
the singularity proof in Section~\ref{subsec:crit-singularity}.

Let $X$ be the mean-zero Gaussian generalized function on the circle
with logarithmic covariance, and let $X_1,X_2$ be independent copies:
\[
 E[X_i(x)X_i(y)]=-2\log|e^{2\pi\mathrm{i}x}-e^{2\pi\mathrm{i}y}|.
\]
For a smooth local mollification $X_\epsilon$, define the subcritical
and critical GMC measures, with normalized circle Lebesgue measure, by
\begin{align*}
 \mc M_X^\gamma(\dd x)
 &=\lim_{\epsilon\downarrow0}
 \exp\!\left(\frac\gamma2X_\epsilon(x)
 -\frac{\gamma^2}{8}E[X_\epsilon(x)^2]\right)\dd x,
 \qquad 0<\gamma<2,\\
 \mc M_X&=\lim_{\gamma\uparrow2}\frac{\mc M_X^\gamma}{2(2-\gamma)}.
\end{align*}
These are weak limits in probability. Throughout this section,
$\mc M_X$ without a superscript denotes the critical measure in the
second line. For a different field $Y$ in a fixed circle chart, $\mc M_Y^\gamma$ uses the same deterministic variance
subtraction $E[X_\epsilon(x)^2]$, rather than the variance of $Y$.
The critical limit uses the same denominator. The local constructions and
their independence of mollification are standard GMC results; see
\cite[Theorems~2.5 and~3.1]{PowellReview} and
\cite[Remark~3.5]{F1}. Set $\mc M_i=\mc M_{X_i}$. These measures are finite,
nonzero, non-atomic and have full support almost surely. Put
\[
 \mu_i=\frac{\mc M_i}{\mc M_i(\m S^1)},\qquad
 F_i(x)=\mu_i([0,x])\quad(0\le x\le1),\qquad
 F_i(x+1)=F_i(x)+1.
\]
The functions $F_i$ are increasing degree-one homeomorphisms of $\m R$.
Let $U$ be uniform on $\m R/\m Z$ and independent of the fields. For
$\varphi\in\Homeo_+(\m S^1)$, define
\[
 H_4=F_2^{-1}\circ R_U\circ F_1,\qquad
 \m P_4=\operatorname{Law}(H_4),\qquad
 \m P_{4,\varphi}=(R_\varphi)_*\m P_4.
\]
The composition defining $H_4$ is read modulo $1$, so it does not depend
on the representative of $U$. Write $\m E_4$ for expectation under
$\m P_4$; $P$ and $E$ refer to the underlying fields and auxiliary
random variables in each construction. The variable $U$ is
uniform in the \emph{mass coordinate}. In particular, conditional on
$\mu_2$, the law of $H_4(0)$ is $\mu_2$.

Let $\mu^0$ be the oriented unrooted SLE$_4$ loop measure of
\cite[Theorem~4.2]{zhan2020sleloop}, and let $\mu_4^\#$ be the shape
probability of \cite[Definition~4.7]{F1}. We use its independent
left-Haar welding lift $\m P_4^{\mathrm{weld}}$ from
Definition~\ref{introdef:sle-haar-lift}: sample the normalized Riemann
maps $f,g$ with $f(1)=g(1)$ and then rotate $g^{-1}\circ f$ on the left
by an independent uniform angle.
Write
\[
 \m P_{4,\varphi}^{\mathrm{weld}}
 =(R_\varphi)_*\m P_4^{\mathrm{weld}},\qquad
 P_I=(\operatorname{Res}_I)_*\m P_4,\qquad
 P_I^{\mathrm{weld}}
 =(\operatorname{Res}_I)_*\m P_4^{\mathrm{weld}}.
\]
\begin{thm}[Mutual absolute continuity on arcs]\label{thm:crit-local-bridge}
For every closed proper arc $I$,
\begin{equation}\label{eq:crit-local-bridge}
 P_I^{\mathrm{weld}}\sim P_I.
\tag{LB}
\end{equation}
\end{thm}
The main conclusion of this section is the following local singularity
theorem and its consequence for SLE$_4$.

\begin{thm}[Critical singularity]\label{thm:crit-all-nonwp-total}
Let $\varphi\in\Homeo_+(\m S^1)\setminus\WP(\m S^1)$. There is a
closed proper arc $J$ such that
\[
 (\operatorname{Res}_J)_*\m P_{4,\varphi}\perp P_J.
\]
For the SLE$_4$ welding law this implies both
\[
 (h\mapsto h\circ\varphi)_*\m P_4^{\mathrm{weld}}
 \perp\m P_4^{\mathrm{weld}},\qquad
 (h\mapsto\varphi\circ h)_*\m P_4^{\mathrm{weld}}
 \perp\m P_4^{\mathrm{weld}}.
\]
\end{thm}

The model law is invariant under left and right rotations and conjugation
by reflection: rotating the source or target field changes the cumulative
origin by a shift absorbed into the uniform phase, while reflection
replaces $U$ by $-U$. It is also invariant under inversion, since
\[
 H_4^{-1}=F_1^{-1}\circ R_{-U}\circ F_2
\]
and $X_1,X_2$ are independent with the same law.

Right composition commutes with left rotations and hence acts on their
quotient. Each orbit has a unique representative fixing $0$. Composing
this representative with an independent uniform rotation gives the
unique left-rotation invariant lift of any quotient law.

For a quotient law $\overline P$ and its left-Haar lift $P$,
\begin{equation}\label{eq:crit-haar-lift}
 (R_\varphi)_*\overline P\ll\overline P
 \quad\Longleftrightarrow\quad (R_\varphi)_*P\ll P.
\end{equation}
Projection proves one direction. For the converse, right composition
preserves the Haar conditional law on each orbit. Integrating this kernel
against the two quotient laws preserves absolute continuity.

\begin{prop}[Local transfer]\label{prop:crit-local-sing}\label{lem:crit-local-transfer}\label{lem:crit-local-inversion}
Assume \eqref{eq:crit-local-bridge}. For each closed proper arc $I$,
\[
 (\operatorname{Res}_I)_*\m P_{4,\varphi}^{\mathrm{weld}}
 \sim(\operatorname{Res}_I)_*\m P_{4,\varphi}.
\]
Consequently, $\m P_{4,\varphi}^{\mathrm{weld}}\ll\m P_4^{\mathrm{weld}}$
implies $(\operatorname{Res}_I)_*\m P_{4,\varphi}\ll P_I$ for every $I$,
whereas singularity of $(\operatorname{Res}_I)_*\m P_{4,\varphi}$ and $P_I$
for one such arc implies
$\m P_{4,\varphi}^{\mathrm{weld}}\perp\m P_4^{\mathrm{weld}}$.
More generally, if two laws on $\Homeo_+(\m S^1)$ have mutually
absolutely continuous restrictions to every closed proper arc, their
inverse-homeomorphism laws have the same property.
\end{prop}
\begin{proof}
Since $\varphi(I)$ is a closed proper arc,
\[
 \operatorname{Res}_I(h\circ\varphi)
 =\bigl(\operatorname{Res}_{\varphi(I)}h\bigr)\circ\varphi|_I.
\]
Pushing forward \eqref{eq:crit-local-bridge} proves the first assertion.
Absolute continuity follows by restriction; singularity follows by
pulling back a separating Borel set under $\operatorname{Res}_I$.

For inversion, let $P,Q$ have mutually absolutely continuous restrictions
to every closed proper arc. Fix such an arc $J$ and enumerate the proper
arcs $K_n$ with rational endpoints. Define
\[
 \mathcal U_n=\{u\in\mathcal X_{K_n}:J\subset\operatorname{int}u(K_n)\},
 \qquad E_n=\{h:h|_{K_n}\in\mathcal U_n\}.
\]
The sets $E_n$ cover the homeomorphism space. On the open set
$\mathcal U_n$, the inverse extraction map
\[
 \iota_n:\mathcal U_n\longrightarrow C(J,\m S^1),\qquad
 \iota_n(u)=u^{-1}|_J
\]
is continuous: uniform convergence of increasing interval lifts gives
uniform convergence of inverses on compact subsets of the limiting
image. Thus, for Borel $A\subset C(J,\m S^1)$,
\[
 E_n\cap\{h:h^{-1}|_J\in A\}
 =\operatorname{Res}_{K_n}^{-1}
     \{u\in\mathcal U_n:\iota_n(u)\in A\}.
\]
Mutual absolute continuity on $K_n$ makes this event null under $P$
exactly when it is null under $Q$. Taking the countable union in $n$
proves mutual absolute continuity of the inverse restriction laws.
\end{proof}

\subsection{Local absolute continuity of critical weldings}\label{subsec:crit-comparison}

Fix closed proper arcs $I=[a_0,a_1]$ and $B=[b_0,b_1]$, written in
oriented interval coordinates. We shall compare both restriction laws
with the same measure on
\[
 O_B=\{u\in\mathcal X_I:u(I)\subset\operatorname{int}B\}.
\]
For a full homeomorphism, $O_B$ denotes the event $h|_I\in O_B$.

\begin{definition*}[Local matching]
Let $\mathfrak M(J)$ be the space of finite, nonzero, atomless Borel
measures with full support on a closed arc $J$, with the weak topology.
For $\mc M_1\in\mathfrak M(I)$ and $\mc M_2\in\mathfrak M(B)$, set
\[
 C_1(x)=\mc M_1([a_0,x]),\quad C_2(y)=\mc M_2([b_0,y]),\qquad
 m_1=\mc M_1(I),\quad m_2=\mc M_2(B).
\]
Define
\begin{align*}
 \mathcal D&=\{(\mc M_1,\mc M_2,s)\in
      \mathfrak M(I)\times\mathfrak M(B)\times\mathbb R:
                  0<s<s+m_1<m_2\},\\
 \mathcal W(\mc M_1,\mc M_2,s)(x)&=C_2^{-1}(C_1(x)+s),\qquad x\in I.
\end{align*}
The map $\mathcal W:\mathcal D\to O_B$ matches equal mass increments.
If $u=\mathcal W(\mc M_1,\mc M_2,s)$, then $s=C_2(u(a_0))$ is the target
mass before the first matched point.

Let $\lambda_J=\operatorname{Law}((X+C)|_{J^\circ})$, where
$C\sim N(0,1)$ is independent of $X$, and let $\gamma_0=N(0,1)$ on
$\mathbb R$. For every Borel set $F\subset\mathcal X_I$, define
\begin{equation}\label{eq:crit-common-measure}
 \begin{aligned}
 \Sigma_{I,B}(F)=\int&
 \1_{\mathcal D}(\mc M_{Y_1},\mc M_{Y_2},s)
 \1_{\{\mathcal W(\mc M_{Y_1},\mc M_{Y_2},s)\in F\}}\\
 &\hspace{4mm}\times\lambda_I(\dd Y_1)\lambda_B(\dd Y_2)\gamma_0(\dd s).
 \end{aligned}
\end{equation}
Here $\mc M_{Y_i}$ is the local critical GMC with the normalization
of Section~\ref{subsec:crit-laws}; the integrand is set to zero outside
$\mathcal D$.
\end{definition*}

The map $\mathcal W$ is continuous: weak convergence to an atomless
measure gives uniform convergence of cumulative functions, and their
inverses converge uniformly on compact subintervals of $(0,m_2)$.
The local GMC is a measurable function of the field, so
\eqref{eq:crit-common-measure} defines a finite Borel measure on $O_B$.
We also use the continuous-shift rule
\begin{equation}\label{eq:crit-continuous-shift}
 \mc M_\epsilon^G\xrightarrow{P}\mc M^G,\quad
 q_\epsilon\xrightarrow{P}q\text{ uniformly}
 \quad\Longrightarrow\quad
 e^{q_\epsilon}\mc M_\epsilon^G\xrightarrow{P}e^q\mc M^G,
\end{equation}
where measure convergence is weak and $q$ is continuous. No independence
is required. See \cite[Theorem~2.5]{PowellReview} for the local critical
construction and its positive approximations.

\begin{prop}[Local comparison]\label{prop:crit-geometric-branch-proved}
For every pair $I,B$ as above,
\begin{equation}\label{eq:crit-geometric-branch}
 P_I|_{O_B}\sim\Sigma_{I,B}\sim P_I^{\mathrm{weld}}|_{O_B}.
\end{equation}
\end{prop}

\begin{lem}\label{lem:crit-common-model}
For every pair of closed proper arcs $I,B$,
$P_I|_{O_B}\sim\Sigma_{I,B}$.
\end{lem}
\begin{proof}
We need the conditional law of a total mass. Let $G=X+C$,
$Y=G|_{J^\circ}$ and $m=\mc M_G(J)$. Choose a smooth nonnegative
Cameron--Martin function $e$, zero on $J$ and positive off $J$.
After normalizing $e$, write $G=\mathcal Y+ae$, with $a\sim N(0,1)$
independent of $\mathcal Y$. The field $Y$ is determined by
$\mathcal Y$. Conditional on $\mathcal Y$,
\[
 T(a)=m+\int_{\m S^1\setminus J}e^{ae(x)}\mc M_{\mathcal Y}(\dd x),
 \qquad T'(a)>0,\qquad T(\mathbb R)=(m,\infty).
\]
Thus $T=\mc M_G(\m S^1)$ has conditional density
$\phi(T^{-1}(t))/T'(T^{-1}(t))$ on $(m,\infty)$, where $\phi$ is the
standard Gaussian density. Averaging over $\mathcal Y$ given $Y$ gives
a density $f(t\mid Y)>0$ for almost every $t>m$.

Apply this to $G_i=X_i+C_i^0$, where $C_1^0,C_2^0$ are independent
$N(0,1)$ constants independent of $X_1,X_2,U$. They cancel in the
normalized GMCs. Write $T_i=\mc M_{G_i}(\m S^1)$ and
$Y_1=G_1|_{I^\circ}$, $Y_2=G_2|_{B^\circ}$, and let $C_i,m_i$ be
the cumulative and total masses of $\mc M_{Y_i}$. On $O_B$,
\[
 C_2(H_4(x))=rC_1(x)+s,\qquad
 r=T_2/T_1,\qquad 0<s<s+rm_1<m_2.
\]
The mass phase has density $1/T_2$. The substitution
$(T_1,T_2)=(t,rt)$ has Jacobian $t$, and hence the conditional density
of $(r,s)$ on this region is
\[
 \frac1r\int_{\max(m_1,m_2/r)}^\infty
           f_1(t\mid Y_1)f_2(rt\mid Y_2)\,\dd t>0
 \quad\text{a.e.}
\]
For fixed $r>0$, translation $Y_1\mapsto Y_1+\log r$ preserves the
law's null sets and gives $\mc M_{Y_1+\log r}=r\mc M_{Y_1}$.
Fubini's theorem therefore identifies the null sets of the matching
data with those in \eqref{eq:crit-common-measure}. Applying
$\mathcal W$ proves the claim.
\end{proof}

\begin{definition*}[Rooted loops and stopping]
Write $\sigma_\eta$ for the $3/2$-dimensional Minkowski content measure
of an SLE$_4$ loop, in the normalization
\[
 \sigma_\eta(A)=\lim_{\epsilon\downarrow0}\epsilon^{-1/2}
   \operatorname{Area}\{z:\operatorname{dist}(z,\eta\cap A)<\epsilon\},
 \qquad L(\eta)=\sigma_\eta(\eta).
\]
For compact $A$, this is the convention of
\cite[Section~2.3]{zhan2020sleloop}. Let
$\nu_{z\rightleftarrows w}^{\#}$ be the two-sided whole-plane SLE$_4$
law: whole-plane SLE$_4(2)$ from $z$ to $w$, followed by chordal
SLE$_4$ back to $z$ in the complement. Zhan's rooted measure is
\[
 \mu_z^1(\dd\eta)=L(\eta)^{-1}\int_{\mathbb C\setminus\{z\}}
             |w-z|^{-1}\nu_{z\rightleftarrows w}^{\#}(\dd\eta)\,\dd^2w.
\]
Forgetting the starting point retains orientation. With this convention,
\cite[Theorems~4.1--4.2, (4.2) and (4.6)]{zhan2020sleloop} gives
\begin{equation}\label{eq:crit-campbell}
 \mu^0(\dd\eta)\sigma_\eta(\dd z)
   =\dd^2z\,(\operatorname{unroot}_*\mu_z^1)(\dd\eta).
\end{equation}
Both sides are measures on oriented unrooted loops with a marked point.

Set $r_n=2^{-n}$, $n\ge1$. For a loop rooted at $z$, let $\tau_n$ be
its first exit time from $B(z,r_n)$ and $T_z$ its return time to $z$.
Put
\[
 E_n=\{\tau_n<T_z\},\qquad
 c_n=\mu_z^1(E_n)=c_*r_n^{-1/2}\in(0,\infty).
\]
The last equality follows from
\cite[Theorem~4.1(v)--(vi)]{zhan2020sleloop}. Fix a bounded, strictly
positive probability density $\rho$ on $\mathbb C$ and $a_n>0$ with
$\sum_na_n=1$. Define
\[
 \mathbf Q(\dd z,\{n\},\dd\eta)
  =\rho(z)\,\dd^2z\,a_nc_n^{-1}\1_{E_n}\mu_z^1(\dd\eta),
 \qquad K=\eta[0,\tau_n].
\]
The arc $K$ retains its root, tip, and orientation.
\end{definition*}

\begin{lem}[Stopped loops]\label{lem:crit-stopped-root}\label{lem:crit-root-avoid}
The measure $\mathbf Q$ is a probability. Its unrooted loop marginal
$Q_{\mathrm{loop}}$ is mutually absolutely continuous with respect to
$\mu^0$. Conditional on $K$, the remaining arc has chordal SLE$_4$
law $\Pi_K$ in $\widehat{\mathbb C}\setminus K$, from the tip to the root.

Suppose additional marks $\mathbf m$ are conditionally independent of
$(z,n)$ given the oriented unrooted loop. For every closed proper subarc
$C=C(\eta,\mathbf m)$,
\[
 \mathbf Q(K\cap C=\varnothing\mid\eta,\mathbf m)>0\quad\text{a.s.}
\]
\end{lem}
\begin{proof}
The normalization in the definition gives
\[
 \mathbf Q(1)=\int_{\mathbb C}\rho(z)\sum_n
 a_nc_n^{-1}\mu_z^1(E_n)\,\dd^2z=1.
\]
By \eqref{eq:crit-campbell},
\[
 \frac{\dd Q_{\mathrm{loop}}}{\dd\mu^0}(\eta)=D(\eta)
 :=\int_\eta\rho(z)\sum_n a_nc_n^{-1}\1_{E_n(\eta,z)}
                                      \,\sigma_\eta(\dd z).
\]
Every root on a nondegenerate loop satisfies $E_n$ for all sufficiently
large $n$. Since natural length is finite, positive, and has full support,
\[
 0<D(\eta)\le\|\rho\|_\infty L(\eta)
            c_*^{-1}\sum_na_n2^{-n/2}<\infty\quad\mu^0\text{-a.e.}
\]
This proves mutual absolute continuity. A single fixed $n$ would exclude
loops of diameter less than $r_n$.

Let $q_{z,n}$ be the law of $K$ under $c_n^{-1}\1_{E_n}\mu_z^1$.
Zhan's stopping formula \cite[Theorem~4.1(ii), (4.3)]{zhan2020sleloop}
gives, for the remaining arc $\eta^K$,
\[
 \mathbf Q(\dd z,\{n\},\dd K,\dd\eta^K)
 =\rho(z)\,\dd^2z\,a_nq_{z,n}(\dd K)\Pi_K(\dd\eta^K).
\]
The same kernel $\Pi_K$ occurs in every summand, so integrating out
$(z,n)$ preserves the conditional continuation law.

Finally, the conditional law of the root and index given $\eta,\mathbf m$
is
\[
 \frac{\rho(z)a_nc_n^{-1}\1_{E_n(\eta,z)}}{D(\eta)}
                                      \,\sigma_\eta(\dd z).
\]
The open arc $\eta\setminus C$ has positive natural length. For each
root there, all sufficiently large $n$ satisfy
$E_n$ and $r_n<\operatorname{dist}(z,C)$, so
$K\subset\overline{B(z,r_n)}$ avoids $C$. Integrating the positive
conditional density over these pairs proves the last assertion.
\end{proof}

\begin{definition*}[Critical wedge fields]
Let $\mathcal S=\mathbb R\times(0,\pi)$ and $\phi(w)=-e^{-w}$.
Let $h_{\mathrm{lat}}$ be the centered Gaussian generalized function
with covariance
\[
 E[h_{\mathrm{lat}}(s+\mathrm{i}\theta)h_{\mathrm{lat}}(t+\mathrm{i}\vartheta)]
 =2\sum_{n\ge1}\frac{e^{-n|s-t|}}n\cos(n\theta)\cos(n\vartheta).
\]
For $\alpha\in\{1,2\}$, independently choose a continuous radial process
$R_\alpha$ as follows:
\[
 \begin{array}{c|cc}
       &s\ge0&s<0\\ \hline
 R_2(s)&-\mathcal R_{2s}&\widehat B_{-2s}\\
 R_1(s)&B_{2s}-s\text{ conditioned to stay negative}&\widehat B_{-2s}-s.
 \end{array}
\]
Here $\mathcal R$ is zero-start BES$(3)$; $B,\widehat B$ are standard
Brownian motions; the two half-line processes are independent. The
conditioning is the entrance law obtained from negative initial values
tending to zero. Set
\[
 \widehat h_\alpha=R_\alpha+h_{\mathrm{lat}},\qquad
 h_\alpha=\widehat h_\alpha\circ\phi^{-1}
                      +2\log| (\phi^{-1})'|,\qquad
 \mathsf W_\alpha=\operatorname{Law}(h_\alpha),\quad\mathsf W=\mathsf W_2.
\]
Thus $R_\alpha$ is the cross-section average:
\[
 \frac1\pi\int_0^\pi\widehat h_\alpha(s+\mathrm{i}\theta)\,\dd\theta
       =R_\alpha(s).
\]
The integral is defined by regularization. The field
$h_{\mathrm{lat}}$ is the zero-cross-section-average component of a
Neumann GFF, hence the term lateral field. The last zero of $R_\alpha$
is at $0$; this fixes the last-exit representative.

The marked surface $(\mathbb H,h_\alpha;0,\infty)$ represents the
$(2,\alpha)$ critical wedge. Under a change of coordinates,
\[
 h'=h\circ A+2\log|A'|,
\]
where the conformal map $A$ preserves the two marked points.
This is the convention of
\cite[Definitions~2.2--2.5 and Remark~2.7]{HP21_critical}.
Its critical boundary length is denoted by $\nu_h$. On a compact interval
away from the marked ends,
\[
 \nu_h(\dd x)=\lim_{\epsilon\downarrow0}
   \left(\log(1/\epsilon)-\tfrac12h_\epsilon(x)\right)
                 \epsilon e^{h_\epsilon(x)}\,\dd x,
 \qquad A^*\nu_h=\nu_{h\circ A+2\log|A'|},
\]
where $h_\epsilon(x)$ is the semicircle average and convergence is weak
in probability. See \cite[Section~2.3]{HP21_critical}.
\end{definition*}

\begin{lem}[Boundary fields of a wedge]\label{lem:crit-wedge-boundary}
Let $N$ be the restriction of $h_{\mathrm{lat}}$ to either strip boundary.
Let $B$ be two-sided Brownian motion with variance parameter $2$,
independent of $N$. On every bounded interval, the laws of $N+R_2$ and
$N+B$ are mutually absolutely continuous.
\end{lem}
\begin{proof}
The covariance $K_N(s)=2\sum_{n\ge1}e^{-n|s|}/n$ has Fourier transform
\[
 \widehat K_N(\xi)=4\sum_{n\ge1}(n^2+\xi^2)^{-1}
                         \asymp(1+\xi^2)^{-1/2}.
\]
Its Cameron--Martin space is therefore $H^{1/2}(\mathbb R)$.
Couple $\mathcal R$ to a standard Brownian motion $W$ by
$\mathcal R_t=W_t+\int_0^t\mathcal R_u^{-1}\,\dd u$. Then
$R_2=B+V$, where
\[
 B(s)=\begin{cases}-W_{2s},&s\ge0,\\\widehat B_{-2s},&s<0,\end{cases}
 \qquad
 V(s)=\begin{cases}-\int_0^{2s}\mathcal R_u^{-1}\,\dd u,&s\ge0,\\0,&s<0.
 \end{cases}
\]
The estimate $E[\mathcal R_u^{-1}]=Cu^{-1/2}$ gives local bounded
variation of $V$. Also $V$ is locally $\alpha$-H\"older for every
$\alpha<1/2$. If $\chi$ is a smooth compact cutoff and $v=\chi V$, then
\[
 \|v(\cdot+t)-v\|_2^2
 \le\|v(\cdot+t)-v\|_1\|v(\cdot+t)-v\|_\infty
 \le C\operatorname{Var}(v)t^{1+\alpha},
\]
so $v\in H^{1/2}(\mathbb R)$. Conditional Cameron--Martin, using the
independence of $N$ from $(B,V)$, gives
\[
 \operatorname{Law}(N+v,B,V)\sim\operatorname{Law}(N,B,V).
\]
Choose $\chi=1$ on the interval in question, add $B$, and restrict.
\end{proof}

\begin{definition*}[Conformal coordinates and local boundary data]
For $u\in\mathbb H$, put
\[
 \mathcal C_u^{\mathrm{in}}(z)=\operatorname{Re}u+
       \mathrm{i}\operatorname{Im}u\frac{1+z}{1-z},\qquad
 \mathcal C_u^{\mathrm{out}}(z)=\operatorname{Re}u+
       \mathrm{i}\operatorname{Im}u\frac{z+1}{z-1}.
\]
Every conformal map $\mathbb D\to\mathbb H$ is uniquely
$\mathcal C_u^{\mathrm{in}}\circ R_\theta$; every conformal map
$\mathbb D^*\to\mathbb H$ is uniquely
$\mathcal C_u^{\mathrm{out}}\circ R_\theta$. Fix probability laws
\[
 \Lambda_1=\operatorname{Law}(\mathcal C_{U_1}^{\mathrm{in}}\circ R_{\Theta_1}),
 \qquad
 \Lambda_2=\operatorname{Law}(\mathcal C_{U_2}^{\mathrm{out}}\circ R_{\Theta_2}),
\]
where $U_i$ has a strictly positive area density on $\mathbb H$ and
$\Theta_i$ is independent Haar. Only their null sets will be used.
Set
\[
 \mathcal A_{I,B}=\{(A_1,A_2):
        A_1(I)\Subset(0,\infty),\ A_2(B)\Subset(-\infty,0)\}.
\]
For such a pair and half-plane fields $h_1,h_2$, define
\begin{align*}
 Y_1&=(h_1\circ A_1+2\log|A_1'|)|_{I^\circ},&
 Y_2&=(h_2\circ A_2+2\log|A_2'|)|_{B^\circ},\\
 L_1(t)&=(2\pi)^{-1}\nu_{h_1}([0,t]),\quad t>0,&
 L_2(t)&=(2\pi)^{-1}\nu_{h_2}([t,0]),\quad t<0,\\
 s&=L_1(A_1(a_0))-L_2(A_2(b_0)).
\end{align*}
The factor $(2\pi)^{-1}$ expresses lengths in the units of circle GMC.
Boundary restrictions are generalized functions defined by local
regularization.
\end{definition*}

\begin{definition*}[The joint probability law]
Extend $\mathbf Q$ once, and denote the resulting law by $\mathbf Q^*$,
using the following conditional kernels.
\begin{enumerate}[label=(\roman*)]
\item Given $K$, choose a measurable conformal map
$\Psi_K:\mathbb H\to\widehat{\mathbb C}\setminus K$ taking $0$ to its
tip and $\infty$ to its root. Set
$\zeta=\Psi_K^{-1}(\eta\setminus K)$. Independently sample
$\mathcal H\sim\mathsf W_1$ and cut along $\zeta$. Let $j_i$ be the
maps from the last-exit representatives of the two resulting surfaces
onto the two components of $\mathbb H\setminus\zeta$. Set
\[
 h_i=\mathcal H\circ j_i+2\log|j_i'|,\qquad
 \Phi_i=\Psi_K\circ j_i:\mathbb H\to D_i.
\]
The domains $D_1,D_2$ lie to the left and right of the oriented loop.
The rays $(0,\infty)$ on side $1$ and $(-\infty,0)$ on side $2$ map
to the continuation; the other rays map to $K$. In particular,
$\Phi_i(0)=K_{\mathrm{tip}}$ and $\Phi_i(\infty)=z$.
\item Conditional on $\eta$, choose $p\in D_1$ and $q\in D_2$ from
normalized spherical area, independently of each other and of the
fields and stopping data given $\eta$. Choose independent Haar
$\alpha,\beta$. Let $M_0$ take $p,q$ to $0,\infty$ and normalize the
inner conformal radius to one. One measurable choice is
\[
 B_{p,q}(z)=\frac{z-p}{z-q},\qquad
 M_0=\frac{B_{p,q}}{\operatorname{CR}(B_{p,q}(D_1),0)},
 \qquad M=R_\alpha\circ M_0.
\]
The marks are finite almost surely. Let
$f:\mathbb D\to M(D_1)$ and $g:\mathbb D^*\to M(D_2)$ satisfy
\[
 f(0)=0,\quad f'(0)=1,\quad g(\infty)=\infty,\quad g(1)=f(1).
\]
\item Define
\[
 g_\beta=g\circ R_{-\beta},\quad h=g_\beta^{-1}\circ f,\qquad
 A_1=(M\circ\Phi_1)^{-1}\circ f,\quad
 A_2=(M\circ\Phi_2)^{-1}\circ g_\beta.
\]
\end{enumerate}
Writing $\omega=(z,n,\eta)$ and $\mathsf a_D$ for normalized spherical
area on $D$, the joint law is
\begin{equation}\label{eq:crit-joint-law}
\begin{aligned}
 &\mathbf Q^*(\dd\omega,\dd\mathcal H,\dd p,\dd q,\dd\alpha,\dd\beta)\\
 &\quad=\mathbf Q(\dd\omega)\,
 \mathsf W_1(\dd\mathcal H)\,
 \mathsf a_{D_1(\eta)}(\dd p)\,
 \mathsf a_{D_2(\eta)}(\dd q)\,
 \dd\alpha\,\dd\beta,
\end{aligned}
\end{equation}
where $\dd\alpha$ and $\dd\beta$ are Haar probability measures on
$\mathbb R/\mathbb Z$. The other variables above are measurable functions
of these variables, with the conformal maps chosen as specified.
Let $\pi$ denote the projection onto the original variables $(z,n,\eta)$,
forgetting the added field and marks. Since all added conditional kernels
are probability kernels, the original marginal is preserved:
\begin{equation}\label{eq:crit-joint-projection}
 \pi_*\mathbf Q^*=\mathbf Q.
\end{equation}
\end{definition*}

By \cite[Theorem~1.2 and Corollary~1.4]{HP21_critical}, independent
chordal SLE$_4$ cuts a $(2,1)$ wedge into independent $(2,2)$ surfaces
whose boundary lengths agree on the interface. Thus their last-exit
representatives satisfy
\begin{equation}\label{eq:crit-cut-field-law}
 \operatorname{Law}_{\mathbf Q^*}(h_1,h_2)=\mathsf W^{\otimes2}.
\end{equation}

\begin{lem}\label{lem:crit-marked-loop-class}
Under $\mathbf Q^*$, the law of $h$ is mutually absolutely continuous
with respect to $\m P_4^{\mathrm{weld}}$.
\end{lem}
\begin{proof}
Independently add a real variable $T$ with strictly positive density
and set $\widetilde M=e^T M$. For fixed $\eta$, the parameters
$(p,q,\alpha,T)$ give a strictly positive density relative to Haar
measure $\dd A$ on the sphere's M\"obius transformations satisfying
$A^{-1}(0)\in D_1(\eta)$ and $A^{-1}(\infty)\in D_2(\eta)$.
Together with $Q_{\mathrm{loop}}\sim\mu^0$, this gives
\[
 \operatorname{Law}(\eta,\widetilde M)
 \sim\1_{\{A^{-1}(0)\in D_1(\eta),\ A^{-1}(\infty)\in D_2(\eta)\}}
                 \mu^0(\dd\eta)\,\dd A.
\]
Set $\eta'=\widetilde M\eta$. Under the change of variables
$(\eta,A)\mapsto(A\eta,A)$, M\"obius invariance turns the reference measure into
\[
 \1_{\{0\in D_1(\eta'),\ \infty\in D_2(\eta')\}}
                 \mu^0(\dd\eta')\,\dd A.
\]
Hence the marginal law of $\eta'$ is mutually absolutely continuous
with respect to the oriented separating-loop measure. Forgetting
orientation and using its scale decomposition $\dd r/r$ times
$\mu_4^\#$ shows that, regarding the normalized loop as unoriented,
\[
 \operatorname{Law}\!\left(
 \frac{\eta'}{\operatorname{CR}(\eta',0)}\right)\sim\mu_4^\#.
\]
Conformal-radius normalization cancels the added factor $e^T$, so this
normalized loop is $M\eta$. Its normalized welding, followed by the
independent left rotation $\beta$, is precisely $h$. The definition
of $\m P_4^{\mathrm{weld}}$ proves the assertion.
\end{proof}

\begin{lem}\label{lem:crit-joint-charts}\label{lem:crit-chart-density}
Under $\mathbf Q^*$,
\[
 \operatorname{Law}(h_1,h_2,A_1,A_2)
       \sim\mathsf W^{\otimes2}\otimes\Lambda_1\otimes\Lambda_2.
\]
\end{lem}
\begin{proof}
Let $\mathcal G=\sigma(z,n,K,\zeta,\mathcal H,j_1,j_2)$. Conditional
on $\mathcal G$, the centers $u_1=\Phi_1^{-1}(p)$ and
$u_2=\Phi_2^{-1}(q)$ have law
\[
 (\Phi_1^{-1})_*\mathsf a_{D_1}\otimes
 (\Phi_2^{-1})_*\mathsf a_{D_2}
 \sim\dd^2u_1\,\dd^2u_2.
\]
Fix these centers. Let $f_0,g_0$ be the normalized maps for $M_0\eta$,
and let $H_0$ be the angular lift of $g_0^{-1}\circ f_0$ fixing zero.
The normalizations give
\begin{equation}\label{eq:crit-rotated-maps}
 f_\alpha=R_\alpha\circ f_0\circ R_{-\alpha},\qquad
 g_\alpha=R_\alpha\circ g_0\circ R_{H_0(-\alpha)}.
\end{equation}
Indeed, these maps satisfy $f_\alpha'(0)=1$ and
$g_\alpha(1)=f_\alpha(1)$. Define $c_1,c_2$ by
\[
 (M_0\circ\Phi_1)^{-1}\circ f_0
       =\mathcal C_{u_1}^{\mathrm{in}}\circ R_{c_1},\qquad
 (M_0\circ\Phi_2)^{-1}\circ g_0
       =\mathcal C_{u_2}^{\mathrm{out}}\circ R_{c_2}.
\]
The chart phases therefore satisfy
\[
 \theta_1=c_1-\alpha,\qquad
 \theta_2=c_2+H_0(-\alpha)-\beta,
\]
where $c_i$ are fixed given $\mathcal G,u_1,u_2$. Integrating first in
$\beta$ and then in $\alpha$ shows that
\[
 \operatorname{Law}(\theta_1,\theta_2\mid\mathcal G,u_1,u_2)
       =\dd\theta_1\,\dd\theta_2.
\]
Thus the conditional chart law is mutually absolutely continuous with
respect to $\Lambda_1\otimes\Lambda_2$. Integrating over $\mathcal G$
and using \eqref{eq:crit-cut-field-law} proves the claim.
\end{proof}

\begin{lem}[Local fields and lengths]\label{lem:crit-wedge-local-fields}\label{lem:crit-wedge-offset}
Fix $(A_1,A_2)\in\mathcal A_{I,B}$ and sample
$(h_1,h_2)\sim\mathsf W^{\otimes2}$. For the local data defined above,
\begin{align*}
 \operatorname{Law}(Y_1,Y_2)&\sim\lambda_I\otimes\lambda_B,\\
 A_i^*\nu_{h_i}&=2\pi\mc M_{Y_i}\quad\text{on the observed arc},\\
 \operatorname{Law}(Y_1,Y_2,s)&\sim\lambda_I\otimes\lambda_B\otimes\dd s.
\end{align*}
\end{lem}
\begin{proof}
For the field law, apply Lemma~\ref{lem:crit-wedge-boundary} and the
fixed strip-to-half-plane coordinate change. In disk coordinates,
write a free GFF as $\Gamma_0+PX$, with independent zero-boundary
$\Gamma_0$ and harmonic extension $PX$. For the fixed chart $A$, set
\[
 \operatorname{Av}_A(F)=\frac1\pi\int_0^\pi
                         F(A^{-1}(e^{\mathrm{i}\theta}))\,\dd\theta.
\]
After fixing this semicircle average to zero, its boundary trace is
\[
 X-\operatorname{Av}_A(PX)-\operatorname{Av}_A(\Gamma_0).
\]
The last term is an independent nondegenerate Gaussian. Conditional
on $X$, the additive constant has a positive density on $\mathbb R$.
Smooth local coordinate terms are Cameron--Martin shifts. This gives
the first assertion.

For the length identity, let $G$ be the semicircle-gauge Neumann GFF,
$g$ its trace, and $Z=g\circ A_i$. Its variance satisfies
\[
 \operatorname{Var}(G_\epsilon(x))=2\log(1/\epsilon)+\ell(x)+o(1),
 \qquad\ell(x)=4\log\max(1,|x|).
\]
The finite variance correction in the circle chart is
$d_A=-2\log|A_i'|+\ell\circ A_i$. For the subcritical boundary length
$\nu_G^\gamma(\dd x)=\lim_{\epsilon\downarrow0}
\epsilon^{\gamma^2/4}e^{\gamma G_\epsilon(x)/2}\,\dd x$, this gives
\begin{align*}
 A_i^*\nu_G^\gamma
 &=2\pi|A_i'|e^{\gamma^2(\ell\circ A_i-d_A)/8}\mc M_Z^\gamma
  =2\pi|A_i'|^{1+\gamma^2/4}\mc M_Z^\gamma,\\
 A_i^*\nu_G&=2\pi|A_i'|^2\mc M_Z
                     =2\pi\mc M_{Z+2\log|A_i'|}.
\end{align*}
The last line follows by division by $2(2-\gamma)$ and the critical
limit. The deterministic weights converge uniformly on the observed
arc. See \cite[Remarks~3.5--3.6 and Lemma~3.8]{F1} and
\cite[Section~2.1 and Definition~2.6]{PowellSepulveda}.
In the coupling of Lemma~\ref{lem:crit-wedge-boundary}, the strip field
$B+h_{\mathrm{lat}}$ gives $G$ in semicircle gauge, and the wedge differs
from $G$ by $q=V\circ\operatorname{Re}\phi^{-1}+2\log|(\phi^{-1})'|$.
On the observed arc the regularizations of $q$ converge uniformly.
Both lengths acquire the factor $e^{q\circ A_i}$ by
\eqref{eq:crit-continuous-shift}; the additional term in the derivative
regularization vanishes by \cite[Remark~3.6]{HP21_critical}.

For the offset law, choose on each welding ray a nonzero smooth nonnegative
radial cutoff $\psi_i$ with compact support strictly between zero and the observed
arc. The finite-energy function $e_i(r,\theta)=\pm\psi_i(r)\cos\theta$
has zero semicircle averages; choose the sign so that its trace on the
welding ray is $\psi_i$. Rescale $\psi_i$ so that $e_i$ has Dirichlet norm one, and write
\[
 h_i=\widetilde h_i+a_ie_i,\qquad
 a_1,a_2\sim N(0,1)
\]
independently of the residual fields and radial processes.
The fields $Y_i$ depend only on the residual fields. The lengths from
zero to the initial observed points have the form
\[
 \ell_i(a_i)=b_i+\frac1{2\pi}
           \int_{\{\psi_i>0\}}e^{a_i\psi_i}\,\nu_{\widetilde h_i}.
\]
Here $b_i$ is the normalized residual length from zero to the initial
observed point on the part where $\psi_i=0$.
These lengths are finite: near zero the length is a constant times
$\int_0^\infty e^{-\mathcal R_{2u}}\mathcal N(\dd u)$, with stationary
lateral critical measure $\mathcal N$ independent of $\mathcal R$.
Its finiteness follows from the block estimate in
Lemma~\ref{lem:crit-random-start-bessel}. Thus
\[
 \ell_i'(a_i)>0,\qquad\ell_i(\mathbb R)=(b_i,\infty).
\]
Conditional on the residual fields, each $\ell_i$ has positive density
$p_i$ on $(b_i,\infty)$. Since $s=\ell_1-\ell_2$, its density is
\begin{equation}\label{eq:crit-offset-density}
 p(s)=\int_{\max(b_2,b_1-s)}^\infty p_1(s+t)p_2(t)\,\dd t>0
                                      \quad\text{for a.e. }s\in\mathbb R.
\end{equation}
Condition on the residual fields and radial processes, then average
over these variables given $Y_1,Y_2$. The conditional law therefore has
the form
\[
 \operatorname{Law}(s\in\dd t\mid Y_1,Y_2)
       =p_{Y_1,Y_2}(t)\,\dd t,\qquad
 0<p_{Y_1,Y_2}(t)<\infty
\]
for $\operatorname{Law}(Y_1,Y_2)\otimes\dd t$-almost every pair.
Positivity follows from \eqref{eq:crit-offset-density}; finiteness follows
because each conditional density integrates to one. Combining this
kernel with the first assertion gives the stated joint measure class.
\end{proof}

\begin{proof}[Proof of Proposition~\ref{prop:crit-geometric-branch-proved}]
The GMC comparison is Lemma~\ref{lem:crit-common-model}. For the SLE
comparison, work under $\mathbf Q^*$ and set
\[
 C_B=M^{-1}(g_\beta(B)),\qquad E_B=\{K\cap C_B=\varnothing\},
 \qquad\mathcal M=\sigma(\eta,p,q,\alpha,\beta).
\]
The marks are independent of $(z,n)$ given $\eta$ by
\eqref{eq:crit-joint-law}. Lemma~\ref{lem:crit-root-avoid} gives
$\mathbf Q^*(E_B\mid\mathcal M)>0$ almost surely. Since $h$ is
$\mathcal M$-measurable, Lemma~\ref{lem:crit-marked-loop-class} implies
that, for every Borel $F\subset\mathcal X_I$,
\[
 P_I^{\mathrm{weld}}(F\cap O_B)=0
 \quad\Longleftrightarrow\quad
 \mathbf Q^*(h|_I\in F,\ O_B\cap E_B)=0.
\]

On $O_B\cap E_B$, both observed arcs lie in the welding rays. Set
$C_1(x)=\mc M_{Y_1}([a_0,x])$ and
$C_2(y)=\mc M_{Y_2}([b_0,y])$. Equality of interface lengths and
Lemma~\ref{lem:crit-wedge-local-fields} give
\[
 \begin{aligned}
 L_1(A_1(x))-L_1(A_1(a_0))&=C_1(x),\\
 L_2(A_2(y))-L_2(A_2(b_0))&=C_2(y),\\
 y=h(x)&\quad\Longleftrightarrow\quad C_2(y)=C_1(x)+s.
 \end{aligned}
\]
The exterior coordinate runs outward along the negative welding ray.
Consequently, the event above is exactly
\[
 (A_1,A_2)\in\mathcal A_{I,B},\qquad
 (\mc M_{Y_1},\mc M_{Y_2},s)\in\mathcal D,\qquad
 \mathcal W(\mc M_{Y_1},\mc M_{Y_2},s)\in F.
\]
By Lemmas~\ref{lem:crit-joint-charts} and
\ref{lem:crit-wedge-local-fields}, the law of the local data, restricted
to $\mathcal A_{I,B}$, is mutually absolutely continuous with respect to
\[
 (\Lambda_1\otimes\Lambda_2)|_{\mathcal A_{I,B}}
       \otimes\lambda_I\otimes\lambda_B\otimes\dd s.
\]
The open set $\mathcal A_{I,B}$ has positive
$\Lambda_1\otimes\Lambda_2$ measure, since the endpoints of the welding
rays can be placed outside the observed arcs. Fubini's theorem and
$\gamma_0(\dd s)\sim\dd s$ identify the resulting null sets with
$\Sigma_{I,B}(F)=0$. All pushforwards are measurable by the construction
of the normalized Riemann maps, the local GMCs and the cumulative
inverses. This proves the SLE comparison.
\end{proof}

\begin{proof}[Proof of Theorem~\ref{thm:crit-local-bridge}]
Choose countably many proper arcs $B$ with rational endpoints. Since
every proper image arc $h(I)$ lies in the interior of one of them,
$\mathcal X_I=\bigcup_BO_B$. Proposition~\ref{prop:crit-geometric-branch-proved}
gives $P_I^{\mathrm{weld}}|_{O_B}\sim P_I|_{O_B}$ for each $B$; taking
the countable union proves $P_I^{\mathrm{weld}}\sim P_I$.
\end{proof}

\subsection{Local shape statistics and smooth deformations}\label{subsec:crit-shape}

We use bounded statistics of normalized restrictions to detect
distortion. We first bound their change under composition. We then
prove density estimates for smooth source deformations; these estimates
will control the means and covariances of the statistics.

\begin{definition*}[Normalized restrictions]
For $h\in\Homeo_+(\m S^1)$ and $I=[a,a+r]$, $0<r<1$, define
\[
 h_I(u)=\frac{h(a+ru)-h(a)}{h(a+r)-h(a)},\qquad 0\le u\le1.
\]
Differences are taken in degree-one increasing lifts. The resulting
homeomorphism of $[0,1]$ is independent of the lift; for the random
welding $H_4$ we write $H_{4,I}=(H_4)_I$.

Fix $0<c<1/10$ and a real nonzero function
$\chi\in C_c^\infty((0.4,0.6))$. Set $\vartheta=\chi''$ and
$\Delta_c f(u)=f(u+c)-2f(u)+f(u-c)$. For $0<s<1/2$, define
\[
 \xi_{x,s}(h)=\frac{h(x+s)-2h(x)+h(x-s)}{h(x+s)-h(x-s)},\qquad
 F_I(h)=\int_{\mathbb R}\vartheta(u)\xi_{a+ru,cr}(h)\,\dd u.
\]
The variable $\xi_{x,s}$ compares the lengths of the two adjacent image
arcs. The signed average $F_I$ is bounded by $\|\vartheta\|_1$ and
uses only $h|_I$. We will use the cancellation identities
$\int\vartheta=\int u\vartheta(u)\,\dd u=0$.

For a fixed $\varphi\in\Homeo_+(\m S^1)$ with lift $\Phi$, put
\[
 \ell_I=\Phi(a+r)-\Phi(a),\qquad k_I=\varphi_I,\qquad
 \varepsilon_I=\|k_I-\operatorname{id}\|_{L^\infty([0,1])},
\]
\[
 w_I=\frac1c\int_{\mathbb R}\vartheta(u)\Delta_c k_I(u)\,\dd u.
\]
Thus $\ell_I$ is the image length, $\varepsilon_I$ measures the departure
of the normalized restriction from the identity, and $w_I$ is a signed
coefficient of this distortion. The support of $\vartheta$ ensures that
all arguments $u,u\pm c$ lie in $[0,1]$.
We integrate over arc parameters $(a,r)\in\m S^1\times(0,1)$ with
\[
 \dd\Lambda(I)=\dd a\,\frac{\dd r}{r^2},\qquad I=[a,a+r].
\]
In $F_{\varphi(I)}$, the image arc is parametrized affinely from
$\Phi(a)$ to $\Phi(a+r)$.
\end{definition*}

\begin{lem}[Stability under composition]\label{lem:crit-packet-bv}
For all $h,\varphi\in\Homeo_+(\m S^1)$ and every proper arc $I$,
\[
 |w_I|\le C\varepsilon_I,\qquad
 |F_I(h\circ\varphi)-F_{\varphi(I)}(h)|\le C\varepsilon_I,
\]
where $C$ depends only on $c$ and $\vartheta$.
\end{lem}
\begin{proof}
The first inequality follows from $\Delta_c k_I=\Delta_c(k_I-\mathrm{id})$.
For the second, put $M(t)=h(\Phi(a)+\ell_I t)$, $\nu=\dd M$ and
\[
 f(u)=\frac{M(u+c)-2M(u)+M(u-c)}{M(u+c)-M(u-c)}.
\]
All expressions use increasing lifts. We first bound the variation of
$f$, uniformly in the increasing function $M$. If $d\ge c/2$ and
$E\subset[0,1]$, partition $[0,1]$ into $N_c=\lceil2/c\rceil$ intervals $J$ of length at
most $c/2$. For $t\in J$, $[t-d,t+d]\supset J$, and hence
\begin{equation}\label{eq:crit-bv-count}
 \int_E\frac{\nu(\dd t)}{\nu([t-d,t+d])}\le N_c.
\end{equation}
First suppose $M$ is smooth with $M'>0$. With
$L=M(u)-M(u-c)$, $R=M(u+c)-M(u)$ and $D=L+R$,
\[
 f'(u)=\frac{2LM'(u+c)+2RM'(u-c)}{D^2}
 -\frac{2M'(u)}D.
\]
The first term is nonnegative. Since $|f|\le1$, \eqref{eq:crit-bv-count} gives
\[
 \operatorname{Var}_{[1/4,3/4]}f
 \le2+4\int_{1/4}^{3/4}\frac{M'(u)}{M(u+c)-M(u-c)}\,\dd u
 \le2+4N_c.
\]
For general $M$, convolve with positive smooth kernels and add a positive
linear term tending to zero. The approximations converge uniformly on
compact sets. On the interval considered the denominators have a positive
minimum, so the same bound follows by lower semicontinuity of variation.

We now compare $f$ with the quotient after composition. Write $k=k_I$
and $\epsilon=\varepsilon_I$, and let $f_k$ be the same quotient
as $f$ with $M(u-c),M(u),M(u+c)$ replaced by
$M(k(u-c)),M(k(u)),M(k(u+c))$. If $\epsilon<c/8$, monotonicity of
$(C-2B+A)/(C-A)$ in the three arguments gives
\[
 f(u-\epsilon)-A_-(u)\le f_k(u)\le f(u+\epsilon)+A_+(u),
\quad
 A_\pm(u)=\frac{2\nu([u-\epsilon,u+\epsilon])}
 {\nu([u-c\pm\epsilon,u+c\pm\epsilon])}.
\]
For $u\in[0.4,0.6]$ the BV translation estimate bounds the integrals of
$|f(u\pm\epsilon)-f(u)|$ by $C\epsilon$. If $|t-u|\le\epsilon$, each
denominator contains $[t-(c-2\epsilon),t+(c-2\epsilon)]$. Fubini and
\eqref{eq:crit-bv-count} therefore give
\[
 \int_{0.4}^{0.6}A_\pm(u)\,\dd u
 \le4\epsilon\int_{0.4-\epsilon}^{0.6+\epsilon}
 \frac{\nu(\dd t)}{\nu([t-(c-2\epsilon),t+(c-2\epsilon)])}
 \le4N_c\epsilon.
\]
Thus $\int_{0.4}^{0.6}|f_k-f|\,\dd u\le C\epsilon$. Multiplying by
$\vartheta$ gives the second inequality. If $\epsilon\ge c/8$, it follows
instead from $|f|,|f_k|\le1$. No derivative of $h$ or $\varphi$ is used.
\end{proof}

The next lemma gives uniform bounds for changes of the Gaussian law
under smooth deformations supported at an arbitrary small scale.

\begin{definition*}[Smooth source deformations]
In the coordinate $x\in\m R/\m Z$, fix a real
$\zeta\in C_c^\infty(\mathbb R)$ and choose $d_0>0$ so that
$d\operatorname{supp}\zeta$ fits inside one period for $0<d\le d_0$.
For a center $x_0$, let $\phi_t$ be the flow of the periodically
extended vector field $v_d(x)=d\zeta((x-x_0)/d)$, with
$\phi_0=\operatorname{id}$. Let $\pi_0$ subtract the circle mean, and
write $\mathsf G$ for the law of $X$. Define
\[
 \mathsf G_t=\operatorname{Law}\bigl(\pi_0(X\circ\phi_t+2\log\phi_t')\bigr).
\]
\end{definition*}

\begin{lem}[Smooth deformations]\label{lem:crit-gaussian-flow}\label{lem:crit-smooth-channel}
The law $\mathsf G_t$ is mutually absolutely continuous with respect to
$\mathsf G$. Write $D_t=\dd\mathsf G_t/\dd\mathsf G$.
For every fixed $T<\infty$,
$D_t$ has a jointly measurable version, is $C^2$ as an $L^1$-valued
function of $t$, and
\begin{equation}\label{eq:crit-density-derivatives}
 \sup_{d,\,|t|\le T}\sum_{j=0}^2
 \|\partial_t^jD_t\|_{L^1}<\infty,
 \qquad
 \sup_d E\!\left[\sup_{|t|\le T}|\log D_t|^2\right]<\infty.
\end{equation}
All suprema in $d$ range over $(0,d_0]$, and all $L^1$ norms are taken
with respect to $\mathsf G$. The constants depend on $T$ and $\zeta$,
but not on $d$ or $x_0$.
For any bounded measurable $F$ on $\Homeo_+(\m S^1)$,
\[
 \m E_4[F(H_4\circ\phi_t)]=E[D_tF(H_4)].
\]
On the right, $D_t$ is evaluated at the source field $X_1$, and the
expectation also integrates over $X_2$ and the independent phase $U$.
\end{lem}
Let $\mathcal H$ be the real mean-zero Cameron--Martin space, with norm
$\|f\|_{\mathcal H}^2=\sum_{n\ne0}|n||\widehat f(n)|^2$, and set
\[
 C_tf=\pi_0(f\circ\phi_t),\qquad L_vf=\pi_0(v_df'),\qquad
 S_t=C_tC_t^*,\qquad b_t=\pi_0(2\log\phi_t').
\]
The estimates behind the lemma are
\[
 \begin{aligned}
 \|L_v+L_v^*\|_{\mathrm{HS}}+\|L_v(L_v+L_v^*)\|_{\mathrm{HS}}&\le C,\\
 \sup_{d,\,|t|\le T}\sum_{j=0}^2
 \bigl(\|\partial_t^j(S_t-I)\|_{\mathrm{HS}}
       +\|\partial_t^jb_t\|_{\mathcal H}\bigr)&\le C_T.
 \end{aligned}
\]
They give uniform $L^p$ bounds for the Gaussian densities for some $p>1$;
differentiating their log densities then gives
\eqref{eq:crit-density-derivatives}. The Fourier calculation and the
passage from finite-dimensional densities to the field law are proved
in Appendix~\ref{subsec:crit-gaussian-proof}.

\subsection{Mass estimates and weighted means}\label{subsec:crit-mass}

The source and target measures require different estimates. For the
source, the quarter masses divided by the total mass stay away from
zero in probability, uniformly over small arcs. Writing $I_j$,
$1\le j\le4$, for the four equal subarcs of $I$, we prove
\begin{equation}\label{eq:crit-quarter-input}
 \lim_{b\downarrow0}\sup_{|I|\le r_0}
 P\!\left(\min_{1\le j\le4}\frac{\mu_1(I_j)}{\mu_1(I)}<b\right)=0
\end{equation}
for some $r_0>0$. For the target, sample $Y$ conditionally with law
$\mu_2$, and let $B_\pm(q)$ be the lengths of the arcs from $Y$ in
the two directions with $\mu_2$-mass $q$. We show that their length
ratios remain tight when the masses differ by bounded factors: for
each fixed $M<\infty$,
\begin{equation}\label{eq:crit-root-input}
 \lim_{K\to\infty}\limsup_{q\downarrow0}
 P\!\left(\sup_{\substack{\sigma,\tau\in\{-,+\}\\|u|,|v|\le M}}
 \left|\log\frac{B_\sigma(qe^u)}{B_\tau(qe^v)}\right|>K\right)=0.
\end{equation}
The ratios are defined for $q<e^{-M}$. The signs $+$ and $-$ denote
the positive and negative circle directions, respectively.
Lemma~\ref{lem:crit-quarter} proves the source estimate. The target
estimate follows from a Bessel description of critical boundary length
at a sampled point. We then combine the two estimates to control
normalized image lengths and prove that their weighted means have a
uniformly positive derivative near the midpoint.

\begin{lem}[Masses of equal subarcs]\label{lem:crit-quarter}
For the critical circle GMC, \eqref{eq:crit-quarter-input} holds.
\end{lem}
\begin{proof}
By rotation invariance take $I=[-r/2,r/2]$ in $\mathbb R/\mathbb Z$.
Put $\tau_r=\tan(\pi r/2)$ and
\[
 x_r(u)=\pi^{-1}\arctan(\tau_ru).
\]
For $r$ sufficiently small this is a chart on a fixed interval
$(-L,L)$ with $L>2$, and $I$ is the image of $[-1,1]$. Set
$Y_r=X_1\circ x_r$. Its covariance is
\[
 -2\log|u-v|-2\log(2\tau_r)
 +\log(1+\tau_r^2u^2)+\log(1+\tau_r^2v^2).
\]
Choose $\rho\in C_c^\infty((-L,L))$ with $\int\rho=1$, and let
$A_r=\langle Y_r,\rho\rangle$, $Z_r=Y_r-A_r$. For
$k(u,v)=-2\log|u-v|$, put
\[
 p_\rho(u)=\int k(u,t)\rho(t)\,\dd t,\qquad
 b_\rho=\iint k(s,t)\rho(s)\rho(t)\,\dd s\,\dd t.
\]
The covariance of $Z_r$ is
$k(u,v)-p_\rho(u)-p_\rho(v)+b_\rho$, independent of $r$. No independence
of $A_r$ and $Z_r$ is needed.

Use a common regularization and the same critical normalization. The
coordinate Jacobian and the deterministic variance difference give
\begin{equation}\label{eq:crit-quarter-pullback}
 x_r^*\mc M_{X_1}(\dd u)
 =\frac{2\tau_r^2}{\pi} e^{A_r}
 (1+\tau_r^2u^2)^{-2}\nu_{Z_r}(\dd u),
\end{equation}
where $\nu_Z=e^{-p_\rho+b_\rho/2}\mc M_Z^{\mathrm{sw}}$ and
$\mc M_Z^{\mathrm{sw}}$ is the self-Wick critical GMC with reference measure
$\dd u$. To check the constant and the weight, the subcritical identity is
\begin{align*}
 x_r^*\mc M_{X_1}^\gamma(\dd u)
 &=\frac{\tau_r(2\tau_r)^{\gamma^2/4}}\pi
 e^{\gamma A_r/2}(1+\tau_r^2u^2)^{-1-\gamma^2/4}\\
 &\quad\times e^{-\gamma^2p_\rho(u)/4+\gamma^2b_\rho/8}
 \mc M_{Z_r}^{\gamma,\mathrm{sw}}(\dd u).
\end{align*}
Indeed, $x_r'(u)=\tau_r/[\pi(1+\tau_r^2u^2)]$ and the finite difference
of the regularized variances of $Y_r$ and $Z_r$ is
\[
 -2\log(2\tau_r)+2\log(1+\tau_r^2u^2)+2p_\rho(u)-b_\rho.
\]
Subtracting this difference in the Wick exponent gives the identity.
Divide by $2(2-\gamma)$ and let $\gamma\uparrow2$ at fixed $r$.
The weights converge uniformly on $[-1,1]$, which proves
\eqref{eq:crit-quarter-pullback}.

The quarter points have preimages $-b_r,0,b_r$, where
\[
 b_r=\frac1{\sqrt{1+\tau_r^2}+1}.
\]
For all small $r$, each of the four preimage subarcs contains a fixed
nondegenerate closed interval $B_j$, $1\le j\le4$. On $[-1,1]$, the
weight $(1+\tau_r^2u^2)^{-2}$ has upper and lower bounds independent of
$r$. Since the laws of $Z_r$ are identical,
\[
 \min_j\frac{\mu_1(I_j)}{\mu_1(I)}
 \quad\text{stochastically dominates}\quad
 C^{-1}\frac{\min_j\nu_Z(B_j)}{\nu_Z([-1,1])}
\]
for a single reference field $Z$. The random variable on the right is
strictly positive almost surely, since the critical GMC is finite and has
full support. Its lower-tail probability tends to zero as $b\downarrow0$,
uniformly in $r$. This proves \eqref{eq:crit-quarter-input}.
\end{proof}

\begin{lem}[Bessel tail integrals]\label{lem:crit-random-start-bessel}
Let $R$ be BES$(3)$ with an almost surely finite random initial value
$R_0\ge0$. Let $(\mathcal N_-,\mathcal N_+)$ be jointly stationary,
non-atomic random measures on $\m R$, independent of $R$, with full support.
Assume that the mass of a unit interval has a finite $p$th moment for every
$0<p<1$ and a finite negative moment of every order. Define
\[
 J_\sigma(s)=\int_s^\infty e^{-R_{2u}}\mathcal N_\sigma(\dd u),\qquad
 S_\sigma(h)=\inf\{s\ge0:J_\sigma(s)\le e^{-h}\},
 \quad \sigma\in\{-,+\}.
\]
The tails $J_\sigma$ are almost surely finite, continuous, and strictly
decreasing to zero. For every $M<\infty$,
\[
 \lim_{T\to\infty}\limsup_{h\to\infty}
 P\!\left(\max_\sigma S_\sigma(h+M)
 -\min_\sigma S_\sigma(h-M)>T\right)=0.
\]
\end{lem}
The proof compares the inverse tails with last-exit times
$L_a=\sup\{s:R_{2s}=a\}$. On $\{R_0<a\}$, the process after $L_a$,
shifted down by $a$, is BES$(3)$ from zero. Stationarity of the lateral
measures makes $(e^aJ_\sigma(L_a))_{\sigma=\pm}$ have a law independent
of $a$, with finite moments of order $p\in(0,1)$ and of every negative
order. These bounds place all the
inverse tails in the window between $L_{h-M-K}$ and $L_{h+M+K}$,
apart from an error tending to zero as $K\to\infty$. The difference
of these last-exit times has a distribution independent of $h$.
The details, including the last-exit argument, are in
Appendix~\ref{subsec:crit-bessel-proof}.

We now turn to the point sampled from the target measure. Let $h$ be
the Neumann GFF on $\mathbb H$, normalized to have zero average on the
unit semicircle, and write $\nu_h$ for its critical boundary length. Fix a compact interval
$J\Subset(-1/4,1/4)$, put $r_c=1/4$ and $L=\log(1/r_c)$. For $y\in J$
use
\[
 \psi_y(s+\mathrm{i}\theta)
 =y-r_ce^{-(s+\mathrm{i}\theta)},
 \qquad s\ge0,\quad0<\theta<\pi.
\]
The Neumann covariance in this gauge is
\[
 G_{\mathbb H}(z,w)=-\log|z-w|-\log|z-\overline w|
 +2\log|z|_++2\log|w|_+.
\]
The maps $\psi_y$ parametrize half-disks of radius $r_c$ centered at $y$,
all contained in the unit disk. Expanding the two
logarithms gives
\begin{align}\label{eq:crit-polar-covariance}
 &G_{\mathbb H}(\psi_y(s+\mathrm{i}\theta),
                \psi_y(t+\mathrm{i}\vartheta))\notag\\
 &\quad=2L+2(s\wedge t)+2\sum_{n\ge1}
 \frac{e^{-n|s-t|}}n\cos(n\theta)\cos(n\vartheta).
\end{align}
Thus, for fixed $y$, the field in these coordinates is $A_s+N(s,\theta)$,
where
\[
 A_s=A_0+B_{2s},\qquad A_0\sim N(0,2L),
\]
and $A_0$, $B$ and the stationary lateral field $N$ are independent.

Let $h_\epsilon(y)$ denote the semicircle average and define
\[
 C_b=\{h_{r_ce^{-u}}(y)<b+2u\text{ for every }y\in J,\ u\ge0\}.
\]
\begin{lem}[A truncated rooted law]\label{lem:crit-barrier-root}
The events $C_b$ increase to a full-probability event. For each sufficiently
large $b$, there is a finite random measure $\overline d_h^b$ on $J$, with
$Z_b=E[\overline d_h^b(J)]\in(0,\infty)$, such that
\[
 \1_{C_b}\overline d_h^b=\1_{C_b}\nu_h|_J.
\]
Under the probability
\[
 \Theta_b(\dd h,\dd y)=Z_b^{-1}P(\dd h)\overline d_h^b(\dd y),
\]
$y$ is uniform on $J$ and the conditional radial process is
\[
 A_s=b+2s-R_{2s},\qquad R_0=b-A_0>0,
\]
where, conditional on $R_0$, $R$ is BES$(3)$. Its initial distribution is
\[
 \operatorname{Law}(A_0\in\dd a)
 =\frac{(b-a)_+e^a}{m_b}N(0,2L)(\dd a),\qquad
 m_b=E[(b-A_0)_+e^{A_0}].
\]
The lateral field $N$ is independent of $R$ under this law. Moreover, the
field and its critical boundary measure on every compact root annulus are
constructed together by the same local GMC functional.
\end{lem}
For $\epsilon=r_ce^{-s}$, the truncated derivative measure is
\[
 d_\epsilon^b(\dd y)=
 \frac{b+2s-h_\epsilon(y)}2\,\epsilon e^{h_\epsilon(y)}
 \1_{\{h_{r_ce^{-u}}(y)<b+2u,\ 0\le u\le s\}}\1_J(y)\,\dd y.
\]
Brownian exponential tilting followed by the positive-distance Doob
transform gives the displayed Bessel law and
$E[d_\epsilon^b(J)]=|J|r_cm_b/2$. Uniform integrability permits a joint
limit of the field, its boundary measure and the sampled point.
On $C_b$ this limit agrees with $\nu_h|_J$. Appendix~\ref{subsec:crit-root-proof}
proves these assertions and checks that the limiting field and boundary
measure are still related by the local GMC construction.

Write $\psi_y^-(s)=y-r_ce^{-s}$ and $\psi_y^+(s)=y+r_ce^{-s}$
for the two boundary parametrizations. Let $\mathcal N_-$ and
$\mathcal N_+$ be the critical measures of $N$ on the corresponding
strip boundaries $\theta=0$ and $\theta=\pi$. With $t=|s-u|$, the boundary covariances are
\[
 K_{++}(t)=K_{--}(t)=-2\log(1-e^{-t}),\qquad
 K_{+-}(t)=-2\log(1+e^{-t}).
\]
The two measures are jointly stationary and independent of $R$; the two
sides need not be independent. Their one-block positive
moments of order $p<1$ and negative moments follow from the critical GMC
moment theorem \cite[Section~2.6, Theorem~2.11]{PowellReview}. The boundary coordinate rule adds the term
$2\log|\psi_y'|=2\log r_c-2s$ to the field.
Consequently, under $\Theta_b$ and with compatible units,
\begin{equation}\label{eq:crit-root-length}
 (\psi_y^\pm)^*\nu_h(\dd s)
 =r_c^2e^b e^{-R_{2s}}\mathcal N_\pm(\dd s),\qquad s>0.
\end{equation}

\begin{lem}[Rooted mass ratios]\label{lem:crit-root-input-proved}\label{lem:crit-root-transfer}
The marked law obtained by sampling $Y$ from normalized critical circle GMC
satisfies \eqref{eq:crit-root-input}.
\end{lem}
\begin{proof}
We first identify the circle root in a Cayley chart. Let
$A:\mathbb D\to\mathbb H$ be a fixed Cayley map, $X$ a centered circle
log field, and $\mu_X=\mc M_X/\mc M_X(\m S^1)$ its normalized critical GMC.
Write $\Gamma=\Gamma_0+PX$, where $PX$ is the harmonic extension of
$X$ and $\Gamma_0$ is an independent zero-boundary GFF. Subtract the
unit-semicircle mean
$c_A$ from $\Gamma\circ A^{-1}$ to obtain a semicircle-gauge Neumann GFF
$h$. Conformal invariance modulo constants and the critical boundary-length
coordinate rule give
\[
 A^*\nu_h=2\pi e^{-c_A}|A'|^2\mc M_X.
\]
This uses a common regularization and does not require $c_A$ to be
independent of $X$. Integrating away from the pole and using non-atomicity,
then normalizing, yields
\begin{equation}\label{eq:crit-cayley-root}
 A_*\mu_X(\dd y)=\frac{a_A(y)\nu_h(\dd y)}{\int a_A\,\dd\nu_h},
 \qquad a_A(y)=|A'(A^{-1}y)|^{-2}.
\end{equation}
The boundary trace and coordinate formulas are those of
\cite[Remarks~3.5--3.6 and Lemma~3.8]{F1} and
\cite[Section~2.1]{PowellSepulveda}.

Choose finitely many Cayley charts $A_j$ and a nonnegative partition of
unity $\chi_j$ subordinate to compact arcs away from their poles. Rescale
each chart by a fixed half-plane automorphism so that its observed real
interval is contained in $(-1/4,1/4)$. Work on the disjoint union of
these chart spaces. On the $j$th component, let $P_j$ be the law of its
semicircle-gauge field and define
\[
 Q_j(\dd h,\dd y)=P_j(\dd h)
 \frac{\chi_j(A_j^{-1}y)a_{A_j}(y)\nu_h(\dd y)}{T_j(h)},
 \qquad T_j(h)=\int a_{A_j}\,\dd\nu_h.
\]
The $j$th component carries the circle measure
$(A_j^{-1})_*(a_{A_j}\nu_h/T_j)$ and marked point $A_j^{-1}y$.
Equation~\eqref{eq:crit-cayley-root} and $\sum_j\chi_j=1$ show that
$Q=\sum_jQ_j$ is a probability whose image is the required marked circle
law. Write $C_{j,b}$ and $\Theta_{j,b}$ for the barrier event and root law
in this chart. The restriction $Q_{j,b}=\1_{C_{j,b}}Q_j$ has density
\begin{equation}\label{eq:crit-root-density}
\begin{aligned}
 D_{j,b}(h,y)&=Z_{j,b}\chi_j(A_j^{-1}y)\1_{C_{j,b}}
 \frac{a_{A_j}(y)}{T_j(h)},\\
 E_{\Theta_{j,b}}[D_{j,b}]&\le1.
\end{aligned}
\end{equation}
This follows by testing the two measures against a nonnegative measurable
function. The density is fixed as the mass level tends to zero.

On each compact chart, $a_{A_j}$ has positive upper and lower bounds
and $T_j(h)\in(0,\infty)$. Thus \eqref{eq:crit-root-length} gives
\[
 C^{-1}J_\pm(s)\le A_\pm(s)\le CJ_\pm(s),\qquad s\ge s_0,
\]
for finite random $C\ge1,s_0$, where $A_\pm(s)$ are the normalized
masses of $[y-r_ce^{-s},y]$ and $[y,y+r_ce^{-s}]$. Let
$S_\pm^A(\lambda)=\inf\{s\ge0:A_\pm(s)\le e^{-\lambda}\}$ and let $S_\pm$ invert
$J_\pm$ as in Lemma~\ref{lem:crit-random-start-bessel}. Comparison of
chart and circle distances gives
\[
 |S_\pm^A(\lambda)+\log B_\pm(e^{-\lambda})|\le C_0
\]
for all sufficiently large $\lambda$, with $C_0<\infty$ almost surely.
Localize $C,s_0,C_0$, the level from which the preceding comparison
holds, and $\min_\pm J_\pm(s_0)$ to deterministic bounds.
Outside an event of arbitrarily small probability, for some fixed $c$,
\[
 S_\pm(\lambda-c)\le S_\pm^A(\lambda)\le S_\pm(\lambda+c)
\]
at all sufficiently large levels. Lemma~\ref{lem:crit-random-start-bessel}
therefore proves \eqref{eq:crit-root-input} under each $\Theta_{j,b}$.

Let $F_{q,K}$ denote the event inside the probability in
\eqref{eq:crit-root-input}, and put
$\delta_b=\|Q-\sum_jQ_{j,b}\|_{\rm TV}$. For every $L_0>0$,
\[
 Q(F_{q,K})\le\delta_b+\sum_j\left(
 L_0\Theta_{j,b}(F_{q,K})
 +E_{\Theta_{j,b}}[D_{j,b}\1_{\{D_{j,b}>L_0\}}]\right).
\]
Moreover,
\[
 \delta_b=\sum_jQ_j(C_{j,b}^{\,c})
             \le\sum_jP_j(C_{j,b}^{\,c})\longrightarrow0.
\]
For fixed $b$, the finitely many densities $D_{j,b}$ are integrable.
Choose first $b$, then $L_0$, and finally $K$; take the upper limit as
$q\downarrow0$. The tightness under $\Theta_{j,b}$ makes the right-hand
side arbitrarily small. Projection to the marked circle law proves
\eqref{eq:crit-root-input}.
\end{proof}

\begin{lem}[Normalized image lengths]\label{lem:crit-balanced-tight}\label{lem:crit-small-end-arc}
For every $\eta>0$ there are $b_\eta\in(0,1/2)$ and $r_\eta>0$ such
that, for every $I=[a,a+r]$ with $r\le r_\eta$,
\[
 \m P_4\bigl(H_{4,I}(1/2)\notin[b_\eta,1-b_\eta]\bigr)\le\eta.
\]
Moreover, for every $\eta>0$ there are $d_\eta,r_\eta'>0$ such that
\[
 \m P_4\bigl(H_{4,I}(3/4)-H_{4,I}(1/4)<d_\eta\bigr)\le\eta
 \qquad (|I|\le r_\eta').
\]

For every $\epsilon>0$,
\begin{align*}
 &\lim_{u\downarrow0}\limsup_{r\downarrow0}\sup_a
 \m P_4\bigl(H_{4,[a,a+r]}(u)>\epsilon\bigr)=0,\\
 &\lim_{u\downarrow0}\limsup_{r\downarrow0}\sup_a
 \m P_4\bigl(1-H_{4,[a,a+r]}(1-u)>\epsilon\bigr)=0.
\end{align*}
Consequently, either probability tends to zero along any deterministic
sequences $r_n\to0$ and $u_n\to0$.
\end{lem}
\begin{proof}
It suffices to take $0<\eta<1$. Let $x=a+r/2$, $q=\mu_1(I)$ and
\[
 t_-=\mu_1([a,x]),\qquad t_+=\mu_1([x,a+r]).
\]
By \eqref{eq:crit-quarter-input}, choose $b>0$ and $r_1>0$ such that
\[
 P(bq\le t_-\le q,\ bq\le t_+\le q)\ge1-\eta/3
 \qquad(r\le r_1).
\]
Conditionally on the source field, $H_4(x)$ is a $\mu_2$-sample. Its
joint law with $\mu_2$ does not depend on the source field. By
\eqref{eq:crit-root-input}, choose $K<\infty$ and $q_0>0$ such that
\[
 P\!\left(
 \frac{\max_{\sigma=\pm,\,v\in\{b,1\}}B_\sigma(vq)}
 {\min_{\sigma=\pm,\,v\in\{b,1\}}B_\sigma(vq)}>e^K
 \right)\le\eta/3
 \qquad(0<q\le q_0).
\]
The same bound holds conditionally on any source field with $q\le q_0$.
On the source event above, if $q\le q_0$ and the target ratio is at
most $e^K$, monotonicity gives
\[
 H_{4,I}(1/2)=\frac{B_-(t_-)}{B_-(t_-)+B_+(t_+)}
 \in\left[\frac1{1+e^K},\frac{e^K}{1+e^K}\right].
\]
By rotation invariance, $E[q]=r$. Thus $P(q>q_0)\le r/q_0$.
Take $r_\eta\le\min(r_1,\eta q_0/3)$ and
$b_\eta=(1+e^K)^{-1}$. The union bound proves the midpoint bound.

For the central-quarter bound, apply the midpoint bound with error $\eta/2$ to $I$
and its left half $I^-=[a,a+r/2]$. On an event of probability at least
$1-\eta$, both normalized midpoint values lie in $[b,1-b]$. Thus
\begin{align*}
 H_{4,I}(3/4)-H_{4,I}(1/4)
 &\ge H_{4,I}(1/2)-H_{4,I}(1/4)\\
 &=H_{4,I}(1/2)\bigl(1-H_{4,I^-}(1/2)\bigr)\ge b^2.
\end{align*}
Take $d_\eta=b^2$. No independence of the two midpoint values is needed.

For the endpoint estimates, fix $\eta>0$ and choose $b,r_\eta$ from the midpoint bound. For $r\le r_\eta$, put
$I_j=[a,a+2^{-j}r]$ and
\[
 V_j=\frac{H_4(a+2^{-j-1}r)-H_4(a)}
 {H_4(a+2^{-j}r)-H_4(a)}.
\]
For every $j$, $\m P_4(V_j>1-b)\le\eta$. Therefore, for
$N_n=\sum_{j=0}^{n-1}\1_{\{V_j>1-b\}}$,
\[
 \m E_4[N_n]\le n\eta,\qquad
 \m P_4(N_n>n/2)\le2\eta.
\]
On $\{N_n\le n/2\}$, at least $n/2$ factors in the following product
are at most $1-b$, and all factors are at most one:
\[
 H_{4,[a,a+r]}(2^{-n})=\prod_{j=0}^{n-1}V_j
 \le(1-b)^{n/2}.
\]
Choose $n$ with $(1-b)^{n/2}\le\epsilon$. Monotonicity gives
\[
 \sup_{a,\,0<r\le r_\eta,\,0<u\le2^{-n}}
 \m P_4\bigl(H_{4,[a,a+r]}(u)>\epsilon\bigr)\le2\eta.
\]
Since $\eta$ is arbitrary, this proves the left-endpoint estimate and its
sequential form. Nested right halves give the right-endpoint estimate by the
same argument. The factors $V_j$ need not be independent.
\end{proof}

We test normalized image lengths against bounded functions of finitely
many welding values. Flows avoiding the observation points fix the weight.

\begin{definition*}[Cylinder functions]
Let $x_1,\ldots,x_m$ be fixed circle points and let
\[
 W(h)=f(h(x_1),\ldots,h(x_m)),\qquad 0\le f\le1,
 \qquad w=\m E_4[W],
\]
where $f:(\m S^1)^m\to[0,1]$ is Borel measurable. A function of this form is a bounded cylinder function. For an arc $I$ whose interior contains
none of the $x_i$, define
\[
 M_{W,I}(u)=\m E_4[W H_{4,I}(u)].
\]
\end{definition*}

\begin{prop}[Weighted midpoint means]\label{prop:crit-weighted-central}
Let $W$ be a bounded cylinder function with $w>0$, and let $I$
contain none of its observation points in its interior.
There are $r_w,c_w,\eta_0>0$ such that
\[
 M_{W,I}'(u)\ge c_w
 \qquad(|I|\le r_w,\ |u-1/2|\le\eta_0).
\]
These functions have uniformly bounded $C^{1,1}$ norms near $1/2$.
The constants $r_w,c_w$ depend on $w$, but not on $m$, $f$, or the
locations of the observation points. The interval for $u$ is fixed.
\end{prop}
\begin{proof}
Put $C_I=H_{4,I}(3/4)-H_{4,I}(1/4)$. Lemma~\ref{lem:crit-balanced-tight}
gives $d,r_w>0$, depending only on $w$, such that
$\m P_4(C_I<d)\le w/2$ for $|I|\le r_w$. Hence
\begin{equation}\label{eq:crit-weighted-central-mass}
 \m E_4[W C_I]\ge d\bigl(w-\m P_4(C_I<d)\bigr)
 \ge dw/2=:a_w>0.
\end{equation}
Every smooth flow supported in $I$ fixes the observations and therefore
fixes $W$. Choose such a flow which, in normalized coordinates,
translates a neighborhood of $u$ by $t$. Since it fixes the endpoints,
Lemma~\ref{lem:crit-smooth-channel} gives
\[
 M_{W,I}(u+t)=E[D_t W H_{4,I}(u)].
\]
The density derivative bounds give a uniform $C^{1,1}$ bound near $1/2$.

Let $\Theta$ be independent and uniform on $[0,1]$, and set
$Y_I=H_{4,I}^{-1}(\Theta)$. Expectations involving $Y_I$ also integrate
over $\Theta$. Conditional on $H_4$, the event $Y_I\le u$ is
$\Theta\le H_{4,I}(u)$, so
\[
 E[W\1_{\{Y_I\le u\}}]=M_{W,I}(u),\qquad
 E[W\1_{\{Y_I\in[1/4,3/4]\}}]\ge a_w.
\]
Choose a smooth vector field $v$, positive on $(0,1)$, flat at the
endpoints, and equal to one on $[1/4,3/4]$. Let $\kappa_t$ be its flow
and $\widehat\kappa_t$ its rescaling to $I$, extended by the identity.
Define
\[
 \zeta'(u)=1/v(u),\qquad \zeta(1/2)=0,\qquad Z=\zeta(Y_I).
\]
Since $(h\circ\widehat\kappa_t)_I=h_I\circ\kappa_t$, composition sends
$Y_I$ to $\kappa_{-t}(Y_I)$ and $Z$ to $Z-t$. Write $D_t$ for the
source-field density associated with $\widehat\kappa_t$.
Changing variables by the inverse flow and integrating in time gives,
for every interval $B\subset(-1/8,1/8)$,
\[
 2E[W\1_{\{Z\in B\}}]
 =E\!\left[W\int_{-1}^1D_{-t}\1_{\{Z-t\in B\}}\,\dd t\right].
\]
Choose $L=L(w)$ using \eqref{eq:crit-density-derivatives} so that
\[
 P\!\left(\sup_{|t|\le1}|\log D_t|>L\right)\le a_w/2.
\]
On $\{Y_I\in[1/4,3/4]\}$, we have $Z\in[-1/4,1/4]$, and the set
of $t\in[-1,1]$ with $Z-t\in B$ has length $|B|$. Since $W\le1$,
the weighted mass where this holds and all densities exceed $e^{-L}$
is at least $a_w/2$. Therefore
\[
 E[W\1_{\{Z\in B\}}]\ge\frac{a_we^{-L}}4|B|.
\]
The coordinate $\zeta$ has derivative one near $1/2$. Thus the weighted
law of $Y_I$ dominates $c_w\,\dd u$ there. The distribution-function
identity and continuity of $M_{W,I}'$ give the required pointwise bound.
\end{proof}

For the same weight $W$, a starting point $x\in\m S^1$, a length
$0<s<1$, and a fraction $v\in[0,1]$, define
\[
 G_W(x,s,v)=\m E_4\!\left[
 W\frac{H_4(x+sv)-H_4(x)}{H_4(x+s)-H_4(x)}\right].
\]
For $I=[a,a+r]$, put
\begin{align*}
 x(u)&=\Phi(a+r(u-c)),&
 s(u)&=\ell_I\bigl(k_I(u+c)-k_I(u-c)\bigr),\\
 x_0(u)&=\Phi(a)+\ell_I(u-c),& s_0&=2c\ell_I,\\
 \rho(u)&=\frac{k_I(u)-k_I(u-c)}{k_I(u+c)-k_I(u-c)}.
\end{align*}
Substitution into the definition of $F_I$ gives
\begin{align}\label{eq:crit-weighted-mean-identity}
 &\m E_4\bigl[W\{F_I(H_4\circ\varphi)-F_{\varphi(I)}(H_4)\}\bigr]\notag\\
 &\quad=-2\int\vartheta(u)
 \bigl\{G_W(x(u),s(u),\rho(u))-G_W(x_0(u),s_0,1/2)\bigr\}\,\dd u.
\end{align}
The two terms differ in their starting point, length, and split
fraction. Proposition~\ref{prop:crit-weighted-central} controls the last
variable. When $W=1$, reflection symmetry gives $G_1(x,s,1/2)=1/2$.
For a general weight, the dependence on $x$ and $s$ also matters.
The following bounds control these changes in the expansion used in
Proposition~\ref{prop:crit-weighted-cutoff}.

\begin{lem}[Regularity of weighted means]\label{lem:crit-weighted-regularity}
Let $W$ be as in Proposition~\ref{prop:crit-weighted-central}, with
observation set $\mathcal O=\{x_1,\ldots,x_m\}$. There are universal
constants $c_0,\eta_0>0$ such that the following estimates hold whenever
\[
 0<d\le1/8,\qquad 0<s\le c_0d,\qquad
 \operatorname{dist}([x,x+s],\mathcal O)\ge d,
 \qquad |v-1/2|\le\eta_0.
\]
Define
\[
 D_W(x,s)=\partial_{\log s}G_W(x,s,1/2),\qquad
 \partial_{\log s}=s\partial_s.
\]
All derivatives are taken with $W$ fixed. In particular,
$\partial_vG_W$ and $\partial_v^2G_W$ denote the first and second
derivatives with respect to the third argument $v$. Then
\[
 \begin{aligned}
 &|\partial_vG_W(x,s,v)|+|\partial_v^2G_W(x,s,v)|\\
 &\qquad+|D_W(x,s)|+|\partial_{\log s}D_W(x,s)|\le C.
 \end{aligned}
\]
For two such intervals $[x,x+s]$ and $[y,y+t]$, assume also that
$|x-y|\le d/4$ and $1/2\le s/t\le2$. In compatible interval lifts,
assume that the smallest interval containing both arcs stays at distance
at least $d$ from $\mathcal O$. Then
\begin{align*}
 |G_W(x,s,v)-G_W(y,s,v)|&\le C|x-y|/d,\\
 |\partial_vG_W(x,s,1/2)-\partial_vG_W(y,t,1/2)|
 &\le C\bigl(|x-y|/d+|\log(s/t)|\bigr)^{1/2},\\
 |D_W(x,s)-D_W(y,t)|
 &\le C\bigl((|x-y|/d)^{1/2}+|\log(s/t)|\bigr).
\end{align*}
The constant $C$ is independent of $W$, the number and locations of
its observations, and the parameters satisfying these conditions.
\end{lem}
\begin{proof}
Choose $c_0,\eta_0$ small enough for the fixed flow profiles below.
The separation condition ensures that their supports avoid $\mathcal O$,
so every flow fixes $W$. A flow supported near $[x,x+s]$, equal to
$z\mapsto x+e^\lambda(z-x)$ on that interval, gives
\[
 G_W(x,se^\lambda,v)=E[D_\lambda W H_{4,[x,x+s]}(v)].
\]
Lemma~\ref{lem:crit-gaussian-flow}, with a fixed profile scaled by $s$,
bounds the first two derivatives in $\lambda$. The flows fixing the
endpoints and moving the middle point give the derivatives in $v$.

A translation flow on an interval of length comparable to $d$ moves
all three evaluation points by $\tau d$ and fixes the observations.
Its first density derivative gives, uniformly on fixed neighborhoods
of $v=1/2$ and $\lambda=0$,
\[
 |G_W(x,se^\lambda,v)-G_W(y,se^\lambda,v)|\le C|x-y|/d.
\]
Together with the scale derivative bound, this gives
\[
 \sup_{|v-1/2|\le\eta_0}|G_W(x,s,v)-G_W(y,t,v)|
 \le C\bigl(|x-y|/d+|\log(s/t)|\bigr).
\]
For a function with uniformly bounded second derivative, the finite
difference estimate
\[
 |f'(u)|\le 2\|f\|_\infty/h+Ch
\]
on a smaller interval bounds $|f'(u)|$ by $C\|f\|_\infty^{1/2}$
when $\|f\|_\infty$ is small. Apply this in $v$ to obtain the bound for
$\partial_vG_W$, and in $\lambda$ to the translation difference to obtain
\[
 |D_W(x,s)-D_W(y,s)|\le C(|x-y|/d)^{1/2}.
\]
The bound on $\partial_{\log s}D_W$ supplies the remaining scale term.
For larger differences the uniform derivative bounds suffice.
\end{proof}

\subsection{Local distortion and symmetry}\label{subsec:crit-local-regularity}

The normalized image-length estimates rule out non-quasisymmetric maps.
For quasisymmetric maps, local affine deformations show that nonzero
common parts of the local laws force symmetry.

\begin{thm}[Singularity outside $\QS(\m S^1)$]\label{thm:crit-nonqs-total}
If $\varphi\in\Homeo_+(\m S^1)\setminus\QS(\m S^1)$, then there is
a closed proper arc $J$ such that
\[
 (\operatorname{Res}_J)_*\m P_{4,\varphi}\perp P_J.
\]
\end{thm}
\begin{proof}
By the quasisymmetry criterion in Lemma~\ref{lem:discrete WP}, choose
intervals $I_n=[a_n,a_n+r_n]$, with midpoints $x_n=a_n+r_n/2$, such that
\[
 v_n=\frac{\Phi(x_n)-\Phi(a_n)}{\Phi(a_n+r_n)-\Phi(a_n)}
 \longrightarrow0\quad\text{or}\quad1.
\]
After passing to a subsequence, $r_n\to0$. Indeed, for $r$ bounded
away from zero, the two half-interval increments have positive minima
by compactness and strict monotonicity. Pass to a further subsequence
with $a_n$ convergent. All remaining $I_n$ lie in the interior of one
closed proper arc $J$.

Suppose first that $v_n\to0$, and define the restriction-measurable
statistic
\[
 Z_n(h)=\frac{h(x_n)-h(a_n)}{h(a_n+r_n)-h(a_n)}.
\]
Lemma~\ref{lem:crit-balanced-tight} gives
\begin{equation}\label{eq:crit-nonqs-base-tail}
 \lim_{b\downarrow0}\limsup_n\m P_4(Z_n\le b)=0.
\end{equation}
On the other hand, with $J_n=\varphi(I_n)$,
\[
 Z_n(H_4\circ\varphi)=H_{4,J_n}(v_n)
 \longrightarrow0\quad\text{in probability},
\]
by Lemma~\ref{lem:crit-small-end-arc}, since $|J_n|\to0$ by uniform
continuity of $\varphi$.

Choose $b_k\downarrow0$ and then increasing indices $n_k$ such that
\[
 \m P_4(Z_{n_k}\le b_k)\le2^{-k},\qquad
 \m P_{4,\varphi}(Z_{n_k}>b_k)\le2^{-k}.
\]
Indeed, \eqref{eq:crit-nonqs-base-tail} supplies each $b_k$ with the
first bound for all sufficiently large indices; convergence under the
transformed law then supplies $n_k$. The Borel set
\[
 A=\{h:Z_{n_k}(h)\le b_k\text{ for all sufficiently large }k\}
\]
is measurable from $h|_J$. The two summable bounds and the first
Borel--Cantelli lemma give
\[
 \m P_4(A)=0,\qquad \m P_{4,\varphi}(A)=1.
\]
If $v_n\to1$, use $1-Z_n$ and the right-endpoint estimate of
Lemma~\ref{lem:crit-small-end-arc}. This proves local singularity.
\end{proof}

We next show that common parts of the local laws force a quasisymmetric
map to be symmetric. Here a \emph{nonzero common part} is a nonzero
finite measure dominated by both laws. The proof uses small affine
deformations to show that the relevant limits of normalized maps
preserve midpoints.

We first record continuity in total variation under right composition
by rotations. Recall that $R_t(x)=x+t$ on $\m R/\m Z$, so this action
sends $h$ to $h\circ R_t$. If $P$ is invariant under this action on
$\Homeo_+(\m S^1)$ with its uniform topology, and $\rho\ll P$ is a
probability, then
\begin{equation}\label{eq:crit-density-rotation}
 \|\bigl(h\mapsto h\circ R_t\bigr)_*\rho-\rho\|_{\mathrm{TV}}
 \longrightarrow0\qquad(t\to0).
\end{equation}
Here and below, $\|P-Q\|_{\mathrm{TV}}=\sup_A|P(A)-Q(A)|$ for finite
measures, with the supremum over Borel sets.
Indeed, writing $q=\dd\rho/\dd P$, the transformed density is
$q(h\circ R_{-t})$. It converges to $q$ in $L^1(P)$ first for
$q\in C_b$ by bounded convergence, and then for all $q\in L^1(P)$
by density and the isometry of right rotation.

\begin{lem}[Local affine approximations]\label{lem:crit-small-affine}
Let $a_n\in\m S^1$, $r_n\downarrow0$, and let $A>0$, $B\in\m R$
be fixed. For each fixed compact interval of $u$-values there are smooth
circle diffeomorphisms $g_n$ which, in lifts, satisfy
\[
 g_n(a_n+r_nu)=a_n+r_n(Au+B)
\]
on that interval, such that
\[
 g_n\longrightarrow\operatorname{id}\text{ uniformly},\qquad
 \|(R_{g_n})_*\m P_4-\m P_4\|_{\mathrm{TV}}\longrightarrow0.
\]
For every finite measure $\rho\ll\m P_4$, also
$\|(R_{g_n})_*\rho-\rho\|_{\mathrm{TV}}\to0$.
\end{lem}
\begin{proof}
We construct the affine map by composing a dilation and a translation,
and prove that each changes the law by a quantity tending to zero in
total variation. For the dilation, fix a small coordinate neighborhood
$(-2d,2d)$ of zero. For $R\downarrow0$, put $L=\log(d/R)$.
There are smooth even functions $a_R$, supported in $(-d,d)$, equal to
one on $[-R,R]$, with
\[
 0\le a_R\le1,\qquad
 \|a_R\|_{H^{1/2}(\m S^1)}^2\le C_d/L,\qquad
 \int a_R\,\dd x\le C_d/L.
\]
To construct $a_R$, start with a radial function on the upper half-plane
which is one for $|z|\le2R$, zero for $|z|\ge d/2$, and equals
$\log(d/(2|z|))/\widetilde L$ in between, where
$\widetilde L=\log(d/(4R))$. Its Dirichlet integral is
\[
 \int_0^\pi\int_{2R}^{d/2}\frac1{\widetilde L^2t^2}
 \,t\,\dd t\,\dd\theta
 =\frac\pi{\widetilde L}\le\frac C L.
\]
Smooth the two junctions on fixed relative annuli, keeping the function
one on the half-disk of radius $R$ and zero outside that of radius $d$. This
changes the bound by a fixed factor. The trace inequality gives the
$\dot H^{1/2}$ bound. Its trace has integral $O_d(L^{-1})$ and squared
$L^2$ norm $O_d(L^{-2}+R)$, giving the stated circle norm.

Choose a fixed smooth even nonnegative $b$, supported in
$(-2d,-d)\cup(d,2d)$, with $\int b\,\dd x=1$. Set
\[
 c_R=\int a_R\,\dd x,\qquad
 v_R(x)=\int_0^x(a_R(t)-c_Rb(t))\,\dd t.
\]
Evenness and zero total integral of $v_R'$ make $v_R$ vanish outside
$(-2d,2d)$. Moreover,
\[
 v_R(x)=x\ (|x|\le R),\quad
 \|v_R'\|_\infty\le C_d,\quad
 \|v_R\|_\infty\to0,\quad
 \|v_R'\|_{H^{1/2}}\to0.
\]
Translate this vector field to $a_n$ and choose $R$ to be a fixed large
multiple of $r_n$. Choose the multiple large enough to contain the trajectories from the
prescribed compact set up to time $\log A$. On this set, the flow is
exactly the required dilation. All these
flows and their inverses have uniformly bounded derivatives.

Use the operators $C_t,S_t,b_t$ defined after
Lemma~\ref{lem:crit-gaussian-flow}. The Fourier calculation in
Appendix~\ref{subsec:crit-gaussian-proof} gives
\[
 \|L_v+L_v^*\|_{\mathrm{HS}}\le C\|v'\|_{H^{1/2}}.
\]
The derivative bounds on the flows give uniform bounds for $C_t,C_t^{-1}$
and, on each fixed time interval,
\[
 \begin{aligned}
 S_t-I&=\int_0^t C_s(L_v+L_v^*)C_s^*\,\dd s,
 &b_t&=2\pi_0\int_0^t v'\circ g_s\,\dd s,\\
 \|S_t-I\|_{\mathrm{HS}}+\|b_t\|_{\mathcal H}
 &\le C\|v'\|_{H^{1/2}}\longrightarrow0.
 \end{aligned}
\]
Let $\gamma_t,\gamma$ be the transformed and reference field laws, and
$\gamma_{t,N},\gamma_N$ their first $N$ real Fourier marginals. The
compressed covariance $S_{t,N}$ has spectrum in a fixed compact
subinterval of $(0,\infty)$. Gaussian relative entropy and Pinsker's
inequality give
\[
 \|\gamma_{t,N}-\gamma_N\|_{\mathrm{TV}}^2
 \le C\left(\|P_Nb_t\|_{\mathcal H}^2+
       \operatorname{Tr}(S_{t,N}-I_N-\log S_{t,N})\right)
 \le C\left(\|b_t\|_{\mathcal H}^2+
            \|S_t-I\|_{\mathrm{HS}}^2\right).
\]
Conditional densities relative to $(\gamma_t+\gamma)/2$ converge along
the increasing Fourier sigma-fields. Hence
\[
 \|\gamma_t-\gamma\|_{\mathrm{TV}}
 =\sup_N\|\gamma_{t,N}-\gamma_N\|_{\mathrm{TV}}
 \le C\|v'\|_{H^{1/2}}\longrightarrow0.
\]
The welding coordinate rule and contraction under measurable maps give
the same conclusion for $\m P_4$. Only the first covariance derivative
is used, so the cutoff profile may vary with $R$.

A translation by $Br_n$ on the same compact set is generated by
$Br_n\chi_0(x-a_n)$, where $\chi_0$ is a fixed smooth cutoff equal to one
on a neighborhood of zero. Its derivative tends to zero in $H^{1/2}$,
so the same argument applies. Compose the dilation and translation.
The triangle inequality and invariance of total variation under a
measurable bijection prove the assertion for $g_n$.

Finally write $\rho=f\m P_4$. Put $T_n=R_{g_n}$ and
$P_n=(T_n)_*\m P_4$. For bounded continuous $f$,
$(T_n)_*(f\m P_4)=(f\circ T_n^{-1})P_n$. Since $g_n^{-1}$ also
converges uniformly to the identity,
\begin{align*}
 &\|(T_n)_*(f\m P_4)-f\m P_4\|_{\mathrm{TV}}\\
 &\quad\le 2\|f\|_\infty\|P_n-\m P_4\|_{\mathrm{TV}}
       +\m E_4[|f\circ T_n^{-1}-f|]\longrightarrow0.
\end{align*}
The second term tends to zero by dominated convergence. Bounded continuous
functions are dense in $L^1(\m P_4)$ on the homeomorphism space. If $f_j$
approximates a nonnegative $f\in L^1(\m P_4)$, the two approximation
errors are at most $2\m E_4[|f-f_j|]$, since pushforward preserves their
total variation. Letting first $n\to\infty$ and then $j\to\infty$
proves the assertion for every finite $\rho\ll\m P_4$.
\end{proof}

\begin{lem}[Midpoint rigidity]\label{lem:crit-weighted-tv-rigidity}
Let $W$ be a Borel function on $\Homeo_+(\m S^1)$ with $0\le W\le1$,
$w=\m E_4[W]>0$, and let $\rho=W\m P_4$.
Let $O$ be an open proper arc. Suppose $\psi_n$ are circle
homeomorphisms with a common quasisymmetry bound, $\psi_n\to\operatorname{id}$
uniformly, and, for every closed $J\Subset O$,
\[
 \|(\operatorname{Res}_J)_*(R_{\psi_n})_*\rho
       -(\operatorname{Res}_J)_*\rho\|_{\mathrm{TV}}\longrightarrow0.
\]
Then, for every closed arc $J\Subset O$,
\[
 \sup_{[x-r,x+r]\subset J}
 \left|\frac{\psi_n(x+r)-\psi_n(x)}{\psi_n(x)-\psi_n(x-r)}-1\right|
 \longrightarrow0,
\]
where $r>0$ and the differences are computed in increasing lifts.
\end{lem}
\begin{proof}
Suppose the conclusion fails. We will obtain a nontrivial limit of the
normalized maps and then show that it must be the identity. Choose
$I_n=[a_n,a_n+r_n]\Subset O$ in one compact subarc, with normalized
midpoint distortion bounded away from zero. Uniform
convergence of $\psi_n$ implies $r_n\to0$. Pass to a subsequence with
$a_n\to a\in O$. With compatible lifts, put
\[
 s_n=\psi_n(a_n+r_n)-\psi_n(a_n),\qquad
 k_n(u)=\frac{\psi_n(a_n+r_nu)-\psi_n(a_n)}{s_n}.
\]
Here $s_n\to0$. The uniform quasisymmetry bound and $k_n(0)=0$,
$k_n(1)=1$ give locally uniform bounds and a common H\"older modulus
for $k_n$ and their inverses on compact intervals. These bounds follow
by repeated bisection of adjacent intervals. A subsequence therefore
converges locally uniformly to a homeomorphism $k$ of $\m R$, with
\[
 k(0)=0,\qquad k(1)=1,\qquad k(1/2)\ne1/2.
\]

For $u<v<z$, define
\[
 M_n(u,v,z)=\m E_4\!\left[
 W\frac{H_4(a_n+r_nv)-H_4(a_n+r_nu)}
 {H_4(a_n+r_nz)-H_4(a_n+r_nu)}\right].
\]
To obtain a locally uniform limit of these means, approximate $W$ in
$L^1$ by conditional expectations $W_m$ on increasing finite subsets
of a fixed dense set avoiding $a$. Then $0\le W_m\le1$ and
$\m E_4[W_m]=w$. At large $n$, smooth flows moving any one of the
three displayed points and fixing the other two also fix $W_m$.
On each compact set of ordered triples with positive separation, the
scaled density derivative bounds give a common Lipschitz bound for
$M_n$ with $W$ replaced by $W_m$. The bound is independent of $m$, although the index from which it
holds may depend on $m$. Also,
\[
 \sup_{u<v<z}|M_n(u,v,z)-M_n^{(m)}(u,v,z)|
 \le\m E_4[|W-W_m|].
\]
It follows that $(M_n)$ is asymptotically equicontinuous on those
compact sets. A diagonal subsequence converges locally uniformly to
a continuous function $M$.

Lemma~\ref{lem:crit-small-affine} applied at scale $r_n$, including its
assertion for $\rho$, gives for every $A>0,B\in\m R$,
\[
 M(Au+B,Av+B,Az+B)=M(u,v,z).
\]
Hence $M(u,v,z)=m((v-u)/(z-u))$ for a continuous nondecreasing $m$.
Proposition~\ref{prop:crit-weighted-central}, applied first to $W_m$
and then passed to $W$, gives
\[
 m(v)-m(u)\ge c_w(v-u)
 \quad\text{if }1/2-\eta_0\le u<v\le1/2+\eta_0.
\]
In particular, monotonicity implies
\begin{equation}\label{eq:crit-midpoint-unique}
 m(v)=m(1/2)\quad\Longrightarrow\quad v=1/2.
\end{equation}

Fix $t\in\m R$. Lemma~\ref{lem:crit-small-affine}, at scale $s_n$
and center $\psi_n(a_n)$, gives a smooth map $\tau_n$ equal to
translation by $s_nt$ on the required compact set, with
\[
 \epsilon_n:=\|(R_{\tau_n})_*\rho-\rho\|_{\mathrm{TV}}\longrightarrow0.
\]
Choose fixed arcs
\[
 a\in\operatorname{int}J_0,\qquad
 J_0\Subset\operatorname{int}J_1,\qquad
 J_1\Subset\operatorname{int}J_2\Subset O,
\]
and put
\[
 d_n:=\|(\operatorname{Res}_{J_2})_*(R_{\psi_n})_*\rho
              -(\operatorname{Res}_{J_2})_*\rho\|_{\mathrm{TV}}
 \longrightarrow0.
\]
For large $n$,
\[
 \psi_n^{-1}(J_1)\subset J_2,\qquad
 \tau_n(\psi_n(J_0))\subset J_1.
\]
Set $q_n=\psi_n^{-1}\circ\tau_n\circ\psi_n$. Pullback and contraction
of total variation give
\[
 \|(\operatorname{Res}_{J_1})_*(R_{\psi_n^{-1}})_*\rho
                 -(\operatorname{Res}_{J_1})_*\rho\|_{\mathrm{TV}}
 \le d_n.
\]
Since $R_{q_n}=R_{\psi_n}\circ R_{\tau_n}\circ R_{\psi_n^{-1}}$,
\begin{align*}
 (R_{q_n})_*\rho-\rho
 &=(R_{\psi_n}\circ R_{\tau_n})_*
       \bigl((R_{\psi_n^{-1}})_*\rho-\rho\bigr)\\
 &\quad+(R_{\psi_n})_*\bigl((R_{\tau_n})_*\rho-\rho\bigr)
       +\bigl((R_{\psi_n})_*\rho-\rho\bigr).
\end{align*}
The two inclusions above therefore imply
\[
 \|(\operatorname{Res}_{J_0})_*(R_{q_n})_*\rho
                 -(\operatorname{Res}_{J_0})_*\rho\|_{\mathrm{TV}}
 \le 2d_n+\epsilon_n\longrightarrow0.
\]

In normalized coordinates the conjugate is
\[
 q_{n,t}(u)=k_n^{-1}(k_n(u)+t)
 \longrightarrow q_t(u)=k^{-1}(k(u)+t)
\]
locally uniformly. Equicontinuity of $M_n$ now gives
\[
 M(q_t(u),q_t(v),q_t(z))=M(u,v,z).
\]
Take $v=(u+z)/2$. By \eqref{eq:crit-midpoint-unique},
\[
 q_t((u+z)/2)=\tfrac12(q_t(u)+q_t(z)).
\]
The continuous map $q_t$ preserves midpoints and is therefore affine.
For $t\ne0$ it has no fixed point, so its slope is one. Write
$q_t(u)=u+b(t)$. The group law makes $b$ additive, and continuity gives
$b(t)=bt$. Since $q_1(0)=k^{-1}(1)=1$, we have $b=1$. Thus
$k(u+t)=k(u)+t$, and $k(0)=0$ implies $k=\operatorname{id}$, contrary
to $k(1/2)\ne1/2$.
\end{proof}

\begin{lem}[Symmetry from conjugated rotations]\label{lem:crit-conjugated-rotations}
For a circle homeomorphism $g$, let
\[
 K(g)=\sup_{x,\,0<r\le1/2}
 \max\left\{\frac{g(x+r)-g(x)}{g(x)-g(x-r)},
            \frac{g(x)-g(x-r)}{g(x+r)-g(x)}\right\},
\]
using degree-one lifts. If $\varphi\in\Homeo_+(\m S^1)$ satisfies
\[
 K(\varphi\circ R_t\circ\varphi^{-1})\longrightarrow1
 \qquad(t\to0),
\]
then $\varphi\in\Sy(\m S^1)$.
\end{lem}
\begin{proof}
Let $G$ be a degree-one lift of $\varphi$, and smooth it by a nonnegative
approximate identity $\eta_\epsilon$ to obtain $G_\epsilon$. Since $\dd G$
has full support, $G_\epsilon$ is a smooth increasing diffeomorphism.
For $k_\epsilon=G_\epsilon\circ G^{-1}$,
\[
 k_\epsilon(x)=\int G(G^{-1}(x)+t)\eta_\epsilon(t)\,\dd t.
\]
Averaging the adjacent-interval inequalities of
$\varphi\circ R_t\circ\varphi^{-1}$ gives $K(k_\epsilon)\to1$.
Repeated bisection bounds the distance of any normalized restriction
of a $K$-quasisymmetric map from the identity by $\omega(K)$, where
$\omega(K)\to0$ as $K\downarrow1$. The inverse restriction has the
same bound, since $\|f^{-1}-\operatorname{id}\|_\infty
=\|f-\operatorname{id}\|_\infty$ for increasing homeomorphisms of $[0,1]$.

Using the normalized restrictions of Section~\ref{subsec:crit-shape},
the identity $G=k_\epsilon^{-1}\circ G_\epsilon$ becomes
\[
 G_J=(k_\epsilon^{-1})_{G_\epsilon(J)}\circ(G_\epsilon)_J.
\]
Consequently,
\[
 \|G_J-\operatorname{id}\|_\infty
 \le\omega(K(k_\epsilon))
       +\|(G_\epsilon)_J-\operatorname{id}\|_\infty.
\]
For fixed $\epsilon$, the second term tends uniformly to zero as
$|J|\to0$, because $G_\epsilon'$ is continuous and positive. Taking
this limit and then $\epsilon\downarrow0$ proves symmetry.
\end{proof}

\begin{lem}[Common parts imply symmetry]\label{lem:crit-common-sym}
Let $\varphi\in\QS(\m S^1)$. Suppose the local laws
$P_I$ and $(\operatorname{Res}_I)_*\m P_{4,\varphi}$ have a nonzero
common part for each arc in a finite family whose interiors cover the
circle. Then $\varphi\in\Sy(\m S^1)$.
\end{lem}
\begin{proof}
Denote the arcs by $I_j$. Choose nonzero finite measures
\[
 \lambda_j\le P_{I_j},\qquad
 \lambda_j\le(\operatorname{Res}_{I_j})_*\m P_{4,\varphi}.
\]
Write $Q_j=(\operatorname{Res}_{I_j})_*\m P_{4,\varphi}$ and choose
the Radon--Nikodym densities of $\lambda_j$ relative to $P_{I_j}$ and
$Q_j$ with values in $[0,1]$. Define
\[
 V_j(h)=\frac{\dd\lambda_j}{\dd P_{I_j}}(h|_{I_j}),\qquad
 W_j(h)=\frac{\dd\lambda_j}{\dd Q_j}((h\circ\varphi)|_{I_j}),
\]
and set $\rho_j=V_j\m P_4$, $\mu_j=W_j\m P_4$.
The first weight depends only on $h|_{I_j}$ and the second only on
$h|_{\varphi(I_j)}$. Both measures have the same positive total mass
as $\lambda_j$, and by construction
\[
 (\operatorname{Res}_{I_j})_*\rho_j
 =(\operatorname{Res}_{I_j})_*(R_\varphi)_*\mu_j=\lambda_j.
\]

Set $c_t=\varphi\circ R_t\circ\varphi^{-1}$. For
$J\Subset\varphi(\operatorname{int}I_j)$ and sufficiently small $t$,
\[
 R_t(\varphi^{-1}(J))\subset\operatorname{int}I_j.
\]
The shared restriction law $\lambda_j$ can therefore be pushed forward:
\[
 (\operatorname{Res}_J)_*(R_{c_t})_*\mu_j
 =\bigl(u\mapsto(u\circ R_t\circ\varphi^{-1})|_J\bigr)_*\lambda_j.
\]
Using $\lambda_j=(\operatorname{Res}_{I_j})_*\rho_j$ gives
\begin{align*}
 &\|(\operatorname{Res}_J)_*(R_{c_t})_*\mu_j
       -(\operatorname{Res}_J)_*\mu_j\|_{\mathrm{TV}}\\
 &\quad\le\|\bigl(h\mapsto h\circ R_t\bigr)_*\rho_j-\rho_j\|_{\mathrm{TV}}
 \longrightarrow0
\end{align*}
by \eqref{eq:crit-density-rotation}, applied after normalization.

The maps $c_t$ tend uniformly to the identity and have a common
quasisymmetry bound, since $\varphi$ is quasisymmetric. Apply
Lemma~\ref{lem:crit-weighted-tv-rigidity} to each $\mu_j$ and any
sequence $t\to0$. On compact subarcs covering the circle, their
adjacent-interval ratios tend uniformly to one. Intervals of length
bounded below are handled by uniform convergence to the identity.
Consequently $K(c_t)\to1$.

Lemma~\ref{lem:crit-conjugated-rotations} now gives
$\varphi\in\Sy(\m S^1)$.
\end{proof}

\subsection{Covariance and accumulated distortion}\label{subsec:crit-energy}

Covariance decay gives $L^2$ bounds for sums of $F_I$ and
$F_{\varphi(I)}$. A reconstruction formula then compares $w_I^2$ with
$\varepsilon_I^2$. For symmetric non-WP maps, this comparison supplies
the divergent normalization used in the singularity proof.

\begin{prop}[Covariance and multiscale bounds]\label{prop:crit-two-frames}
For proper arcs $I=[a,a+r]$ and $J=[b,b+R]$ with $0<r\le R<1$,
the critical model satisfies
$\m E_4[F_I(H_4)]=0$ and
\begin{equation}\label{eq:crit-packet-covariance}
 |\m E_4[F_I(H_4)F_J(H_4)]|
 \le C\begin{cases}r/R,&\operatorname{dist}(a,b)\le AR,\\
 r^2/\operatorname{dist}(a,b)^2,&\operatorname{dist}(a,b)>AR,
 \end{cases}
\end{equation}
where $A,C$ are fixed constants.

Moreover, if $\varphi\in\Sy(\m S^1)$, then for every deterministic
$f\in L^2(\Lambda)$ supported in $r\in[\delta,r_0]$ for some
$0<\delta<r_0<1$,
\begin{align*}
 &\m E_4\!\left[\left|\int f(I)F_I(H_4)\,\dd\Lambda(I)\right|^2\right]\\
 &\quad+\m E_4\!\left[\left|\int f(I)F_{\varphi(I)}(H_4)\,\dd\Lambda(I)\right|^2\right]\\
 &\qquad\le C_\varphi\int |f(I)|^2\,\dd\Lambda(I).
\end{align*}
\end{prop}
\begin{proof}
We first estimate the covariance by translating one triple of points
while fixing the other. Let
$z_-=z_0-S,z_0,z_+=z_0+S$ be at distance at least $d\ge8s$ from $x$.
Choose a smooth profile equal to one near zero and supported in
$(-1/2,1/2)$, and scale its vector field by $d$. For small time the flow
translates $x-s,x,x+s$ by $td$ and fixes the three points $z_-,z_0,z_+$.
Lemma~\ref{lem:crit-gaussian-flow},
applied to the product of the two bounded split variables, gives
\begin{equation}\label{eq:crit-six-point}
 \left|\partial_x^2 E[\xi_{x,s}(H_4)\xi_{z_0,S}(H_4)]\right|
 \le Cd^{-2}.
\end{equation}
The two factors $d^{-1}$ come from changing the flow time to the
translation coordinate. The estimate also holds for the distributional derivative.

Reflection symmetry gives $E[F_I(H_4)]=0$. Suppose $r\le R$ and put
$x_*=a+r/2$. For the parameter $v$ in $F_J$, denote the three evaluation points by
$z_j(v)=b+Rv+jcR$, $j=-1,0,1$, and set
$d(v)=\min_j\operatorname{dist}(x_*,z_j(v))$. The set of $v$ with
$d(v)\le16r$ has length at most $Cr/R$; there we use $|\xi|\le1$.

On its complement, \eqref{eq:crit-six-point} applies to every
evaluation point in $F_I$. Since $\int\vartheta=\int(u-1/2)\vartheta(u)\,\dd u=0$,
Taylor's theorem gives a contribution bounded by
\[
 Cr^2\int_{\{d(v)>16r\}}\frac{\dd v}{d(v)^2}
 \le C\frac rR.
\]
For the last bound, sum over the three marked points and use
$\dd z_j=R\,\dd v$; each integral is at most
$R^{-1}\int_{|z-x_*|>16r}|z-x_*|^{-2}\,\dd z\le C/(Rr)$.
If $\operatorname{dist}(a,b)>AR$ for a sufficiently large $A$, all the
separations are comparable to $\operatorname{dist}(a,b)$ and the same
argument gives $Cr^2/\operatorname{dist}(a,b)^2$.

We next turn the covariance estimate into an $L^2$ bound. Let $\sigma$
be a positive measure on arc parameters assigning mass at most $K$ to
each box $[b,b+R]\times[R,2R]$. We claim that
\[
 \m E_4\left[\left|\int f(I)F_I(H_4)\,\dd\sigma(I)\right|^2\right]
 \le CK\int|f(I)|^2\,\dd\sigma(I).
\]
Initially $f$ is bounded and supported on finitely many scale layers.
Take dyadic scale layers $r\asymp2^{-j}$ and $R\asymp2^{-k}$, with
$r\le R$. For a fixed fine arc, $O(1)$ coarse boxes meet the near region;
by \eqref{eq:crit-packet-covariance}, their contribution to the integral
of the absolute covariance is at most $CKr/R$. At distance
$\asymp2^nR$ there are $O(2^n)$ coarse boxes, each contributing at most
$CKr^2/(2^nR)^2$. Summing in $n$ gives the same bound.

For a fixed coarse arc there are $O(R/r)$ fine boxes in the near region,
so the corresponding integral is at most $CK$. At distance comparable to
$2^nR$, there are $O(2^nR/r)$ fine boxes. Their total contribution is
at most $CK(r/R)2^{-n}$. Summing in $n$ gives a bounded far contribution;
periodic lifts handle boxes meeting the circle seam. The two-sided Schur
estimate therefore bounds the operator
between these scale layers by $CK\sqrt{r/R}$. If $f_j$ is the restriction
of $f$ to the $j$th layer, the full covariance form is at most
\[
 CK\sum_{j,k}2^{-|j-k|/2}\|f_j\|_{L^2(\sigma)}
 \|f_k\|_{L^2(\sigma)}
 \le CK\sum_j\|f_j\|_{L^2(\sigma)}^2.
\]
The last inequality follows from Cauchy--Schwarz and summability of
$2^{-|j-k|/2}$. Approximation by bounded functions extends the estimate to $L^2(\sigma)$.

To obtain the bound for the image arcs, push $\Lambda$ forward under
\[
 T_\varphi(a,r)=(\Phi(a),\Phi(a+r)-\Phi(a)).
\]
This measure has uniformly bounded mass on the same boxes.
Let $q=\Phi^{-1}$ and let $M\ge1$ be an adjacent-interval quasisymmetry
constant for $q$. Consider the target box
$[b,b+R]\times[R,2R]$ and set $s=q(b+R)-q(b)$. A preimage start $a$
lies in an interval of length $s$. If $\Phi(a)\le b+R/2$, the target
arc contains $[b+R/2,b+R]$; if $\Phi(a)\ge b+R/2$, it contains
$[b+R,b+3R/2]$. Quasisymmetry bounds the preimage length of either
interval below by $s/[M(M+1)]$. Hence each preimage arc has
$r\ge s/[M(M+1)]$, and its parameter mass is at most
\[
 s\int_{s/[M(M+1)]}^\infty\frac{\dd r}{r^2}=M(M+1).
\]
Compatible lifts handle boxes crossing the circle seam.
Applying the preceding $L^2$ bound first to $\sigma=\Lambda$, and then
to $\sigma=(T_\varphi)_*\Lambda$ using the measurable inverse of
$T_\varphi$ on its image, proves both asserted estimates.
\end{proof}

We next compare the coefficients $w_I$ with the full distortion
$\varepsilon_I$. For the remainder of this subsection, fix a symmetric
map $\varphi$ and choose $r_0>0$ sufficiently small. For $j\ge0$, set
\begin{align*}
 e_j^{(4)}&=\int_{2^{-j-1}r_0\le r\le2^{-j}r_0}w_I^2\,\dd\Lambda(I),\\
 s_j^{(4)}&=\int_{2^{-j-1}r_0\le r\le2^{-j}r_0}\varepsilon_I^2\,\dd\Lambda(I).
\end{align*}
\begin{lem}[Reconstruction from distortion coefficients]\label{lem:crit-calderon-reconstruction}
For the fixed $c$ and $\chi\in C_c^\infty((0.4,0.6))$ chosen above,
set $\psi=\Delta_c\chi''/c$.
Then
\[
 C_\psi=\int_0^\infty|\widehat\psi(t)|^2\,\frac{\dd t}{t}
 \in(0,\infty).
\]
Let $G$ be a degree-one lift of a quasisymmetric circle homeomorphism,
and define $\ell_I,w_I$ from this lift as above. Put
\[
 A(y,t)=\int\psi(u)G(y+tu)\,\dd u.
\]
For each $r_0>0$ there is a smooth periodic function $H_0$ such that
\[
 G(z)=z+H_0(z)+\frac1{C_\psi}\int_0^{r_0}\int_{\mathbb R}
 A(y,t)\psi\!\left(\frac{z-y}{t}\right)\frac{\dd y\,\dd t}{t^2}.
\]
The small-scale integral converges absolutely and uniformly in $z$.
Moreover, $A(a,r)=\ell_Iw_I$ for $I=[a,a+r]$.
\end{lem}
\begin{proof}
With the Fourier convention $\widehat f(\xi)=\int e^{-\mathrm{i}\xi u}
 f(u)\,\dd u$,
\[
 \widehat\psi(\xi)=\frac{4\xi^2\sin^2(c\xi/2)}c\widehat\chi(\xi).
\]
It vanishes to order at least four at zero and is rapidly decreasing at
infinity. Since $\chi$ is nonzero, $\psi$ is nonzero and $C_\psi>0$.
The vanishing moments also eliminate the affine part of $G$ and give
$A(a,r)=\ell_Iw_I$.

A quasisymmetric circle homeomorphism is H\"older continuous. Thus, for some
$\alpha>0$, $|A(y,t)|\le Ct^\alpha$ at small scales. The $y$-support of
$\psi((z-y)/t)$ has length $O(t)$, so the absolute small-scale integral is
bounded by $C\int_0^{r_0}t^{\alpha-1}\,\dd t$, uniformly in $z$.

Write $G(z)=z+v(z)$, where $v$ is continuous and periodic. For a nonzero
Fourier mode $e^{2\pi\mathrm{i}nz}$, the full-scale analysis and synthesis
multiplier is
\[
 \frac1{C_\psi}\int_0^\infty
 \widehat\psi(-2\pi nt)\widehat\psi(2\pi nt)\,\frac{\dd t}{t}=1,
\]
because $\psi$ is real. The contribution of $t\ge r_0$ has a Fourier
multiplier which decays faster than any power of $|n|$; together with the
constant mode it defines a smooth periodic $H_0$. The equality follows
first for trigonometric polynomials and then for $v$ by smooth periodic
mollification, which preserves the H\"older bound, and uniform absolute
convergence.
\end{proof}

\begin{prop}[Comparison of distortion energies]\label{prop:crit-shape-layers}
For $\varphi\in\Sy(\m S^1)$ and a sufficiently small fixed $r_0>0$,
\begin{equation}\label{eq:crit-shape-layers}
 s_j^{(4)}\le C_\varphi\left(
 \sum_{n\ge0}2^{-3n/2}e_{j+n}^{(4)}
 +\sum_{n=0}^j2^{-n/2}e_{j-n}^{(4)}+2^{-j/2}\right).
\end{equation}
\end{prop}
\begin{proof}
Use $a_*=7/8$ and $b_*=9/8$ in Lemma~\ref{lem:sym ab}, decreasing
$r_0$ so that its bounds apply to all common parent intervals below.
Let $G$ be a lift of $\varphi$, and let $D_If$ denote $f$ minus its
affine interpolant on $I$. Thus $\|D_IG\|_\infty=\ell_I\varepsilon_I$.
Lemma~\ref{lem:crit-calderon-reconstruction} gives the uniformly
convergent decomposition
\[
 G(z)=z+H_0(z)+\sum_{k\ge0}\mathcal W_k(z),\qquad h=2^{-k}r_0,
\]
\[
 \mathcal W_k(z)=\frac1{C_\psi}
 \int_{h/2}^h\int_{\mathbb R}A(y,t)
 \psi\!\left(\frac{z-y}{t}\right)\frac{\dd y\,\dd t}{t^2}.
\]
Set $q_k(y)=\int_{h/2}^h|w_{[y,y+t]}|^2\,\dd t/t$. Then
\begin{equation}\label{eq:crit-qk}
 \frac h2e_k^{(4)}\le\int_{\m S^1}q_k(y)\,\dd y\le he_k^{(4)}.
\end{equation}
Fix a parent layer $R/2\le r\le R$, where $R=2^{-j}r_0$, and put
$U_k(I)=\|D_I\mathcal W_k\|_\infty/\ell_I$.
The support of $\psi$, Taylor's theorem for $D_I\psi$, and
Lemma~\ref{lem:sym ab} give
\[
 U_k(I)^2\le C_\varphi\begin{cases}
 (h/R)^{2a_*}h^{-1}
   \displaystyle\int_{a-Ch}^{a+r+Ch}q_k(y)\,\dd y,&h\le R,\\[2mm]
 (R/h)^{4-2b_*}h^{-1}
   \displaystyle\int_{a-Ch}^{a+Ch}q_k(y)\,\dd y,&h>R.
 \end{cases}
\]
Indeed, in the first case the image-length ratio is at most
$C_\varphi(h/R)^{a_*}$, while in the second it is at most
$C_\varphi(h/R)^{b_*}$ and subtracting the affine interpolant supplies
$(R/h)^2$. Cauchy--Schwarz in $\dd y\,\dd t/t^2$ gives the displayed
squared bounds. For fixed $y$, the allowed starts $a$ have length
$O(R)$ in the first case and $O(h)$ in the second. Integrating and
using \eqref{eq:crit-qk} therefore gives
\begin{align}
 \int_{R/2}^R\int_{\m S^1}U_k(I)^2\,\dd a\,
       \frac{\dd r}{r^2}
 &\le C_\varphi(h/R)^{2a_*}e_k^{(4)} &&(h\le R),
       \label{eq:crit-fine-layer}\\
 &\le C_\varphi(R/h)^{3-2b_*}e_k^{(4)} &&(h>R).
       \label{eq:crit-coarse-layer}
\end{align}
For finite sums, weighted Cauchy--Schwarz yields
\[
 \left(\sum_kU_k(I)\right)^2
 \le\left(\sum_k2^{-|k-j|/4}\right)
       \sum_k2^{|k-j|/4}U_k(I)^2.
\]
Since $2a_*=7/4$ and $3-2b_*=3/4$, integration gives the two
convolutions in \eqref{eq:crit-shape-layers}. The smooth term contributes
at most $C_\varphi R^{3-2b_*}\le C_\varphi2^{-j/2}$, because
$\|D_IH_0\|_\infty\le C_\varphi r^2$ and
$\ell_I\ge C_\varphi^{-1}r^{b_*}$ by \eqref{eq:nonwp-sym-powers}.
The affine term vanishes.
Uniform convergence of the reconstruction and Fatou's lemma complete
the proof.
\end{proof}

The preceding estimates now give the normalization and the local
statistics used in the singularity proof.

\begin{definition*}[Accumulated distortion and statistics]
For the fixed symmetric map $\varphi$ and the layers $e_j^{(4)},s_j^{(4)}$ above,
set
\[
 \mathcal E_N=1+\sum_{j=0}^Ne_j^{(4)},\qquad
 \mathcal S_N=\sum_{j=0}^Ns_j^{(4)}.
\]
For $N\ge0$, put
\[
 \Omega_N^{(4)}=\m S^1\times[2^{-N-1}r_0,r_0].
\]
For any deterministic Borel function $\beta:\m S^1\to[0,1]$, define
\begin{align*}
 \mathcal T_N^\beta(h)&=\frac1{\mathcal E_N}
  \int_{\Omega_N^{(4)}}\beta(a)w_IF_I(h)\,\dd\Lambda(I),\\
 \mathcal B_N^\beta(h)&=\frac1{\mathcal E_N}
  \int_{\Omega_N^{(4)}}\beta(a)w_IF_{\varphi(I)}(h)\,\dd\Lambda(I),\\
 \mathcal A_N^\beta(h)&=\mathcal T_N^\beta(h\circ\varphi)
                         -\mathcal B_N^\beta(h).
\end{align*}
\end{definition*}

\begin{lem}[Accumulated energy and local statistics]\label{lem:crit-record}\label{lem:crit-statistics}
If $\varphi\in\Sy(\m S^1)\setminus\WP(\m S^1)$, then
$\mathcal E_N\to\infty$, and there is a deterministic subsequence
$N_k\uparrow\infty$ such that
$\mathcal S_{N_k}\le C_\varphi\mathcal E_{N_k}$.

Uniformly over Borel $\beta:\m S^1\to[0,1]$ and along this subsequence,
\begin{equation}\label{eq:crit-global-frame}
 \begin{aligned}
 \|\mathcal T_N^\beta\|_{L^2(\m P_4)}
 +\|\mathcal B_N^\beta\|_{L^2(\m P_4)}
 &\le C_\varphi\mathcal E_N^{-1/2},\\
 \sup_h|\mathcal A_N^\beta(h)|&\le C_\varphi.
 \end{aligned}
\end{equation}
The same bounds hold after omitting a fixed range of coarse scales.
If $\beta$ is supported in a compact subset of the interior of a proper
arc $J$, then, at sufficiently small scales, $\mathcal T_N^\beta$ depends
only on $h|_J$, and $\mathcal B_N^\beta,\mathcal A_N^\beta$ depend only
on $h|_{\varphi(J)}$.
\end{lem}

\begin{proof}
For $I=[x-y,x+y]$, symmetry and the definitions give, at sufficiently
small scales,
\[
 \delta_\varphi(x,y)
 =\log\frac{1-k_I(1/2)}{k_I(1/2)},\qquad
 |\delta_\varphi(x,y)|\le C\varepsilon_I.
\]
The change of variables $a=x-y$, $r=2y$ therefore yields
\[
 E(\varphi)\le C_\varphi
       \left(1+\sum_{j\ge0}s_j^{(4)}\right).
\]
The coarse scales contribute a finite constant because $\varphi$ is
quasisymmetric. If $\sum_j e_j^{(4)}<\infty$,
Proposition~\ref{prop:crit-shape-layers} gives
$\sum_j s_j^{(4)}<\infty$, and Lemma~\ref{lem:discrete WP} would imply
$\varphi\in\WP(\m S^1)$. Thus $\mathcal E_N\uparrow\infty$.
Since $|w_I|\le C$ and the $j$th layer has $\Lambda$-mass $O(2^j)$,
$\mathcal E_N=O(2^N)$. Hence
$\mathcal E_N2^{-5N/4}\to0$. This positive sequence attains its maximum
on every tail. We may therefore choose infinitely many indices $N$ for which
\[
 \mathcal E_{N+n}2^{-5(N+n)/4}
 \le\mathcal E_N2^{-5N/4}\quad(n\ge0).
\]
Summing Proposition~\ref{prop:crit-shape-layers} over $j\le N$, the terms with index at most $N$ are bounded
by $C\mathcal E_N$. The remaining terms are bounded by
\[
 C\sum_{n\ge1}2^{-3n/2}\mathcal E_{N+n}
 \le C\mathcal E_N\sum_{n\ge1}2^{-n/4}.
\]
This proves the claim along the selected indices.

For the statistics, Proposition~\ref{prop:crit-two-frames} and $\beta^2\le1$ give
\[
 \begin{aligned}
 \|\mathcal T_N^\beta\|_2^2+\|\mathcal B_N^\beta\|_2^2
 &\le\frac{C_\varphi}{\mathcal E_N^2}
              \int_{\Omega_N^{(4)}}w_I^2\,\dd\Lambda(I)\\
 &\le\frac{C_\varphi}{\mathcal E_N}.
 \end{aligned}
\]
By Lemma~\ref{lem:crit-packet-bv} and Cauchy--Schwarz,
\begin{align*}
 |\mathcal A_N^\beta(h)|
 &\le\frac C{\mathcal E_N}\int_{\Omega_N^{(4)}}|w_I|\varepsilon_I\,\dd\Lambda(I)\\
 &\le C\sqrt{\frac{\mathcal E_N-1}{\mathcal E_N}}
          \sqrt{\frac{\mathcal S_N}{\mathcal E_N}}\le C_\varphi.
\end{align*}
Restricting the integrals preserves these estimates. A fixed range
$r\ge\delta>0$ contributes $O_\delta(\mathcal E_N^{-1})$, since the
integrands are bounded and that range has finite $\Lambda$-mass.
Finally, for $r<\operatorname{dist}(\operatorname{supp}\beta,\partial J)$,
every contributing $I=[a,a+r]$ is contained in $J$. The stated
measurability follows from the definitions.
\end{proof}

\subsection{Proof of critical singularity}\label{subsec:crit-singularity}

It remains to treat $\varphi\in\Sy(\m S^1)\setminus\WP(\m S^1)$. In this case
$\varepsilon_I$ tends uniformly to zero at small scales, but the
accumulated distortion diverges. The statistics
$\mathcal T_N^\beta,\mathcal B_N^\beta$ have vanishing $L^2$ norms,
whereas their difference after composition has a positive weighted mean.
The identity below cancels the extra term caused by changing the image
length in the weighted mean.

\begin{lem}[Cancellation of the length term]\label{lem:crit-square-cancellation}
Let $\Phi$ be a degree-one lift of $\varphi\in\Homeo_+(\m S^1)$,
fix $0<r<1$, and put $I=[a,a+r]$. Define
\[
 U_r(a)=\int\vartheta(u)\Phi(a+ru)\,\dd u,\qquad
 L_I=\frac1{2c}\int\vartheta(u)
       \bigl(k_I(u+c)-k_I(u-c)\bigr)\,\dd u.
\]
Then $U_r$ is periodic and
\begin{align*}
 w_I&=\frac{U_r(a+cr)-2U_r(a)+U_r(a-cr)}{c\ell_I},\\
 L_I&=\frac{U_r(a+cr)-U_r(a-cr)}{2c\ell_I}.
\end{align*}
For every bounded measurable periodic $b(a)$, integration over one
period gives
\begin{align*}
 &\int b(a)w_IL_I\,\dd a\\
 &\quad=\frac1{2c^2}\int
 \left\{\frac{b(a)}{\ell_{[a,a+r]}^2}
 -\frac{b(a+cr)}{\ell_{[a+cr,a+(1+c)r]}^2}\right\}\\
 &\hspace{45mm}\times\bigl(U_r(a+cr)-U_r(a)\bigr)^2\,\dd a.
\end{align*}
Also,
\[
 |U_r(a+cr)-U_r(a)|\le C\ell_I\varepsilon_I.
\]
\end{lem}
\begin{proof}
The identities for $w_I,L_I$ follow by linearity; the constant terms
in $k_I$ vanish since $\int\vartheta=0$. The identity
\[
 (u_+-2u+u_-)(u_+-u_-)
 =(u_+-u)^2-(u-u_-)^2
\]
and the substitution $a\mapsto a+cr$ in the second integral give the
claimed formula. All integrals are finite at fixed $r$.

Both $\int\vartheta$ and $\int u\vartheta(u)\,\dd u$ vanish. Subtract
$\Phi(a)+\ell_I u$ inside $U_r(a)$ and
$\Phi(a)+\ell_I(u+c)$ inside $U_r(a+cr)$. All the arguments lie in
$I$ because $c<1/10$ and $\operatorname{supp}\vartheta\subset(0.4,0.6)$.
Each remainder is bounded by
$\|\vartheta\|_{L^1}\ell_I\varepsilon_I$.
\end{proof}

\begin{prop}[A weighted lower bound]\label{prop:crit-weighted-cutoff}
Let $\varphi\in\Sy(\m S^1)\setminus\WP(\m S^1)$ and use the
subsequence of Lemma~\ref{lem:crit-record}. Let $W$ be a bounded cylinder
function as in Proposition~\ref{prop:crit-weighted-central}, with
$\m E_4[W]=w>0$. Let $0\le\zeta\le1$ be a smooth
circle function and let $0\le\eta\le1$ be Lipschitz and vanish on a
neighborhood of every observation point. Put
\[
 \beta(a)=\zeta(a)\eta(\varphi(a)).
\]
For the statistics of Lemma~\ref{lem:crit-statistics}, there is a
constant $c_w>0$, depending on $w$ but not on $\zeta,\eta$ or the
number of observations, such that
\[
 \liminf_k\left\{\m E_4[W\mathcal A_{N_k}^\beta]
 -\frac{c_w}{2\mathcal E_{N_k}}
   \int_{\Omega_{N_k}^{(4)}}\beta(a)w_I^2\,\dd\Lambda(I)\right\}\ge0.
\]
The same assertion holds if a fixed range of coarse scales is omitted.
\end{prop}
\begin{proof}
Choose $0<d\le1/8$ so that $\eta$ vanishes within distance $4d$
of every observation point, and put
\[
 \omega(r)=\sup_{|J|\le(1+c)r}\varepsilon_J,\qquad
 t(r)=\sup_{|J|\le(1+c)r}|\varphi(J)|.
\]
Both tend to zero. If $\beta(a)$ or $\beta(a+cr)$ is nonzero,
the image of $[a,a+(1+c)r]$ stays at distance at least $4d-t(r)$
from the observations. For $t(r)<\min(d/4,c_0d)$, all translation
and dilation neighborhoods used below therefore satisfy
Lemma~\ref{lem:crit-weighted-regularity}, even when one cutoff value
is zero. On this region define
\[
 x_I^*=\Phi(a)+\ell_I(1/2-c),\quad s_I^*=2c\ell_I,\quad
 a_I=\partial_vG_W(x_I^*,s_I^*,1/2),\quad d_I=D_W(x_I^*,s_I^*).
\]
Proposition~\ref{prop:crit-weighted-central} gives $a_I\ge c_w$.
Lemma~\ref{lem:crit-weighted-regularity} bounds $a_I,d_I$ uniformly.

Write $e(u)=k_I(u)-u$. The arguments in
\eqref{eq:crit-weighted-mean-identity} satisfy
\begin{align*}
 \rho(u)-1/2&=-\frac{e(u+c)-2e(u)+e(u-c)}{4c}
                 +O(\varepsilon_I^2),\\
 \log\frac{s(u)}{s_I^*}
 &=\frac{e(u+c)-e(u-c)}{2c}+O(\varepsilon_I^2),\\
 x(u)-x_0(u)&=\ell_I e(u-c).
\end{align*}
First expand in $v$, then in $\log s$, and use the translation estimate
for $x$. Replacing the resulting derivatives by $a_I,d_I$ gives
\begin{equation}\label{eq:crit-weighted-packet-expanded}
 \m E_4\bigl[W\{F_I(H_4\circ\varphi)-F_{\varphi(I)}(H_4)\}\bigr]
 =\tfrac12a_Iw_I-2d_IL_I+R_I,
\end{equation}
where, for small $I$,
\[
 |R_I|\le C\varepsilon_I
 \left(\varepsilon_I^{1/2}+\varepsilon_I+
       (\ell_I/d)^{1/2}+\ell_I/d\right).
\]
The two derivative errors are $O(\varepsilon_I^{1/2}+(\ell_I/d)^{1/2})$
and $O((\ell_I/d)^{1/2})$, respectively. The Taylor remainders are
$O(\varepsilon_I^2)$, and the endpoint translation costs
$O(\varepsilon_I\ell_I/d)$.

We treat the term $d_IL_I$ before estimating $R_I$. For fixed $r$, apply
Lemma~\ref{lem:crit-square-cancellation} with $b(a)=\beta(a)d_I$.
Set $b=0$ wherever $\beta=0$; no derivative is needed outside the
neighborhood just described. Write $I'=[a+cr,a+(1+c)r]$ and
$J=[a,a+(1+c)r]$. From $\|k_J-\operatorname{id}\|_\infty\le\omega(r)$,
\[
 \frac{\ell_I}{\ell_J}=\frac1{1+c}+O(\omega(r)),\qquad
 \frac{\ell_{I'}}{\ell_J}=\frac1{1+c}+O(\omega(r)).
\]
In particular, for sufficiently small $r$,
\[
 \ell_{I'}/\ell_I=1+O(\omega(r)),\qquad
 |x_{I'}^*-x_I^*|\le C\ell_I.
\]
The regularity estimate and the Lipschitz bounds on the cutoffs imply
\begin{align*}
 |d_{I'}-d_I|&\le C\bigl(\omega(r)+(t(r)/d)^{1/2}\bigr),\\
 |\beta(a+cr)-\beta(a)|&\le C_{\zeta,\eta}(r+t(r)).
\end{align*}
Subtract the coefficients by changing $\beta$, $d_I$, and
$\ell_I^{-2}$ in turn. The preceding bounds give
\begin{align*}
 \left|\frac{\beta(a)d_I}{\ell_I^2}
 -\frac{\beta(a+cr)d_{I'}}{\ell_{I'}^2}\right|
 &\le\frac{C}{\ell_I^2}
       \left(\omega(r)+\sqrt{t(r)/d}+r+t(r)\right)\\
 &=:\frac{q(r)}{\ell_I^2},\qquad q(r)\longrightarrow0.
\end{align*}
Here $C$ may depend on the cutoffs and $d$, but not on $r$. The bound
also holds when both coefficients vanish.
Lemma~\ref{lem:crit-square-cancellation} therefore gives
\begin{equation}\label{eq:crit-weighted-cross-small}
 \left|\int\beta(a)d_Iw_IL_I\,\dd a\right|
 \le Cq(r)\int\varepsilon_I^2\,\dd a.
\end{equation}

Write $|R_I|\le q_1(r)\varepsilon_I$, where $q_1(r)\to0$ by the
remainder bound above. The coarse scales $r\ge\delta>0$ contribute
$O_\delta(\mathcal E_N^{-1})$. At smaller scales, Cauchy--Schwarz and
$\mathcal S_N\le C_\varphi\mathcal E_N$ give
\[
 \begin{aligned}
 \frac1{\mathcal E_N}\int_{\Omega_N^{(4)}\cap\{r<\delta\}}
 q(r)\varepsilon_I^2\,\dd\Lambda(I)
 &\le C_\varphi\sup_{r<\delta}q(r),\\
 \frac1{\mathcal E_N}\int_{\Omega_N^{(4)}\cap\{r<\delta\}}
 |w_IR_I|\,\dd\Lambda(I)
 &\le C_\varphi\sup_{r<\delta}q_1(r).
 \end{aligned}
\]
Let $N=N_k\to\infty$ and then $\delta\downarrow0$. The length term
vanishes by \eqref{eq:crit-weighted-cross-small}, and the remainder
vanishes by the second bound. Since $a_I\ge c_w$ and $\beta\ge0$,
\eqref{eq:crit-weighted-packet-expanded} proves the assertion.
\end{proof}

\begin{thm}[Singularity in $\Sy(\m S^1)\setminus\WP(\m S^1)$]\label{thm:crit-sym-total}
If $\varphi\in\Sy(\m S^1)\setminus\WP(\m S^1)$, there is a closed
proper arc $J$ such that
\[
 (\operatorname{Res}_J)_*\m P_{4,\varphi}\perp P_J.
\]
\end{thm}
\begin{proof}
Choose a finite smooth partition of unity subordinate to proper arcs,
and omit a fixed range of coarse scales so that each contributing
interval stays inside the corresponding arc. Since the omitted energy
is bounded and $\mathcal E_N\to\infty$, one member $\zeta$, supported
in the interior of an arc $J$, has the following property along a
subsequence of the indices in Lemma~\ref{lem:crit-statistics}. With
$K=\varphi(J)$, the measures
\[
 \nu_N(B)=\frac1{\mathcal E_N}\int
 \zeta(a)w_I^2\1_{\{\varphi(a)\in B\}}\,\dd\Lambda(I)
\]
converge weakly to a measure $\nu$ of mass $p>0$, supported in a
compact subset of $K^\circ$. All integrals below use these truncated
$\Omega_N^{(4)}$.

Suppose the two local laws on $J$ have a nonzero common part. Choose
\[
 0\ne\lambda\le P_J,\qquad
 \lambda\le(\operatorname{Res}_J)_*\m P_{4,\varphi}.
\]
Pulling $\lambda$ back under $u\mapsto u\circ\varphi|_J$ gives
$W P_K$, where $0\le W\le1$ and
$w=\m E_4[W]=\lambda(\mathcal X_J)>0$. Choose a countable dense set
$(x_i)$ in $K^\circ$ containing no atom of $\nu$, and set
\[
 W_m=\m E_4[W\mid H_4(x_1),\ldots,H_4(x_m)].
\]
Then $W_m\to W$ in $L^1$, $0\le W_m\le1$, and
$\m E_4[W_m]=w$, because evaluation on the dense set generates the
Borel sigma-field of the continuous restriction.

Fix $m$. Let $\eta_\delta$ be Lipschitz, zero within distance $\delta$
of $x_1,\ldots,x_m$, and one beyond distance $2\delta$. By
Lemma~\ref{lem:crit-packet-bv}, Cauchy--Schwarz, and
$\mathcal S_N\le C_\varphi\mathcal E_N$,
\[
 \sup_h|\mathcal A_N^\zeta(h)
       -\mathcal A_N^{\zeta\eta_\delta\circ\varphi}(h)|
 \le C_\varphi\left(\int(1-\eta_\delta)\,\dd\nu_N\right)^{1/2}.
\]
Proposition~\ref{prop:crit-weighted-cutoff} therefore gives
\[
 \liminf_N\m E_4[W_m\mathcal A_N^\zeta]
 \ge\frac{c_w}2\int\eta_\delta\,\dd\nu
       -C_\varphi\left(\int(1-\eta_\delta)\,\dd\nu\right)^{1/2}.
\]
The selected points are not atoms of $\nu$. Letting
$\delta\downarrow0$ makes the right side tend to $c_wp/2$.
This bound is independent of $m$, since $c_w$ depends only on $w$.
The uniform bound for $\mathcal A_N^\zeta$ in
Lemma~\ref{lem:crit-statistics} and $W_m\to W$ in $L^1$ then imply
\begin{equation}\label{eq:crit-weighted-packet-needed}
 \liminf_N\m E_4[W\mathcal A_N^\zeta]\ge c_wp/2>0.
\end{equation}

On the other hand, the definition of the pullback measure gives
\[
 \m E_4[W\mathcal A_N^\zeta]
 =\int\mathcal T_N^\zeta\,\dd\lambda
   -\int W\mathcal B_N^\zeta\,\dd P_K.
\]
Since $\lambda\le P_J$ and $W P_K\le P_K$, Cauchy--Schwarz and
\eqref{eq:crit-global-frame} yield
\[
 |\m E_4[W\mathcal A_N^\zeta]|
 \le\sqrt w\left(
 \|\mathcal T_N^\zeta\|_{L^2(P_J)}
 +\|\mathcal B_N^\zeta\|_{L^2(P_K)}\right)\longrightarrow0,
\]
contradicting \eqref{eq:crit-weighted-packet-needed}. Thus the local
laws are mutually singular.
\end{proof}

The three cases now prove the main theorem.

\begin{proof}[Proof of Theorem~\ref{thm:crit-all-nonwp-total}]
If $\varphi\notin\QS(\m S^1)$, the conclusion is
Theorem~\ref{thm:crit-nonqs-total}. If
$\varphi\in\QS(\m S^1)\setminus\Sy(\m S^1)$, choose a finite cover
by interiors of proper arcs. By Lemma~\ref{lem:crit-common-sym}, at
least one of the corresponding pairs of local laws has no nonzero
common part. For finite measures this is precisely mutual singularity.
If $\varphi\in\Sy(\m S^1)\setminus\WP(\m S^1)$,
Theorem~\ref{thm:crit-sym-total} gives the required arc.

Proposition~\ref{prop:crit-local-sing} transfers the local singularity
to right composition of the SLE$_4$ welding law. For left composition, apply the same local
conclusion to $\varphi^{-1}\notin\WP(\m S^1)$. By
Proposition~\ref{lem:crit-local-inversion}, the inverse SLE$_4$ welding law
has the same local measure classes as $\m P_4$, because the model law
is invariant under inversion. The right-composition result therefore
applies to this inverse law. Inverting the two laws turns right
composition by $\varphi^{-1}$ into left composition by $\varphi$,
which proves the second assertion.
\end{proof}

The local proof also gives an absolute-continuity criterion on a finite
cover and a singularity criterion for Stieltjes measures.

\begin{cor}[Consequences of local singularity]\label{thm:crit-local-converse}\label{cor:crit-sle-converse}\label{cor:crit-sing-stieltjes}
\begin{enumerate}[label=(\roman*)]
\item Suppose the interiors of the closed proper arcs $I_1,\ldots,I_M$
cover the circle. If
\begin{equation}\label{eq:crit-local-ac}
 (\operatorname{Res}_{I_j})_*\m P_{4,\varphi}\ll P_{I_j}
 \qquad(1\le j\le M),
\end{equation}
then $\varphi\in\WP(\m S^1)$.
\item If either $\m P_{4,\varphi}^{\mathrm{weld}}\ll\m P_4^{\mathrm{weld}}$
or $\m P_4^{\mathrm{weld}}\ll\m P_{4,\varphi}^{\mathrm{weld}}$, then
$\varphi\in\WP(\m S^1)$. The same conclusion holds on the left-rotation
quotient.
\item If $\dd\varphi|_I\perp\dd x|_I$ for some closed proper arc $I$, or the
same condition holds for $\varphi^{-1}$, then right and left composition
by $\varphi$ make the SLE$_4$ welding law mutually singular with its
original law. Here $\dd\varphi(A)=|\varphi(A)|$.
\end{enumerate}
\end{cor}
\begin{proof}
For (i), if $\varphi$ were not quasisymmetric, the intervals in the proof of
Theorem~\ref{thm:crit-nonqs-total} would eventually lie in one $I_j$
whose interior contains their limiting base point. The separating event
on that arc contradicts \eqref{eq:crit-local-ac}. Thus $\varphi$ is
quasisymmetric. Absolute continuity gives a nonzero common part on each
$I_j$, so Lemma~\ref{lem:crit-common-sym} gives symmetry. If $\varphi$
were not WP, use a partition of unity subordinate to these same $I_j$
in the proof of Theorem~\ref{thm:crit-sym-total}. That proof gives
singularity on one of them, again a contradiction.

For (ii), two probability measures that are mutually singular cannot be absolutely
continuous in either direction. Apply Theorem~\ref{thm:crit-all-nonwp-total}.
For the quotient, use \eqref{eq:crit-haar-lift}; reverse absolute
continuity reduces to the forward direction by $R_{\varphi^{-1}}$.

For (iii), a WP homeomorphism is absolutely continuous with positive derivative
almost everywhere; see \cite[Definition~2.3]{F1}. Its Stieltjes measure
cannot be singular on a nondegenerate arc. The same applies to its inverse,
since $\WP(\m S^1)$ is a group. Apply Theorem~\ref{thm:crit-all-nonwp-total}.
\end{proof}

\begin{remark}[The subcritical parameters]\label{rem:subcritical-shape-proof}
The local singularity argument also applies to each fixed
$0<\kappa<4$. Let $\gamma=\sqrt\kappa$, and let $H$, $\m P$ and
$\mu_1,\mu_2$ be the subcritical GMC welding, its law and its two
normalized measures from Lemma~\ref{introlem:sle-gmc-lift}. Then, for
every $\varphi\notin\WP(\m S^1)$, there is a closed proper arc $J$
such that
\[
 (\operatorname{Res}_J)_*(R_\varphi)_*\m P
 \perp(\operatorname{Res}_J)_*\m P.
\]
Here are the changes needed in the proof.

For smooth source deformations, the coordinate rule in
\cite[Lemma~3.8]{F1} replaces the transformed Gaussian law by
\[
 \operatorname{Law}\bigl(\pi_0(X\circ\phi_t+
 Q_\gamma\log\phi_t')\bigr),\qquad
 Q_\gamma=2/\gamma+\gamma/2.
\]
The covariance operator $S_t$ is unchanged, and the mean is
$b_t^\gamma=(Q_\gamma/2)b_t$, with $b_t$ as in
Lemma~\ref{lem:crit-gaussian-flow}. Thus its density derivative bounds
and the entropy estimates in Lemma~\ref{lem:crit-small-affine} hold
with constants depending on $\gamma$.

For the source estimate, divide $I$ into four equal subarcs and put
$m_j=\mu_1(I_j)$ and
$Z_I=\max_{1\le j\le3}|\log(m_{j+1}/m_j)|$. Then
\[
 \min_j\frac{m_j}{\sum_i m_i}\ge\frac14e^{-3Z_I},\qquad
 \sup_I E[Z_I^p]<\infty\quad(1\le p<\infty)
\]
by Lemma~\ref{lem:log GMC moment}. Markov's inequality gives
\eqref{eq:crit-quarter-input} for $\mu_1$.

For the target estimate, let $Q=F_{X_2}^{-1}$ and let $V$ be Haar and
independent of the fields. The point $Y=Q(V)$ has conditional law
$\mu_2$, and the two arc lengths are
\[
 B_+(d)=Q(V+d)-Q(V),\qquad B_-(d)=Q(V)-Q(V-d).
\]
For $M\ge0$ and $0<q<e^{-M}/2$, set
\[
 K_{q,M}=\sup_{\substack{\sigma,\tau\in\{-,+\}\\|u|,|v|\le M}}
 \left|\log\frac{B_\sigma(qe^u)}{B_\tau(qe^v)}\right|,
 \qquad
 L_\pm=\log\frac{B_\pm(qe^M)}{B_\mp(qe^{-M})}.
\]
Monotonicity gives
\[
 K_{q,M}=\log\frac{\max_\sigma B_\sigma(qe^M)}
                         {\min_\sigma B_\sigma(qe^{-M})}
 \le |L_+|+|L_-|.
\]
Indeed, each same-side ratio is bounded by $e^{L_++L_-}$, and the
two opposite-side ratios are controlled by $L_+$ and $L_-$.
Lemma~\ref{lem:inverse rotation} applies since
$qe^M+qe^{-M}<1$, and yields
\[
 \sup_{0<q<e^{-M}/2}E[K_{q,M}^p]
 \le C_{\gamma,p}(1+2M)^p.
\]
This proves \eqref{eq:crit-root-input} for the subcritical target.
Conditional on the source field, $F_{X_1}(x)+U$ is still Haar and
independent of the target field. The estimate therefore applies at
the random source mass $q=\mu_1(I)$ in the proof of
Lemma~\ref{lem:crit-balanced-tight}.

The model has the same rotation, reflection and inversion invariances
as $H_4$. The weighted mean estimates and the proofs in
Sections~\ref{subsec:crit-local-regularity}--\ref{subsec:crit-singularity}
now apply with $H_4,\m P_4$ replaced by $H,\m P$. The reconstruction
and energy comparison are deterministic. Finally,
Lemma~\ref{introlem:sle-gmc-lift}, based on \cite[Lemma~4.8]{F1} and
the Haar lift, transfers singularity to right composition of the
SLE$_\kappa$ welding law; inversion gives left composition.
The transfer below criticality uses the subcritical mass estimates
in place of the critical Bessel construction. The proof in
Sections~\ref{sec:subcritical-estimates}--\ref{sec:subcritical-singularity}
is retained as a separate argument. Constants may depend on the fixed
$\gamma$; no uniform estimate as $\gamma\uparrow2$ is claimed.
\end{remark}

\section{The WP normalizer and two-sided singularity}\label{sec:normalizer}

The one-sided dichotomy identifies WP as the set of non-singular
composition directions. To pass to simultaneous left and right
composition, we prove that the full WP group is its own normalizer in
$\Homeo_+(\m S^1)$. The proof first establishes quasisymmetry of a
normalizing homeomorphism and then WP regularity through its composition
operator. Sections~\ref{subsec:two-sided-quasi-invariance}
and~\ref{subsec:two-sided-singularity} apply this result to the welding law.

Throughout this section we use angular coordinates
$\m S^1=\m R/(2\pi\m Z)$, with $\theta\mapsto e^{\ii\theta}$,
and write $m(\dd\theta)=\dd\theta/(2\pi)$ for normalized Haar measure.
The rotation parameter has the same normalization, so all averages over
$\W$ are taken against $m(\dd t)$. This convention is local to the
present section. It differs from the coordinate $x\in\m R/\m Z$ used
in the preceding welding estimates only by $\theta=2\pi x$; the
normalizer argument is independent of those coordinate choices.

\subsection{The normalizer theorem}\label{subsec:normalizer-theorem}

We work with the full circle group, with no prescribed fixed point.
The distinction from a fixed-point subgroup is explained in
Remark~\ref{introrem:fixed-point-normalizer}.

\begin{df}
    The normalizer of $\WP(\m S^1)$ in $\Homeo_+(\m S^1)$ is
    \begin{equation*}
        G:= \mathcal{N}(\WP(\m S^1)):=\{ \varphi \in \Homeo_+(\m S^1) \vert \varphi\circ\WP(\m S^1) \circ \varphi^{-1}=\WP(\m S^1) \}.
    \end{equation*}
    For a subgroup $W \subset \Homeo_+(\m S^1)$, for example $W=\QS(\m S^1)$, write
    \begin{equation*}
        G(W):= \mathcal{N}(\WP(\m S^1),W):=\{ \varphi \in W \vert \varphi\circ\WP(\m S^1) \circ \varphi^{-1}=\WP(\m S^1) \}.
    \end{equation*}
    
\end{df}

\begin{remark}
    The group property gives $\WP(\m S^1) \subset G(\QS(\m S^1))\subset G$. 
\end{remark}
\begin{thm}\label{thm:main}
    $ G = \WP(\m S^1)$.
\end{thm}
\begin{proof}
Theorems~\ref{thm:qs regularity} and~\ref{thm:from qs to wp} give the reverse inclusions.
\end{proof}
\subsection{Conjugated rotations and an invariant measure}\label{subsec:normalizer-rotations}

Conjugating rotations by a normalizing homeomorphism produces WP maps.
We first make their quasisymmetry bounds uniform in the rotation
parameter, then average their pullbacks of Lebesgue measure. The resulting
invariant measure will identify the regularity of the conjugating map.

Fix $\varphi\in G$.
For $t \in \m R / 2\pi \m Z$, let $R_t(z) = e^{\ii t}z$ and define
\[
    h_t = \varphi^{-1} \circ R_t \circ \varphi.
\]
Since $R_t \in \WP(\m S^1)$, we have 
\begin{equation}\label{eq:assumption}
    h_t \in \WP(\m S^1)  \qquad \text{for every}~t. 
\end{equation}
For $M \geq 1$, let $\mc Q_M \subset \Homeo_+(\m S^1)$ consist of those $h$ for which
\[
    M^{-1} \leq  \frac{m(h(I_1))}{m(h(I_2))}\leq M,
\]
whenever $I_1$, $I_2$ are adjacent arcs of equal length at most $\pi$.
Define
\[
    E_M = \{t|h_t \in \mc Q_M\}.
\]
\begin{lem}\label{lem:qc M}
    \begin{itemize}
        \item Every $h \in \QS(\m S^1)$ belongs to $\mc Q_M$ for some integer $M \geq 1$.
        \item For every $M \geq 1$, there exists $\Phi(M) \geq 1$ such that
                \[
                    h_1,h_2 \in \mc Q_M \implies h_1^{-1} \circ h_2 \in \mc Q_{\Phi(M)}.
                \]
        \item For fixed $M$, $\mc Q_M$ is closed in the relative uniform topology of $\Homeo_+(\m S^1)$.
    \end{itemize}
\end{lem}
\begin{proof}
The first assertion is the classical adjacent-interval characterization
of quasisymmetric boundary maps; see \cite[Chapter~II]{lehto1973quasiconformal}.
The quantitative form of this characterization supplies a distortion
function $\eta_M$ depending only on $M$. The inverse of an
$\eta$-quasisymmetry has distortion function
$t\mapsto1/\eta^{-1}(1/t)$, and distortion functions compose under
composition. Applied to $h_1^{-1}\circ h_2$, these facts give a bound
$\Phi(M)$ depending only on $M$; these are the standard inversion and
composition properties in \cite{lehto1973quasiconformal}.
Finally, if $h_j\to h$ uniformly and $h$ is a homeomorphism, the image
lengths of each fixed pair of adjacent arcs converge to those under $h$.
Their denominators remain positive in the limit, so the defining
inequalities for $\mc Q_M$ pass to $h$.
\end{proof}

\begin{lem}\label{lem:continuity}
    The map
    \[
        t \mapsto h_t
    \]
    is continuous from the parameter circle to $\Homeo_+(\m S^1)$ equipped
    with the uniform group metric $d_H$ defined in
    Section~\ref{introsec:subcritical}.
\end{lem}
\begin{proof}
We give a direct verification in the uniform topology used for boundary
correspondences in \cite{lehto1973quasiconformal}. Let $d$ be the circle
distance used in $d_H$, and let $\omega$ be the modulus of continuity of
$\varphi^{-1}$ with respect to $d$. Since rotations are isometries,
\[
 \sup_z d(h_t(z),h_s(z))
 \le\omega\bigl(d(e^{\ii t},e^{\ii s})\bigr).
\]
Since $h_t^{-1}=h_{-t}$, the same bound holds for the inverses. Hence
\[
 d_H(h_t,h_s)\le2\omega\bigl(d(e^{\ii t},e^{\ii s})\bigr)
 \longrightarrow0\qquad(t\to s).
\]
\end{proof}

\begin{lem}
    There exists $M_1<\infty$ such that 
    \[
        h_t \in \mc Q_{M_1}  \qquad \text{for every}~t.
    \]
\end{lem}
\begin{proof}
    By Lemma~\ref{lem:qc M} and~\ref{lem:continuity}, each $E_M$ is closed. From assumption~\eqref{eq:assumption}, we have
    \[
        \m R / 2\pi\m Z = \bigcup_{M=1}^\infty E_M.
    \]
    Then from Baire-category theorem, some $E_M$ has nonempty interior. Choose a nonempty open arc $I_0 \subset E_M$. The difference set $I_1 :=I_0 -I_0$ contains an open neighborhood $U$ of $0$. If $u=b-a$ with $a$, $b \in I_0$, then the group relation $h_{s+t} = h_s \circ h_t$ gives
    \[
        h_u = h_{-a} \circ h_b = h_a^{-1} \circ h_b.
    \]
    By Lemma~\ref{lem:qc M}, there is $M_0 = \Phi(M)$, independent of $a$, $b$, such that
    \[
        h_t \in \mc Q_{M_0}  \qquad \text{for every}~t \in U.
    \]
    Choose $\d>0$ with $(-\d,\d)\subset U$ and an integer $N$ with
    $\pi/N<\d$. For each parameter, take its representative
    $t\in[-\pi,\pi]$ and put $\tilde t=t/N\in U$. Then
    \[
        h_t=(h_{\tilde t})^N.
    \]
    The distortion function of an $N$-fold composition of $M_0$-controlled quasisymmetries depends only on $M_0$ and $N$. The proof is finished.
\end{proof}
For each $t$, define the probability measure
\[
    \mu_t = (h_t^{-1})_* m, \qquad \mu_t(A) = m(h_t(A)).
\]
\begin{lem}
    The family $t \mapsto \mu_t$ is weak-star continuous. Consequently there exists a Borel probability measure $\nu$ satisfying
    \begin{equation}\label{eq:nu}
        \int_{\m S^1} F \dd \nu = \int_{\m R / 2\pi \m Z} \left(\int_{\m S^1} F \dd \mu_t\right) \dd m(t),
    \end{equation}
    for every $F \in C(\m S^1)$.
\end{lem}
\begin{proof}
    For $F \in C(\m S^1)$,
    \[
        \int_{\m S^1} F \dd \mu_t = \int_{\m S^1} F\circ h_t^{-1} \dd m.
    \]
    By Lemma~\ref{lem:continuity}, $t \mapsto h_t^{-1}$ is continuous. Hence the above integrand is jointly continuous and bounded, so the right-hand side depends continuously on $t$. The functional on the right-hand side of \eqref{eq:nu} is linear, bounded and takes the value $1$ at the constant function $1$. The Riesz representation theorem gives a unique Borel probability measure $\nu$.
\end{proof}
\begin{lem}\label{lem:nu property}
     The measure $\nu$ satisfies
     \begin{itemize}
         \item $M_1^{-1}\nu(I_2)\le\nu(I_1)\le M_1\nu(I_2)$ for all adjacent equal-length arcs $I_1$, $I_2$.
         \item $(h_s)_* \nu =\nu$ for every $s$.
     \end{itemize}
\end{lem}
\begin{proof}
    For the first one, we have
    \begin{align*}
        \nu(I_1) = &\int_{\m S^1} \mathbf{1}_{I_1} \dd \nu \\
        =& \int_{\m R / 2\pi \m Z} \left(\int_{\m S^1} \mathbf{1}_{I_1} \dd \mu_t\right) \dd m(t) \\
        =&\int_{\m R / 2\pi \m Z}  \mu_t(I_1) \dd m(t) \\
        \leq&  M_1\int_{\m R / 2\pi \m Z}  \mu_t(I_2) \dd m(t) \\
        =&  M_1\nu(I_2).
    \end{align*}
    Interchanging $I_1$ and $I_2$ gives the reverse bound.
    For the second one, for every Borel set $A$ and every $s$,
    \begin{align*}
        \nu(h_s(A)) =&  \int_{\m R / 2\pi \m Z}  \mu_t(h_s(A)) \dd m(t) \\
        =& \int_{\m R / 2\pi \m Z}  m(h_t(h_s(A))) \dd m(t) \\
        =&\int_{\m R / 2\pi \m Z}  m(h_u(A)) \dd m(u) \\
        =& \nu(A),
    \end{align*}
    where the third equality is Haar invariance of the parameter measure. Replacing $A$ by $h_s^{-1}(A)$ finishes the proof.
    
\end{proof}
\subsection{Quasisymmetry of a normalizing homeomorphism}\label{subsec:normalizer-qs}

The averaged measure is invariant under the conjugated rotations and has
uniformly comparable masses on adjacent arcs. Pushing it forward by
$\varphi$ gives Haar measure; the adjacent-arc bounds therefore become
quasisymmetry bounds for $\varphi$.

\begin{thm}\label{thm:qs regularity}
    $G \subset \QS(\m S^1)$.
\end{thm}
\begin{proof}
    For $\varphi \in G$, define $\l = \varphi_* \nu$. Since $\varphi \circ h_s \circ \varphi^{-1} = R_s$, $(h_s)_* \nu =\nu$ implies 
    \[
        (R_s)_* \l =  (\varphi \circ h_s \circ \varphi^{-1})_*(\varphi_* \nu)=\varphi_* (h_s)_* \nu = \varphi_*  \nu=\l,
    \]
    for every $s$. By uniqueness of normalized Haar measure, $\l = m=\varphi_* \nu$, that is,
    \[
        \nu(A) = m(\varphi(A)) \qquad  \text{for every Borel set}~A.
    \]
    Therefore, from Lemma~\ref{lem:nu property}, we have $\varphi \in \mc Q_{M_1} \subset \QS(\m S^1)$.
\end{proof}
\subsection{WP regularity from the composition operator}\label{subsec:normalizer-wp}

For a quasisymmetric map, WP regularity is characterized by a
Hilbert--Schmidt condition on its composition operator. Normalization
of the WP group imposes this condition modulo operators invariant under
rotations. The operator lemmas below remove the remaining diagonal term.

\begin{thm}\label{thm:from qs to wp}
    $G(\QS(\m S^1))\subset \WP(\m S^1)$.
\end{thm}
Let
\[
    \mc H = H^{1/2}_\m C(\m S^1) = \mc H_+ \oplus\mc H_-,
\]
where
\[
    \mc H_+ = \overline{\text{span}}\{e_n=\frac{z^n}{\sqrt{\abs{n}}}: n \geq 1\}, \qquad \mc H_- = \overline{\text{span}}\{e_{-n}=\frac{z^{-n}}{\sqrt{\abs{n}}}: n \geq 1\}.
\]
For $h \in \QS(\m S^1)$, let $C_h$ denote the normalized composition operator on $\mc H$ and let $\L_h = (C_hC_h^*-I)/2$ denote its welding Grunsky operator. We have
\[
    \L_{h \circ \varphi}=C_\varphi\L_h C_\varphi^*+\L_\varphi
\]
for $h$, $\varphi \in \QS(\m S^1)$.
We use the composition-operator formulation of \cite{F3}. The next two lemmas express WP regularity and the normalizer condition in these coordinates.
\begin{lem}\label{lem:df of WP}
    For $\varphi\in\QS(\m S^1)$, the following are equivalent:
    \begin{itemize}
        \item $\varphi \in \WP(\m S^1)$;
        \item $\L_\varphi \in \HS$;
        \item $\L_\varphi-D\in\HS$ for some bounded operator $D$ diagonal in the
        Fourier basis $(e_n)_{n\ne0}$.
    \end{itemize}
\end{lem}
\begin{proof}
The first two conditions are equivalent by \cite{F3}; the second
implies the third by taking $D=0$. For the remaining implication, write
\[
 C_\varphi=\begin{pmatrix}a&b\\\bar b&\bar a\end{pmatrix}
 \quad\text{on }\mc H_+\oplus\mc H_-.
\]
The symplectic identities are
\[
 aa^*-bb^*=I,\qquad a^*a-b^T\bar b=I,\qquad ab^T=ba^T.
\]
The first makes $a$ surjective and the second makes it bounded below,
so $a$ is invertible. Consequently,
\[
 (\L_\varphi)_{+-}=\tfrac12(ab^T+ba^T)=ab^T.
\]
A Fourier-diagonal $D$ has zero $+-$ block. Thus
$\L_\varphi-D\in\HS$ implies $ab^T\in\HS$, and multiplication by
$a^{-1}$ gives $b\in\HS$. The composition-operator WP criterion
\cite{TT06,F3} now yields $\varphi\in\WP(\m S^1)$.
\end{proof}
\begin{lem}
If $\varphi\in G(\QS(\m S^1))$, then $\varphi\in\QS(\m S^1)$ and
\[
 C_\psi\L_\varphi C_\psi^*-\L_\varphi\in\HS
 \qquad\text{for every }\psi\in\WP(\m S^1).
\]
In particular, this holds for every rotation $\psi$.
\end{lem}\begin{proof}
By the normalizer condition, for every $\psi\in\WP(\m S^1)$ there is
$h=h(\psi)\in\WP(\m S^1)$ such that
    \[
        h \circ \varphi =\varphi \circ \psi
    \]
    and hence
    \[
        C_\varphi\L_h C_\varphi^*+\L_\varphi=\L_{h \circ \varphi} = \L_{\varphi \circ \psi}=C_\psi \L_\varphi C_\psi^*+\L_\psi,
    \]
    therefore,
    \[
        C_\psi \L_\varphi C_\psi^* - \L_\varphi =C_\varphi\L_h C_\varphi^*-\L_\psi\in\HS.
    \]
\end{proof}
Let $\rho_t = C_{R_t^{-1}}$ for $t \in \W$.
\begin{lem}\label{lem:representation continuous}
The representation $\rho:\W\to\U$ is strongly continuous. For every
fixed bounded linear operator $T:\mc H\to\mc H$, the map
\[
 \W\ni t\longmapsto\alpha_t(T)=\rho(t)T\rho(t)^*
\]
is continuous in the strong operator topology.
\end{lem}
\begin{proof}
For a trigonometric polynomial $p$ with zero constant term and frequencies
$0<|j|\le m$, we have
\[
 \norm{\rho(t_n)p-\rho(t)p}_{\mc H}
 \le\max_{0<|j|\le m}|e^{\ii jt_n}-e^{\ii jt}|
       \norm{p}_{\mc H}\longrightarrow0\qquad(t_n\to t).
\]
Such polynomials are dense in $\mc H$. Since the operators $\rho(t)$
are unitary, for $u\in\mc H$,
\[
 \norm{\rho(t_n)u-\rho(t)u}_{\mc H}
 \le2\norm{u-p}_{\mc H}+\norm{\rho(t_n)p-\rho(t)p}_{\mc H}.
\]
First let $n\to\infty$, then let $p\to u$. This proves strong continuity.
Also $\rho(t_n)^*=\rho(-t_n)\to\rho(-t)=\rho(t)^*$ strongly.
For each fixed $u\in\mc H$, boundedness of $T$ gives
\begin{align*}
 \norm{(\alpha_{t_n}(T)-\alpha_t(T))u}_{\mc H}
 &\le\norm{T}\,\norm{(\rho(t_n)^*-\rho(t)^*)u}_{\mc H}\\
 &\quad+\norm{(\rho(t_n)-\rho(t))T\rho(t)^*u}_{\mc H}
 \longrightarrow0.
\end{align*}
Both terms tend to zero by strong continuity, since $u$ and
$T\rho(t)^*u$ are fixed vectors.
\end{proof}
For the remainder of the proof, set
\[
    T = \L_\varphi, \qquad \alpha_t(A) = \rho(t) A \rho(t)^*, \qquad b(t) = \alpha_t(T)-T.
\]
Hence $b(t) \in \HS$ for every $t \in \W$. Moreover,
\[
    b(s+t)=b(s)+\alpha_s(b(t)).
\]
Hence
\[
    \norm{b(s+t)}_2 \leq \norm{b(s)}_2 + \norm{b(t)}_2,
\]
as $\alpha_s$ preserves the $\HS$-norm.
\begin{lem}\label{lem:lower-semicontinuity of HS}
    Suppose $A_j \in \HS$, $A_j \to A$ strongly, and $\sup_j \norm{A_j}_2 < \infty$. Then $A \in \HS$ and 
    \[
        \norm{A}_2 \leq \liminf_{j \to \infty} \norm{A_j}_2.
    \]
\end{lem}
\begin{proof}
    For every positive integer $N$, strong convergence gives $A_j e_n \to A e_n$ for each $|n|\leq N$, $n \neq 0$. Therefore,
    \[
        \sum_{0<|n|\le N}\norm{Ae_n}_{\mc H}^2
        =\lim_{j\to\infty}\sum_{0<|n|\le N}\norm{A_j e_n}_{\mc H}^2
        \le\liminf_{j\to\infty}\norm{A_j}_2^2.
    \]
    Letting $N \to \infty$ we complete the proof.
\end{proof}
\begin{lem}\label{lem:cocycle bounded}
    The cocycle $b$ is bounded in $\HS$-norm on $\W$.
\end{lem}
\begin{proof}
    Define
    \[
        F_M = \{t \in \W|\norm{b(t)}_2 \leq M\}.
    \]
    By Lemma~\ref{lem:representation continuous}, $t_j \to t$ implies $\alpha_{t_j}(T) \to \alpha_t(T)$ strongly and hence $b(t_j) \to b(t)$ strongly. Lemma~\ref{lem:lower-semicontinuity of HS} then shows that each $F_M$ is closed. Since $\W=\bigcup_{M=1}^\infty F_M$ is a locally compact Lie group hence is a Baire space, some $F_M$ contains a nonempty open set. Choose a nonempty open arc $I_0 \subset F_M$. The difference set $I_1 :=I_0 -I_0$ contains an open neighborhood $U$ of $0$. There is $M_0$ such that
    \[
        \norm{b(t)}_2 \leq M_0  \qquad \text{for every}~t \in U.
    \]
    Choose $\d>0$ with $(-\d,\d)\subset U$ and an integer $N$ with
    $\pi/N<\d$. For each parameter, take its representative
    $t\in[-\pi,\pi]$ and put $\tilde t=t/N\in U$. Then
    \[
        \norm{b(t)}_2\leq N\norm{b(\tilde t)}_2 \leq NM_0 .
    \]
   The proof is finished.
\end{proof}
\begin{lem}
    The map $b:\W \to \HS$ is strongly Borel measurable. Furthermore, the Bochner integral
     \[
        q = \int_{\W} b(t) \dd m(t)
    \]
    is well-defined and belongs to $\HS$.
\end{lem}
\begin{proof}
    Firstly, let $F$ be a finite-rank operator. The scalar function
    \[
        t \mapsto \brac{b(t),F}_2 = \tr(F^*b(t))
    \]
    is a finite linear combination of matrix coefficients $\brac{b(t)u,v}$, hence is continuous by Lemma~\ref{lem:representation continuous}. Let $S \in \HS$, and choose finite-rank operators $F_j \to S$ in Hilbert--Schmidt norm. By Lemma~\ref{lem:cocycle bounded}, with $M_0=\sup_{t \in \W} \norm{b(t)}_2$,
    \[
        \sup_{t \in \W} |\brac{b(t),S-F_j}| \leq M_0\norm{S-F_j}_2 \to 0.
    \]
    Thus $t \mapsto \brac{b(t),S}_2$ is a uniform limit of continuous functions and is therefore continuous. Hence $b$ is weakly continuous as an $\HS$-valued map. The Hilbert space $\HS$ is separable. Pettis' measurability theorem therefore implies that weak measurability is equivalent to strong Borel measurability. By Lemma~\ref{lem:cocycle bounded},
    \[
        \int_{t \in \W} \norm{b(t)}_2\,\dd m(t) \leq \sup_{t \in \W} \norm{b(t)}_2 < \infty,
    \]
    which is exactly the Bochner integrability criterion.
\end{proof}
\begin{lem}
    The operator $D=T+q$ is invariant under rotation conjugation:
    \[
        \a_t(D) = D, \qquad t \in \W.
    \]
    Furthermore, there exists a bounded sequence $(d_n)_{n \neq 0}$ such that
    \[
        De_n = d_n e_n, \qquad n \neq 0.
    \]
\end{lem}
\begin{remark}
    Formally,
    \[
        D = \int_{\W} \a_t(T) \dd m(t).
    \]
\end{remark}
\begin{proof}
    Since $\a_t : \HS \to \HS$ is bounded and linear, it commutes with the Bochner integral. We have
    \begin{align*}
        \a_t(q) =& \int_{\W} \a_t(b(s))\dd m(s) \\
        =& \int_{\W} b(s+t)\dd m(s) - \int_{\W} b(t)\dd m(s)\\
        =&q-b(t),
    \end{align*}
    where the second equality uses cocycle identity and the last equality uses invariance of Haar measure. Therefore,
    \[
        \a_t(D)=\a_t(T)+\a_t(q) = T+b(t)+q-b(t)=D.
    \]
    It follows that
    \[
        (e^{-\ii m \t} - \e^{-\ii n \t})\brac{De_n,e_m} = 0,\qquad \text{for every}~\t, n, m.
    \]
    So $D$ is diagonal and is also bounded as $D=T+q$ is bounded.
\end{proof}
We summarize the above argument as a theorem.
\begin{thm}\label{thm:generalized Schur lemma}
    Let $T $ be a bounded linear operator such that
    \[
        b(t) \in \HS, \qquad \text{for every}~t \in \W.
    \]
    Then the $\HS$-valued Bochner integral
    \[
        q = \int_{\W} b(t) \dd m(t)
    \]
    is well-defined and belongs to $\HS$. The operator
    \[
        D =T +q
    \]
    is diagonal in the basis $\{e_n\}_{n \neq 0}$.
\end{thm}
\begin{proof}[Proof of Theorem~\ref{thm:from qs to wp}]
    Let $T = \L_\varphi$, then from Theorem~\ref{thm:generalized Schur lemma}, we have there exist a diagonal matrix $D$ and a Hilbert--Schmidt matrix $q$ such that $T = D-q$, hence from Lemma~\ref{lem:df of WP} we have $\varphi \in \WP(\m S^1)$.
\end{proof}

\subsection{Two-sided translates of the WP group}\label{subsec:normalizer-translates}
The normalizer theorem also determines which pairs of homeomorphisms
preserve WP under left and right composition. The reduction is algebraic:
apply the proposed identity to the identity map, then remove the resulting
WP factor.
\begin{thm}\label{thm:wp-two-sided-translates}
Let $\varphi_1,\varphi_2\in\Homeo_+(\m S^1)$. Then
\[
 \varphi_1\circ\WP(\m S^1)\circ\varphi_2^{-1}=\WP(\m S^1)
 \quad\Longleftrightarrow\quad
 \varphi_1,\varphi_2\in\WP(\m S^1).
\]
\end{thm}
\begin{proof}
Suppose the equality holds. Applying it to the identity map gives
$\alpha=\varphi_1\circ\varphi_2^{-1}\in\WP(\m S^1)$. Hence
\[
 \varphi_1\circ\WP(\m S^1)\circ\varphi_1^{-1}
 =\WP(\m S^1)\circ\alpha^{-1}=\WP(\m S^1).
\]
Theorem~\ref{thm:main} yields $\varphi_1\in\WP(\m S^1)$, and then
$\varphi_2=\alpha^{-1}\circ\varphi_1\in\WP(\m S^1)$.
The converse follows from the group property.
\end{proof}

\subsection{Two-sided quasi-invariance}\label{subsec:two-sided-quasi-invariance}
The preceding theorem identifies the normalizer of the full WP group.
We now fix $\k\in(0,4]$ and apply the normalizer theorem to the two
composition actions. The one-sided results in
Section~\ref{sec:critical-local} and
Remark~\ref{rem:subcritical-shape-proof} identify $\QIR$ with WP;
the following definitions record its left and simultaneous analogues.
Define
\[
    \QI = \{ (\varphi,\psi) \in \Homeo_+(\m S^1) \times \Homeo_+(\m S^1)| \varphi \circ \SLE_\k^\text{weld}\circ \psi^{-1}  \sim \SLE_\k^\text{weld}\}.
\]
Define the quasi-invariance group:
\[
    \QIL = \{ \varphi \in \Homeo_+(\m S^1)| \varphi \circ \SLE_\k^\text{weld}  \sim \SLE_\k^\text{weld}\},
\]
\[
    \QIR = \{ \varphi \in \Homeo_+(\m S^1)|  \SLE_\k^\text{weld} \circ \varphi^{-1}  \sim \SLE_\k^\text{weld}\},
\]
\[
    \QIC = \{ \varphi \in \Homeo_+(\m S^1)| \varphi \circ \SLE_\k^\text{weld}\circ \varphi^{-1}  \sim \SLE_\k^\text{weld}\},
\]
\begin{lem}\label{lem:containing relation}
    \[
        \WP(\m S^1) \subset \QIL=\QIR \subset \QIC \subset   \mc N(\QIL) .
    \]
    \[
        \WP(\m S^1) \times \WP(\m S^1) \subset \QI.
    \]
\end{lem}
\begin{proof}
    The WP inclusions follow from Theorems~\ref{body:thm:one}
and~\ref{body:thm:two}; for $\k<4$ they also follow from
\cite[Corollary~4.12]{F1}, with the normalization of
Definition~\ref{introdef:sle-haar-lift}. The first identity follows from the inverse-invariance of  $\SLE_\k^\text{weld}$. The second containing relation comes from the expression
    \[
        \varphi \circ \SLE_\k^\text{weld}\circ \varphi^{-1} =( \varphi \circ \SLE_\k^\text{weld})\circ \varphi^{-1}.
    \]
    The third containing relation comes from the expression,  for $\varphi \in \QIC$ and $\psi \in \QIL$,
    \[
        \varphi \circ  \psi \circ \varphi^{-1} \circ \SLE_\k^\text{weld}= \varphi\circ (\psi \circ  ( \varphi^{-1} \circ \SLE_\k^\text{weld} \circ \varphi) ) \circ \varphi^{-1}.
    \]
    The last containing relation comes from the expression
    \[
        \varphi \circ \SLE_\k^\text{weld}\circ \psi =( \varphi \circ \SLE_\k^\text{weld})\circ \psi.
    \]
\end{proof}
\begin{cor}\label{cor:two-sided-quasi-invariance}
    We have $\QI \subset \mc N(\QIR)\times \mc N(\QIR)$. Furthermore, if $\QIR = \WP(\m S^1)$, then $\QI = \WP(\m S^1) \times \WP(\m S^1)$.
\end{cor}
\begin{proof}
    Note that $ (\varphi,\psi) \in \QI$ implies that, for each $\d \in \QIL$,
    \[
        \varphi \circ \d \circ \varphi^{-1} \circ  \SLE_\k^\text{weld}= \varphi\circ (\delta \circ  ( \varphi^{-1} \circ \SLE_\k^\text{weld} \circ \psi) ) \circ \psi^{-1} \sim  \SLE_\k^\text{weld},
    \]
    it follows that $\varphi \in \mc N(\QIL)= \mc N(\QIR)$, similarly, $\psi \in \mc N(\QIR)$.  By Theorem~\ref{thm:main}, $\mc N(\WP(\m S^1)) =\WP(\m S^1)$.
\end{proof}
\subsection{From one-sided to two-sided singularity}\label{subsec:two-sided-singularity}

Fix $\kappa\in(0,4]$ and write $P=\m P^{\mathrm{weld}}$. We use the
non-singularity classes $\tilde\QIL,\tilde\QIR,\tilde\QI$ defined in
Section~\ref{introsec:two-sided}. Theorems~\ref{thm:sing-conditional}
and~\ref{thm:crit-all-nonwp-total}, together with WP quasi-invariance
and inversion invariance, give
\[
 \tilde\QIL=\QIL=\tilde\QIR=\QIR=\WP(\m S^1).
\]
In particular, for a one-sided change, a nonzero common part implies
mutual absolute continuity.

\begin{thm}[Two-sided singularity]\label{thm:two-sided-singularity}
We have
\[
 \tilde\QI\subset\mc N(\QIR)\times\mc N(\QIR)
 =\WP(\m S^1)\times\WP(\m S^1).
\]
The first inclusion is also an equality.
\end{thm}
\begin{proof}
Let $(\varphi,\psi)\in\tilde\QI$, put $T(h)=\varphi\circ h\circ\psi^{-1}$
and $Q=T_*P$, and take a nonzero common part $\nu\le P,Q$. For
$a\in H:=\QIR$, let $S_a(h)=h\circ(\psi a\psi^{-1})^{-1}$. This is the
conjugate by $T$ of right composition by $a^{-1}$, so $(S_a)_*Q\sim Q$.
Set $K=H\cap\psi^{-1}H\psi$ and let $H/K$ denote the set of left cosets
of $K$ in $H$. If $a,b\in H$ and $aK\ne bK$, then $a^{-1}b\in H\setminus K$,
so $\psi a^{-1}b\psi^{-1}\notin H$. Since $S_a^{-1}S_b=S_{a^{-1}b}$,
the identity $\tilde\QIR=H$ gives $(S_a)_*P\perp(S_b)_*P$.
The measures $(S_a)_*\nu$ indexed by coset
representatives are therefore pairwise singular, and each is absolutely
continuous with respect to $Q$. Their densities have pairwise disjoint
positive sets up to $Q$-null sets, so $H/K$ is at most countable.

The subgroup $K$ is Borel in $H$ with its relative uniform topology, since
\[
 K=\{a\in H:\|(S_a)_*P-P\|_{\mathrm{TV}}<1\}.
\]
Indeed, $a\mapsto(S_a)_*P$ is weakly continuous and total variation
is lower semicontinuous under weak convergence. Write $H$ as the
product of the circle M\"obius group and the subgroup fixing three
prescribed points. The latter is identified with $T_0(1)$, a connected
Hilbert manifold and topological group whose topology is finer than
uniform convergence \cite[Chapter~I, Sections~2--3 and Appendix~A]{TT06}.
Both factors are connected Baire groups in their respective topologies.
For either factor $F$, its inclusion into $H$ with the uniform topology
is continuous. Hence $K_F=K\cap F$ is Borel in $F$.
The map $F/K_F\to H/K$ given by $fK_F\mapsto fK$ is injective,
since $f^{-1}g\in K$ for $f,g\in F$ implies $f^{-1}g\in K_F$.
Thus $K_F$ has countable index in $F$. If $K_F$ were meagre, its countably
many translates would make $F$ meagre in itself, contrary to the Baire
category theorem. By Pettis' theorem, $K_F^{-1}K_F=K_F$ therefore
contains a neighborhood of the identity. The subgroup $K_F$ and all
its cosets are open, so connectedness gives $K_F=F$. Thus $K$ contains
both factors and hence $K=H$.

This proves $\psi H\psi^{-1}\subset H$. Applying the same argument to
$T^{-1}$ gives the reverse inclusion. Inverting the welding interchanges
the two factors and gives $\varphi\in\mc N(H)$ as well.
Now $H=\WP(\m S^1)$, so Theorem~\ref{thm:main} identifies the normalizer
with $\WP(\m S^1)$. Theorem~\ref{body:thm:two} gives the reverse inclusion
for $\tilde\QI$.
\end{proof}

\subsection{Further consequences for group isomorphisms}\label{subsec:normalizer-isomorphisms}
The normalizer theorem determines the abstract automorphisms of the WP
group. Orientation-reversing conjugacies must also be included. Let
$j(z)=\bar z$ and set
\[
 \WP^\pm(\m S^1)=\WP(\m S^1)\cup
 \{\chi\circ j:\chi\in\WP(\m S^1)\}.
\]
Reflection preserves the WP class, so this is
$\WP(\m S^1)\rtimes\langle j\rangle$.

\begin{thm}[Automorphisms of the WP group]\label{thm:wp-automorphisms}
Every group automorphism $\mathcal A$ of $\WP(\m S^1)$ has the form
\[
 \mathcal A(\varphi)=\Gamma\circ\varphi\circ\Gamma^{-1}
 \qquad(\varphi\in\WP(\m S^1))
\]
for a unique $\Gamma\in\WP^\pm(\m S^1)$. Conversely, every such
$\Gamma$ defines an automorphism. In particular,
\[
 \operatorname{Aut}(\WP(\m S^1))
 \cong\WP(\m S^1)\rtimes\m Z_2.
\]
\end{thm}
\begin{proof}
Rubin's reconstruction theorem \cite[Theorem~2]{Rub89} supplies a
conjugating homeomorphism $\Gamma\in\Homeo(\m S^1)$. Its local-density
hypothesis holds because $\Diff_+(\m S^1)\subset\WP(\m S^1)$: for an
open arc $J$ and $x\in J$, smooth diffeomorphisms supported in $J$ move
$x$ throughout a neighborhood of $x$. The circle is locally compact
and has no isolated points, as required.

The map $\Gamma$ normalizes WP. If it preserves orientation,
Theorem~\ref{thm:main} gives $\Gamma\in\WP(\m S^1)$. If it reverses orientation,
then $\Gamma\circ j$ preserves orientation and normalizes WP, so
$\Gamma\circ j\in\WP(\m S^1)$. Conversely, both WP maps and reflection
normalize the WP group.

For uniqueness, it suffices to show that the centralizer of
$\Diff_+(\m S^1)$ in $\Homeo(\m S^1)$ is trivial. Suppose a commuting
homeomorphism $q$ has $q(x)\ne x$. Choose a small open arc containing
$x$ but not $q(x)$ and a smooth diffeomorphism $u$ supported there
with $u(x)\ne x$. Then $q(u(x))=u(q(x))=q(x)$, contradicting
injectivity of $q$. Thus $q=\operatorname{id}$.
\end{proof}

\section{Radon--Nikodym derivatives for WP actions}\label{sec:rn-wp-sle}

We prove WP quasi-invariance for SLE welding throughout $0<\kappa\le4$,
adding the critical case to the subcritical result of \cite{F1}.
We extend the analytic
density formula of \cite[Corollary~1.3]{BJ25} to fixed WP deformations.
Theorems~\ref{body:thm:one} and~\ref{body:thm:two} supply the density
assertion of Theorem~\ref{introthm:main}. For a finite-energy welding,
the relative action is a difference of Loewner energies plus a
difference of logarithmic conformal radii. We give their sum a meaning
when the original welding is sampled from SLE.

The proof has two steps. Integration by parts gives the density for
analytic deformations. We then separate its linear part from a bounded
quadratic remainder. The linear part extends in $L^2$, while resolvent
estimates control the remainder. A moment bound makes the analytic
densities uniformly integrable and identifies their limit with the WP
pushforward. Inversion gives the two-sided formula.

The finite-mode right integration by parts formula, its horizontal
$H^{3/2}$ covariance, and the density for finite compositions of
horizontal finite Fourier flows also appear in \cite[Theorems~1.2--1.3,
Proposition~4.4, Proposition~6.2]{Cai26}. The horizontal space there has
zero modes $n=0,\pm1$. The formulas below also include $n=\pm1$; these
modes require separate verification. The difference in notation does not
change these overlapping results. Our main target is instead the endpoint
formula for a general fixed WP map and the convergence needed to identify
it as a Radon--Nikodym derivative.

\subsection{The law and the inputs}\label{subsec:rn-assumptions}

We fix the loop and welding normalizations, state the integration-by-parts
and inversion inputs, and specify the conventions for the densities.

Fix $\kappa\in(0,4]$ and put
\begin{equation}\label{eq:central}
 c_{\mathrm m}=1-6\left(\frac2{\sqrt\kappa}-\frac{\sqrt\kappa}{2}\right)^2,
 \qquad c_{\mathrm L}=26-c_{\mathrm m},\qquad c=\frac{c_{\mathrm L}}{12}>2.
\end{equation}
Let $M$ be the SLE$_\kappa$ loop measure restricted to loops separating $0$
and $\infty$. We use the normalization and the scale decomposition
\begin{equation}\label{eq:scale}
 M=\pi^{-1/2}\dd s\otimes\sigma,\qquad
 F=e^s f,\qquad f(z)=z+\sum_{j\ge2}a_jz^j,
\end{equation}
where $\sigma$ is a probability measure; see \cite[Section~1.3]{BJ24} for this normalization. Let $g$ map $\m D^*=\{z:|z|>1\}$
onto the exterior of $f(\m S^1)$, with $g(\infty)=\infty$ and $g'(\infty)>0$.
Let $\alpha$ be an independent uniform angle, and write
\begin{equation}\label{eq:law}
 \widehat H=U\circ g^{-1}\circ f,\qquad U=r_{-\alpha},\qquad
 \m P^{\mathrm{weld}}=\operatorname{Law}(\widehat H),\qquad r_\alpha(z)=e^{\ii\alpha}z.
\end{equation}
We write $\m E^{\mathrm{weld}}$ for expectation under
$\m P^{\mathrm{weld}}$. For $\kappa<4$, the GMC comparison in
Lemma~\ref{introlem:sle-gmc-lift} identifies this measure class;
the Radon--Nikodym derivatives computed here are those of
$\m P^{\mathrm{weld}}$ itself.
We use the Borel realization in which $f=f_{\widehat H}$ is normalized by $f(0)=0$
and $f'(0)=1$. The SLE loops in this range are almost surely Jordan
and conformally removable. For removability, use the local SLE
construction of \cite{zhan2020sleloop}, the subcritical regularity
results of \cite{RS05}, and the endpoint theorem
\cite[Theorem~1.1]{kavvadias2022conformal}. Thus the normalized welding
pair is determined by $\widehat H$ almost surely. The integration-by-parts
and reflection identities stated next hold for all $\kappa\in(0,4]$;
they are cited inputs, not additional hypotheses on the parameter.

The loop integration-by-parts formula is
\cite[Theorem~4.1]{BJ24}, whose parameter range is $\kappa\in(0,4]$;
see also \cite[Corollary~3.11]{GQW24}.
We use the following specialization, verified in the proof of
Lemma~\ref{lem:mixed-ibp}. Let $w_j(\xi)=-\xi^{j+1}$, and let $\mathcal L_j$ be the complex
variation corresponding to the conformal deformation of the loop
$\xi\mapsto\xi+\varepsilon w_j(\xi)+o(\varepsilon)$. This acts on the loop
in the plane. It differs from a variation by composition on $\m S^1$.
For $j=-2,-1,0$,
\begin{equation}\label{eq:BJ}
 \begin{gathered}
 \int (\mathcal L_jU)V\dd M
 =-\int U\mathcal L_jV\dd M+
 \frac{c_{\mathrm m}}{12}\int\theta_jUV\dd M,\\
 \theta_0=\theta_{-1}=0,\qquad \theta_{-2}=6e^{-2s}(a_3-a_2^2).
 \end{gathered}
\end{equation}
The initial test functions are finite sums of
$\chi(s)P(a,\bar a)\zeta(K)$, where $\chi,\zeta\in C_c^\infty(\m R)$,
$P$ is a polynomial in finitely many coefficients, and
\begin{equation}\label{eq:K}
 K(h)=\log|g'(\infty)|-\log|f'(0)|
\end{equation}
for a welding pair $h=g^{-1}f$. We also use rotation invariance and
the standard differentiability of local conformal deformations.

For the two-sided result we also use the reflection identity, valid
throughout the same parameter range; see
\cite[Lemma~2.6]{GQW24} and \cite[Section~1.5]{BJ24}:
\begin{equation}\label{eq:inversion-input}
 \iota_*M=M,\qquad \iota(z)=1/\bar z.
\end{equation}
The one-sided result does not use \eqref{eq:inversion-input}. We use the
measurable Riemann mapping theorem, the Grunsky inequalities, standard
properties of the Beurling transform and Schatten ideals, the group property
of $\WP(\m S^1)$, and its characterization by $\log\varphi'\in H^{1/2}(\m S^1)$.

For $\varphi\in\WP(\m S^1)$ put $R_\varphi(h)=h\circ\varphi$ and, when it
exists, set
\begin{equation}\label{eq:rho}
 \rho_\varphi=\frac{\dd(R_\varphi)_*\m P^{\mathrm{weld}}}{\dd\m P^{\mathrm{weld}}},\qquad
 p_\chi=\rho_{\chi^{-1}}.
\end{equation}
Here $p_\chi$ is the density for $h\mapsto h\chi^{-1}$; densities are
evaluated at the target point. When they exist, change of variables and
the chain rule give
\begin{equation}\label{eq:pinv}
 \m E^{\mathrm{weld}}[p_\chi(\widehat H)F(\widehat H\chi)]=\m E^{\mathrm{weld}}[F(\widehat H)],\qquad
 p_{\chi\omega}(h)=p_\chi(h)p_\omega(h\chi).
\end{equation}
All statements below hold for each fixed choice of the parameters, outside
an exceptional set that may depend on them.

\subsection{Energy and the two contributions to the density}\label{subsec:rn-contributions}

We use $h\chi$ for $h\circ\chi$. Since densities are evaluated at the
target welding, the density $\rho_\varphi$ compares $h$ with
$h\varphi^{-1}$. To describe this comparison, we first recall the
Grunsky representation of Loewner energy.

Let $(f,g)$ be a normalized disjoint pair, with $f(0)=0$, $f'(0)=1$,
$g(\infty)=\infty$, and $g'(\infty)=b\ne0$. Following \cite{TT06}, define
$b_{mn}$ by
\begin{align}\label{eq:grunsky-coefficients}
 \log\frac{g(z)-g(w)}{z-w}
 &=\log b-\sum_{m,n\ge1}b_{mn}z^{-m}w^{-n},\\
 \log\frac{g(z)-f(w)}{bz}
 &=-\sum_{m\ge1,n\ge0}b_{m,-n}z^{-m}w^n,\\
 \log\frac{f(z)-f(w)}{z-w}
 &=\log b-\sum_{m,n\ge0}b_{-m,-n}z^mw^n.
\end{align}
The logarithms are normalized by their expansions at the indicated points;
$b_{0,0}=\log b$. Put $b_{-n,m}=b_{m,-n}$. On $\ell^2(\m N)$ define the four operators, denoted here by $G_1,G_2,G_3,G_4$, by
\begin{equation}\label{eq:TT}
 \begin{aligned}
 (G_1)_{mn}&=\sqrt{mn}\,b_{-m,-n},&
 (G_2)_{mn}&=\sqrt{mn}\,b_{-m,n},\\
 (G_3)_{mn}&=\sqrt{mn}\,b_{m,-n},&
 (G_4)_{mn}&=\sqrt{mn}\,b_{m,n}.
 \end{aligned}
\end{equation}
For a Jordan welding $h$, fix a normalized welding pair and use it to
define $G_j(h)$ and $K(h)$. For a composed welding, use the pair obtained
from this pair by quasiconformal surgery. For a conformally removable
curve, the normalized pair is uniquely determined by $h$.
We also write $G_1(f)$ for the inner operator defined by the coefficients of $f$.
In particular, $G_1(h)$ depends only on the normalized inner map. The
Grunsky inequality gives $\|G_1(h)\|\le1$ even when $h$ is not
quasisymmetric. If $\varphi\in\WP(\m S^1)$, then $G_1(\varphi)$ and
$G_4(\varphi)$ are Hilbert--Schmidt. We do not assume that $G_1(\widehat H)$ is
Hilbert--Schmidt. Entrywise complex conjugation of a matrix is denoted by a
bar, while $*$ denotes the Hilbert space adjoint.

For $h\in\WP(\m S^1)$, let $I^{\mathrm L}(h)$ be the Loewner energy
of its welded curve. The energy identity of \cite[Theorem~8.1]{W2}
and the Fredholm determinant formula of \cite{TT06} give
\begin{equation}\label{eq:rn-loewner-energy}
 \begin{aligned}
 I^{\mathrm L}(h)
 &=\frac1\pi\int_{\m D}\left|\frac{f''}{f'}\right|^2\dd A
   +\frac1\pi\int_{\m D^*}\left|\frac{g''}{g'}\right|^2\dd A-4K(h)\\
 &=-12\log\det\bigl(I-G_1(h)^*G_1(h)\bigr).
 \end{aligned}
\end{equation}
Both expressions are independent of the normalization of the welding
pair. The finite-energy form of our action is
\begin{equation}\label{eq:rn-geometric-action}
 \frac{c_{\mathrm L}}{24}
 \bigl[I^{\mathrm L}(h\varphi^{-1})-I^{\mathrm L}(h)\bigr]
 +2\bigl[K(h\varphi^{-1})-K(h)\bigr],
 \qquad h,\varphi\in\WP(\m S^1).
\end{equation}
Proposition~\ref{prop:rn-finite-energy} below identifies this expression
with the explicit action. In particular, the second contribution is the
change in the ratio of the exterior and interior conformal radii, not a
boundary derivative evaluated at an unspecified point.

For analytic $\chi$ and a Jordan welding $h$, the determinant cocycle
$D(h,\chi)$ in \eqref{eq:geom} is defined even when $I^{\mathrm L}(h)$
is infinite. Lemma~\ref{lem:analytic-cocycle} proves
\begin{equation}\label{eq:rn-two-contributions}
 p_\chi(h)=\exp\{-c_{\mathrm L}D(h,\chi)-2\Delta_\chi K(h)\},
 \qquad \Delta_\chi K(h)=K(h\chi)-K(h).
\end{equation}
When $h,\chi\in\WP(\m S^1)$, it also proves
$24D(h,\chi)=I^{\mathrm L}(h\chi)-I^{\mathrm L}(h)$.
For a general WP deformation and an SLE welding, the ordinary mixed
Fredholm determinant need not exist. We therefore combine the linear
parts of the two contributions before taking the WP limit. Their sum is
the logarithmic derivative defined next.

\begin{df}\label{def:rn-divergence}
For a real vector field $v(\theta)=\sum_{n\in\m Z}v_ne^{\ii n\theta}$ let
\begin{equation}\label{eq:sobolev}
 [v]_{3/2}^2=2\sum_{n\ge1}n^3|v_n|^2,\qquad v_{-n}=\bar v_n.
\end{equation}
For a normalized univalent $f$ define
\begin{equation}\label{eq:stress}
 \begin{aligned}
 \mathcal T_f(z)&=cS_f(z)+\left(\frac{f'(z)}{f(z)}\right)^2-z^{-2}
 =\sum_{n\ge1}B_n(f)z^{n-2},\\
 S_f&=\frac{f'''}{f'}-\frac32\left(\frac{f''}{f'}\right)^2.
 \end{aligned}
\end{equation}
For a finite Fourier series set
\begin{equation}\label{eq:divergence}
 \mathfrak d_v(\widehat H)=2\operatorname{Im}\sum_{n\ge1}\bar v_n B_n(f_{\widehat H}).
\end{equation}
For finite Fourier $v$, the integration by parts formula below identifies
$\mathfrak d_v$ as the logarithmic derivative in the direction $v$.
\end{df}
\begin{lem}[The logarithmic derivative]\label{lem:divergence-estimate}
The map $v\mapsto \mathfrak d_v$ has a unique continuous extension from real finite
Fourier series to $H^{3/2}(\m S^1)/\m R$ with values in $L^2(\m P^{\mathrm{weld}})$. Moreover,
\begin{equation}\label{eq:divergence-moment}
 \|\mathfrak d_v\|_{L^2(\m P^{\mathrm{weld}})}^2
 =2\sum_{n\ge1}\bigl[c(n^3-n)+2n\bigr]|v_n|^2,
 \quad 2[v]_{3/2}^2\le\|\mathfrak d_v\|_2^2\le c[v]_{3/2}^2.
\end{equation}
The series in \eqref{eq:divergence} is interpreted as an $L^2(\m P^{\mathrm{weld}})$ limit for
general $v\in H^{3/2}$. It is not asserted to converge at every welding.
\end{lem}

\subsection{The explicit WP action and the one-sided theorem}\label{subsec:rn-formula}

The formula for a fixed WP deformation uses only its Grunsky operator,
a Beltrami representative, and the normalized inner map of the random
welding. We define these terms before stating the density theorem.
The representative will disappear from the final random variable.

Let $\varphi\in\WP(\m S^1)$. Choose a quasiconformal extension
$a_\mu:\m D\to\m D$ with $a_\mu(0)=0$, $a_\mu|_{\m S^1}=\varphi$ and
\begin{equation}\label{eq:mu}
 \begin{aligned}
 \mu&=0\quad\text{on }|z|<r_0,\qquad k=\|\mu\|_\infty<1,\\
 N(\mu)^2&=\frac1\pi\int_{\m D}
 \frac{|\mu(z)|^2}{(1-|z|^2)^2}\dd A(z)<\infty.
 \end{aligned}
\end{equation}
Such a representative exists: modify an arbitrary WP extension on a compact
subdisk, make it conformal near $0$, and join the two boundary values by a
quasiconformal annulus map. The hyperbolic $L^2$ condition is unchanged near
$\m S^1$.

Extend $\mu$ by zero outside $\m D$ and put
\begin{equation}\label{eq:eta}
 \eta(z)=\left(\frac z{\bar z}\right)^2\overline{\mu(1/\bar z)},
 \qquad |z|>1.
\end{equation}
Let
\begin{equation}\label{eq:beurling}
 \mathcal Cu(z)=\frac1\pi\int_{\m C}\frac{u(w)}{z-w}\dd A(w),\quad
 \mathcal B=\partial\mathcal C,\quad
 \mathcal B_*=\mathbf 1_{\m D^*}\mathcal B\mathbf 1_{\m D^*}.
\end{equation}
Write $M_\lambda$ for multiplication by $\lambda$. Let $J,\bar J:\ell^2\to
L^2(\m D^*)$ be the complex-linear isometries
\begin{equation}\label{eq:J}
 Jx=\sum_{n\ge1}x_n\sqrt{\frac n\pi}z^{-n-1},\qquad
 \bar Jx=\sum_{n\ge1}x_n\sqrt{\frac n\pi}\bar z^{-n-1}.
\end{equation}
Both series converge in $L^2(\m D^*)$, and both maps are complex
linear. Their ranges are the holomorphic and anti-holomorphic Bergman
subspaces, respectively. Set
\begin{equation}\label{eq:Q}
 Q_\eta=(\bar J^*M_\eta)\mathcal B_*(I-M_\eta\mathcal B_*)^{-1}(M_\eta J).
\end{equation}
The outer factors map $\ell^2$ to $L^2(\m D^*)$ and back; the middle
factors, including the resolvent, act on the full $L^2(\m D^*)$ space.
Thus $Q_\eta$ is an operator on $\ell^2$.
Lemma~\ref{lem:resolvents} relates it to $G_1(\varphi)$ and proves
that it is trace class. Matrix realizations and the relation to Fourier
coordinates are given in Section~\ref{subsec:operator-realization}.

For a normalized univalent $f$ put
\begin{equation}\label{eq:lambdaf}
 \lambda_f(f(z))=\mu(z)\frac{f'(z)}{\overline{f'(z)}}\quad(z\in\m D),
 \qquad \lambda_f=0\quad\text{on }\m C\setminus f(\m D).
\end{equation}
\begin{df}\label{def:rn-resolvent}
For $\lambda=\mu$ or $\lambda=\lambda_f$ define $q_\lambda=\lambda/z$ and
\begin{equation}\label{eq:resolvent}
 \mathcal R(\lambda)=\frac1\pi\int_{\m C} q_\lambda\,
 \mathcal B(I-M_\lambda\mathcal B)^{-1}q_\lambda\dd A.
\end{equation}
The pairing in \eqref{eq:resolvent} is bilinear, without complex conjugation.
Both factors belong to $L^2$, so the integral is absolutely convergent. We call $\mathcal R(\lambda)$ the resolvent term.
\end{df}
For $n\ge1$ put
\begin{equation}\label{eq:vmu}
 \begin{aligned}
 m_n(\mu)&=\frac1\pi\int_{\m D}\mu(z)z^{n-2}\dd A(z),\\
 (v_\mu)_n&=-\ii\overline{m_n(\mu)},\qquad
 (v_\mu)_{-n}=\overline{(v_\mu)_n},\qquad (v_\mu)_0=0.
 \end{aligned}
\end{equation}
We show below that $v_\mu\in H^{3/2}(\m S^1)$.

We write $\HS$ for the Hilbert--Schmidt class and $\|\cdot\|_{\HS}$ for
its norm. For a trace-class operator $A$, write
$\|A\|_{\tr}=\tr\bigl((A^*A)^{1/2}\bigr)$ for its trace norm. Here $\tr A$
is the trace of $A$; it does not denote the trace-class ideal.
For $L\in\HS(\ell^2)$ with $\|L\|<1$, set
\begin{equation}\label{eq:det2}
 \log\det{}_2(I-L)=-\sum_{j\ge2}\frac{\tr(L^j)}j.
\end{equation}
When $L$ is trace class, $\log\det(I-L)=\log\det{}_2(I-L)-\tr L$.
For an analytic deformation, expanding the ordinary determinant gives
$\log\det(I-L)=\log\det{}_2(I-L)-\tr L$. The resolvent identity
\eqref{eq:B1-resolvent} splits $G_1(\varphi)$ into a term linear in
$\mu$ and the trace-class remainder $Q_\eta$. The linear trace and the
linear change of conformal radius combine into $\mathfrak d_{v_\mu}$.
This cancellation leaves the formula below. Lemma~\ref{lem:analytic-action}
verifies the equality term by term.

\begin{df}\label{body:def:rn-action}
The action is
\begin{equation}\label{body:eq:action}
 \begin{aligned}
 \mathscr I_\mu(\varphi,h)={}&-\frac{c_{\mathrm L}}2
 \log\det\bigl(I-G_1(\varphi)^*G_1(\varphi)\bigr)
 +2K(\varphi^{-1})+\mathfrak d_{v_\mu}(h)\\
 &+c_{\mathrm L}\operatorname{Re}\log\det{}_2\bigl(I-G_1(\varphi)\overline{G_1(h)}\bigr)
 -c_{\mathrm L}\operatorname{Re}\tr\bigl(Q_\eta\overline{G_1(h)}\bigr)\\
 &+2\operatorname{Re}\bigl[\mathcal R(\lambda_{f_h})-\mathcal R(\mu)\bigr].
 \end{aligned}
\end{equation}
The first determinant is an ordinary positive Fredholm determinant.
The mixed product $G_1(\varphi)\overline{G_1(h)}$ is Hilbert--Schmidt
with operator norm strictly below one, so its regularized logarithm is defined
by \eqref{eq:det2}. The trace involving $Q_\eta$ and the two scalar
resolvents are finite. Lemma~\ref{lem:bounds} proves that all these
nonlinear terms are bounded uniformly in the normalized univalent $f_h$.
Only $\mathfrak d_{v_\mu}$ requires an $L^2(\m P^{\mathrm{weld}})$
interpretation. Thus \eqref{body:eq:action} defines a real random
variable modulo null sets without subtracting two infinite energies.
For deterministic WP weldings, Proposition~\ref{prop:rn-finite-energy}
proves absolute convergence of the linear series and identifies the
action with \eqref{eq:rn-geometric-action}.
\end{df}

\begin{thm}[One-sided WP density]\label{body:thm:one}
For every $\kappa\in(0,4]$ and fixed
$\varphi\in\WP(\m S^1)$, the expression \eqref{body:eq:action} is a finite real
$L^2(\m P^{\mathrm{weld}})$ random variable and
\begin{equation}\label{body:eq:one}
 \frac{\dd(R_\varphi)_*\m P^{\mathrm{weld}}}{\dd\m P^{\mathrm{weld}}}(h)
 =e^{-\mathscr I_\mu(\varphi,h)}>0\quad\m P^{\mathrm{weld}}\text{-a.s.}
\end{equation}
The action $\mathscr I_\mu(\varphi,\cdot)$ is independent of the choice
of $\mu$ in \eqref{eq:mu} as an $L^2(\m P^{\mathrm{weld}})$ equivalence
class; we denote it by $\mathscr I(\varphi,\cdot)$. If $k>0$, then
\begin{equation}\label{body:eq:one-moments}
 \m E^{\mathrm{weld}}[\rho_\varphi(\widehat H)^r]<\infty\qquad(|r|k<1).
\end{equation}
For each fixed $\varphi$, $\rho_\varphi\in L^{2+\epsilon}(\m P^{\mathrm{weld}})$ for some $\epsilon>0$, and $\mathscr I(\varphi,\cdot)\in L^p(\m P^{\mathrm{weld}})$ for every finite $p$.
If $\varphi$ is a rotation, $\rho_\varphi=1$.
\end{thm}

\subsection{Circle derivations and integration by parts}\label{subsec:rn-infinitesimal}

We pass from the loop integration-by-parts formula to differentiation
under right composition. The circle derivations identify this variation,
and their Lie bracket reduces the calculation to the first two modes.
The resulting covariance identity extends the logarithmic derivative
from finite Fourier fields to $H^{3/2}$.

\begin{df}\label{def:circle-derivations}
Let $\mathscr P$ be the algebra of polynomials in finitely many $a_j,\bar a_j$.
All its elements are bounded by the univalent coefficient estimates. Define
$Z_nf=z^{n+1}f'$ for $n\ge1$. If $\omega=f^{-1}$ is the inverse germ at $0$, define
\begin{equation}\label{eq:Yn}
 \begin{aligned}
 (V_nf)(\xi)&=\bigl[\omega(\xi)^{1-n}f'(\omega(\xi))\bigr]_{<0},\\
 (W_nf)(z)&=(V_nf)(f(z))-z^{1-n}f'(z),\\
 (Y_nf)(z)&=(W_nf)(z)-(W_nf)(0)-(W_nf)'(0)f(z),\\
 D_n&=Y_n+\bar Z_n,\qquad \bar D_n=Z_n+\bar Y_n.
 \end{aligned}
\end{equation}
Here $[\cdot]_{<0}$ denotes the sum of the strictly negative powers in the
Laurent expansion at $0$. For fixed $f$, $(V_nf)(f(z))$ means that the function $V_nf$ is
evaluated at $f(z)$. The singular terms cancel in $W_nf$, so it is
holomorphic near $0$. The maps $V_n$ and $W_n$ assign functions to $f$;
they are not asserted to be linear operators or derivations. The derivations
$Z_n,Y_n$ act on the holomorphic coefficients and vanish on their
conjugates. Their conjugates act on the conjugate coefficients.
\end{df}
\begin{lem}[Circle variations and their bracket]\label{lem:circle-derivations}
The derivations $D_n$ and $\bar D_n$ in
Definition~\ref{def:circle-derivations} preserve $\mathscr P$. Let
$\chi_\varepsilon$ be a real two-parameter family of analytic circle
diffeomorphisms with
\[
 \chi_\varepsilon(z)=z-\varepsilon z^{1-n}
 +\bar\varepsilon z^{n+1}+O(|\varepsilon|^2)\qquad(z\in\m S^1).
\]
For the normalized inner map $f_\varepsilon$ of $h\circ\chi_\varepsilon$,
write $P_\varepsilon=P(f_\varepsilon,\bar f_\varepsilon)$. Then, for
$P\in\mathscr P$,
\[
 P_\varepsilon=P+\varepsilon D_nP+\bar\varepsilon\,\bar D_nP
 +o(|\varepsilon|).
\]
Thus $D_n$ and $\bar D_n$ are the Wirtinger components of the real
variation. Moreover,
\begin{equation}\label{eq:V}
 (V_nf)(\xi)=\sum_{j=0}^{n-2}[z^{n-2}](f'(z)^2f(z)^j)\xi^{-j-1},
 \qquad [D_1,D_n]=(1-n)D_{n+1}\quad(n\ge2).
\end{equation}
\end{lem}
\begin{proof}
We derive the variation by solving the Beltrami equation. Choose a
quasiconformal $a_\varepsilon:\m D\to\m D$ which is the identity near $0$
and agrees with $\chi_\varepsilon^{-1}$ on a collar of $\m S^1$. Let
$\mu_\varepsilon$ be its Beltrami coefficient and put
$\lambda_\varepsilon=f_*\mu_\varepsilon$, extended by zero outside
$f(\m D)$. Let $W^\varepsilon$ solve the Beltrami equation with coefficient
$\lambda_\varepsilon$, normalized by $W^\varepsilon(0)=0$ and
$(W^\varepsilon)'(\infty)=1$. If $h=g^{-1}\circ f$, the chain rule gives
\[
 \widetilde f_\varepsilon=W^\varepsilon\circ f\circ a_\varepsilon^{-1},
 \qquad \widetilde g_\varepsilon=W^\varepsilon\circ g,\qquad
 \widetilde g_\varepsilon^{-1}\circ\widetilde f_\varepsilon
 =h\circ\chi_\varepsilon.
\]
Both maps are conformal in their respective disks. The deformed curve is
$W^\varepsilon(f(\m S^1))$. Thus right composition by
$\chi_\varepsilon$ is realized by this deformation of the curve.

Set $u_n=(W_nf)(0)$ and $\ell_n=(W_nf)'(0)$. Write
$\dot W=\partial_\varepsilon W^\varepsilon|_0$ and
$\dot W^-=\partial_{\bar\varepsilon}W^\varepsilon|_0$. Since the support
of $\lambda_\varepsilon$ stays in a compact subset of $f(\m D)$, these
fields are holomorphic near and outside the loop and are bounded at
infinity. On a collar inside $\m S^1$, differentiation gives
\[
 \partial_\varepsilon\widetilde f_\varepsilon|_0
 =\dot W(f(z))-z^{1-n}f'(z),\qquad
 \partial_{\bar\varepsilon}\widetilde f_\varepsilon|_0
 =\dot W^-(f(z))+z^{n+1}f'(z).
\]
The two left sides are holomorphic in $\m D$. In the coordinate
$\xi=f(z)$, the principal part of the first right side is cancelled by
$(V_nf)(\xi)$, and the second has no principal part. Gluing the resulting
holomorphic functions across the collar and using boundedness at infinity
shows that $\dot W-V_nf$ and $\dot W^-$ are constants on the exterior.
Since $\widetilde f_\varepsilon(0)=0$, these constants are $-u_n$ and $0$.
Hence, as holomorphic germs at $0$,
\[
 \widetilde f_\varepsilon=f+\varepsilon(W_nf-u_n)
 +\bar\varepsilon z^{n+1}f'+o(|\varepsilon|),\qquad
 \widetilde f_\varepsilon'(0)=1+\varepsilon\ell_n+o(|\varepsilon|).
\]
Normalize by $f_\varepsilon=\widetilde f_\varepsilon/
\widetilde f_\varepsilon'(0)$. Then
\[
 f_\varepsilon=f+\varepsilon Y_nf+\bar\varepsilon Z_nf+o(|\varepsilon|).
\]
Taking conjugate coefficients and applying the chain rule to $P$ proves
its stated variation. Equivalently, after the affine normalization, the
$\varepsilon$ component of the loop deformation is
$(V_nf)(\xi)-u_n-\ell_n\xi$.

For the coefficient of $\xi^{-j-1}$ in $V_nf$, change variables
$\xi=f(z)$ in the residue. This gives
\[
 [\xi^{-j-1}]V_nf
 =\operatorname{Res}_{z=0}z^{1-n}f'(z)^2f(z)^j\,\dd z
 =[z^{n-2}]\bigl(f'(z)^2f(z)^j\bigr).
\]
Only $0\le j\le n-2$ occurs. This proves the first formula in
\eqref{eq:V}. Each coefficient of $Y_nf$ is a polynomial in finitely
many $a_j$. With $a_1=1$, the other component is explicitly
\[
 Z_na_k=\begin{cases}(k-n)a_{k-n},&k\ge n+1,\\0,&k\le n,
 \end{cases}
 \qquad
 D_na_k=[z^k]Y_nf,\quad D_n\bar a_k=\overline{Z_na_k}.
\]
The conjugate formulas give $\bar D_na_k=Z_na_k$ and
$\bar D_n\bar a_k=\overline{D_na_k}$. For instance, writing
$a=a_2$, $b=a_3$, and $d=a_4$,
\[
 \begin{aligned}
 D_1a&=2a^2-3b,&D_1\bar a&=1,\\
 D_2a&=-2a^3+6ab-5d,&D_2\bar a&=0.
 \end{aligned}
\]
For the bracket, use the angle coordinate $z=e^{\ii\theta}$. The
$\varepsilon$ component of the displayed circle variation has angle
velocity $v_n(\theta)=\ii e^{-\ii n\theta}$. Thus
\[
 v_1v_n'-v_nv_1'=(n-1)v_{n+1}.
\]
To fix the sign, let $\phi_s^v,\phi_t^w$ be the circle flows of smooth
real fields $v,w$, and define $T_v(s)P(h)=P(h\circ\phi_s^v)$. Then
\[
 \left.\partial_s\partial_t\right|_{s=t=0}
 \bigl(\phi_s^v\circ\phi_t^w-\phi_t^w\circ\phi_s^v\bigr)(x)
 =v'(x)w(x)-w'(x)v(x).
\]
The surgery construction applies to both parameters. Differentiating its
coefficient functions gives
\[
 [D_v,D_w]P=D_{v'w-w'v}P=-D_{vw'-wv'}P.
\]
Extend the identity complex linearly. Taking $v=v_1$ and $w=v_n$ gives
$[D_1,D_n]=(1-n)D_{n+1}$.
\end{proof}

For a real finite Fourier field
$v(\theta)=v_0+\sum_{n\ge1}(v_ne^{\ii n\theta}+\bar v_ne^{-\ii n\theta})$,
define the derivation $A_v$ on $\mathscr P$ by
\[
 A_v=v_0A_0+\sum_{n\ge1}(-\ii\bar v_nD_n+\ii v_n\bar D_n),\qquad
 A_0=\ii\sum_{j\ge2}(j-1)
 \left(a_j\frac{\partial}{\partial a_j}
 -\bar a_j\frac{\partial}{\partial\bar a_j}\right).
\]
The sums are finite on each polynomial. If $e_t^v$ is the circle flow of
$v$, Lemma~\ref{lem:circle-derivations} gives
\[
 A_vP(h)=\left.\frac{\dd}{\dd t}\right|_{t=0}P(h e_t^v).
\]
For the constant field, this follows from the normalized inner map
$f_t(z)=e^{-\ii t}f(e^{\ii t}z)$ of $h e_t^1$.

\begin{lem}[Right-composition integration by parts]\label{lem:mixed-ibp}
For $P\in\mathscr P$ and $n\ge1$,
\begin{equation}\label{eq:ibp-all}
 \m E^{\mathrm{weld}}[D_nP]=\m E^{\mathrm{weld}}[B_nP],\qquad
 \m E^{\mathrm{weld}}[\bar D_nP]=\m E^{\mathrm{weld}}[\bar B_nP].
\end{equation}
\end{lem}
\begin{proof}
We first verify \eqref{eq:BJ} with the stated test domain. In \cite[(2.5)--(2.6)]{BJ24}, the test space is generated
by $\chi(s)P(a,\bar a)\zeta(\widetilde c_0)$, where
$\widetilde c_0=-K$. Replacing $\zeta(x)$ by $\zeta(-x)$ therefore gives the test domain
specified in Section~\ref{subsec:rn-assumptions}. The adjoint formula of
\cite[Theorem~4.1]{BJ24}, applied to $U$ and $\bar V$, gives
\[
 \int (\mathcal L_jU)V\,\dd M
 =-\int U\mathcal L_jV\,\dd M
 +\frac{c_{\mathrm m}}{12}\int(\vartheta_j-\widetilde\vartheta_j)UV\,\dd M.
\]
For $j=-2,-1,0$, the exterior Schwarzian has no residue against
$w_j(\xi)=-\xi^{j+1}$. At $0$, $S(f^{-1})(0)=-6(a_3-a_2^2)$.
The residue formula in \cite[(1.16)]{BJ24} therefore gives
\[
 \widetilde\vartheta_j=0,\quad
 \vartheta_0=\vartheta_{-1}=0,\quad
 \vartheta_{-2}=-e^{-2s}S(f^{-1})(0)
 =6e^{-2s}(a_3-a_2^2).
\]
This proves \eqref{eq:BJ} on the stated test domain.

We first remove the cutoff in $K$. Take $\zeta=1$ on $[-1,1]$ and
$\operatorname{supp}\zeta\subset[-2,2]$. All terms have scale support in a
fixed compact interval $J$, and $M\{s\in J\}=\pi^{-1/2}|J|<\infty$.
The coefficient bounds $|a|\le2$ and $|b|\le3$ give the needed
uniform estimates. Indeed, a translation and the local deformation with
vector field $-1/\xi$ do not change the exterior conformal radius to first
order, while the vector field $-\xi$ contributes $-1/2$ to its logarithm.
Using the scale derivatives computed below gives
\[
 \mathcal L_{-1}K=-e^{-s}a,\qquad
 \mathcal L_{-2}K=-e^{-2s}(2b-\tfrac72a^2),\qquad
 \mathcal L_0K=0.
\]
Thus these three derivatives are bounded on the fixed scale interval $J$.
The remaining coefficients and $\theta_j$ are bounded there as well.
Differentiating $\zeta(K/R)$ produces a factor
$R^{-1}\zeta'(K/R)\mathcal L_jK$, so its $L^1(M)$ norm is $O(R^{-1})$.
Dominated convergence removes the remaining cutoffs in \eqref{eq:BJ}.

At scale $s=0$, for a vector field acting on the loop $w$ write
$q=w(f)/f'=\sum q_jz^j$ and $f^\#(z)=z+\sum_{j\ge2}\bar a_jz^j$.
The normalized variations are
\begin{equation}\label{eq:physical}
 \begin{aligned}
 \mathcal L_ws&=q_1/2,\\
 \mathcal L_wf&=w(f)-f'\sum_{j\le0}q_jz^j-\tfrac12q_1(zf'+f),\\
 \mathcal L_wf^\#&=(f^\#)'\sum_{j\le0}q_jz^{2-j}
 +\tfrac12q_1(z(f^\#)'-f^\#).
 \end{aligned}
\end{equation}
To check the first two identities, let the infinitesimal change of the
interior coordinate have meromorphic component $u(z)$, so that
$\mathcal L_w(e^sf)=w(e^sf)-e^sf'u$. The boundary condition for a real
circle deformation determines the conjugate component of $u$ by reflection
in $\m S^1$. Cancellation of the nonpositive powers in $q$ and normalization
of $f(0)$ and $f'(0)$ give
\[
 u(z)=\sum_{j\le0}q_jz^j+\frac{q_1}{2}z,
 \qquad \mathcal L_ws=\frac{q_1}{2}.
\]
More explicitly, the first-order change of the disk coordinate has the
form $tu+\bar t\widetilde u$, where the tangency condition is
$\widetilde u(z)=-z^2\overline{u(1/\bar z)}$ on the circle. Hence
\[
 \widetilde u(z)=-\sum_{j\le0}\bar q_jz^{2-j}
 -\frac{\bar q_1}{2}z.
\]
The coefficient of $\bar t$ in the scale is $\bar q_1/2$. Taking the
coefficient of $\bar t$ in the variation of $f$, then conjugating its
Taylor coefficients, gives the last line of \eqref{eq:physical}. In this
operation $f^\#$ has the conjugate coefficients, while the complex parameter
in $\mathcal L_w$ is not conjugated.

The calculations for $n=1,2$ provide the two initial cases of the
induction below. The case $n=1$ alone does not produce $n=2$, since
$[D_1,D_1]=0$. Put $a=a_2$, $b=a_3$, $d=a_4$, $e=a_5$. At $s=0$,
\[
 \begin{aligned}
 q_{-1}(z)&=-\frac1{f'(z)}
 =-1+2az+(3b-4a^2)z^2+O(z^3),\\
 q_{-2}(z)&=-\frac1{f(z)f'(z)}
 =-z^{-1}+3a+(4b-7a^2)z+O(z^2),\\
 q_0(z)&=-\frac{f(z)}{f'(z)}=-z+az^2+O(z^3).
 \end{aligned}
\]
Substitution into \eqref{eq:physical} yields the following values; bars
indicate the independent conjugate coefficient variables:
\[
 \begin{array}{c|ccc}
 &\mathcal L_js&\mathcal L_ja&\mathcal L_jb\\ \hline
 j=-1&a&3b-3a^2&4d-4ab\\
 j=-2&2b-\frac72a^2&\frac{23}{2}a^3-17ab+5d&-a^4+17a^2b-14ad-9b^2+6e\\
 j=0&-\frac12&\frac a2&b
 \end{array}
\]
\[
 \begin{array}{c|c}
 &\mathcal L_j\bar a\\ \hline
 j=-1&a\bar a-1\\
 j=-2&(2b-\frac72a^2)\bar a+3a\\
 j=0&-\frac12\bar a.
 \end{array}
\]
The same substitution for each
coefficient and scale covariance give, on $\mathscr P$,
\begin{equation}\label{eq:basefields}
 \begin{aligned}
 D_1&=-e^s\mathcal L_{-1}-2a\mathcal L_0,\\
 D_2&=-e^{2s}\mathcal L_{-2}-3e^sa\mathcal L_{-1}
 -(4b-a^2)\mathcal L_0.
 \end{aligned}
\end{equation}
For example, at scale zero the first vector field is $1+2a\xi$ and the
second is $\xi^{-1}+3a+(4b-a^2)\xi$. These equal
$(V_nf)(\xi)-(W_nf)(0)-(W_nf)'(0)\xi$ for $n=1,2$, respectively. Hence the
identities hold on every coefficient, not only on $a$ and $b$.
The scale components vanish: in the first line they are $-a+a$, and in
the second they are
\[
 -\left(2b-\tfrac72a^2\right)-3a^2+
 \tfrac12(4b-a^2)=0.
\]
For a test function $P$, apply \eqref{eq:BJ} to $C_jP$ and expand
$\mathcal L_j(C_jP)$. This gives
\[
 \int C_j\mathcal L_jP\,\dd M
 =\int P\left(\frac{c_{\mathrm m}}{12}C_j\theta_j
 -\mathcal L_jC_j\right)\dd M.
\]
Thus the divergence of $\sum_jC_j\mathcal L_j$ is the sum of the displayed
terms. For $D_1$ the coefficient derivatives are
$\mathcal L_{-1}(-e^s)=-a$ and $\mathcal L_0(-2a)=-a$, giving $2a$.
For $D_2$, at $s=0$ the table gives
\[
 \begin{aligned}
 \mathcal L_{-2}(-e^{2s})&=7a^2-4b,\\
 \mathcal L_{-1}(-3e^sa)&=-3(a^2+3b-3a^2)=6a^2-9b,\\
 \mathcal L_0(-(4b-a^2))&=a^2-4b.
 \end{aligned}
\]
The same result holds at arbitrary $s$ by scale covariance;
the loop term is $-c_{\mathrm m}(b-a^2)/2$. Their sum gives
\[
 -\bigl[(7a^2-4b)+(6a^2-9b)+(a^2-4b)\bigr]
 -\frac{c_{\mathrm m}}2(b-a^2)
 =(6c+4)b-(6c+1)a^2.
\]
Here we used $12c=26-c_{\mathrm m}$. These are precisely
\begin{equation}\label{eq:base-t}
 B_1=2a_2,\qquad B_2=(6c+4)a_3-(6c+1)a_2^2.
\end{equation}
This proves \eqref{eq:ibp-all} for $n=1,2$.

For the induction, the variation of \eqref{eq:stress} gives
\begin{equation}\label{eq:stress-rec}
 D_1B_n-D_nB_1=(1-n)B_{n+1}.
\end{equation}
Since $V_1f=0$ and $W_1f=-f'$, the normalized variation is
\[
 Y_1f=-f'+1+2a_2f.
\]
The variation $-f'$ is a translation of the argument of $f$, and
$1+2a_2f$ is an infinitesimal affine change of its value. The Schwarzian is
invariant under affine changes, so $Y_1S_f=-S_f'$. With $p=f'/f$,
\[
 Y_1p=\left(\frac{Y_1f}{f}\right)'
 =\left(-p+\frac1f+2a_2\right)'=-p'-\frac{f'}{f^2}.
\]
Consequently
\[
 \begin{aligned}
 Y_1\mathcal T_f
 &=-cS_f'-2pp'-2p\frac{f'}{f^2}\\
 &=-\mathcal T_f'+2z^{-3}-2\frac{f'^2}{f^3}.
 \end{aligned}
\]
We now derive the coefficient used in the recurrence. Write
$A(\xi)=\omega(\xi)^{-3}\omega'(\xi)
-a\omega(\xi)^{-2}\omega'(\xi)$. Since
$\omega(\xi)=\xi-a\xi^2+O(\xi^3)$, differentiation of the Laurent
series gives $A(\xi)=\xi^{-3}+O(1)$. Write
$\omega^{1-n}f'(\omega)=V_nf+\sum_{j\ge0}r_j\xi^j$. Changing variables
$\xi=f(z)$ and using the definition of $Y_n$, we obtain
\[
 \begin{aligned}
 Y_na_2
 &=\operatorname{Res}_{z=0}
 \bigl((V_nf)(f(z))-z^{1-n}f'(z)\bigr)(z^{-3}-az^{-2})\,\dd z\\
 &=-\operatorname{Res}_{\xi=0}\left(\sum_{j\ge0}r_j\xi^j\right)
 A(\xi)\,\dd\xi=-r_2.
 \end{aligned}
\]
The last equality uses $A(\xi)=\xi^{-3}+O(1)$. Since $V_nf$ has only
negative powers, changing variables again gives
\[
 r_2=\operatorname{Res}_{\xi=0}
 \frac{\omega(\xi)^{1-n}f'(\omega(\xi))}{\xi^3}\,\dd\xi
 =\operatorname{Res}_{z=0}\frac{z^{1-n}f'(z)^2}{f(z)^3}\,\dd z
 =[z^{n-2}]\frac{f'(z)^2}{f(z)^3}.
\]
Thus $Y_na_2=-[z^{n-2}]f'^2/f^3$. Let
$R_n=[z^{n-2}]f'^2/f^3$. Since $B_1=2a_2$, it follows that
$D_nB_1=-2R_n$. Taking the coefficient of $z^{n-2}$ in the preceding
identity for $Y_1\mathcal T_f$ gives
\[
 D_1B_n=-(n-1)B_{n+1}-2R_n.
\]
Subtracting $D_nB_1$ proves \eqref{eq:stress-rec}.

For the induction, assume the identities for $1$ and $n$. By the product
rule and the identity for $D_n$ applied to $B_1P$,
\[
 \begin{aligned}
 \m E^{\mathrm{weld}}[D_1D_nP]
 &=\m E^{\mathrm{weld}}[B_1D_nP]\\
 &=\m E^{\mathrm{weld}}[B_nB_1P-(D_nB_1)P].
 \end{aligned}
\]
Interchanging $1$ and $n$ gives
$\m E^{\mathrm{weld}}[D_nD_1P]
=\m E^{\mathrm{weld}}[B_1B_nP-(D_1B_n)P]$. Subtraction yields
\[
 \m E^{\mathrm{weld}}[[D_1,D_n]P]
 =\m E^{\mathrm{weld}}[(D_1B_n-D_nB_1)P].
\]
Substituting \eqref{eq:V} and \eqref{eq:stress-rec} gives
\[
 (1-n)\m E^{\mathrm{weld}}[D_{n+1}P]
 =(1-n)\m E^{\mathrm{weld}}[B_{n+1}P].
\]
For $n\ge2$, division by $1-n$ proves the first identity at $n+1$.
Starting from $n=1,2$ proves it for all $n\ge1$. Complex conjugation gives
the second identity. The input \eqref{eq:BJ} was specialized from \cite{BJ24} above;
the cutoff and coefficient calculations are part of this proof.
\end{proof}

\begin{proof}[Proof of Lemma~\ref{lem:divergence-estimate}]
For a holomorphic coefficient polynomial $F$, the derivation $\bar Y_n$
vanishes on $F$, so $\bar D_nF=Z_nF$. The second integration-by-parts
identity gives
\begin{equation}\label{eq:ward}
 \m E^{\mathrm{weld}}[Z_nF]=\m E^{\mathrm{weld}}[F\bar B_n].
\end{equation}
For $\delta f=z^{n+1}f'$ put $p=f'/f$. Direct differentiation gives
\[
 Z_nS_f=z^{n+1}S_f'+2(n+1)z^nS_f+(n^3-n)z^{n-2},
 \quad Z_np=(z^{n+1}p)'.
\]
Since $S_f$ is analytic at $0$ and $p^2=z^{-2}+O(z^{-1})$, the
coefficient of $z^{n-2}$ in $Z_n\mathcal T_f$ is
\begin{equation}\label{eq:diag}
 Z_nB_n=c(n^3-n)+2n.
\end{equation}
Each $B_n$ is a bounded coefficient polynomial. Taking $F=B_n$ in
\eqref{eq:ward} gives
$\m E^{\mathrm{weld}}[|B_n|^2]=c(n^3-n)+2n$.
For $f_\theta(z)=e^{-\ii\theta}f(e^{\ii\theta}z)$, the chain rule gives
$B_n(f_\theta)=e^{\ii n\theta}B_n(f)$. The law of the normalized shape is
rotation invariant. Therefore, for every $\theta$,
\[
 \begin{aligned}
 \m E^{\mathrm{weld}}[B_m\bar B_n]
 &=e^{\ii(m-n)\theta}\m E^{\mathrm{weld}}[B_m\bar B_n],\\
 \m E^{\mathrm{weld}}[B_mB_n]
 &=e^{\ii(m+n)\theta}\m E^{\mathrm{weld}}[B_mB_n].
 \end{aligned}
\]
The first expectation vanishes when $m\ne n$, and the second vanishes for
all $m,n\ge1$. Set $X=\sum_{n\ge1}\bar v_nB_n$ for a finite Fourier
series. The identities above imply
\[
 \m E^{\mathrm{weld}}[X^2]=0,\qquad
 \m E^{\mathrm{weld}}[|X|^2]
 =\sum_{n\ge1}\bigl[c(n^3-n)+2n\bigr]|v_n|^2.
\]
Since $\mathfrak d_v=2\operatorname{Im}X$, it follows that
$\m E^{\mathrm{weld}}[\mathfrak d_v^2]=2\m E^{\mathrm{weld}}[|X|^2]$. This proves
the identity in \eqref{eq:divergence-moment}; the bounds follow from
$2n^3\le c(n^3-n)+2n\le cn^3$. Fourier truncation proves the unique
$L^2$ extension. For a real finite Fourier field $v$, the
integration by parts identities in Lemma~\ref{lem:mixed-ibp}, together
with rotation invariance, also give
\[
 \begin{aligned}
 \m E^{\mathrm{weld}}[A_vP]
 &=\m E^{\mathrm{weld}}\!\left[
 \sum_{n\ge1}(-\ii\bar v_nB_n+\ii v_n\bar B_n)P\right]\\
 &=\m E^{\mathrm{weld}}[\mathfrak d_vP].
 \end{aligned}
\]
\end{proof}

\subsection{Finite flows and the analytic density}\label{subsec:rn-analytic}

We first integrate the logarithmic derivative along finite Fourier
flows. We then identify the integral with a determinant cocycle and a
change of conformal radius. The determinant calculation also fixes the
factor $1/24$ in the finite-energy formula.

Let $e_t^v$ be the circle flow of a real finite Fourier field $v$. Recall
that $A_vP=\left.\frac{\dd}{\dd t}\right|_{t=0}P(\,\cdot\,e_t^v)$, and that
$\m E^{\mathrm{weld}}[A_vP]=\m E^{\mathrm{weld}}[\mathfrak d_vP]$.

\begin{lem}[Density along a Fourier flow]\label{lem:finite-flow}
For every finite $T$,
\begin{equation}\label{eq:finite-flow}
 p_{e_T^v}(h)=\exp\left\{-\int_0^T \mathfrak d_v(he_t^v)\dd t\right\}.
\end{equation}
Thus $p_{e_T^v}=\rho_{e_{-T}^v}$, while the density for right
composition by $e_T^v$ is
\[
 \rho_{e_T^v}(h)=\exp\left\{\int_0^T \mathfrak d_v(he_{-t}^v)\dd t\right\}.
\]
The density and its reciprocal are bounded. In particular,
$(\rho_{e_t^v}-1)/t\to \mathfrak d_v$ in $L^1(\m P^{\mathrm{weld}})$ as $t\to0$.
\end{lem}
\begin{proof}
The inverse circle flow has a quasiconformal extension $a_t$ that is the
identity near $0$, agrees with the inverse flow on a fixed collar of
$\m S^1$, and has Beltrami coefficient supported in a fixed annulus
$r<|z|<R<1$. For small complex $t$, the coefficients are holomorphic in
$t$ in $L^\infty$ and have norm $O(|t|)$. After pushforward by any normalized
univalent $f$, Koebe's theorem places their support in a common bounded
annulus away from $0$. Let this annulus be $A$. Choose $p>2$ and $\delta>0$
so that, uniformly in $f$ and $|t|<\delta$,
\[
 \|M_{\lambda_{t,f}}\mathcal B\|_{L^p\to L^p}\le\tfrac12,\qquad
 \|\lambda_{t,f}\|_p\le C|t|.
\]
For $q_{t,f}=(I-M_{\lambda_{t,f}}\mathcal B)^{-1}\lambda_{t,f}$,
\[
 \|q_{t,f}\|_p\le2C|t|,\qquad
 W_{t,f}(z)=z+\mathcal Cq_{t,f}(z)-\mathcal Cq_{t,f}(0).
\]
The coefficient $q_{t,f}$ is supported in $A$. Differentiating the Cauchy
integral at $0$ and using H\"older's inequality gives, for every $j\ge1$,
\[
 \left|W_{t,f}^{(j)}(0)-\mathbf1_{\{j=1\}}\right|
 \le C_j\|q_{t,f}\|_p\le 2CC_j|t|.
\]
After decreasing $\delta$, $|W_{t,f}'(0)|\ge1/2$. Since $a_t$ is the
identity near $0$, the normalized transformed inner map has the same germ
as $W_{t,f}(f(z))/W_{t,f}'(0)$. The Neumann series is holomorphic in $t$.
The last estimate and the univalent coefficient bounds give a common
complex time disk on which each transformed coefficient is holomorphic
and uniformly bounded. The conjugate coefficient is complexified by
$t\mapsto\overline{a_j(\bar t)}$. The same construction works for finitely
many parameters.

For $P\in\mathscr P$ define, for small real $t$,
\[
 Q_tP(h)=e^{-\int_0^t \mathfrak d_v(he_s^v)\dd s}P(he_t^v).
\]
The same time disk applies to $\mathfrak d_v$, a finite coefficient polynomial. On a
smaller disk, $|Q_tP|\le C_Pe^{rM_v}$ uniformly over normalized initial
maps. Cauchy's estimate gives a sample-independent bound on each Taylor
coefficient. The surgery construction is independent of the extension:
comparing two extensions with the same analytic boundary values gives
conformal maps on both complementary domains that agree on a boundary
collar; after normalization the comparison is the identity. It also gives
the composition law. The composition law gives $Q_{t+s}P=Q_t(Q_sP)$. Here $A_v-\mathfrak d_v$
means $P\mapsto A_vP-\mathfrak d_vP$, with $\mathfrak d_v$ acting by multiplication. In
particular,
\[
 \begin{aligned}
 \left.\partial_tQ_tP\right|_{t=0}
 &=A_vP-\mathfrak d_vP,\\
 \left.\partial_t^2Q_tP\right|_{t=0}
 &=A_v^2P-2\mathfrak d_vA_vP+(\mathfrak d_v^2-A_v\mathfrak d_v)P.
 \end{aligned}
\]
The second identity uses the product rule for $A_v$. Differentiating the
composition law repeatedly gives, for every $j\ge0$,
\[
 \left.\partial_t^jQ_tP\right|_{t=0}=(A_v-\mathfrak d_v)^jP.
\]
Set $L=A_v-\mathfrak d_v$. We use the integration by parts identities in
Lemma~\ref{lem:mixed-ibp}. For $R\in\mathscr P$ they give
\[
 \m E^{\mathrm{weld}}[D_nR]=\m E^{\mathrm{weld}}[B_nR],\qquad
 \m E^{\mathrm{weld}}[\bar D_nR]=\m E^{\mathrm{weld}}[\bar B_nR].
\]
Also, rotation invariance gives $\m E^{\mathrm{weld}}[A_0R]=0$. Using the
definitions of $A_v$ and $\mathfrak d_v$, and the fact that $v$ is a finite Fourier
series, we obtain
\begin{equation}\label{eq:flow-ibp}
 \begin{aligned}
 \m E^{\mathrm{weld}}[A_vR]
 &=\m E^{\mathrm{weld}}\!\left[
 \sum_{n\ge1}(-\ii\bar v_nB_n+\ii v_n\bar B_n)R\right]\\
 &=\m E^{\mathrm{weld}}[\mathfrak d_vR],\qquad
 \m E^{\mathrm{weld}}[LR]=0.
 \end{aligned}
\end{equation}
Each $D_n$ and $\bar D_n$ maps $\mathscr P$ to itself by the coefficient
formulas above; the same holds for $A_0$.
Moreover $\mathfrak d_v\in\mathscr P$. Hence $L$ maps $\mathscr P$ to itself. For
$j\ge1$ we may take $R=L^{j-1}P$ in \eqref{eq:flow-ibp}, which gives
\[
 \m E^{\mathrm{weld}}[L^jP]
 =\m E^{\mathrm{weld}}[A_v(L^{j-1}P)]
 -\m E^{\mathrm{weld}}[\mathfrak d_vL^{j-1}P]=0.
\]
The Taylor expansion and the sample-independent Cauchy bound above then give
\[
 \m E^{\mathrm{weld}}[Q_tP]
 =\sum_{j\ge0}\frac{t^j}{j!}\m E^{\mathrm{weld}}[L^jP]
 =\m E^{\mathrm{weld}}[P].
\]
The Cauchy bound, rather than the IBP identity, justifies the exchange of
expectation and sum. Mixed coefficient polynomials determine finite measures on
the compact normalized univalent class. The independent uniform phase in
\eqref{eq:law} then gives the same identity for the full welding law.
Iteration proves \eqref{eq:finite-flow} for finite $T$. The divergence is a
bounded polynomial for finite Fourier $v$, giving the two-sided bounds and
the stated derivative at zero.
\end{proof}

For a quasisymmetric circle homeomorphism $\chi$, let $C_\chi u=u\circ\chi$
on $\dot H^{1/2}(\m S^1;\m C)$, with constants factored out. Denote its positive and negative frequency subspaces by $\mc H_+$ and $\mc H_-$. In the bases
$z^n/\sqrt n$ and $z^{-n}/\sqrt n$, write
\begin{equation}\label{eq:block}
 C_\chi=\begin{pmatrix}a&b\\\bar b&\bar a\end{pmatrix},\qquad
 X_\chi=a^{-1}b.
\end{equation}
These $a,b$ are composition blocks and are not the Grunsky operators in
\eqref{eq:TT}. The indefinite form $\|u_+\|^2-\|u_-\|^2$ is preserved. Therefore
\[
 aa^*-bb^*=I,\qquad a^*a-\bar b^*\bar b=I.
\]
The first identity makes $a$ surjective; the second gives
$\|ax\|\ge\|x\|$. Hence $a$ is invertible and
\[
 I-X_\chi X_\chi^*=a^{-1}(a^*)^{-1}
 \ge\|a\|^{-2}I,\qquad \|X_\chi\|<1.
\]
The off-diagonal block $b$ is Hilbert--Schmidt precisely when
$\chi\in\WP(\m S^1)$; see \cite[Chapter~II, Section~5]{TT06}.
For analytic $\chi$ it is trace class, as the following Fourier estimate
shows. Let $\widetilde\chi$ be an angle lift. In the bases above,
\[
 b_{mn}=\sqrt{\frac mn}\frac1{2\pi}\int_0^{2\pi}
 e^{-\ii n\widetilde\chi(\theta)}e^{-\ii m\theta}\,\dd\theta.
\]
For a sufficiently small $\delta>0$, the analytic extension satisfies
$\operatorname{Im}\widetilde\chi(\theta-\ii\delta)\le-\delta'$ for
some $\delta'>0$. Shifting the contour to $\theta-\ii\delta$ gives
\[
 |b_{mn}|\le\sqrt{\frac mn}e^{-m\delta-n\delta'},\qquad
 \sum_{m,n\ge1}|b_{mn}|<\infty.
\]
Thus $b$ and $X_\chi=a^{-1}b$ are trace class.

\begin{lem}\label{lem:grunsky-graph}
For quasisymmetric $h$, the composition and Grunsky operators satisfy
the classical graph identity \cite[Chapter~II, Section~5]{TT06}:
\begin{equation}\label{eq:graph}
 G_1(h)=b_h\bar a_h^{-1},\qquad
 C_h\mc H_- =\{(G_1(h)x,x):x\in\ell^2\}.
\end{equation}
Let $Z_h=\overline{G_1(h)}$. If $h,\chi$ are quasisymmetric, or if
$\chi$ is analytic and $h$ is a Jordan welding, then
\begin{equation}\label{eq:mobius}
 Z_{h\chi}=(\bar b+\bar a Z_h)(a+bZ_h)^{-1}.
\end{equation}
\end{lem}
\begin{proof}
For quasisymmetric weldings, \eqref{eq:graph} is the period-mapping
formula cited above. We recall its coefficient form for analytic $h$
to fix the conventions; the extension to an arbitrary Jordan welding
will not use a composition operator for that welding.
Define $P_n$ by
\[
 -\log(1-f(z)/w)+\log(f(z)/z)
 =\sum_{n\ge1}P_n(1/w)z^n/n.
\]
Then
\[
 P_n(1/f(\zeta))=\zeta^{-n}+n\sum_{m\ge1}b_{-n,-m}\zeta^m.
\]
To prove density, truncate the exterior Dirichlet expansion in the exterior
map coordinate. Each truncated function extends across the analytic
boundary. After $w\mapsto1/w$, approximate it by polynomials uniformly on
a slightly larger bounded Jordan domain; Cauchy estimates give uniform
convergence of derivatives. Thus the Dirichlet energies converge, and the
Faber functions span a dense set.

For quasisymmetric $h,\chi$, apply $C_\chi$ to the conjugate graph
$\{(x,Z_hx):x\in\ell^2\}$. The image is
\[
 \left\{\bigl((a+bZ_h)x,(\bar b+\bar aZ_h)x\bigr):x\in\ell^2\right\}.
\]
Since $a+bZ_h=a(I+X_\chi Z_h)$ and $\|X_\chi\|<1$, the first component is
invertible. This proves \eqref{eq:mobius} for quasisymmetric $h,\chi$.

For a Jordan welding, let $f_r(z)=f(rz)/r$ and $D_r e_n=r^ne_n$.
The Grunsky expansion gives
\[
 G_1(f_r)=D_rG_1(f)D_r,\qquad
 G_1(f_r)^*=D_rG_1(f)^*D_r.
\]
As $D_r\to I$ strongly, $Z_r=\overline{G_1(f_r)}$ and $Z_r^*$ converge
strongly to $Z=\overline{G_1(f)}$ and $Z^*$, respectively. The Grunsky
inequality gives $\|Z_r\|,\|Z\|\le1$. Put
$M_r=a+bZ_r$, $M=a+bZ$, and $X=a^{-1}b$. Then
\[
 M_r^{-1}=(I+XZ_r)^{-1}a^{-1},\qquad
 \sup_r\|M_r^{-1}\|\le\frac{\|a^{-1}\|}{1-\|X\|}.
\]
The resolvent identity
$M_r^{-1}-M^{-1}=M_r^{-1}b(Z-Z_r)M^{-1}$ shows that
$M_r^{-1}\to M^{-1}$ strongly. Hence
\[
 (\bar b+\bar aZ_r)M_r^{-1}
 \longrightarrow(\bar b+\bar aZ)M^{-1}\qquad\text{strongly}.
\]

It remains to identify this limit with the transformed inner map. Choose a
smooth quasiconformal extension $a_\mu$ of $\chi^{-1}$ whose coefficient
$\mu$ is compactly supported in $\m D\setminus\{0\}$. Write
\[
 \lambda_r(f_r(z))=\mu(z)\frac{f_r'(z)}{\overline{f_r'(z)}},\qquad
 \lambda(f(z))=\mu(z)\frac{f'(z)}{\overline{f'(z)}}.
\]
Extend both coefficients by zero. On a compact neighborhood of
$\operatorname{supp}\mu$, the maps $f_r,f_r'$ converge uniformly to
$f,f'$, and $f'$ is bounded away from zero. Their inverse maps converge on
compact subsets of $f(\m D)$. Koebe bounds put the supports of $\lambda_r$
in a fixed compact subset of $\m C\setminus\{0\}$. Consequently, for every
finite $p$,
\[
 \lambda_r\to\lambda\quad\text{in }L^p(\m C),\qquad
 \sup_r\|\lambda_r\|_\infty\le\|\mu\|_\infty<1.
\]
Choose $p>2$ with $\|\mu\|_\infty\|\mathcal B\|_{L^p\to L^p}<1$, and set
$q_r=(I-M_{\lambda_r}\mathcal B)^{-1}\lambda_r$,
$q=(I-M_\lambda\mathcal B)^{-1}\lambda$. Then
\[
 q_r-q=(I-M_{\lambda_r}\mathcal B)^{-1}
 \bigl[(\lambda_r-\lambda)+(\lambda_r-\lambda)\mathcal Bq\bigr]
 \longrightarrow0\quad\text{in }L^p.
\]
The normalized solutions $W_r(z)=z+\mathcal Cq_r(z)-\mathcal Cq_r(0)$
therefore converge locally uniformly to $W$. The inner surgery maps
$W_rf_ra_\mu^{-1}$, normalized at $0$, converge locally uniformly to
$Wfa_\mu^{-1}$. Their Grunsky coefficients converge entrywise. Applying
the analytic formula to $f_r$ and taking the strong limit proves
\eqref{eq:mobius}. No composition operator for the rough $h$ is used.
\end{proof}

\begin{lem}[The relative determinant and the analytic density]\label{lem:analytic-cocycle}
For analytic $\chi$ and any Jordan welding $h$, or for $h,\chi\in\WP(\m S^1)$,
set
\begin{equation}\label{eq:geom}
 D(h,\chi)=\tfrac12\log\det(I+bb^*)+
 \operatorname{Re}\log\det(I+X_\chi\overline{G_1(h)}).
\end{equation}
If $\chi,\omega$ are analytic or $h,\chi,\omega\in\WP(\m S^1)$, then
$D(h,\chi\omega)=D(h,\chi)+D(h\chi,\omega)$. If
$\varphi=\chi^{-1}$, then $X_\chi=-G_1(\varphi)$ and
\begin{equation}\label{eq:geom-inverse}
 D(h,\chi)=-\tfrac12\log\det(I-G_1(\varphi)^*G_1(\varphi))
 +\operatorname{Re}\log\det(I-G_1(\varphi)\overline{G_1(h)}).
\end{equation}
For $h,\chi\in\WP(\m S^1)$,
\begin{equation}\label{eq:rn-energy-cocycle}
 D(h,\chi)=\frac{I^{\mathrm L}(h\chi)-I^{\mathrm L}(h)}{24}.
\end{equation}
For analytic $\chi$, with $\Delta_\chi K(h)=K(h\chi)-K(h)$,
\begin{equation}\label{eq:analytic-rn}
 p_\chi(h)=\exp\{-c_{\mathrm L}D(h,\chi)-2\Delta_\chi K(h)\}.
\end{equation}
\end{lem}
\begin{proof}
First let $h,\chi\in\WP(\m S^1)$, put $Z=\overline{G_1(h)}$,
$Z'=Z_{h\chi}$ and $M=a+bZ=a(I+X_\chi Z)$. Both $b$ and $Z$ are
Hilbert--Schmidt. Thus $X_\chi Z$ is trace class, and all the
Fredholm determinants below are defined. The graph transformation and
symplectic block relations give
\[
 I-Z'^*Z'=M^{-*}(I-Z^*Z)M^{-1},\qquad
 \det(M^*M)=\det(I+bb^*)\,|\det(I+X_\chi Z)|^2.
\]
For the second identity, use $a^*a-I$ trace class and
$\det(a^*a)=\det(aa^*)=\det(I+bb^*)$. Taking determinants in the first
identity therefore gives
\[
 D(h,\chi)=-\tfrac12\log\det(I-Z'^*Z')+
 \tfrac12\log\det(I-Z^*Z).
\]
Equation~\eqref{eq:rn-loewner-energy} proves
\eqref{eq:rn-energy-cocycle}, including its sign and factor $1/24$.
It also proves the cocycle on WP weldings.

For analytic $\chi$ and a general Jordan welding $h$, use the radial
approximation in Lemma~\ref{lem:grunsky-graph}.
If $T$ is trace class and $Z_r,Z_r^*$ converge strongly to $Z,Z^*$
with uniformly bounded norms, finite-rank approximation gives
$\|T(Z_r-Z)\|_{\tr}\to0$. In particular,
\[
 \|X_\chi(Z_r-Z)\|_{\tr}\longrightarrow0,\qquad
 \det(I+X_\chi Z_r)\longrightarrow\det(I+X_\chi Z).
\]
The same argument applies after a fixed analytic composition. Since
$\|X_\chi Z_r\|\le\|X_\chi\|<1$, the logarithm is the branch continuous
from $I$. Thus the determinant cocycle passes to the limit. Taking
$h=\varphi$ in
\eqref{eq:mobius} gives $X_\chi=-G_1(\varphi)$; the block identity gives
the determinant term in \eqref{eq:geom-inverse}.

Let $\chi=e_t^v$, with angle flow $\theta\mapsto\theta+tv(\theta)+O(t^2)$.
For the negative basis element $z^{-n}/\sqrt n$, the first-order change is
$-\ii\sqrt n\,z^{-n}v(\theta)$. Its coefficient in the positive basis
$z^m/\sqrt m$ is $-\ii\sqrt{mn}v_{m+n}$. Hence
$(b'_0)_{mn}=-\ii\sqrt{mn}v_{m+n}$.

Differentiate the inner Grunsky expansion in both variables and then take
the diagonal. Using
\[
 \lim_{w\to z}\left(
 \frac{f'(z)f'(w)}{(f(z)-f(w))^2}-\frac1{(z-w)^2}\right)
 =\frac16S_f(z)
\]
gives
\begin{equation}\label{eq:schwarzian-b1}
 S_f(z)=-6\sum_{m,n\ge1}\sqrt{mn}\,(G_1(h))_{mn}z^{m+n-2}.
\end{equation}
At $t=0$, $a=I$ and $b=X=0$, so the derivative of the first determinant
in \eqref{eq:geom} is zero and $X'_0=b'_0$. Therefore
\[
 \left.\frac\dd{\dd t}\right|_{t=0}c_{\mathrm L}D(h,e_t^v)
 =c_{\mathrm L}\operatorname{Im}
 \sum_{m,n\ge1}\sqrt{mn}v_{m+n}\overline{(G_1(h))_{mn}}.
\]
For finite Fourier $v$ the sum is finite. Since $c_{\mathrm L}=12c$,
\eqref{eq:schwarzian-b1} shows that this is exactly the contribution of
$cS_f$ to $\mathfrak d_v$. For the capacity part, let
$\omega(\xi)^{1-n}f'(\omega(\xi))=(V_nf)(\xi)+\sum_{j\ge0}r_{n,j}\xi^j$.
Since $V_nf$ has only negative powers, changing variables $\xi=f(z)$ gives
\[
 r_{n,1}=\operatorname{Res}_{\xi=0}
 \frac{\omega(\xi)^{1-n}f'(\omega(\xi))}{\xi^2}\,\dd\xi
 =\operatorname{Res}_{z=0}\frac{z^{1-n}f'(z)^2}{f(z)^2}\,\dd z
 =[z^{n-2}]\left(\frac{f'}f\right)^2.
\]
Also $(W_nf)(z)=-\sum_{j\ge0}r_{n,j}f(z)^j$, so
$-(W_nf)'(0)=r_{n,1}$. The loop vector field has this linear coefficient at
infinity. Its complex derivative of $\log|g'(\infty)|$ is half that
coefficient; the normalized inner scale does not change. Thus
\[
 2D_nK=[z^{n-2}]\left(\frac{f'}f\right)^2,\qquad
 2A_vK=2\operatorname{Im}\sum_{n\ge1}\bar v_n
 [z^{n-2}]\left(\frac{f'}f\right)^2.
\]
The term $z^{-2}$ has no coefficient of $z^{n-2}$ for $n\ge1$. Thus
\[
 \left.\frac\dd{\dd t}\right|_{t=0}
 [c_{\mathrm L}D(h,e_t^v)+2\Delta_{e_t^v}K(h)]=\mathfrak d_v(h).
\]
Together with the cocycle and Lemma~\ref{lem:finite-flow}, this proves
\eqref{eq:analytic-rn} for finite Fourier flows and their finite products.

For an arbitrary analytic $\chi$, connect its angle lift to the identity
by a linear path. Approximate its right velocity in a common narrower
analytic strip by finite Fourier fields and time steps. Complex ODE
stability gives convergence of the maps and inverses in that strip. The
contour estimate is uniform on a smaller strip, so entrywise convergence
gives $b_j\to b$ in trace norm. By the difference-quotient norm and
density of smooth functions, $a_j\to a$ strongly. Since $\|a_j^{-1}\|\le1$,
also $a_j^{-1}\to a^{-1}$ strongly. Finite-rank approximation then gives
$X_j\to X$ in trace norm, hence $\sup_j\|X_j\|<1$.
Interpolate each analytic error $\chi_j^{-1}\chi$ with the identity
using a fixed radial cutoff. These extensions converge in $C^1$ to the
identity. Composing them with a fixed smooth extension of $\chi^{-1}$
gives Beltrami coefficients converging in $L^\infty$, supported in a
common annulus away from $0$ and $\m S^1$. Koebe bounds and the
logarithmic derivative identity below give pointwise convergence and a
uniform bound for $\Delta K$.
Thus the geometric action converges pointwise and is uniformly bounded;
dominated convergence identifies the limiting density. In particular its
exponential has expectation one.
\end{proof}

\begin{remark}
    The proof above derives \eqref{eq:analytic-rn} from
Lemmas~\ref{lem:finite-flow} and \ref{lem:grunsky-graph}, without using
\cite[Corollary~1.3]{BJ25}. That corollary also gives a Radon--Nikodym
formula for analytic circle diffeomorphisms, expressed in terms of the
relative universal Liouville action and the capacity. When the two
formulas are applied to the same action with matching conventions, they
agree almost surely by uniqueness of the Radon--Nikodym derivative.
\end{remark}

\subsection{Resolvent identities and operator bounds}\label{subsec:rn-resolvents}

We rewrite the two analytic contributions before taking a WP limit.
The Grunsky operator splits into a linear Hilbert--Schmidt term and a
quadratic trace-class remainder. The logarithmic derivative of the
normalized Beltrami solution has the same linear-plus-remainder form.
Combining the two linear terms gives the $L^2$ variable
$\mathfrak d_{v_\mu}$ in \eqref{body:eq:action}.

\begin{lem}[Grunsky and conformal-radius resolvents]\label{lem:resolvents}
For the representative in \eqref{eq:mu},
\begin{equation}\label{eq:B1-resolvent}
 G_1(\varphi)=\bar J^*(I-M_\eta\mathcal B_*)^{-1}M_\eta J
 =\bar J^*M_\eta J+Q_\eta,
\end{equation}
and
\begin{equation}\label{eq:bounds-B}
 \|G_1(\varphi)\|\le k,\qquad
 \|G_1(\varphi)\|_{\HS}\le\frac{N(\mu)}{1-k},\qquad
 \|Q_\eta\|_{\tr}\le\frac{N(\mu)^2}{1-k}.
\end{equation}
If $\lambda$ is bounded, supported in a bounded set away from $0$, and
$\|\lambda\|_\infty<1$, let $W^\lambda$ be its normalized quasiconformal
solution with $W^\lambda(0)=0$ and $(W^\lambda)'(\infty)=1$. Then
\begin{equation}\label{eq:log-derivative}
 \log(W^\lambda)'(0)=-\frac1\pi\int_{\m C}\frac{\lambda(z)}{z^2}\dd A(z)
 -\mathcal R(\lambda),\qquad
 |\mathcal R(\lambda)|\le
 \frac{\|\lambda/z\|_2^2}{\pi(1-\|\lambda\|_\infty)}.
\end{equation}
The logarithm is continued from $\lambda=0$.
\end{lem}
\begin{proof}
Let $W^\mu$ have coefficient $\mu$ inside $\m D$ and zero outside. The pair
$(W^\mu a_\mu^{-1},W^\mu)$ has welding $\varphi^{-1}$. Reflection
$z\mapsto1/\bar z$ exchanges its sides, and the inner map of the reflected
pair has welding $\varphi$ and exterior coefficient $\eta$. If $V$ is
the reflected quasiconformal solution and $\zeta\in\m D$, set
$u_\zeta(\xi)=V'(\zeta)/(V(\xi)-V(\zeta))$. The $\bar\partial$ derivative
of $u_\zeta-(\xi-\zeta)^{-1}$ is supported in $\m D^*$ and equals
\[
 -(I-M_\eta\mathcal B_*)^{-1}M_\eta(\xi-\zeta)^{-2}.
\]
Put $R_\eta=(I-M_\eta\mathcal B_*)^{-1}$. The function $u_\zeta$
decays at infinity, so the Cauchy transform gives
\[
 \partial_\zeta\log\frac{V(\xi)-V(\zeta)}{\xi-\zeta}
 =\mathcal C\!\left[R_\eta M_\eta(\,\cdot-\zeta)^{-2}\right](\xi),
 \qquad \xi,\zeta\in\m D.
\]
For $|w|>1$ and $\xi,\zeta\in\m D$, expand
\[
 \frac1{\xi-w}=-\sum_{m\ge1}\xi^{m-1}w^{-m},\qquad
 \frac1{(w-\zeta)^2}=\sum_{n\ge1}n\zeta^{n-1}w^{-n-1}.
\]
Comparing the coefficient of $\xi^m\zeta^{n-1}$ with the inner
Grunsky expansion in \eqref{eq:grunsky-coefficients} gives
\[
 \sqrt{mn}\,b_{-m,-n}
 =\left(\bar J^*R_\eta M_\eta J\right)_{mn}.
\]
The expansions hold locally uniformly in the disk; the corresponding
operator identity follows from the $L^2$ resolvent bound. The relation
$R_\eta=I+M_\eta\mathcal B_*R_\eta$ then gives
\[
 G_1(\varphi)=\bar J^*M_\eta J+
 (\bar J^*M_\eta)\mathcal B_*R_\eta(M_\eta J).
\]
The first term is linear in $\eta$ and the second is $Q_\eta$, containing
at least two factors of $\eta$. This proves \eqref{eq:B1-resolvent}.
The identity
$\sum_{n\ge1}|\sqrt{n/\pi}z^{-n-1}|^2=1/[\pi(|z|^2-1)^2]$ gives
$\|M_\eta J\|_{\HS}=\|\bar J^*M_\eta\|_{\HS}=N(\mu)$.
The trace-class and Hilbert--Schmidt estimates follow by the resolvent bound
$(1-k)^{-1}$.

For the operator norm, let $P$ map the $n$-th basis vector to
$\sqrt{n/\pi}\bar z^{n-1}\mathbf1_\m D$. Then $\mathcal B P=-J$ and
$\mathcal B_*^*J=0$. Also
$\|\mathbf1_\m D\mathcal B q\|_2=\|\bar J^*q\|_{\ell^2}$ for $q$ supported in
$\m D^*$. Put
$q=(I-M_\eta\mathcal B_*)^{-1}M_\eta Jx$. Then
$q=M_\eta(Jx+\mathcal B_*q)$ and $G_1(\varphi)x=\bar J^*q$. Since
$Jx\perp\mathcal B_*q$ and $\mathcal B$ is unitary on $L^2(\m C)$,
\[
 \begin{aligned}
 \|Jx+\mathcal B_*q\|_2^2
 &=\|x\|^2+\|\mathcal B_*q\|_2^2\\
 &=\|x\|^2+\|q\|_2^2-\|\bar J^*q\|_{\ell^2}^2.
 \end{aligned}
\]
The multiplication bound $\|M_\eta\|\le k$ therefore gives
\[
 (1-k^2)\|q\|_2^2+k^2\|G_1(\varphi)x\|^2\le k^2\|x\|^2.
\]
As $\|G_1(\varphi)x\|=\|\bar J^*q\|\le\|q\|_2$, the left side is at
least $\|G_1(\varphi)x\|^2$. Thus $\|G_1(\varphi)\|\le k$.

For \eqref{eq:log-derivative}, the nonvanishing function $W^\lambda(z)/z$
has winding number zero. Let $U=\log(W^\lambda/z)$ be the branch tending
to zero at infinity. Near $0$ it is analytic; on compact annuli its
weak gradient is square integrable; at infinity its gradient is
$O(|z|^{-2})$. Hence $\nabla U\in L^2(\m C)$. The Beltrami equation implies
\[
 \bar\partial U=\lambda/z+\lambda\mathcal B(\bar\partial U),\qquad
 U=\mathcal C(\bar\partial U).
\]
Evaluating at $0$ and expanding the resolvent once gives
\eqref{eq:log-derivative}. The remainder bound follows from the $L^2$
isometry of $\mathcal B$ and Cauchy--Schwarz.
\end{proof}

Koebe's estimate gives, for normalized univalent $f$,
\begin{equation}\label{eq:koebe}
 \left|\frac{f'(z)}{f(z)}\right|
 \le\frac{(1+|z|)^2}{|z|(1-|z|^2)},\quad
 \left\|\frac{\lambda_f}{z}\right\|_2^2\le
 \frac{16\pi N(\mu)^2}{r_0^2},\quad
 \left\|\frac\mu z\right\|_2^2\le\frac{\pi N(\mu)^2}{r_0^2}.
\end{equation}
The estimates are uniform in $f$.

\begin{lem}[Identification for analytic deformations]\label{lem:analytic-action}
If $\mu$ is also supported away from $\m S^1$, then, for every Jordan
welding $h$,
\begin{equation}\label{eq:analytic-identification}
 \mathscr I_\mu(\varphi,h)=c_{\mathrm L}D(h,\varphi^{-1})
 +2\bigl[K(h\varphi^{-1})-K(h)\bigr].
\end{equation}
\end{lem}
\begin{proof}
Recall that $\Delta_\chi K(h)=K(h\circ\chi)-K(h)$; this is a
composition difference, not a Laplacian. Write $W=W^{\lambda_f}$.
The surgery pair for $h\circ\varphi^{-1}$ is
$(W\circ f\circ a_\mu^{-1},W\circ g)$. Since $f'(0)=1$ and
$W'(\infty)=1$,
\[
 (W\circ f\circ a_\mu^{-1})'(0)=\frac{W'(0)}{a_\mu'(0)},\qquad
 (W\circ g)'(\infty)=g'(\infty).
\]
Using $K(h)=\log|g'(\infty)|-\log|f'(0)|$ gives
\[
 \Delta_{\varphi^{-1}}K(h)
 =K(h\circ\varphi^{-1})-K(h)
 =\log|a_\mu'(0)|-\log|W'(0)|.
\]
The change of variables $w=f(z)$ gives
\[
 \frac1\pi\int_{\m C}\frac{\lambda_f(w)}{w^2}\,\dd A(w)
 =\frac1\pi\int_{\m D}\mu(z)\left(\frac{f'(z)}{f(z)}\right)^2\dd A(z).
\]
Taking the original pair to be the circle and using
\eqref{eq:log-derivative} gives
\begin{equation}\label{eq:capacity-linear}
 \begin{aligned}
 \Delta_{\varphi^{-1}}K(h)-K(\varphi^{-1})
 ={}&\frac1\pi\operatorname{Re}\int_{\m D}\mu(z)
 \left[\left(\frac{f'(z)}{f(z)}\right)^2-z^{-2}\right]\dd A(z)\\
 &+\operatorname{Re}\bigl[\mathcal R(\lambda_f)-\mathcal R(\mu)\bigr].
 \end{aligned}
\end{equation}
Let $H_\eta=\bar J^*M_\eta J$. Reflection gives
$(H_\eta)_{mn}=\sqrt{mn}\,\overline{m_{m+n}(\mu)}$. Thus
\eqref{eq:schwarzian-b1} implies
\[
 -c_{\mathrm L}\operatorname{Re}\tr(H_\eta\overline{G_1(h)})
 =\frac{2c}{\pi}\operatorname{Re}\int_{\m D}\mu S_f\dd A.
\]
Put $G=G_1(\varphi)$ and $Z=\overline{G_1(h)}$. Since $\mu$ is
supported away from the circle, $G$ is trace class. Equations
\eqref{eq:det2} and \eqref{eq:B1-resolvent} give
\[
 \log\det(I-GZ)
 =\log\det{}_2(I-GZ)-\tr(H_\eta Z)-\tr(Q_\eta Z).
\]
Moreover, by \eqref{eq:stress} and \eqref{eq:vmu},
\[
 \begin{aligned}
 \frac2\pi\operatorname{Re}\int_{\m D}\mu\mathcal T_f\,\dd A
 &=2\operatorname{Re}\sum_{n\ge1}m_n(\mu)B_n(f)\\
 &=2\operatorname{Im}\sum_{n\ge1}\overline{(v_\mu)_n}B_n(f)
 =\mathfrak d_{v_\mu}(h).
 \end{aligned}
\]
Substitute these two identities and \eqref{eq:capacity-linear} into
\eqref{eq:geom-inverse}. The linear terms combine into $\mathfrak d_{v_\mu}$, and
the remaining terms give \eqref{body:eq:action}.
\end{proof}

\subsection{Moment estimates and the WP extension}\label{subsec:rn-extension}

We now carry out the extension to a general fixed WP map. The decisive
estimate is a quadratic bound for the nonlinear remainder, uniform in
the welding. It has two uses: radial cutoffs give convergence of the
actions, and comparison with a real multiple of the Beltrami coefficient
gives exponential moments. These moments prevent loss of mass when we
pass to the limiting density.

\begin{lem}[A linear term and a quadratic remainder]\label{lem:bounds}
For $\mu$ as in \eqref{eq:mu}, $v_\mu\in H^{3/2}(\m S^1)$ and
\begin{equation}\label{eq:v-bound}
 [v_\mu]_{3/2}^2\le
 \bigl(16+4\log(1/r_0)\bigr)N(\mu)^2.
\end{equation}
Furthermore,
\begin{equation}\label{eq:action-split}
 \mathscr I_\mu(\varphi,h)=C_\mu+\mathfrak d_{v_\mu}(h)+E_\mu(h),
\end{equation}
where
\begin{equation}\label{eq:constant}
 C_\mu=-\frac{c_{\mathrm L}}2\log\det(I-G_1(\varphi)^*G_1(\varphi))
 +2K(\varphi^{-1}),
\end{equation}
and $E_\mu$ is the sum of the last three terms of \eqref{body:eq:action}.
The remainder satisfies
\begin{equation}\label{eq:rn-quadratic-remainder}
 \|E_\mu\|_\infty
 \le C_{\kappa,r_0}\frac{N(\mu)^2}{(1-k)^3},
\end{equation}
where the bound is uniform over normalized univalent inner maps.
\end{lem}
\begin{proof}
From $(H_\eta)_{mn}=\sqrt{mn}\,\overline{m_{m+n}}$ and
$\sum_{m=1}^{n-1}m(n-m)=(n^3-n)/6$,
\[
 \|H_\eta\|_{\HS}^2
 =\sum_{m,\ell\ge1}m\ell|m_{m+\ell}|^2
 =\frac16\sum_{n\ge2}(n^3-n)|m_n|^2\le N(\mu)^2.
\]
Cauchy--Schwarz, with $\mu=0$ for $|z|<r_0$, gives
\[
 |m_1|^2\le N(\mu)^2\frac1\pi\int_{r_0<|z|<1}
 \frac{(1-|z|^2)^2}{|z|^2}\,\dd A(z)
 \le2\log(1/r_0)N(\mu)^2.
\]
Since $n^3\le\frac43(n^3-n)$ for $n\ge2$,
\[
 [v_\mu]_{3/2}^2
 =2|m_1|^2+2\sum_{n\ge2}n^3|m_n|^2
 \le\bigl(4\log(1/r_0)+16\bigr)N(\mu)^2.
\]

By \eqref{eq:bounds-B}, $G_1(\varphi)$ is Hilbert--Schmidt and has norm
at most $k<1$. Thus
\[
 0\le-\log\det(I-G_1(\varphi)^*G_1(\varphi))
 \le\frac{\|G_1(\varphi)\|_{\HS}^2}{1-k^2}.
\]
For $L=G_1(\varphi)\overline{G_1(h)}$, the Grunsky inequality gives
$\|L\|\le k$ and $\|L\|_{\HS}\le N(\mu)/(1-k)$. Consequently
\[
 |\log\det{}_2(I-L)|
 \le\frac{\|L\|_{\HS}^2}{2(1-k)}.
\]
The trace term is bounded by $\|Q_\eta\|_{\tr}$. Finally,
\eqref{eq:log-derivative} and \eqref{eq:koebe} bound the resolvent terms
uniformly in $f$. More explicitly, the three remainders in
\eqref{body:eq:action} are bounded respectively by
\[
 \frac{c_{\mathrm L}N(\mu)^2}{2(1-k)^3},\qquad
 \frac{c_{\mathrm L}N(\mu)^2}{1-k},\qquad
 \frac{34N(\mu)^2}{r_0^2(1-k)}.
\]
This proves \eqref{eq:rn-quadratic-remainder}.
Equation~\eqref{eq:divergence-moment} controls the remaining linear term.
\end{proof}

We now approximate the Beltrami coefficient while retaining the same
bound $\|\mu_j\|_\infty\le k<1$. This common bound is essential for the
moment estimate. Choose Borel representatives for all Beltrami coefficients. Let
$0\le\vartheta_j\le1$ be radial cutoffs supported away from $\m S^1$ and
equal to one eventually on each compact subset of $\m D$. Put
$\mu_j=\vartheta_j\mu$ and let $\varphi_j$ be the corresponding boundary
map, using the reflected solution normalized at $0,\infty$ and by
$\varphi_j(1)=\varphi(1)$. Then $\varphi_j$ is analytic.

\begin{lem}[Moments and convergence under radial cutoffs]\label{lem:approx}
For every real $r$ with $|r|k<1$,
\begin{equation}\label{eq:amplification}
 \sup_j\m E^{\mathrm{weld}}\exp\{-r\mathscr I_{\mu_j}(\varphi_j,\widehat H)\}<\infty.
\end{equation}
The maps $\varphi_j$ and their inverses converge uniformly to $\varphi$ and
$\varphi^{-1}$, respectively, and
\begin{equation}\label{eq:action-L2}
 \mathscr I_{\mu_j}(\varphi_j,\cdot)\longrightarrow
 \mathscr I_\mu(\varphi,\cdot)\quad\text{in }L^2(\m P^{\mathrm{weld}}).
\end{equation}
\end{lem}
\begin{proof}
We first prove the moment estimate. Fix a real $r$ with $|r|k<1$.
The reflected coefficients for $\mu_j$ and $r\mu_j$ have common
bounded supports, converge in $L^p$ for $p>2$ sufficiently close to $2$,
and have uniform dilatation bounds below one. Their normalized
quasiconformal solutions therefore converge locally uniformly,
including their conformal derivatives at $0$ and $\infty$.
Thus $\varphi_j,\varphi_j^{-1}$ converge uniformly to
$\varphi,\varphi^{-1}$, and the capacities converge. The same
conclusion holds for the boundary maps determined by $r\mu_j$.

The capacity limits and the determinant bound in
Lemma~\ref{lem:bounds} give
\[
 \sup_j\bigl(|C_{\mu_j}|+|C_{r\mu_j}|\bigr)<\infty.
\]
The remainder bounds in that lemma are also uniform in $j$, since
$N(\mu_j)\le N(\mu)$ and $\|\mu_j\|_\infty\le k$.
For real $r$, $v_{r\mu_j}=rv_{\mu_j}$. Hence
\[
 \begin{aligned}
 \mathscr I_{r\mu_j}-r\mathscr I_{\mu_j}
 &=(C_{r\mu_j}-rC_{\mu_j})+(E_{r\mu_j}-rE_{\mu_j}),\\
 \sup_j\|\mathscr I_{r\mu_j}-r\mathscr I_{\mu_j}\|_{L^\infty(\m P^{\mathrm{weld}})}
 &\le A_r<\infty.
 \end{aligned}
\]
Here $\mathscr I_{r\mu_j}$ uses the boundary map determined by
$r\mu_j$, with the same random welding $h$ in both actions. Since $r\mu_j$ is supported away from the circle,
the analytic density formula gives
$\m E^{\mathrm{weld}}[e^{-\mathscr I_{r\mu_j}}]=1$.
Consequently,
\[
 \m E^{\mathrm{weld}}[e^{-r\mathscr I_{\mu_j}}]
 \le e^{A_r}\m E^{\mathrm{weld}}[e^{-\mathscr I_{r\mu_j}}]
 =e^{A_r}.
\]
This proves \eqref{eq:amplification}; it does not identify the boundary
map of $r\mu_j$ with a power of $\varphi_j$.

We next prove \eqref{eq:action-L2}. We have
$N(\mu_j-\mu)\to0$ and $\|\mu_j\|_\infty\le k$.
Uniformly bounded strong operator convergence becomes
Hilbert--Schmidt convergence after multiplication by a fixed
Hilbert--Schmidt operator; finite-rank approximation proves this fact.
The multiplication operators and their Neumann resolvents converge
strongly, while $M_{\eta_j}J$ and $\bar J^*M_{\eta_j}$ converge in
Hilbert--Schmidt norm. The resolvent formulas therefore give
\[
 Q_{\eta_j}\to Q_\eta\ \text{in trace norm},\qquad
 G_j:=G_1(\varphi_j)\to G:=G_1(\varphi)\ \text{in Hilbert--Schmidt norm}.
\]
Since $\|G_j\|,\|G\|\le k<1$, the two determinant limits follow from
\[
 \begin{aligned}
 \|G_j^*G_j-G^*G\|_{\tr}
 &\le(\|G_j\|_{\HS}+\|G\|_{\HS})\|G_j-G\|_{\HS},\\
 |\log\det{}_2(I-L_j)-\log\det{}_2(I-L)|
 &\le\frac{\|L_j\|_{\HS}+\|L\|_{\HS}}{1-k}
          \|L_j-L\|_{\HS},
 \end{aligned}
\]
where $L_j=G_j\overline{G_1(h)}$ and $L=G\overline{G_1(h)}$.
The second estimate follows by telescoping the powers in
\eqref{eq:det2}; the ordinary logarithm is Lipschitz in trace norm
on positive operators bounded below by $1-k^2$.
The Grunsky inequality $\|G_1(h)\|\le1$ makes both limits uniform in $h$.

For each $f$, \eqref{eq:koebe} gives
$\|(\lambda_{f,j}-\lambda_f)/z\|_2\to0$. Set
\[
 q_j=\lambda_{f,j}/z,\quad q=\lambda_f/z,\qquad
 S_j=(I-M_{\lambda_{f,j}}\mathcal B)^{-1},\quad
 S=(I-M_{\lambda_f}\mathcal B)^{-1}.
\]
The Neumann series gives $S_j\to S$ strongly and
$\sup_j\|S_j\|\le(1-k)^{-1}$. Thus
\[
 \begin{aligned}
 \|S_jq_j-Sq\|_2
 &\le(1-k)^{-1}\|q_j-q\|_2+\|(S_j-S)q\|_2\longrightarrow0,\\
 \pi|\mathcal R(\lambda_{f,j})-\mathcal R(\lambda_f)|
 &\le\|q_j-q\|_2\|\mathcal B S_jq_j\|_2
 +\|q\|_2\|\mathcal B(S_jq_j-Sq)\|_2\longrightarrow0.
 \end{aligned}
\]
The uniform bounds in Lemma~\ref{lem:bounds} give
$L^2(\m P^{\mathrm{weld}})$ convergence by dominated convergence.
The same applies to $\mathcal R(\mu_j)$. Finally \eqref{eq:v-bound} and
\eqref{eq:divergence-moment} give
$\mathfrak d_{v_{\mu_j}}\to \mathfrak d_{v_\mu}$ in $L^2$. This proves
\eqref{eq:action-L2}.
\end{proof}

\begin{prop}[Agreement with the finite-energy action]\label{prop:rn-finite-energy}
For a deterministic $h\in\WP(\m S^1)$,
\begin{equation}\label{eq:rn-stress-dual}
 \sum_{n\ge1}\frac{|B_n(f_h)|^2}{n^3}<\infty.
\end{equation}
Thus the series \eqref{eq:divergence} converges absolutely for every
real $v\in H^{3/2}$ and gives a pointwise interpretation of its linear
term. For $h,\varphi\in\WP(\m S^1)$ and any representative $\mu$ in
\eqref{eq:mu}, the action so interpreted satisfies
\begin{equation}\label{eq:rn-action-energy}
 \mathscr I_\mu(\varphi,h)
 =\frac{c_{\mathrm L}}{24}
   \bigl[I^{\mathrm L}(h\varphi^{-1})-I^{\mathrm L}(h)\bigr]
  +2\bigl[K(h\varphi^{-1})-K(h)\bigr].
\end{equation}
In particular, the complete action is independent of $\mu$ on WP
weldings as well as almost surely under the SLE welding law.
\end{prop}
\begin{proof}
Write $S_n=[z^{n-2}]S_f$. Since $G_1(h)$ is Hilbert--Schmidt,
\eqref{eq:schwarzian-b1} and Cauchy--Schwarz on each anti-diagonal give
\[
 \sum_{n\ge2}\frac{|S_n|^2}{n^3}
 \le6\sum_{m,\ell\ge1}|(G_1(h))_{m\ell}|^2<\infty;
\]
here $\sum_{m=1}^{n-1}m(n-m)=(n^3-n)/6$.
For the remaining part of $\mathcal T_f$, Koebe's bounds on
$1/2<|z|<1$ imply
\[
 (1-|z|^2)^2\left|\left(\frac{f'}f\right)^2-z^{-2}\right|^2
 \le C\bigl(1+|f'(z)|^2\bigr).
\]
The right side is integrable, since $f(\m D)$ is bounded and
$\int_{\m D}|f'|^2\dd A=\operatorname{Area}(f(\m D))<\infty$.
Parseval's identity on this annulus, where the radial integral for the
$n$th Laurent coefficient is comparable to $n^{-3}$, proves
\eqref{eq:rn-stress-dual}. Cauchy--Schwarz now gives absolute convergence
of \eqref{eq:divergence} and continuity in $[v]_{3/2}$ for this fixed $h$.

Take the radial cutoffs $\mu_j$ of Lemma~\ref{lem:approx}, and set
$G_j=G_1(\varphi_j)$, $G=G_1(\varphi)$ and
$Z=\overline{G_1(h)}$. Since $Z$ is Hilbert--Schmidt,
\[
 \|(G_j-G)Z\|_{\tr}\le\|G_j-G\|_{\HS}\|Z\|_{\HS}\longrightarrow0.
\]
Thus \eqref{eq:geom-inverse} converges to the ordinary determinant
formula for $D(h,\varphi^{-1})$. Equation~\eqref{eq:rn-energy-cocycle}
identifies this limit with the energy difference divided by $24$.

For the capacity term, let $W_j=W^{\lambda_{f,j}}$ and
$W=W^{\lambda_f}$. These coefficients converge in $L^p$ for some
$p>2$ sufficiently close to $2$, have a common dilatation bound below
one, and are supported in a common bounded set away from $0$.
Their normalized solutions and conformal derivatives at $0$ and
$\infty$ therefore converge. The same holds for the reflected solutions
defining $a_{\mu_j}$. The surgery identity gives
\[
 K(h\varphi_j^{-1})-K(h)
 =\log|a_{\mu_j}'(0)|-\log|W_j'(0)|.
\]
The right side converges to
$\log|a_\mu'(0)|-\log|W'(0)|=K(h\varphi^{-1})-K(h)$.
Finally, $v_{\mu_j}\to v_\mu$ in $H^{3/2}$ by \eqref{eq:v-bound},
so the pointwise linear terms converge by \eqref{eq:rn-stress-dual}.
All other terms of the explicit action converge pointwise by the
operator and scalar resolvent estimates in Lemma~\ref{lem:approx}.
Passing to the limit in Lemma~\ref{lem:analytic-action} proves
\eqref{eq:rn-action-energy}. Its right side is independent of the
Beltrami representative.
\end{proof}

\begin{proof}[Proof of Theorem~\ref{body:thm:one}]
If $k=0$, then $a_\mu$ is a disk automorphism fixing $0$, so $\varphi$ is
a rotation. The welding law is invariant under right rotations by rotation
invariance of the loop and the independent uniform left rotation. Every
term in \eqref{body:eq:action} then vanishes. The result follows in this case.

Assume $0<k<1$. By Lemma~\ref{lem:bounds}, \eqref{body:eq:action} is a real
$L^2$ random variable. For the analytic maps $\varphi_j$,
Lemmas~\ref{lem:analytic-cocycle} and \ref{lem:analytic-action} give
$\rho_{\varphi_j}=e^{-\mathscr I_{\mu_j}}$ and $\m E^{\mathrm{weld}}[\rho_{\varphi_j}]=1$.
Choose $1<r<1/k$. By \eqref{eq:amplification}, these densities are bounded
in $L^r$ and are uniformly integrable. Equation \eqref{eq:action-L2} and
Vitali's theorem give
\begin{equation}\label{eq:density-L1}
 e^{-\mathscr I_{\mu_j}}\longrightarrow e^{-\mathscr I_\mu}\quad\text{in }L^1(\m P^{\mathrm{weld}}),
 \qquad \m E^{\mathrm{weld}}[e^{-\mathscr I_\mu}]=1.
\end{equation}
For a bounded continuous $F$ on the circle homeomorphism group,
\[
 \m E^{\mathrm{weld}}[F(\widehat H\varphi_j)]=\m E^{\mathrm{weld}}[F(\widehat H)e^{-\mathscr I_{\mu_j}}].
\]
The maps $h\mapsto h\varphi_j$ converge pointwise to
$h\mapsto h\varphi$. Thus, by bounded convergence and \eqref{eq:density-L1},
\[
 \m E^{\mathrm{weld}}[F(\widehat H\varphi)]
 =\m E^{\mathrm{weld}}[F(\widehat H)e^{-\mathscr I_\mu(\varphi,\widehat H)}].
\]
Bounded continuous $F$ determine Borel measures. This identifies the
limit as the actual pushforward and proves \eqref{body:eq:one}. Since the action
is finite, the density is strictly positive. For a second representative,
the same argument gives the same actual density. Uniqueness of the real
logarithm proves representative independence.

Taking a subsequence in \eqref{eq:action-L2} and using Fatou's lemma in
\eqref{eq:amplification} proves \eqref{body:eq:one-moments}. For a fixed WP
map, the inverse belongs to WP. If its representative has norm $k_-<1$,
the inverse-density identity and change of variables give, in the extended
sense,
\[
 \rho_\varphi(x\varphi)\rho_{\varphi^{-1}}(x)=1,\qquad
 \begin{aligned}
 \m E^{\mathrm{weld}}[\rho_\varphi^p]
 &=\m E^{\mathrm{weld}}[\rho_\varphi(x\varphi)^{p-1}]\\
 &=\m E^{\mathrm{weld}}[\rho_{\varphi^{-1}}^{1-p}]
 <\infty\quad(1<p<1+1/k_-).
 \end{aligned}
\]
In particular some $p>2$ is allowed. A small positive and negative moment
of $\rho_\varphi$ gives $\m E^{\mathrm{weld}}[e^{a|\mathscr I(\varphi,\widehat H)|}]<\infty$ for some $a>0$,
and hence all finite moments of $\mathscr I$. For a rotation, the density
is one by the first part of the proof and representative independence.
\end{proof}

The same proof shows that the cutoff actions converge in every finite
$L^p(\m P^{\mathrm{weld}})$: \eqref{eq:amplification} gives a uniform exponential moment
of their absolute values, and \eqref{eq:action-L2} gives convergence in
probability. We also record the one-sided cocycle
\begin{equation}\label{eq:one-cocycle}
 \mathscr I(\varphi\chi,h)=\mathscr I(\chi,h)+\mathscr I(\varphi,h\chi^{-1})
 \quad\m P^{\mathrm{weld}}\text{-a.s.}
\end{equation}
It follows from the actual density and $R_{\varphi\chi}=R_\chi\circ R_\varphi$.

\subsection{Inversion and two-sided composition}\label{subsec:rn-two-sided}

Finally, inversion symmetry converts the one-sided density into a left-composition density. The chain rule then gives the formula for $T_{\psi,\varphi}$, and a change of variables combines the one-sided moment bounds into the asserted two-sided integrability. We also record the corresponding infinitesimal two-sided divergence.

\begin{lem}[Inversion invariance]\label{lem:inversion}
For $\kappa\in(0,4]$, $\widehat H^{-1}$ has law $\m P^{\mathrm{weld}}$.
\end{lem}

\begin{proof}
Before adding the phase let $h_0=g^{-1}f$, $A=g'(\infty)>0$ and
$K(f)=\log A$. The maps for the inverted loop, after normalizing the inner
map, are
\[
 f^\star(z)=\frac A{\overline{g(1/\bar z)}},\qquad
 g^\star(z)=\frac A{\overline{f(1/\bar z)}},\qquad
 s^\star=-s-K(f).
\]
In particular $K(f^\star)=K(f)$ and $h_0(f^\star)=h_0(f)^{-1}$.
The scale decomposition and \eqref{eq:inversion-input} imply, for bounded
$F$ and compactly supported $\chi$,
\[
 \int\chi(s)F(f)\dd s\dd\sigma(f)
 =\int\chi(-s-K(f))F(f^\star)\dd s\dd\sigma(f).
\]
Integrating first in $s$ gives $\operatorname{Law}(f^\star)=\sigma$; no capacity weight
remains. Physical rotation gives invariance of $h_0$ under conjugation by
rotations. Therefore its left and right Haar averages agree. Since
$\widehat H=r_{-\alpha}h_0$, we obtain
$\widehat H^{-1}=h_0^{-1}r_\alpha\overset d=r_{-\alpha}h_0=\widehat H$.
\end{proof}

For $\psi,\varphi\in\WP(\m S^1)$ let
\begin{equation}\label{body:eq:two-action}
 T_{\psi,\varphi}(h)=\psi^{-1}\circ h\circ\varphi.
\end{equation}

\begin{thm}[Two-sided WP density]\label{body:thm:two}
For $\kappa\in(0,4]$ and $\psi,\varphi\in\WP(\m S^1)$,
the pushforward by $T_{\psi,\varphi}$ is mutually absolutely
continuous with respect to $\m P^{\mathrm{weld}}$. Its density satisfies
\begin{equation}\label{body:eq:two}
 \begin{aligned}
 \frac{\dd(T_{\psi,\varphi})_*\m P^{\mathrm{weld}}}{\dd\m P^{\mathrm{weld}}}(h)
 &=\rho_\varphi(h)\rho_\psi(\varphi h^{-1})\\
 &=\rho_\psi(h^{-1})\rho_\varphi(\psi h)\\
 &=e^{-\mathscr J_{\psi,\varphi}(h)},\\
 \mathscr J_{\psi,\varphi}(h)
 &=\mathscr I(\varphi,h)+\mathscr I(\psi,\varphi h^{-1})
 =\mathscr I(\psi,h^{-1})+\mathscr I(\varphi,\psi h).
 \end{aligned}
\end{equation}
The evaluations of $\mathscr I$ are well defined modulo $\m P^{\mathrm{weld}}$-null sets. For some
$\epsilon>0$ depending on $\psi,\varphi$, the density in \eqref{body:eq:two} belongs to $L^{1+\epsilon}(\m P^{\mathrm{weld}})$ and its $\epsilon$-th negative moment is finite. Moreover $\mathscr J_{\psi,\varphi}\in L^p(\m P^{\mathrm{weld}})$ for each finite $p$.
\end{thm}

\begin{proof}[Proof of Theorem~\ref{body:thm:two}]
Write $\jmath(h)=h^{-1}$ and $L_\ell(h)=\ell h$. By
Lemma~\ref{lem:inversion}, $\jmath_*\m P^{\mathrm{weld}}=\m P^{\mathrm{weld}}$, and the group identity
$L_\ell=\jmath R_{\ell^{-1}}\jmath$ gives
\begin{equation}\label{eq:left-density}
 \frac{\dd(L_\ell)_*\m P^{\mathrm{weld}}}{\dd\m P^{\mathrm{weld}}}(h)=\rho_{\ell^{-1}}(h^{-1}).
\end{equation}
For two invertible maps $S,T$ preserving the measure class, their densities
satisfy
\begin{equation}\label{eq:composition-density}
 \rho_{S\circ T}(h)=\rho_S(h)\rho_T(S^{-1}h).
\end{equation}
Indeed, for every bounded nonnegative $F$,
\[
 \begin{aligned}
 \int F\,\dd(S\circ T)_*\m P^{\mathrm{weld}}
 &=\int F(Sx)\rho_T(x)\,\dd\m P^{\mathrm{weld}}(x)\\
 &=\int F(h)\rho_T(S^{-1}h)\rho_S(h)\,
 \dd\m P^{\mathrm{weld}}(h).
 \end{aligned}
\]
Since $T_{\psi,\varphi}=R_\varphi L_{\psi^{-1}}
=L_{\psi^{-1}}R_\varphi$, the two decompositions give the two products in
\eqref{body:eq:two}. The map $h\mapsto\varphi h^{-1}$ sends $\m P^{\mathrm{weld}}$ to a
measure mutually absolutely continuous with respect to $\m P^{\mathrm{weld}}$ by \eqref{eq:left-density}; hence evaluation
of any fixed $\m P^{\mathrm{weld}}$-version of $\mathscr I(\psi,\cdot)$ there is well defined.
The same holds for the other evaluation. Taking negative logarithms proves
\eqref{body:eq:two} and the representative independence.

It remains to prove the moment assertion. Denote the two-sided density by
$R$. For any real $r$, the change of variables $h=x\varphi$ yields, in the
extended nonnegative sense,
\begin{equation}\label{eq:two-moment}
 \m E^{\mathrm{weld}}[R^r]=\m E^{\mathrm{weld}}\bigl[\rho_{\varphi^{-1}}(x)^{1-r}
 \rho_\psi(x^{-1})^r\bigr].
\end{equation}
Put $A=\rho_{\varphi^{-1}}$ and $B=\rho_\psi\circ\jmath$. Theorem
\ref{body:thm:one}, inversion invariance, and the inverse-density change of
variables give $a,b>0$ such that
$\m E^{\mathrm{weld}}[A^{1+a}+A^{-a}+B^{1+b}+B^{-b}]<\infty$.
Choose conjugate exponents $p,q>1$ with $p<1+\min(a,b)$, and then
choose $\epsilon>0$ so that
$p(1+\epsilon)\le1+\min(a,b)$ and $q\epsilon\le\min(a,b)$. H\"older's
inequality gives
\[
 \begin{aligned}
 \m E^{\mathrm{weld}}[R^{1+\epsilon}]
 &\le\m E^{\mathrm{weld}}[A^{-q\epsilon}]^{1/q}
      \m E^{\mathrm{weld}}[B^{p(1+\epsilon)}]^{1/p},\\
 \m E^{\mathrm{weld}}[R^{-\epsilon}]
 &\le\m E^{\mathrm{weld}}[A^{p(1+\epsilon)}]^{1/p}
      \m E^{\mathrm{weld}}[B^{-q\epsilon}]^{1/q}.
 \end{aligned}
\]
Both right sides are finite. It follows that
$\m E^{\mathrm{weld}}[e^{\epsilon|\mathscr J_{\psi,\varphi}|}]<\infty$ after decreasing $\epsilon$,
and hence all finite $L^p$ moments of $\mathscr J_{\psi,\varphi}$.
\end{proof}

For deterministic $h,\varphi,\psi\in\WP(\m S^1)$, the interpretation in
Proposition~\ref{prop:rn-finite-energy} gives
\begin{equation}\label{eq:rn-two-energy}
 \begin{aligned}
 \mathscr J_{\psi,\varphi}(h)
 ={}&\frac{c_{\mathrm L}}{24}
 \bigl[I^{\mathrm L}(\psi h\varphi^{-1})-I^{\mathrm L}(h)\bigr]\\
 &+2\bigl[K(\psi h\varphi^{-1})-K(h)\bigr].
 \end{aligned}
\end{equation}
Indeed, both $I^{\mathrm L}$ and $K$ are invariant under inversion;
substituting \eqref{eq:rn-action-energy} into \eqref{body:eq:two}
makes the intermediate terms cancel. The welding appearing here is
$T_{\psi,\varphi}^{-1}(h)=\psi h\varphi^{-1}$, as required by the
target-point convention for the density.

The second term in \eqref{body:eq:two} can be computed from the original pair.
If $h=g^{-1}f$, $f'(0)=1$, $\mu$ represents $\varphi$, and
$W=W^{\lambda_f}$, then $(Wfa_\mu^{-1},Wg)$ is a pair for
$h\varphi^{-1}$. Its inverse has normalized inner map
\begin{equation}\label{eq:second-inner}
 f_{\varphi h^{-1}}(z)
 =\frac{\overline{g'(\infty)}}{\overline{W(g(1/\bar z))}}.
\end{equation}
Thus every term in the second action is defined by the original welding
pair and the two chosen WP representatives.

\begin{prop}\label{prop:two-divergence}
For real $u,v\in H^{3/2}(\m S^1)$, set
\begin{equation}\label{eq:two-divergence}
 \mathfrak d^{(2)}_{u,v}(h)=\mathfrak d_v(h)+\mathfrak d_u(h^{-1}),\qquad
 \|\mathfrak d^{(2)}_{u,v}\|_2^2
 =2\sum_{n\ge1}\bigl[c(n^3-n)+2n\bigr](|u_n|^2+|v_n|^2).
\end{equation}
For finite Fourier fields this is the logarithmic derivative of the density for
$h\mapsto(e_t^u)^{-1}he_t^v$. For general $H^{3/2}$ fields, the assertion
is the $L^2$ extension; it does not require a classical flow for each field.
\end{prop}

\begin{proof}[Proof of Proposition~\ref{prop:two-divergence}]
For finite Fourier fields the derivative of the actual two-sided density
follows from \eqref{eq:composition-density} and
Lemma~\ref{lem:finite-flow}. Its logarithmic derivative is
$\mathfrak d_v(h)+\mathfrak d_u(h^{-1})$. The inner map is unchanged by a left rotation of $h$,
so $\mathfrak d_v(r_\theta h)=\mathfrak d_v(h)$. On the other hand
$(r_\theta h)^{-1}=h^{-1}r_{-\theta}$, and
$\int_0^{2\pi}\mathfrak d_u(h^{-1}r_{-\theta})\dd\theta=0$ for finite Fourier $u$.
Since $\m P^{\mathrm{weld}}$ is invariant under left rotations,
$\m E^{\mathrm{weld}}[\mathfrak d_v(\widehat H)\mathfrak d_u(\widehat H^{-1})]=0$. This proves orthogonality, not independence.
Now use Lemma~\ref{lem:inversion} and \eqref{eq:divergence-moment}. The extension
to $H^{3/2}$ follows by Fourier truncation in $L^2(\m P^{\mathrm{weld}})$.
\end{proof}

\begin{prop}\label{prop:two-approx}
Let $\varphi_j,\psi_j$ be the analytic maps obtained by the radial cutoffs
in Lemma~\ref{lem:approx}, applied separately to fixed representatives of
$\varphi,\psi$. Then
\begin{equation}\label{eq:two-approx}
 \begin{aligned}
 \rho_{\psi_j,\varphi_j}&\to\rho_{\psi,\varphi}
 \quad\text{in }L^1(\m P^{\mathrm{weld}}),\\
 \mathscr J_{\psi_j,\varphi_j}&\to\mathscr J_{\psi,\varphi}
 \quad\text{in }L^p(\m P^{\mathrm{weld}}),\qquad 1\le p<\infty.
 \end{aligned}
\end{equation}
where $\rho_{\psi,\varphi}$ denotes the density in \eqref{body:eq:two}.
\end{prop}
\begin{proof}
We first use the following continuity fact. If $S_j,S$ are invertible
Borel maps on a Polish probability space, $S_j^{-1}x\to S^{-1}x$ for every
$x$, and their densities $q_j\to q$ in $L^1$, then
$(S_j)_*(f_j\m P^{\mathrm{weld}})\to S_*(f\m P^{\mathrm{weld}})$ in total variation whenever
$f_j\to f$ in $L^1$. For bounded continuous $f$, this follows from
$q_j(f\circ S_j^{-1})\to q(f\circ S^{-1})$ in $L^1$, using dominated
convergence for the term with $q$. Approximation by bounded continuous
functions and the total variation isometry of invertible pushforwards prove
the general case.

By Theorem~\ref{body:thm:one}, the one-sided cutoff densities converge in
$L^1$. Take $S_j=R_{\varphi_j}$ and
$f_j=\rho_{\psi_j}\circ\jmath$ in the continuity fact. Since
$T_{\psi_j,\varphi_j}=R_{\varphi_j}L_{\psi_j^{-1}}$ and
$\jmath_*\m P^{\mathrm{weld}}=\m P^{\mathrm{weld}}$, we obtain the first convergence in
\eqref{eq:two-approx}. The limiting density is positive, so its negative
logarithm converges in probability.

The moment argument in the proof of Theorem~\ref{body:thm:two} is uniform for
these cutoffs: the one-sided positive and negative moment bounds follow
uniformly from \eqref{eq:amplification} and the inverse-density identity.
Thus, for some $\epsilon>0$,
\[
 \sup_j\m E^{\mathrm{weld}}\left[\rho_{\psi_j,\varphi_j}^{1+\epsilon}
 +\rho_{\psi_j,\varphi_j}^{-\epsilon}\right]<\infty.
\]
This gives a uniform exponential moment for
$|\mathscr J_{\psi_j,\varphi_j}|$. Vitali's theorem proves the second convergence.
No pointwise continuity of $\mathfrak d_v$ as a function of a rough welding is used.
\end{proof}

The two-sided action also satisfies, for fixed parameters,
\begin{equation}\label{eq:two-cocycle}
 \mathscr J_{\psi_1\psi_2,\varphi_1\varphi_2}(h)
 =\mathscr J_{\psi_2,\varphi_2}(h)
 +\mathscr J_{\psi_1,\varphi_1}(\psi_2h\varphi_2^{-1}),\qquad
 \mathscr J_{\psi,\varphi}(h)=\mathscr J_{\varphi,\psi}(h^{-1}).
\end{equation}
These identities follow from \eqref{eq:composition-density} and inversion
invariance. For the alternative action $h\mapsto\ell hr$, its density is
$\exp\{-\mathscr I(r,h)-\mathscr I(\ell^{-1},rh^{-1})\}$.
In particular, for the convention $h\mapsto\varphi h\psi^{-1}$ used in
the symbol table, the density with respect to $\m P^{\mathrm{weld}}$ is
\begin{equation}\label{eq:table-convention-sle}
 \exp\{-\mathscr I(\psi^{-1},h)
       -\mathscr I(\varphi^{-1},\psi^{-1}h^{-1})\}.
\end{equation}
\section{Appendix}\label{sec:appendix}

The appendix collects auxiliary material on Weil--Petersson circle homeomorphisms. It includes the finite-increment energy estimate, a discussion of the subgroup fixing three marked points, and a one-variable log-ratio estimate. These results complement the regularity and normalizer arguments developed in the main text.

\subsection{Finite increment energy for WP homeomorphisms}

For $\varphi\in\Homeo_+(\m S^1)$, let $\Phi$ be its degree-one lift and
let $\delta_\varphi$ be as in Section~\ref{subsec:subcritical-wp}. We write
\begin{equation}\label{eq:wp-increment-energy}
 E(\varphi)=\int_0^{1/2}\int_{\m S^1}
 |\delta_\varphi(x,r)|^2\,\dd x\,\frac{\dd r}{r^2}.
\end{equation}
All circle lengths and integrals are normalized as in Section~\ref{subsec:subcritical-wp}. Put
\[
 d(s,t)=\min_{k\in\mathbb Z}|s-t+k|,
 \qquad
 [u]_{1/2}^2=\int_{\m S^1}\int_{\m S^1}
 \frac{|u(s)-u(t)|^2}{d(s,t)^2}\,\dd s\,\dd t.
\]
This is an equivalent seminorm on $H^{1/2}(\m S^1)$. For an arc $I$, set
\[
 u_I=\frac1{|I|}\int_Iu(t)\,\dd t,\qquad
 V_I(u)=\frac1{|I|}\int_I|u(t)-u_I|^2\,\dd t.
\]

\begin{lem}\label{lem:wp-energy-exp}
Let $I$ be an interval, and suppose $v\in L^1(I)$ satisfies
\[
 \frac1{|J|}\int_J|v(t)-v_J|\,\dd t\le\varepsilon
 \quad\text{for every subinterval }J\subset I.
\]
If $0<\varepsilon<\log2/8$, then
\begin{equation}\label{eq:wp-energy-exp}
 \frac1{|I|}\int_I e^{2|v(t)-v_I|}\,\dd t
 \le\frac{e^{4\varepsilon}}{1-e^{8\varepsilon}/2}.
\end{equation}
\end{lem}
\begin{proof}
We give the dyadic stopping argument. For an interval $J$, split it into
its two equal halves and continue recursively. Take the maximal dyadic
subintervals $J'$ of $J$ such that
\[
 \frac1{|J'|}\int_{J'}|v(t)-v_J|\,\dd t>2\varepsilon.
\]
The selected intervals are disjoint. Their total length satisfies
\[
 \sum_{J'}|J'|
 \le\frac1{2\varepsilon}\sum_{J'}\int_{J'}|v-v_J|\,\dd t
 \le\frac1{2\varepsilon}\int_J|v-v_J|\,\dd t
 \le\frac{|J|}{2}.
\]
Their dyadic parents are not selected. Hence, for every selected $J'$,
\[
 |v_{J'}-v_J|
 \le\frac1{|J'|}\int_{J'}|v(t)-v_J|\,\dd t
 \le4\varepsilon.
\]
By dyadic differentiation, outside their union we have
$|v-v_J|\le2\varepsilon$ almost everywhere.

Repeat this selection inside each selected interval, now subtracting its
own average. Let $B_n$ be the union of intervals selected at generation $n$,
with $B_0=I$. Then $|B_n|\le2^{-n}|I|$. On $B_n\setminus B_{n+1}$,
\[
 |v-v_I|\le4n\varepsilon+2\varepsilon
 \quad\text{almost everywhere}.
\]
Since $|B_n|\to0$, these sets cover $I$ up to a null set. Consequently,
\[
 \begin{aligned}
 \frac1{|I|}\int_I e^{2|v-v_I|}\,\dd t
 &\le\sum_{n\ge0}2^{-n}e^{2(4n+2)\varepsilon}\\
 &=e^{4\varepsilon}\sum_{n\ge0}
 \left(\frac{e^{8\varepsilon}}2\right)^n.
 \end{aligned}
\]
The ratio is less than $1$, which proves \eqref{eq:wp-energy-exp}.
\end{proof}

\begin{thm}\label{thm:wp-finite-increment-energy}
If $\varphi\in\WP(\m S^1)$, then $E(\varphi)<\infty$.
In particular, the quantities $E_N$ defined in Section~\ref{subsec:subcritical-wp} are bounded.
\end{thm}
\begin{proof}
We use the standard WP characterization: the lift $\Phi$ is absolutely
continuous, $\Phi'>0$ almost everywhere, and
\[
 u:=\log\Phi'\in H^{1/2}(\m S^1).
\]
See \cite[Definition~2.3]{F1} for this characterization. In particular,
$u$ is a real periodic $L^2$ function and $[u]_{1/2}<\infty$.

For every arc $J$ shorter than the circle, expanding the square gives
\begin{equation}\label{eq:wp-energy-variance}
 \begin{aligned}
 \frac1{2|J|^2}\int_J\int_J|u(s)-u(t)|^2\,\dd s\,\dd t
 &=\frac1{|J|}\int_Ju(t)^2\,\dd t-u_J^2\\
 &=V_J(u).
 \end{aligned}
\end{equation}
Since $d(s,t)\le|J|$ for $s,t\in J$, it follows that
\begin{equation}\label{eq:wp-energy-local}
 V_J(u)\le\frac12\int_J\int_J
 \frac{|u(s)-u(t)|^2}{d(s,t)^2}\,\dd s\,\dd t.
\end{equation}
The integrand on the right is in $L^1(\m S^1\times\m S^1)$. Given
$\eta>0$, absolute continuity of its integral gives $\delta>0$ such that
its integral over any measurable set of area less than $\delta$ is less
than $2\eta$. For $|J|<\sqrt\delta$, the area of $J\times J$ is less than
$\delta$. Hence \eqref{eq:wp-energy-local} gives $V_J(u)<\eta$, and
\begin{equation}\label{eq:wp-energy-local-small}
 \lim_{\ell\downarrow0}\sup_{|J|\le\ell}V_J(u)=0.
\end{equation}
In particular, by Cauchy--Schwarz,
\[
 \frac1{|J|}\int_J|u-u_J|\,\dd t\le V_J(u)^{1/2}.
\]

Fix $\varepsilon_0=\log2/16$. By \eqref{eq:wp-energy-local-small}, choose
$0<r_0<1/8$ so that $V_J(u)^{1/2}\le\varepsilon_0$ for every arc $J$
of length at most $2r_0$. The same bound holds on every subinterval of
such a $J$. Lemma~\ref{lem:wp-energy-exp} gives a finite constant
\[
 C_0:=\frac{e^{4\varepsilon_0}}{1-e^{8\varepsilon_0}/2}
\]
such that, for every arc $I$ of length at most $2r_0$,
\begin{equation}\label{eq:wp-energy-exp-bound}
 \frac1{|I|}\int_I e^{2|u-u_I|}\,\dd t\le C_0.
\end{equation}

Fix $x\in\m S^1$ and $0<r\le r_0$. Put
\[
 I=[x-r,x+r],\quad I_-=[x-r,x],\quad I_+=[x,x+r],\quad
 v=u-u_I,\quad \sigma_I=V_I(u)^{1/2}.
\]
Define the two positive averages
\[
 A_\pm=\frac1r\int_{I_\pm}e^{v(t)}\,\dd t.
\]
Jensen's inequality, followed by Cauchy--Schwarz, gives
\[
 \begin{aligned}
 A_\pm
 &\ge\exp\!\left(\frac1r\int_{I_\pm}v(t)\,\dd t\right)\\
 &\ge\exp\!\left(-\frac1r\int_{I_\pm}|v(t)|\,\dd t\right)
 \ge e^{-\sqrt2\sigma_I}.
 \end{aligned}
\]
Indeed,
\[
 \frac1r\int_{I_\pm}|v(t)|\,\dd t
 \le\left(\frac1r\int_{I_\pm}|v(t)|^2\,\dd t\right)^{1/2}
 \le\sqrt2\sigma_I.
\]
Since $\sigma_I\le\varepsilon_0$, both averages are at least
$e^{-\sqrt2\varepsilon_0}$.

Using $|e^z-1|\le |z|e^{|z|}$ for real $z$, we obtain
\[
 \begin{aligned}
 |A_+-A_-|
 &\le\frac1r\int_I|e^{v(t)}-1|\,\dd t\\
 &\le\frac1r\int_I|v(t)|e^{|v(t)|}\,\dd t\\
 &\le2\sigma_I
 \left(\frac1{|I|}\int_Ie^{2|v(t)|}\,\dd t\right)^{1/2}\\
 &\le2\sqrt{C_0}\,\sigma_I.
 \end{aligned}
\]
The last two steps use $|I|=2r$, Cauchy--Schwarz, and
\eqref{eq:wp-energy-exp-bound}. By the mean value theorem for $\log$,
\begin{equation}\label{eq:wp-energy-delta}
 \begin{aligned}
 |\log A_+-\log A_-|
 &\le e^{\sqrt2\varepsilon_0}|A_+-A_-|\\
 &\le2e^{\sqrt2\varepsilon_0}\sqrt{C_0}\,\sigma_I.
 \end{aligned}
\end{equation}

Absolute continuity of $\Phi$ and $\Phi'=e^u$ give
\[
 \frac{\Phi(x+r)-\Phi(x)}{\Phi(x)-\Phi(x-r)}
 =\frac{\int_{I_+}e^{u(t)}\,\dd t}
 {\int_{I_-}e^{u(t)}\,\dd t}
 =\frac{A_+}{A_-}.
\]
Thus \eqref{eq:wp-energy-delta} yields the pointwise estimate
\begin{equation}\label{eq:wp-energy-pointwise}
 |\delta_\varphi(x,r)|^2
 \le4e^{2\sqrt2\varepsilon_0}C_0\,V_{[x-r,x+r]}(u),
 \qquad 0<r\le r_0.
\end{equation}

It remains to integrate the variances. Put $I(x,r)=[x-r,x+r]$. By
\eqref{eq:wp-energy-variance} and Tonelli's theorem,
\[
 \begin{aligned}
 &\int_0^{r_0}\int_{\m S^1}V_{I(x,r)}(u)\,\dd x\,
 \frac{\dd r}{r^2}\\
 &\quad=\frac18\int_{\m S^1}\int_{\m S^1}|u(s)-u(t)|^2
 \left(\int_0^{r_0}\frac{1}{r^4}
 \int_{\m S^1}\mathbf1_{\{s,t\in I(x,r)\}}\,\dd x\,\dd r\right)
 \dd s\,\dd t.
 \end{aligned}
\]
For fixed $s,t$, the inner set of $x$ has length at most $2r$. It is empty
unless $d(s,t)\le2r$. Therefore, for $s\ne t$, the expression in
parentheses, multiplied by $1/8$, is at most
\[
 \frac14\mathbf1_{\{d(s,t)\le2r_0\}}
 \int_{d(s,t)/2}^{r_0}\frac{\dd r}{r^3}
 \le\frac1{2d(s,t)^2}.
\]
It follows that
\begin{equation}\label{eq:wp-energy-variance-integral}
 \int_0^{r_0}\int_{\m S^1}V_{I(x,r)}(u)\,\dd x\,
 \frac{\dd r}{r^2}
 \le\frac12[u]_{1/2}^2<\infty.
\end{equation}
Combining \eqref{eq:wp-energy-pointwise} and
\eqref{eq:wp-energy-variance-integral},
\[
 \int_0^{r_0}\int_{\m S^1}|\delta_\varphi(x,r)|^2\,\dd x\,
 \frac{\dd r}{r^2}
 \le2e^{2\sqrt2\varepsilon_0}C_0[u]_{1/2}^2<\infty.
\]

Finally, $\varphi\in\WP(\m S^1)$ is quasisymmetric. With
$D=\sup_{x,\,0<r\le1/2}|\delta_\varphi(x,r)|<\infty$,
\[
 \int_{r_0}^{1/2}\int_{\m S^1}|\delta_\varphi(x,r)|^2\,\dd x\,
 \frac{\dd r}{r^2}
 \le D^2\left(\frac1{r_0}-2\right)<\infty.
\]
Adding the two ranges proves $E(\varphi)<\infty$. Since
$\Omega_N\subset\m S^1\times(0,1/2]$,
\[
 E_N=1+\int_{\Omega_N}\delta_\varphi^2\,\dd\eta
 \le1+E(\varphi)<\infty.
\]
The upper bound is independent of $N$.
\end{proof}

\subsection{Normalizers of WP subgroups with fixed points}\label{subsec:fixed-point-normalizers}
For a nonempty finite set $P\subset\m S^1$, define the subgroup
\[
 \mathcal K_P=\{\varphi\in\WP(\m S^1):\varphi(p)=p
                       \text{ for every }p\in P\}.
\]
Unlike the full WP group, this subgroup is not its own normalizer.
The failure persists even among maps fixing $P$ pointwise.

\begin{thm}\label{thm:fixed-point-normalizer}
For every nonempty finite $P$, the normalizer
$\mathcal N(\mathcal K_P)$ in $\Homeo_+(\m S^1)$ contains a
bi-Lipschitz map outside WP which fixes $P$ pointwise. In particular,
$\mathcal K_P\subsetneqq\mathcal N(\mathcal K_P)$.
\end{thm}
\begin{proof}
For this construction identify $\m S^1$ with $\m R/\m Z$.
Choose one marked point as $0$ and write
$0=p_0<p_1<\cdots<p_{m-1}<p_m=1$ for the marked lifts. Put
$C(t)=(e^t-1)/(e-1)$ on $[0,1]$ and define
\[
 A(x)=p_i+(p_{i+1}-p_i)
 C\!\left(\frac{x-p_i}{p_{i+1}-p_i}\right)
 \quad(p_i\le x\le p_{i+1}),\qquad A(x+1)=A(x)+1.
\]
The induced circle map $a$ is bi-Lipschitz and fixes $P$ pointwise.
The function $\log A'$ has a nonzero jump at every marked point, so
it is not in $H^{1/2}$. Indeed, the contribution across a jump bounds
its Douglas integral below by a positive multiple of
$\int_0^\epsilon\int_0^\epsilon(s+t)^{-2}\,\dd s\,\dd t=\infty$.
Hence $a\notin\WP(\m S^1)$.

Let $\varphi\in\mathcal K_P$ and choose its lift $\Phi$ fixing all
$p_i$. The WP characterization gives $\log\Phi'\in H^{1/2}$ and
$\Phi'\in L^2$; the latter follows from exponential integrability,
see \cite[Remark~3.2]{F1}. For either $B=A$ or $B=A^{-1}$, let
$q_B=\log B'$. On each $[p_i,p_{i+1}]$, $q_B$ is $C^1$ with bounded
one-sided derivative. Thus
\[
 r_B=q_B\circ\Phi-q_B\in W^{1,2}([p_i,p_{i+1}]),\qquad
 r_B(p_i+)=r_B(p_{i+1}-)=0.
\]
The pieces join to a periodic $W^{1,2}$ function. Since composition
by a bi-Lipschitz map preserves $H^{1/2}$,
\[
 \log(B\circ\Phi\circ B^{-1})'
 =(\log\Phi'+r_B)\circ B^{-1}\in H^{1/2}.
\]
The conjugate is absolutely continuous, quasisymmetric, and fixes $P$;
hence it belongs to $\mathcal K_P$. Applying this to $A$ and $A^{-1}$
proves $a\mathcal K_Pa^{-1}=\mathcal K_P$.
\end{proof}

We now fix $P=\{1,\ii,-1\}$ and write $\KK=\mathcal K_P$.
The set $\mob(\m S^1)\backslash\WP(\m S^1)$ is identified with
$\KK$ by choosing the representative fixing these three points; this
is an identification of normalized representatives, not a quotient
group by a normal subgroup. The remaining definitions formulate a
conjectural description of the full normalizer.

\begin{df}
For $T\ge1$ and a real $v\in L^2_{\rm loc}([T,\infty))$, set
\[
 k_0(u)=\frac1{4\sinh^2(u/2)}\quad(u\ne0),\qquad
 E_T(v)=\iint_{[T,\infty)^2}|v(t)-v(s)|^2k_0(t-s)\,\dd t\,\dd s.
\]
The diagonal is omitted in this nonnegative integral, which may be
infinite. Define
\[
 \begin{aligned}
 \mc J_T=\{(v_+,v_-)\in L^2_{\rm loc}([T,\infty))^2:
 E_T(v_+)+E_T(v_-)+\|v_+-v_-\|_{L^2([T,\infty))}^2<\infty\}.
 \end{aligned}
\]
For $m=1,2$ and $s_0\in\m R$, let $H^m_{\rm uloc}(s_0)$ be the
space of $a\in W^{m,2}_{\rm loc}([s_0,\infty))$ satisfying
\[
 \|a\|_{H^m_{\rm uloc}(s_0)}^2
 :=\sup_{\tau\ge s_0}\sum_{j=0}^m
       \int_\tau^{\tau+1}|a^{(j)}(s)|^2\,\dd s<\infty.
\]
\end{df}

\begin{df}
A pair of real functions $W=(W_+,W_-)$ belongs to $\mc C$ if, for
some $T\ge1$ and almost every $t\ge T$,
\[
 W_\pm(t)=t\,a(\log t)\pm b(\log t)+e_\pm(t),
\]
where $a\in H^2_{\rm uloc}(\log T)$,
$b\in H^1_{\rm uloc}(\log T)$ and $(e_+,e_-)\in\mc J_T$.
\end{df}
Equivalently, write $W_\pm=F\pm G+e_\pm$ on a sufficiently large tail,
where $F\in W^{2,2}_{\rm loc}$, $G\in W^{1,2}_{\rm loc}$,
$F',G\in L^\infty$, $(e_+,e_-)\in\mc J_T$, and
\[
 \sup_{R\ge T}R\int_R^\infty|F''(t)|^2\,\dd t<\infty,\qquad
 \sup_{R\ge T}R\int_R^\infty|G'(t)|^2\,\dd t<\infty.
\]
To see the equivalence, put $F(t)=ta(\log t)$ and $G(t)=b(\log t)$.
Then
\[
 F'=a+a',\qquad F''=(a'+a'')/t,\qquad G'=b'/t,
\]
with the right sides evaluated at $\log t$. For example,
\[
 R\int_R^\infty|F''(t)|^2\,\dd t
 =\int_{\log R}^\infty|a'(s)+a''(s)|^2e^{\log R-s}\,\dd s.
\]
Uniform bounds on these exponentially weighted integrals are equivalent
to uniform $L^2$ bounds on unit intervals. In the reverse direction,
$F'\in L^\infty$ makes $a=F/t$ and $a'=F'-a$ bounded, giving exactly
the stated $H^2_{\rm uloc}$ condition. The argument for $b$ is the same,
using boundedness of $G$. The forward direction uses the one-dimensional
Sobolev embedding on unit intervals.

\begin{df}
Let $\varphi\in\Homeo_+(\m S^1)$, and suppose that both $\varphi$
and $\varphi^{-1}$ are absolutely continuous. For $p\in P$, choose
a local increasing angular lift
$h_p$ such that
\[
 \varphi(pe^{\ii x})=\varphi(p)e^{\ii h_p(x)},\qquad h_p(0)=0.
\]
For small $r>0$ and sufficiently large $t$, define
\[
 H_{p,\pm}(r)=\pm h_p(\pm r),\qquad
 A_{p,\pm}(t)=-\log H_{p,\pm}(e^{-t}),\qquad
 U_{p,\pm}(t)=\log H_{p,\pm}'(e^{-t}).
\]
Both maps being absolutely continuous implies positive derivatives
almost everywhere. The chain rule gives the exact identity
\[
 U_{p,\pm}(t)=t-A_{p,\pm}(t)+\log A_{p,\pm}'(t)
 \quad\text{for almost every }t.
\]
Set
\[
 V_{p,\pm}(s)=-U_{p,\pm}(A_{p,\pm}^{-1}(s)),\qquad
 U_p=(U_{p,+},U_{p,-}),\quad V_p=(V_{p,+},V_{p,-}).
\]
The pair $V_p$ gives the corresponding logarithmic derivatives of
$\varphi^{-1}$ at $\varphi(p)$, on a sufficiently large tail.
\end{df}

\begin{conj}
A homeomorphism $\varphi\in\Homeo_+(\m S^1)$ belongs to
$\mc N(\KK)$ if and only if all the following conditions hold:
\begin{itemize}
 \item $\varphi(P)=P$;
 \item $\varphi\in\QS(\m S^1)$;
 \item $\varphi$ and $\varphi^{-1}$ are absolutely continuous;
 \item $\log\varphi'\in H^{1/2}_{\rm loc}(\m S^1\setminus P)$;
 \item $U_p,V_p\in\mc C$ for every $p\in P$.
\end{itemize}
\end{conj}
\subsection{Estimate on one-variable log-ratio}
Let $\dot{H}^{1/2}(\m R)$ denote the space of real functions with finite Douglas energy:
\[
    \norm{v}_{\dot{H}^{1/2}(\m R)}^2 = \frac{1}{4\pi^2} \iint_{\m R^2} \left|\frac{v(x)-v(y)}{x-y}\right|^2 \dd x \dd y,
\]
which differs from ${H}^{1/2}(\m R)$ without requiring $v \in L^2(\m R)$. For example, $G(x) = 2x$ belongs to $\WP(\m R)$ but $\log G' \notin L^2(\m R)$.
\begin{lem}\label{lem: one-variable log-ratio}
        For every $z_0 \in \m S^1$ and $h \in \WP(\m S^1)$, the function $F$
    \[
        z \mapsto\log  \left|\frac{h(z)-h(z_0)}{z-z_0}\right |
    \]
    belongs to $\mc H \oplus \m C$.
    
\end{lem}
\begin{remark}
    We strengthen the main theorem of \cite{F1}
\end{remark}
\begin{proof}
    Let $C_a(z)=\ii \frac{z+a}{z-a}$, when $a \in \m S^1$ $C_a$ is a M\"obius map mapping $\m S^1$ to $\hat {\m R}=\m R \cup \infty$.
    \[
        C_a'(z) =  -\frac{2a\ii}{(z-a)^2}
    \]
    Let  $\psi = C_{h(z_0)}\circ h \circ C_{z_0}^{-1} \in \WP(\hat {\m R})$, we have
    \[
        \psi' \circ C_{z_0} C_{z_0}' = C_{h(z_0)}'\circ h h',
    \]
    and hence
    \[
        \log ( \frac{h(z)-h(z_0)}{z-z_0})^2 = \log \frac{C_{z_0}'}{C_{h(z_0)}'\circ h} +\log\frac{h(z_0)}{z_0} = \log h'-\log \psi'\circ C_{z_0} +\log\frac{h(z_0)}{z_0}.
    \]
    By \cite[Theorem~1.3]{Shen-Tang-Wu},
$\psi\in\WP(\hat{\m R})$ implies
$\log\psi'\in\dot H^{1/2}(\m R)$. Conformal invariance of the
Douglas seminorm transfers this statement to the circle. Taking real
parts in the identity above, together with
$\log|h'|\in H^{1/2}(\m S^1)$, proves the assertion.
\end{proof}

\subsection{Operator realizations for the WP action}\label{subsec:operator-realization}

We give the basis calculations underlying \eqref{eq:J}--\eqref{eq:Q}.
They distinguish operators on the full area $L^2$ space from their
compressions to holomorphic or anti-holomorphic subspaces.

We now specify the bases and the meaning of the products. Put
$\phi_n(z)=\sqrt{n/\pi}\,z^{-n-1}$ and
$\psi_n(z)=\sqrt{n/\pi}\,\bar z^{-n-1}$ on $\m D^*$. The identities
\[
 \int_{\m D^*}\phi_m\overline{\phi_n}\,\dd A
 =\int_{\m D^*}\psi_m\overline{\psi_n}\,\dd A=\delta_{mn}
\]
follow from integration in polar coordinates. In particular, $Je_n=\phi_n$
and $\bar Je_n=\psi_n$. The ranges of $J$ and $\bar J$ are the
holomorphic and anti-holomorphic $L^2$ subspaces on $\m D^*$, respectively.
They are subspaces of $L^2(\m D^*)$; neither is all of that space.

The multiplication operator $M_\eta:L^2(\m D^*)\to L^2(\m D^*)$ is
defined by $M_\eta u=\eta u$. Its restriction to the two systems above has
the matrix coefficients
\[
 (\bar J^*M_\eta J)_{mn}
 =\frac{\sqrt{mn}}\pi\int_{\m D^*}\eta(z)z^{-m-n-2}\,\dd A(z).
\]
This is a compression of $M_\eta$, not a matrix for $M_\eta$ on all of
$L^2(\m D^*)$. To write a matrix for $M_\eta$ itself, choose any
orthonormal basis $(u_j)$ of $L^2(\m D^*)$; its entries are
$\int\eta u_j\overline{u_i}\,\dd A$. The basis must include functions
outside the ranges of $J$ and $\bar J$. In particular, the intermediate
$\mathcal B_*$ and resolvent in $Q_\eta$ act on the full $L^2$ space, and
cannot be replaced by products of the compressed matrices without proof.

For an explicit matrix of $\mathcal B_*$ on the full space, let
$P_l$ be the Legendre polynomial of degree $l$, normalized by $P_l(1)=1$,
and set, for $k\in\m Z$ and $l\ge0$,
\[
 p_l(r)=\sqrt{\frac{2l+1}{\pi}}\,r^{-2}P_l(2r^{-2}-1),
 \qquad u_{k,l}(re^{\ii\theta})=p_l(r)e^{\ii k\theta}.
\]
The functions $(u_{k,l})_{k\in\m Z,l\ge0}$ form an orthonormal basis of
$L^2(\m D^*)$. Indeed, angular orthogonality and the substitution
$s=r^{-2}$ give
\[
 2\pi\int_1^\infty p_l(r)p_{l'}(r)r\dd r
 =\sqrt{(2l+1)(2l'+1)}
 \int_0^1P_l(2s-1)P_{l'}(2s-1)\dd s=\delta_{ll'}.
\]
For $k\in\m Z$, define
\[
 (T_kp_l)(r)=
 \begin{cases}
 p_l(r)-2(1-k)r^{k-2}\displaystyle\int_1^r s^{1-k}p_l(s)\dd s,
 &k\le0,\\[1ex]
 p_l(r)-2(k-1)r^{k-2}\displaystyle\int_r^\infty s^{1-k}p_l(s)\dd s,
 &k\ge1.
 \end{cases}
\]
Then the matrix entries of $\mathcal B_*$ in this basis are
\begin{equation}\label{eq:beurling-matrix}
 (\mathcal B_*)_{(j,m),(k,l)}
 =\delta_{j,k-2}\,2\pi\int_1^\infty p_m(r)(T_kp_l)(r)r\dd r.
\end{equation}
To check the formula, extend $u_{k,l}$ by zero to $\m D$ and write
$\mathcal C u_{k,l}=e^{\ii(k-1)\theta}c(r)$ for $r>1$. The relations
\[
 \frac12\left(c'(r)-\frac{k-1}{r}c(r)\right)=p_l(r),\qquad
 \mathcal B_*u_{k,l}
 =\frac12e^{\ii(k-2)\theta}
 \left(c'(r)+\frac{k-1}{r}c(r)\right)
\]
give $\mathcal B_*u_{k,l}=e^{\ii(k-2)\theta}T_kp_l$ after imposing the
Cauchy transform boundary condition. The case $k=1$ is $T_1p_l=p_l$.
The term $p_l(r)$ is part of the principal-value Beurling transform.
In the same basis,
\[
 (M_\eta)_{(j,m),(k,l)}
 =\int_{\m D^*}\eta(z)u_{k,l}(z)\overline{u_{j,m}(z)}\dd A(z).
\]
These matrices represent bounded operators on $\ell^2(\m Z\times\m N_0)$;
thus the inverse in \eqref{eq:Q} is the operator inverse of the full
matrix $I-M_\eta\mathcal B_*$, not the inverse of a compression.

The connection with $\dot H^{1/2}(\m S^1)$ is as follows. Give the negative
frequency subspace the orthonormal basis $z^{-n}/\sqrt n$. If
$u_-(z)=\sum_{n\ge1}x_nz^{-n}/\sqrt n$ is its harmonic extension to
$\m D^*$, then the equality
\[
 -\frac1{\sqrt\pi}\partial_z u_-(z)
 =-\frac1{\sqrt\pi}\sum_{n\ge1}\frac{x_n}{\sqrt n}
 (-n)z^{-n-1}
 =\sum_{n\ge1}x_n\sqrt{\frac n\pi}\,z^{-n-1}=Jx
\]
holds in $L^2(\m D^*)$. It follows first for finite sums and then for
$x\in\ell^2$ by $L^2$ convergence.
The analogous anti-holomorphic extension uses $\bar J$. Thus $J$ is the
Dirichlet-energy realization of the negative Fourier coordinates, whereas
$G_1,\ldots,G_4$ in \eqref{eq:TT} are operators between Fourier
coordinates. The two types of operator are related by
Lemma~\ref{lem:resolvents}.

\subsection{Gaussian estimates for smooth source deformations}\label{subsec:crit-gaussian-proof}

\begin{proof}[Proof of Lemma~\ref{lem:crit-gaussian-flow}]

Let $\mathcal H$ be the real mean-zero Cameron--Martin space, with
$\|f\|_{\mathcal H}^2=\sum_{n\ne0}|n||\widehat f(n)|^2$.
We perform the Gaussian calculations using its associated isonormal
process. Define
\[
 C_tf=\pi_0(f\circ\phi_t),\qquad L_vf=\pi_0(v_df'),\qquad
 S_t=C_tC_t^*,\qquad b_t=\pi_0(2\log\phi_t').
\]
We first estimate the covariance and mean. The derivatives $\phi_t'$
are uniformly bounded above and below for $|t|\le T$, so the
difference-quotient norm gives uniform bounds for $C_t$ and $C_t^{-1}$
on $\mathcal H$. Thus the spectrum of $S_t$ lies in a fixed compact
subinterval of $(0,\infty)$.

In the Fourier basis $|n|^{-1/2}e^{2\pi\mathrm{i}nx}$, the matrix of
$A=L_v+L_v^*$ is
\[
 A_{nm}=2\pi\mathrm{i}(\operatorname{sgn}m-\operatorname{sgn}n)
 \sqrt{|nm|}\,\widehat v_d(n-m).
\]
Consequently,
\[
 \|A\|_{\mathrm{HS}}^2
 =\frac{4(2\pi)^2}{3}\sum_{k\ge2}(k^3-k)|\widehat v_d(k)|^2\le C.
\]
The same scaling gives the additional bound
\begin{equation}\label{eq:crit-flow-LA}
 \|L_vA\|_{\mathrm{HS}}\le C.
\end{equation}
For completeness, let $p=n-\ell$ and $q=\ell-m$. In a nonzero term of
$A_{\ell m}$, $|\ell|+|m|=|q|$. With $p$ fixed, the squared Hilbert--Schmidt
norm of the corresponding contribution is at most
\[
 C|\widehat v_d(p)|^2\sum_q|\widehat v_d(q)|^2
 \bigl(|q|^5+|p||q|^4\bigr).
\]
For every $N$, $|\widehat v_d(k)|\le C_Nd^2(1+d|k|)^{-N}$. Thus
\[
 \sum_p|\widehat v_d(p)|=O(d),\quad
 \sum_p|p|^{1/2}|\widehat v_d(p)|=O(d^{1/2}),
\]
\[
 \sum_q|q|^5|\widehat v_d(q)|^2=O(d^{-2}),\quad
 \sum_q|q|^4|\widehat v_d(q)|^2=O(d^{-1}).
\]
Minkowski's inequality proves \eqref{eq:crit-flow-LA}. The matrix
calculation is first made on Fourier polynomials. Closedness of the flow
generator then identifies the Hilbert--Schmidt extension.

In Hilbert--Schmidt norm,
\[
 S_t'=C_tAC_t^*,\qquad
 S_t''=C_t\bigl(L_vA+(L_vA)^*\bigr)C_t^*.
\]
Also $b_t'=2\pi_0(v_d'\circ\phi_t)$ and
$b_t''=2\pi_0((v_dv_d'')\circ\phi_t)$. Scale invariance of the
$\dot H^{1/2}$ seminorm and the bounds on $C_t$ give
\begin{equation}\label{eq:crit-flow-bounds}
 \sup_{d,\,|t|\le T}\left(
 \|S_t-I\|_{\mathrm{HS}}+\|S_t'\|_{\mathrm{HS}}
 +\|S_t''\|_{\mathrm{HS}}+\sum_{j=0}^2\|b_t^{(j)}\|_{\mathcal H}
 \right)<\infty.
\end{equation}

We next derive the density bounds. Let $W$ be the reference isonormal
Gaussian field and $\mathsf R_t=S_t^{-1}$. Use the Wick quadratic form
$I_2$ normalized by
$E[I_2(K)^2]=2\|K\|_{\mathrm{HS}}^2$. The Gaussian density has
\begin{equation}\label{eq:crit-log-gaussian}
 \log D_t=\tfrac12I_2(I-\mathsf R_t)+W(\mathsf R_tb_t)
 +\tfrac12\log\det{}_2\mathsf R_t-\tfrac12\langle b_t,\mathsf R_tb_t\rangle.
\end{equation}
Let $P_N$ project onto the first $N$ pairs of real Fourier modes, and set
$S_{t,N}=I+P_N(S_t-I)P_N$, $b_{t,N}=P_Nb_t$.
Write $\mathsf R_{t,N}=S_{t,N}^{-1}$. On $P_N\mathcal H$, let $D_{t,N}$ be the
density of $N(b_{t,N},S_{t,N}|_{P_N\mathcal H})$ relative to the standard
Gaussian law. It satisfies the corresponding version of
\eqref{eq:crit-log-gaussian}. The finite-dimensional log densities converge in probability to
\eqref{eq:crit-log-gaussian}. To pass to the densities and their
derivatives, we first prove uniform integrability. For $p>1$ sufficiently
close to $1$, uniformly in $d,t,N$,
\begin{align*}
 E[D_{t,N}^p]
 &=\det(S_{t,N})^{-(p-1)/2}
 \det\!\bigl(pI-(p-1)S_{t,N}\bigr)^{-1/2}\\
 &\quad\times\exp\!\left\{\frac{p(p-1)}2
 \left\langle b_{t,N},
 \bigl(pI-(p-1)S_{t,N}\bigr)^{-1}b_{t,N}\right\rangle\right\}
 \le C.
\end{align*}
The first-order terms of the two log determinants cancel at $S=I$;
\eqref{eq:crit-flow-bounds} controls the remaining terms.

Write $\ell_t=\log D_t$ and $\ell_{t,N}=\log D_{t,N}$. The first two
derivatives of \eqref{eq:crit-log-gaussian} have uniformly bounded
Gaussian chaos coefficients. In particular,
$\mathsf R_t'=-\mathsf R_tS_t'\mathsf R_t$ and
\[
 \frac{\dd}{\dd t}\log\det{}_2\mathsf R_t
 =\operatorname{Tr}[(\mathsf R_t^{-1}-I)\mathsf R_t'],
\]
where the product in the trace is trace class. Its second derivative is
\[
 \operatorname{Tr}\!\left[-\mathsf R_t^{-1}\mathsf R_t'\mathsf R_t^{-1}\mathsf R_t'
 +(\mathsf R_t^{-1}-I)\mathsf R_t''\right].
\]
Both terms are trace class by \eqref{eq:crit-flow-bounds}. The same
bounds hold for the finite-rank approximations. Gaussian chaos moment
bounds therefore give, for every finite $q\ge1$,
\[
 \sup_{d,N,\,|t|\le T}\sum_{j=0}^2
 \|\partial_t^j\ell_{t,N}\|_{L^q}<\infty.
\]
Choose $a$ with $1<a<p$, where $p$ is from the preceding density bound.
For any polynomial $B$ in $\ell_{t,N}'$ and $\ell_{t,N}''$, H\"older's
inequality gives
\[
 E[|D_{t,N}B|^a]
 \le E[D_{t,N}^p]^{a/p}
 E[|B|^{ap/(p-a)}]^{(p-a)/p}.
\]
Apply this to $B=1$, $\ell_{t,N}'$ and
$(\ell_{t,N}')^2+\ell_{t,N}''$. Since
\[
 D_{t,N}'=D_{t,N}\ell_{t,N}',\qquad
 D_{t,N}''=D_{t,N}\bigl((\ell_{t,N}')^2+\ell_{t,N}''\bigr),
\]
the densities and their first two derivatives are uniformly bounded in
$L^a$, hence uniformly integrable.

We can now pass to the limit. For fixed $d$, the coefficients $S_t-I,b_t$
and their first two derivatives range over compact sets in Hilbert--Schmidt or Hilbert norm
as $t$ varies in $[-T,T]$. Hence their Fourier compressions converge
in these norms, uniformly in $t$, together with both derivatives.
The common spectral lower bound and the identities
\[
 \begin{aligned}
 \mathsf R_{t,N}-\mathsf R_t
 &=\mathsf R_{t,N}(S_t-S_{t,N})\mathsf R_t,\\
 \mathsf R_t''&=2\mathsf R_tS_t'\mathsf R_tS_t'\mathsf R_t
                 -\mathsf R_tS_t''\mathsf R_t
 \end{aligned}
\]
give the same convergence for the inverses. The trace formulas above
then give convergence of the determinant terms. Thus, for $j\le2$,
$\partial_t^j\ell_{t,N}\to\partial_t^j\ell_t$ in every finite $L^q$,
uniformly on this time interval.

The uniform $L^a$ bounds proved above now imply, for each $t$,
\[
 (D_{t,N},D_{t,N}',D_{t,N}'')
 \longrightarrow
 \bigl(D_t,D_t\ell_t',D_t((\ell_t')^2+\ell_t'')\bigr)
 \quad\text{in }L^1.
\]
Denote the three limits by $Y_0,Y_1,Y_2$. Continuity of the Gaussian
coefficients and uniform integrability make these limits continuous in
$L^1$. The uniform bounds also allow passage to the limit in the
finite-dimensional fundamental theorem of calculus:
\[
 Y_0(t)-Y_0(s)=\int_s^tY_1(u)\,\dd u,\qquad
 Y_1(t)-Y_1(s)=\int_s^tY_2(u)\,\dd u.
\]
This proves the asserted $C^2(L^1)$ regularity and the first bound in
\eqref{eq:crit-density-derivatives}. Testing against bounded cylindrical
functions identifies $D_t$ as the density of the transformed Gaussian law.

Finally, the coefficient bounds give $\ell\in
L^2(\Omega;W^{1,2}([-T,T]))$. The one-dimensional Sobolev inequality gives
\[
 E\!\left[\sup_{|t|\le T}|\ell_t|^2\right]
 \le C_T\int_{-T}^T E[|\ell_t|^2+|\ell_t'|^2]\,\dd t\le C_T.
\]
Its continuous representative also provides a jointly measurable version
of $D_t=e^{\ell_t}$. This proves the second bound.

To obtain the welding identity, the critical change-of-coordinates formula for a smooth circle
homeomorphism gives
$\mu_{X_1\circ\phi_t+2\log\phi_t'}=\phi_t^*\mu_1$; see
\cite[Lemma~3.8]{F1}. Its cumulative function is
\[
 F_1(\phi_t(x))-F_1(\phi_t(0)).
\]
Replace $U$ by $U+F_1(\phi_t(0))$. Conditional on the two fields, this is
still Haar and hence is independent of the transformed source field and
the target field. The resulting welding is $H_4\circ\phi_t$, and the
change of law of the source field has density $D_t$.
\end{proof}

\subsection{Inverse tails of Bessel-weighted measures}\label{subsec:crit-bessel-proof}

\begin{proof}[Proof of Lemma~\ref{lem:crit-random-start-bessel}]

We locate each inverse tail between two last-exit times of the same
Bessel process. On $\{R_0<a\}$, set
\[
 \widehat L_a=\sup\{t\ge0:R_t=a\},\qquad
 L_a=\widehat L_a/2=\sup\{s\ge0:R_{2s}=a\}.
\]
These times are finite almost surely. We first work in standard BES time.
The BES$(3)$ generator is
$\mathcal L=\frac12\partial_x^2+x^{-1}\partial_x$. For $x>a$ its
probability of hitting $a$ is $a/x$. The Doob transform of the process
killed at $a$ by $q_a(x)=1-a/x$ has generator
\[
 q_a^{-1}\mathcal L(q_af)=\tfrac12f''+(x-a)^{-1}f'.
\]
Thus a BES$(3)$ started at $a+\epsilon$ and conditioned never to return to
$a$, shifted by $a$, is BES$(3)$ started at $\epsilon$. Starting below $a$,
consider successive excursions from $a$ to $a+\epsilon$. Apply the strong
Markov property at their hitting times and sum over the first excursion
which never returns to $a$. As $\epsilon\downarrow0$, its starting time
approaches the last exit $\widehat L_a$, and the shifted conditioned process
converges to BES$(3)$ from $0$. This proves that the post-last-exit process,
shifted by $a$, has that law. The argument uses actual hitting times; it
does not apply the strong Markov property at $\widehat L_a$ itself.

Let $\mathcal R$ be a BES$(3)$ process started at zero, let $g_d$ be
its last visit to $d$, and let $\Phi_{\mathrm N}$ be the standard normal
distribution function. Its radial Gaussian density and its probability
of hitting $d$ from $x>d$ give
\[
 P(g_d>t)=E[\min(1,d/\mathcal R_t)]
 =2\Phi_{\mathrm N}(d/\sqrt t)-1=P(\tau_d>t),
\]
where $\tau_d$ is the first time a standard Brownian motion reaches $d$.
It follows, after the time change $t=2s$, that on $\{R_0<a\}$,
\[
 \operatorname{Law}(L_b-L_a\mid R_0<a)
 =\operatorname{Law}(\tau_{b-a}/2),\qquad b>a.
\]
The post-last-exit law also holds after conditioning on $R_0<a$ because
it does not depend on the initial value below $a$.

Take $\mathcal R$ independent of the two measures and define
\[
 U_\sigma=\int_0^\infty e^{-\mathcal R_{2u}}\mathcal N_\sigma(\dd u).
\]
For $0<p<1$, block decomposition and the inequality
$(\sum_n x_n)^p\le\sum_n x_n^p$ give $E[U_\sigma^p]<\infty$. Indeed,
realize $\mathcal R_t=|W_t|$ with a standard three-dimensional Brownian
motion $W$, and put $V_n=\sup_{0\le v\le2}|W_{2n+v}-W_{2n}|$. Then
\[
 E\!\left[\left(\int_n^{n+1}e^{-\mathcal R_{2u}}
 \mathcal N_\sigma(\dd u)\right)^p\right]
 \le C_p E[e^{-p\mathcal R_{2n}}]
 \le C_p(1+n)^{-3/2}.
\]
Indeed, $\mathcal R_{2u}\ge\mathcal R_{2n}-V_n$ on the block. The
variable $V_n$ is independent of $W_{2n}$ and of $\mathcal N_\sigma$,
and $E[e^{pV_n}]<\infty$. Stationarity and the unit-block moment bound
control the remaining measure factor. The last estimate follows from
the three-dimensional Gaussian density.
Similarly,
\[
 U_\sigma\ge e^{-\sup_{0\le u\le1}\mathcal R_{2u}}
 \mathcal N_\sigma([0,1])
\]
and independence give $E[U_\sigma^{-q}]<\infty$ for every $q>0$.
The same block estimate, with the Gaussian density bound uniform in its
starting point, makes $J_\sigma(0)$ finite almost surely. Full support and
non-atomicity make $J_\sigma$ continuous and strictly decreasing to zero.

On $\{R_0<a\}$, stationarity and independence permit shifting both
measures by the random time $L_a$. Precisely, on that event, for every
bounded measurable functional $F$ of the pair of measures,
\[
 E[F(\theta_{L_a}\mathcal N_-,\theta_{L_a}\mathcal N_+)\mid R]
 =E[F(\mathcal N_-,\mathcal N_+)].
\]
Here $(\theta_s\mathcal N)(B)=\mathcal N(s+B)$ for a Borel set $B$.
Thus, jointly for the two sides,
\[
 \operatorname{Law}\bigl((e^aJ_\sigma(L_a))_{\sigma=\pm}
 \mid R_0<a\bigr)=\operatorname{Law}((U_\sigma)_{\sigma=\pm}).
\]
Fix $K>0$ and set $a=h-M-K$, $b=h+M+K$. On $\{R_0<a\}$, if
$e^aJ_\sigma(L_a)>e^{-K}$ and $e^bJ_\sigma(L_b)\le e^K$ for both
$\sigma$, every $S_\sigma(h+u)$ with $|u|\le M$ lies in $[L_a,L_b]$.
Markov's inequality for the positive and negative moments, followed by the
last-exit increment formula, yields
\begin{align*}
 &P\!\left(\max_\sigma S_\sigma(h+M)
 -\min_\sigma S_\sigma(h-M)>T\right)\\
 &\quad\le P(R_0\ge h-M-K)
 +C_{p,q}(e^{-pK}+e^{-qK})
 +P(\tau_{2(M+K)}>2T).
\end{align*}
After taking the upper limit as $h\to\infty$, the first term vanishes
because $R_0$ is finite almost surely. We can then make the remaining
terms arbitrarily small by first choosing $K$ and then $T$ large.
This proves the assertion.
\end{proof}

\subsection{Construction of the truncated rooted law}\label{subsec:crit-root-proof}

\begin{proof}[Proof of Lemma~\ref{lem:crit-barrier-root}]

For $\epsilon=r_ce^{-s}$, write the derivative approximation as
\[
 \nu_{h,\epsilon}(\dd y)=
 \left(\log(1/\epsilon)-\tfrac12h_\epsilon(y)\right)
 \epsilon e^{h_\epsilon(y)}\,\dd y,
\]
where $h_\epsilon(y)$ is the semicircle average, and set
\[
 d_\epsilon^b(\dd y)=
 \frac{b+2s-h_\epsilon(y)}2\,\epsilon e^{h_\epsilon(y)}
 \1_{\{h_{r_ce^{-u}}(y)<b+2u,\ 0\le u\le s\}}\1_J(y)\,\dd y.
\]
Let $C_b$ impose the same barrier for every $y\in J$ and $u\ge0$.
On $C_b$, as measures on $J$,
\[
 d_\epsilon^b=\nu_{h,\epsilon}
                  +(b/2-L)\epsilon e^{h_\epsilon}\,\dd y.
\]
The derivative approximations converge weakly in probability to $\nu_h$;
see \cite[Definition~2.6]{PowellSepulveda}. The second term tends to zero
in probability by \cite[Remark~3.6]{HP21_critical}.

We need uniform integrability of $d_\epsilon^b(J)$. On a half-disk
containing the relevant semicircles, write $h=g+H$, with mixed-boundary
GFF $g$ and independent harmonic $H$. Even reflection across the free
boundary gives $H_\epsilon(y)=H(y)$, hence
\[
 d_{h,\epsilon}^b(\dd y)
       =e^{H(y)}d_{g,\epsilon}^{\,b-H(y)}(\dd y).
\]
Let $K$ bound $|H|$ on the compact neighborhood. On $\{K\le k\}$,
the barrier for $g$ at height $b-k$ implies $C_b(h)$, so
\cite[Lemma~3.5]{HP21_critical} gives $P(C_b)\to1$ by first sending $b$
and then $k$ to infinity. On $\{K\le k\}$ we also have
$d_{h,\epsilon}^b\le e^k d_{g,\epsilon}^{b+k}$.
The mixed-field bounds of \cite[Lemmas~3.4--3.5]{HP21_critical} give
uniform integrability on $\{K\le k\}$ and
\[
 \sup_\epsilon E[d_{h,\epsilon}^b(J)\1_{\{K>k\}}]
 \le C E[e^K(1+|b|+K)\1_{\{K>k\}}]\longrightarrow0.
\]
Here the estimate uses
$E[d_{g,\epsilon}^a(\dd y)]/\dd y\le C(1+|a|)$ and the Gaussian
exponential moments of $K$.

For fixed $y$, \eqref{eq:crit-polar-covariance} makes the finite-cutoff
density radial. With $A_s=A_0+B_{2s}$, its radial factor, apart from
$r_c/2$, is
\[
 (b+2s-A_s)e^{A_s-s}
 \1_{\{A_u<b+2u,\ 0\le u\le s\}}.
\]
Conditional on $A_0=a<b$, the factor $e^{B_{2s}-s}$ changes
$B_{2u}$ to Brownian motion with drift $2u$. The remaining factor is
the positive-distance Doob transform for Brownian motion killed at zero;
it gives $b+2u-A_u=R_{2u}$ up to time $s$. This is precisely
\cite[Lemma~3.9]{HP21_critical}, with variance rate $2$, exponential
tilt parameter $1$, and barrier $b$; see also its Lemma~3.10 for the
corresponding rooted field law. At time zero the factor is
$(b-A_0)_+e^{A_0}$, which gives the stated tilted initial distribution.
The lateral field stays independent because the density is radial.
In particular,
\[
 E[d_\epsilon^b(J)]=|J|r_cm_b/2=Z_b.
\]
To pass to a joint limit, let
$\omega=(h,\nu_h,(h_{r_ce^{-u}}(\cdot))_{u\ge0})$ and let $\mathsf P$
be its law, using respectively the $H^{-1}_{\rm loc}$, vague, and locally
uniform topologies. These form a Polish product space; the last two
coordinates are measurable functions of $h$. Define
$\Theta_{b,\epsilon}(\dd\omega,\dd y)
=Z_b^{-1}\mathsf P(\dd\omega)d_\epsilon^b(\dd y)$. For Borel $A$,
\[
 \Theta_{b,\epsilon}(A\times J)
 \le Z_b^{-1}\left(
 L_0\mathsf P(A)+
 E[d_\epsilon^b(J)\1_{\{d_\epsilon^b(J)>L_0\}}]\right).
\]
Uniform integrability gives tightness and a subsequential weak limit
$\Theta_b$ whose $\omega$-marginal is $r_b\mathsf P$ for a weak
$L^1(\mathsf P)$ limit $r_b$. Disintegrate and define
\[
 \Theta_b=r_b(\omega)\mathsf P(\dd\omega)K_b(\omega,\dd y),
 \qquad
 \overline d_h^b(\dd y)=Z_b r_b(\omega)K_b(\omega,\dd y).
\]
Absolute continuity of the marginal preserves the graph of the two
measurable functions of $h$, so this defines a random measure indexed
by $h$, with expected mass $Z_b$.

For bounded measurable $F(\omega)$ and continuous $g$ on $J$, uniform
absolute continuity allows continuous approximation of $F$ uniformly
under these marginals. The identity on $C_b$, convergence of
$\nu_{h,\epsilon}$, and uniform integrability yield
\[
 E_{\Theta_b}[\1_{C_b}F(\omega)g(y)]
 =Z_b^{-1}E\!\left[\1_{C_b}F(\omega)
                   \int_Jg(y)\nu_h(\dd y)\right].
\]
A countable separating family of $g$ and a monotone-class argument
prove $\1_{C_b}\overline d_h^b=\1_{C_b}\nu_h|_J$.

Finally fix $0<a<T$ and a buffer $0<\eta<a$. Let $X_y$ be the field
pulled back to $a-\eta\le s\le T+\eta$ and $V_y$ the pulled-back length
on $a\le s\le T$. Common local regularizations admit a deterministic
subsequence defining a jointly Borel GMC functional $\mathcal F$ with
$V_y=\mathcal F(y,X_y)$ for $\mathsf P(\dd\omega)\dd y$-almost every
$(\omega,y)$; this follows from fixed-$y$ convergence and Fubini.
It holds under every $\Theta_{b,\epsilon}$ by absolute continuity.
Once the cutoff exceeds $T+\eta$, the radial Doob-transform law of
$(y,X_y)$, and hence of $(y,X_y,V_y)$, is independent of the cutoff.
Buffered pullbacks are continuous under joint convergence: shifted
field tests converge in $H^1$, and shifted measure tests converge
uniformly on their common compact support. The same triple law
therefore passes to $\Theta_b$. A countable exhaustion by annuli
identifies the field and its GMC jointly, and also passes the radial
and lateral laws to the limit.
\end{proof}

\bibliographystyle{alpha}
\bibliography{ref}
\end{document}